\documentclass[3p,times]{elsarticle}

\usepackage{framed,multirow}

\usepackage{amssymb}
\usepackage{latexsym}
\usepackage{graphicx}	
\usepackage{amsmath}	
\usepackage{subfig}
\usepackage{amsthm}

\usepackage{mathtools}
   \DeclarePairedDelimiter{\norm}{\lVert}{\rVert}

\usepackage{url}
\usepackage[dvipsnames]{xcolor}
\definecolor{newcolor}{rgb}{.8,.349,.1}
\usepackage{multirow}

\newcommand{\mbf}[1]{\mathbf{#1}}			%
\newcommand{\x}{\mbf{x}}
\newcommand{\xcn}{\mbf{x}_{\mbf{c}}^n}               
\newcommand{\xcin}{\mbf{x}_{\mbf{c}_i}^n}    
\newcommand{\xcnp}{\mbf{x}_\mbf{c}^{n+1}}               
\newcommand{\xcinp}{\mbf{x}_{\mbf{c}_i}^{n+1}}    
\newcommand{\xcijnp}{\mbf{x}_{a_{k},\,i_2}^{n+1}} 
\newcommand{\xciknp}{\mbf{x}_{a_{k},\,i_3}^{n+1}}
\newcommand{\xcinpstar}{\hat{\mbf{x}}_{\mbf{c}_{i}}^{n+1}} 
\newcommand{\xcistar}{\mbf{x}_{\mbf{c}_{i}}^{\ast}}

\newcommand{\xbin}{\mbf{x}_{\mbf{b}_i}^n}    
\newcommand{\xbinp}{\mbf{x}_{\mbf{b}_i}^{n+1}}    
\newcommand{\xbixn}{{x}_{{b}_i}^n}    
\newcommand{\xbiyn}{{y}_{{b}_i}^n}

\newcommand{\Q}{\mathbf{Q}}

\newcommand{\0}{\mathbf{0}}
\renewcommand{\u}{\mathbf{u}}
\newcommand{\w}{\mathbf{w}}
\newcommand{\q}{\mathbf{q}}
\newcommand{\F}{\mathbf{F}}

\newcommand{\f}{\mathbf{f}}
\newcommand{\g}{\mathbf{g}}

\newcommand{\halb}{\frac{1}{2}} 
 
\renewcommand{\v}{\mathbf{v}}
\newcommand{\B}{\mathbf{B}}

\newcommand{\de}[2]{\frac {\partial #1}{\partial#2}}

\renewcommand{\S}{\mathbf{S}}
\newcommand{\be}{\begin{equation} \begin{aligned} }
\newcommand{\ee}{\end{aligned} \end{equation}}

\renewcommand{\epsilon}{\varepsilon}
\renewcommand{\phi}{\varphi}
\newcommand{\phitilde}{\tilde{\phi}}

\newcommand{\Popt}{\mathbf{P}_{\mathrm{opt}}}
\DeclareMathOperator{\argmin}{argmin}

\newcommand{\RIcolor}[1]{{\leavevmode\color{black} #1}}
\newcommand{\RIIcolor}[1]{{\leavevmode\color{black} #1}}
\newcommand{\RIIIcolor}[1]{{\leavevmode\color{black} #1}}
\newcommand{\RIVcolor}[1]{{\leavevmode\color{black} #1}}
\newcommand{\RVcolor}[1]{{\leavevmode\color{black} #1}}
\newcommand{\RVIcolor}[1]{{\leavevmode\color{black} #1}}

\newcommand{\RShashkovcolor}[1]{{\leavevmode\color{black} #1}}

\journal{Journal of Computational Physics}

\begin{document}


\begin{frontmatter}

\title{High order direct Arbitrary-Lagrangian-Eulerian schemes on moving Voronoi meshes with topology changes} 


\author[UniTN]{Elena Gaburro$^{*}$}
\ead{elena.gaburro@unitn.it}
\cortext[cor1]{Corresponding author}

\author[UniFE]{Walter Boscheri}
\author[UniTN]{Simone Chiocchetti}
\author[UniWue]{Christian Klingenberg}
\author[MPI]{Volker Springel}
\author[UniTN]{\\Michael Dumbser}
%
\address[UniTN]{Department of Civil, Environmental and Mechanical Engineering, University of Trento, Via Mesiano 77, 38123 Trento, Italy} 
\address[UniFE]{Department of Mathematics and Computer Science, University of Ferrara, via Machiavelli 30, 44121 Ferrara, Italy }
\address[UniWue]{Department of Mathematics at W\"urzburg University, Emil Fischer Str. 40, W\"urzburg, 97074, Germany}
\address[MPI]{Max-Planck-Institut f\"ur Astrophysik, Karl-Schwarzschild-Str 1, D-85748 Garching, Germany}


\begin{abstract}
We present a new family of very high order accurate direct Arbitrary-Lagrangian-Eulerian (ALE) 
Finite Volume (FV) and Discontinuous Galerkin (DG) schemes 
for the solution of nonlinear hyperbolic PDE systems on moving \RIIcolor{two-dimensional} Voronoi meshes 
that are \textit{regenerated} at each time step 
and which explicitly allow \textit{topology changes} in time. 
The Voronoi tessellations are obtained from a set of generator points 
that move with the local fluid velocity. We employ an AREPO-type approach~\cite{Springel},  
which rapidly rebuilds a new high quality mesh exploiting the previous one, 
but \textit{rearranging} the element shapes and neighbors 
in order to guarantee that the mesh evolution is robust 
even for vortex flows and for very long simulation times.
The old and new Voronoi elements associated to the same generator point 
are connected in space--time to construct closed \textit{space--time} control volumes, 
whose bottom and top faces may be polygons with a different number of sides. 
We also need to incorporate some degenerate \textit{space--time sliver elements}, 
which are needed in order to fill the space--time holes that arise because of the topology 
changes in the mesh between time $t^n$ and time $t^{n+1}$.	
The final ALE FV-DG scheme is obtained by a novel redesign of the high order accurate fully discrete 
direct ALE schemes of Boscheri and Dumbser~\cite{Lagrange3D,ALEDG},  which have been extended here 
to general moving Voronoi meshes and space--time sliver elements. Our new numerical scheme is 
based on the integration over  arbitrary shaped closed space--time control volumes combined with a 
fully-discrete space--time  conservation formulation of the governing hyperbolic PDE system. 
In this way the discrete solution is \textit{conservative} and satisfies the geometric 
conservation law (GCL) \textit{by construction}. 
Numerical convergence studies as well as a large set of benchmark problems for hydrodynamics and  
magnetohydrodynamics (MHD) demonstrate the accuracy and robustness of the proposed method. 
Our numerical results clearly show that the new combination of very high order schemes with 
regenerated meshes that allow topology changes in each time step lead to substantial improvements  
compared to direct ALE methods on moving conforming meshes without topology change.  
\end{abstract}


\end{frontmatter}


\section{Introduction}

The aim of this work is to present a novel family of explicit \textit{arbitrary high order accurate} direct ALE  Finite Volume (FV) and Discontinuous Galerkin (DG) schemes on moving Voronoi meshes that are \textit{regenerated} at each time-step and which consequently allow also \textit{topology changes} of the computational grid during the time evolution of the PDE system. 
The \textit{main novelty} lies in the use of a \textit{space--time conservation formulation} of the governing PDE  system over closed, non-overlapping \textit{space-time} control volumes that are constructed from the moving,  regenerated Voronoi meshes between time $t^n$ and time $t^{n+1}$. 
On these closed space--time control volumes the governing equations are then directly integrated by means 
of a high order fully discrete one-step ADER method. 
To the best knowledge of the authors, this is the first time that a unified framework for arbitrary high order accurate explicit non-oscillatory direct ALE FV and DG schemes on moving Voronoi meshes is developed, with an embedded mesh generator that builds a new mesh with a  different topology at each time step.

\subsection{State of the art}
Lagrangian algorithms~\cite{Neumann1950, BENSON1992235, Despres2009, de2017lagrange, gaburro2018Thesis, munz94, Caramana1998, Smith1999, Maire2007} are characterized by a moving computational mesh displaced with a velocity chosen as close as possible to the local fluid velocity. In the Lagrangian description of the fluid, the nonlinear convective terms disappear and, as a consequence, Lagrangian schemes exhibit virtually no numerical dissipation at contact discontinuities and  material interfaces. 
Therefore, the aim of these methods is to \textit{reduce the numerical dissipation errors} due to the convective terms, so that contact discontinuities are sharply captured and material interfaces can be properly identified and tracked.

Lagrangian finite volume schemes~\cite{munz94,despres2005lagrangian,Depres2012,ShashkovCellCentered,Maire2009,Maire2009b,Maire2010,Maire2011,burton2013cell,burton2015reduction} have been developed for the solution of nonlinear hyperbolic systems of PDEs, using the conservation form of the equations based on the physically conserved quantities like mass, momentum and total energy. Higher order Lagrangian-type schemes have been introduced in~\cite{chengshu1,chengshu2,chengshu3}, where high order of accuracy in space is achieved with the aid of an ENO/WENO reconstruction and Runge-Kutta time stepping guarantees high order time discretization as well. Contrarily to the \textit{cell-centered} methods listed so far, where all variables are located at the cell center of the primal mesh, \textit{staggered} Lagrangian schemes~\cite{StagLag,LoubereSedov3D,maire_loubere_vachal10} define the velocity at the grid vertexes and the other variables at the cell center, hence avoiding the need of a nodal solver to compute the mesh velocity of the grid nodes.

Another option for the numerical solution of hyperbolic conservation laws is given by Discontinuous Galerkin~\cite{reed} and Finite Element (FE) schemes~\cite{cremonesi2010lagrangian,cremonesi2011lagrangian,cremonesi2017explicit}, where the numerical solution is approximated by piecewise polynomials within each control volume. Robust Lagrangian DG schemes are presented in~\cite{jia2011new,morgan2018reducing,liu2018lagrangian,liu2019high} and they have been extended to third order for the first time in~\cite{Vilar1,Vilar2,Vilar3,Yuetal}, while high order FE methods applied to Lagrangian hydrodynamics and elasto-plasticity can be found in~\cite{scovazzi1,scovazzi2,Dobrev1,Dobrev2012,Dobrev2013}. 

Although these schemes are widely used, a common problem that affects almost all Lagrangian methods is the severe mesh 
distortion or mesh tangling that happens in the presence of shear flows, which may even cause a breakdown of the computation. This is the reason which led to the development of so-called Arbitrary-Lagrangian-Eulerian (ALE) methods~\cite{ShashkovCellCentered,ShashkovRemap1,ShashkovRemap3,ShashkovRemap4,ShashkovRemap5,MaireMM2,indALE-AWE2016}, where the mesh velocity can be chosen \textit{independently} of the local fluid velocity and thus the grid nodes can be moved at an arbitrary velocity. Cell-centered indirect ALE schemes aim at improving the mesh quality and the overall scheme  robustness by performing a purely Lagrangian phase with subsequent rezoning (mesh optimization) ~\cite{Winslow1997,KnuppRezoning,MaireRezoning} and remapping~\cite{Blanchard2016}, where the numerical solution defined on the old mesh is transferred onto the new grid. To overcome the problem of mesh tangling, sliding line techniques have also been proposed~\cite{Caramana2009, DelPino2010, LoubereSL2013}, which deal with moving \textit{nonconforming}  meshes, whose element sides can slide in order to accommodate the distortion induced by shear flows. 
\RIcolor{A very effective implicit DG method for dealing with weakly compressible Navier-Stokes flows with moving boundaries, using a 
tetrahedralization of space-time, has been
presented in~\cite{wang2015high}.}
For what concerns indirect ALE schemes, interesting techniques for handling the mesh motion have been introduced by the so-called Reconnection ALE (ReALE) algorithms~\cite{ReALE2010,ReALE2011,ReALE2015,ShashkovMultiMat3}, where the rezoning phase allows for topology changes at each time step of the computation. There, moving Voronoi tessellations have been employed and the obtained numerical results demonstrate that the flow features that have been computed in the Lagrangian phase can be better preserved compared to standard indirect ALE methods. \RVcolor{ReALE schemes have also proven to be particularly well suited for dealing with multimaterial fluid flows, in order to sharply capture the interfaces across different materials thanks to a conservative remapping phase which transfers the information from the old mesh to the newly generated one, which keeps tracking the interface. ReALE schemes can therefore be considered as the seminal works concerning moving mesh methods with topology changes. }

Among the different approaches that have been presented in the literature (pure Lagrangian, indirect ALE based on rezoning and remapping, ReALE as well as special nonconforming slide line treatments), a novel family of methods has been proposed, so-called \textit{direct} Arbitrary-Lagrangian-Eulerian (ALE) schemes. 
Also in the framework of direct ALE the mesh velocity can be chosen in an arbitrary way. Usually, it is chosen close to the local fluid velocity. However, the mesh quality can be optimized by a rezoning phase which takes place \textit{before} the computation of the numerical fluxes, hence allowing the space-time control volumes to be defined for each computational cell by connecting the element configuration at the current time level $t^n$ to the next time level $t^{n+1}$. Next, the mesh motion is taken into account \textit{directly} in the numerical flux computation of the FV or DG scheme, without needing any remeshing plus remapping strategy. Furthermore, such approaches naturally extend to unstructured meshes in multiple space dimensions~\cite{ALEMQF} and to slide line treatment with nonconforming meshes~\cite{gaburro2016direct, gaburro2018well}. Direct ALE schemes have been recently presented in~\cite{boscheri2014high,  Lagrange2D,Lagrange3D, ALEDG, boscheri2017high} by employing either very high order FV and DG schemes, also in combination with time-accurate local time stepping (LTS), see~\cite{ALELTS1D, ALELTS2D}. 
These works are characterized by a \textit{fixed} mesh topology, which makes it impossible to study phenomena affected by strong shear motion and vortex flows for very long simulation times, since mesh tangling would \textit{inevitably} occur and lead to a breakdown of the simulation before the final time is reached, \RIIcolor{unless strong mesh smoothing or relaxation procedures are introduced}. 
\RIIcolor{Also, it should be remarked that other ways to prevent, or at least to remarkably postpone, the breakdown of simulations consist in employing high order curvilinear meshes, see the results published in \cite{Dobrev2012,indALE-AWE2016,anderson2015monotonicity, anderson2018high}}. 
Moreover, notice that direct ALE schemes, even when constrained to a fixed connectivity, already ameliorate standard Lagrangian results for complex flow patterns.

From what was observed so far, the idea of \textit{allowing a change of topology} at each time step within the direct ALE framework arises.  A seminal work along this direction is represented by the \textit{AREPO} code of Springel and collaborators~\cite{Springel, springel2010moving, pakmor2014magnetic, pakmor2015improving}. \textit{AREPO} is a massively parallel second order accurate two- and three-dimensional direct ALE finite volume scheme on moving Voronoi tessellations that are rebuilt at each time step from a set of generator points which are moving with the local fluid velocity. The documented results obtained with the \textit{AREPO} technique clearly highlight the robustness and potential of that approach. Similar work in the context of finite element schemes can be found in the well-known particle finite element method of  O\~{n}ate and Idelsohn \textit{et al.}, see~\cite{PFEM1,PFEM2,PFEM3,PFEM4,PFEM5,PFEM6}. In the above-mentioned references, the mesh is completely regenerated at each time step, thus naturally allowing for large deformations and strong shear flows without causing mesh tangling and highly distorted elements. 

\subsection{Challenges of this work}

Up to now the \textit{AREPO} algorithm~\cite{Springel, springel2010moving} is at most second order accurate in space and time. We therefore believe that its results can still be improved by (i) increasing the order of accuracy of the underlying FV scheme in both space and time and by (ii) introducing a higher order DG method into the AREPO framework. 
However, above all, the main difficulty arises from the fact that high order direct ALE schemes need a complete knowledge of the \textit{space--time connectivity} between two consecutive time steps $t^n$ and $t^{n+1}$, and not only of the  \textit{spatial} connectivity at each time level.  
Furthermore, if a change of connectivity is allowed, the space-time connectivity does not coincide neither with the connectivity at time $t^n$, nor with the one at time $t^{n+1}$. Hence, an automatic way to construct the missing space-time connectivity from the available spatial connectivities at $t^n$ and $t^{n+1}$ must be found. 
In addition, the space--time control volumes should be allowed to have as bottom and top faces polygons with a different number of edges, and, moreover, even degenerate \textit{space--time sliver elements} must be incorporated in order to  fill the space-time holes that are caused by the changing topology.
With sliver elements we refer to space--time elements whose areas at time $t^n$ and $t^{n+1}$ are null, but whose space--time volume is not zero, see Sections~\ref{ssec.SpaceTimeConnection}~and~\ref{ssec.Sliver}. In other words, sliver elements exist only in the space-time volume strictly bounded between two consecutive time levels, therefore they must be taken into account only if the numerical scheme requires the full space-time connectivity. 

Finally, this kind of elements should be not only built, but also the one-step ADER finite volume and DG schemes must be substantially modified to handle the integration of the PDE over these new types of space-time control volumes.  
A proof of concept that direct ALE methods can work even on degenerate space-time elements was already given in ~\cite{gaburro2016direct} for second order FV schemes on moving nonconforming meshes, but a much greater effort is  necessary for dealing with such a general situation as the one treated in this work.

\subsection{Structure of the paper}

The rest of the paper is organized as follows. 
In Section~\ref{sec.NumMethod_mesh} we first introduce our \textit{moving computational mesh}, \RIIIcolor{the data representation over it, and the reconstruction procedure needed to obtain high order in space.
Next, in Section~\ref{sec.MeshEvolution} we introduce the mesh motion strategy which is obtained by computing the new coordinates of a set of points \RIIIcolor{(eventually with high order of accuracy, see Section~\ref{ssec.HighOrderTraj}) and \textit{re-drawing} around them a \textit{new} Voronoi tessellation, whose topology could differ from the previous one. 
Mesh optimization techniques (see Section~\ref{ssec.smoothing}) can be employed as well in order to improve the quality of the new tessellations.}

Sections~\ref{ssec.SpaceTimeConnection} and~\ref{ssec.Sliver} represent the first key ingredient of our algorithm:} we explain how to deal with the \textit{topology changes} that are caused by the \textit{regeneration} of the \textit{Voronoi tessellation} at each time step. Then, we explain how to automatically construct the space--time connectivity and the space--time sliver elements. 

Once this has been set up, in Section~\ref{sec.NumMethod_ALE_FV-DGscheme} we describe our \textit{direct} ALE FV-DG scheme, namely an algorithm belonging to the class of direct ALE $P_NP_M$ schemes~\cite{Dumbser2008}, which allows us to formulate a Finite Volume (FV) and a Discontinuous Galerkin (DG) scheme within a \textit{unique} framework. 
The method is first presented for {standard moving Voronoi elements}, i.e. Voronoi elements that are displaced without modifying their shape, i.e. the number of their nodes remains the same at each time level. 
Then, the method is extended to Voronoi elements with different bottom and top faces and finally to sliver elements in Sections~\ref{sssec.NumScheme_sliver_Predictor} and~\ref{sssec.NumScheme_sliver_Flux}, which is the second key ingredient of our scheme.
\RIIIcolor{In particular, for both types of elements we have detailed: i) the \textit{predictor} step,
 which is essential for obtaining high order in time in a fully-discrete one-step procedure, ii) the \textit{corrector} step, which allows the solution to be updated, and iii) the 
 \textit{limiter} that prevents spurious oscillations in the DG scheme.}

\RVcolor{Concerning a potential extension to compressible multi-material flows, the numerical approach proposed in this paper is obviously best suited for so-called \textit{diffuse interface models}, such as the Baer-Nunziato model of compressible two-phase flows  \cite{BaerNunziato1986,AndrianovWarnecke}, the Saurel-Abgrall model \cite{SaurelAbgrall,AbgrallSaurel,SaurelPetitpas} or the compressible multi-phase model recently proposed by Romenski et al. \cite{Rom1998,RomenskiTwoPhase2007,RomenskiTwoPhase2010}. }

In Section~\ref{sec.Results} we show a large set of numerical results, including convergence rates up to fifth order of 
accuracy in space and time for smooth problems,  
as well as a wide set of benchmark test cases solved with our ALE FV-DG schemes on moving Voronoi meshes with topology change for different systems of hyperbolic equations, namely the Euler equations of compressible gas dynamics, including the gravity source term, and the ideal MHD equations. The numerical results are commented and compared with available reference solutions wherever possible.

The paper is closed by some conclusive remarks and an outlook to future work in Section~\ref{sec.conclusion}.

\section{Numerical method I: handling a moving Voronoi tessellation with topology changes and data reconstruction} 
\label{sec.NumMethod_mesh}

We consider a very general formulation of the governing equations in order to model a wide class of physical phenomena, namely all those which are described by equations that can be cast into the following form, 
\be
\label{eq.generalform}
\de{\Q}{t} + \nabla \cdot \F(\Q) = \mathbf{S}(\Q), \qquad \x \in \Omega(t) \subset \mathbb{R}^2, 
\qquad \Q \in \Omega_{\Q} \subset \mathbb{R}^{\nu},
\ee
where $\x = (x,y)$ is the spatial position vector, $t$ represents the time, $\Q = (q_1,q_2, \dots, q_{\nu})$ is the vector of conserved variables defined in the space of the admissible states $\Omega_{\Q} \subset \mathbb{R}^{\nu}$, $ \F(\Q) = (\,\f(\Q), \g(\Q)\,) $ is the non linear flux tensor, and $\mathbf{S}(\Q)$ represents a non linear algebraic source term.

To discretize the moving two-dimensional domain $\Omega(t)$ we employ a centroid based Voronoi-type tessellation made of $N_P$ non overlapping polygons $P_i, i=1, \dots N_P$. The tessellation is first built at time $t=0$ and then it is regenerated at each time step~$t^n$. Data are represented via high order polynomials in each Voronoi polygon, which are either given by a (C)WENO reconstruction procedure for FV schemes, or are directly available from the numerical solution when a DG method is considered.

\subsection{Computational grid} 
At time $t^n=0$ we fix the position of $N_P$ points, called generator points: their coordinates are denoted as $\xcin, i=1,\dots, N_P$ and they are uniformly distributed inside the rectangular domain $\Omega(0)$ as well as on its boundary.
Next, we build a Delaunay triangulation having these generators $\xcn$ as vertexes of the triangles. 
The  defining  property  of  the  Delaunay  triangulation  is that the circumcircle of each triangle is not allowed to contain any of the other generator points in its interior. 
This \textit{empty circumcircle property} distinguishes the Delaunay triangulation from the many other triangulations of the
plane that are possible for the point set. 
Furthermore, this condition uniquely determines the triangulation for points in general position (except for circles with more than three generator points on them for which the Delaunay triangulation contains degenerate cases where it may flip by an infinitesimal motion of one of the points). 
For this step we follow the Delaunay algorithm presented in~\cite{bowyer1981computing, watson1981computing}, where the point location phase is efficiently performed by employing a \textit{jump-and-walk} method~\cite{mucke1999fast}.

\begin{figure}[!bp]
	\centering
	\includegraphics[width=0.33\linewidth]{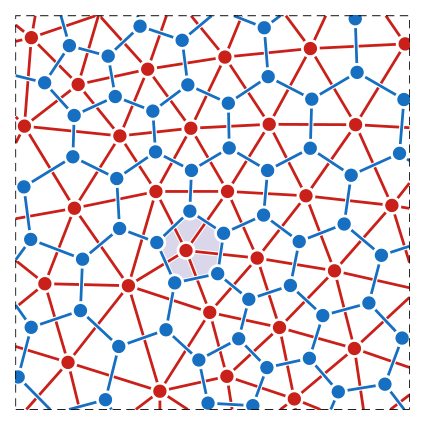}%
	\includegraphics[width=0.33\linewidth]{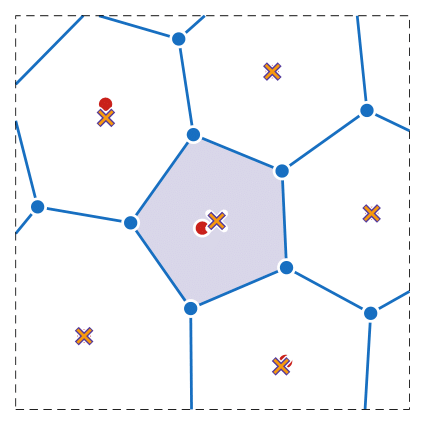}%
	\includegraphics[width=0.33\linewidth]{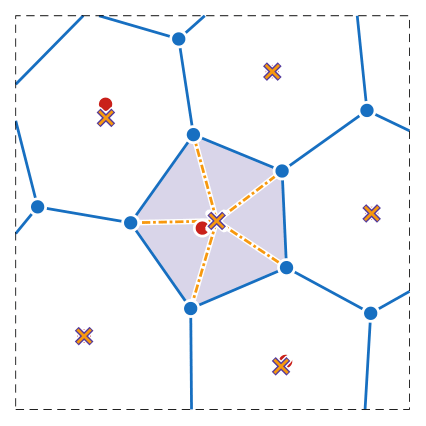}%
	\caption{In these three panels we report the Delaunay triangulation and the generator points in red. 
		The barycenters of the Delaunay triangles and the Voronoi tessellation are represented in blue. 
		Finally, the barycenters of the Voronoi polygons are represented with orange crosses.
	Note that to each generator point corresponds a Voronoi polygon which is obtained by connecting the barycenters of the triangles having this generator point as a vertex. Note also that we employ its barycenter to construct the sub--triangulation  of each Voronoi element (orange dotted line in the right panel).
	}
	\label{fig.illustrationDelvsVor}
\end{figure}

Each generator point $\xcin$ is then associated to a centroid based Voronoi element $P_i^n$ by connecting counterclockwise the \textit{barycenters} of all the Delaunay triangles having this generator point as a vertex. 
Note that the use of barycenters (instead of circumcenters) to construct these Voronoi-type elements avoids degenerate situations caused by the violation of the empty circumcircle property, thus it does not need to be resolved.
We refer to Figure~\ref{fig.illustrationDelvsVor} for a graphical interpretation (generator points are always plotted in red and Voronoi vertexes in blue). 
In particular, given a Voronoi polygon $P_i^n$ we denote 
by $\mathcal{V}(P_i^n) = \{v_{i_1}^n, \dots, v_{i_j}^n, \dots, v_{i_{N_{V_i}^n}}^n \}$ the set of its $N_{V_i}^n$ Voronoi neighbors, 
by $\mathcal{E}(P_i^n) = \{e_{i_1}^n, \dots, e_{i_j}^n, \dots, e_{i_{N_{V_i}^n}}^n \}$ the set of its $N_{V_i}^n$ edges, 
and 
by $\mathcal{D}(P_i^n) = \{d_{i_1}^n, \dots, d_{i_j}^n, \dots, d_{i_{N_{V_i}^n}}^n \}$ the set of its $N_{V_i}^n$ vertexes, consistently \textit{ordered counterclockwise}.
Finally, the barycenter of a Voronoi polygon $P_i^n$ is noted as $\xbin = (\xbixn, \xbiyn)$ (note that usually it does not coincide with the generator point, and it is always plotted in orange). By connecting $\xbin$ with each vertex of $\mathcal{D}(P_i)$ we subdivide the Voronoi polygon $P_i^n$ in $N_{V_i}^n$ subtriangles denoted as $\mathcal{T}(P_i^n) = \{T_{i_1}^n, \dots, T_{i_j}^n, \dots, T_{i_{N_{V_i}^n}}^n \}$.

\subsection{Spatial representation of the numerical solution}
The numerical solution for the conserved quantities $\Q$ in~\eqref{eq.generalform} is represented via a cell-centered approach 
inside each Voronoi polygon $P_i^n$ at the current time $t^n$ by piecewise polynomials of degree $N \geq 0$ denoted by $\mathbf{u}_h^n(\x,t^n)$ and defined in the space $\mathcal{U}_h$,
\begin{equation}
\mathbf{u}_h^n(\x,t^n) = \sum \limits_{\ell=0}^{\mathcal{N}-1} \phi_\ell(\x,t^n) \, \hat{\mathbf{u}}^{n}_{\ell,i} 
:= \phi_\ell(\x,t^n) \, \hat{\mathbf{u}}^{n}_{\ell,i} , \qquad \x \in P_i^n,
\label{eqn.uh}
\end{equation}
where $\phi_\ell(\x,t^n)$ are \textit{modal} spatial basis functions used to span the space of polynomials $\mathcal{U}_h$ up to degree $N$.
In the rest of the paper we will use classical tensor index notation based on the Einstein summation convention, which implies summation over two equal indices. The total number $\mathcal{N}$ of expansion coefficients (degrees of freedom, DOFs) $\hat{\mathbf{u}}^{n}_{l}$ for the basis functions depends on the polynomial degree $N$ and is given by $\mathcal{N} = \mathcal{L}(N,d)$, with 
\begin{equation}
\mathcal{L}(N,d) = \frac{1}{d!} \prod \limits_{m=1}^{d} (N+m),
\label{eqn.nDOF}
\end{equation}
where $d=2$ in this paper, since we are dealing only with two-dimensional domains.
As basis functions $\phi_\ell$ in~\eqref{eqn.uh} we employ a Taylor series of degree $N$ in the variables $\mathbf{x}=(x,y)$ directly defined on the \textit{physical element} $P_i^n$, expanded about its current barycenter $\xbin$ and normalized by its current characteristic length $h_i$
\begin{equation} 
\label{eq.Dubiner_phi_spatial}
\phi_\ell(\x,t^n) |_{P_i^n} = \frac{(x - \xbixn)^{p_\ell}}{p_\ell! \, h_i^{p_\ell}} \, \frac{(y - \xbiyn)^{q_\ell}}{q_\ell! \, h_i^{q_\ell}}, \qquad 
\ell = 0, \dots, \mathcal{N}-1, \quad \ 0 \leq p_\ell + q_\ell \leq N,
\end{equation} 
 $h_i$ being the radius of the circumcircle of $P_i^n$. The unknown expansion coefficients $\hat{\mathbf{u}}^{n}_{\ell,i}$ in~\eqref{eqn.uh} are the rescaled  derivatives $h_i^{p_\ell} h_i^{q_\ell}  \frac{\partial^{p_\ell+q_\ell}}{\partial x^{p_\ell} \partial y^p{_\ell}} \mathbf{Q}\left(\xbin\right)$ 
 of the Taylor expansion about $\xbin$. The time dependence of $\phi(\x,t^n)$ 
 derives from the time-dependence of the cell barycenter $\xbin$.  

The discontinuous finite element data representation~\eqref{eqn.uh} leads naturally to both a Discontinuous Galerkin (DG) scheme if $N>0$, but also to a Finite Volume (FV) scheme in the case $N=0$. This indeed means that for $N=0$ we have $\phi_\ell(\x) = 1$, with $\ell=0$ and~\eqref{eqn.uh} reduces to the classical piecewise constant data representation that is typical of finite volume schemes: 
\be
\mathbf{u}_h^n(\x,t^n) = 1 \cdot  \hat{\mathbf{u}}^{n}_{0,i} \RVIcolor{= \hat{\mathbf{u}}^{n}_{0,i}},  
\qquad \x \in P_i^n, \qquad 
\hat{\mathbf{u}}^{n}_{0,i}  =  \frac{1}{|P_i^n|} \int \limits_{P^n_i} \Q(\mathbf{x},t^n) d\x. 
\ee 
Here, the only degree of freedom per element $P_i^n$ is the usual cell average $\hat{\mathbf{u}}^{n}_{0,i}$. 
Note also that in the case $N>0$ the representation given by~\eqref{eqn.uh} already provides a spatially high order accurate data representation with accuracy $N+1$, which is not the case when $N=0$. If we are interested in increasing the spatial order of accuracy of a finite volume scheme, up to $M+1$ for example, we need to perform a spatial \textit{reconstruction} that generates a spatially high order accurate reconstruction polynomial  $\mathbf{w}_h^n(\x,t^n)$ of degree $M>N$ (see the CWENO procedure presented in~\ref{ssec.CWENO}) that reads 
\begin{equation}
\mathbf{w}_h^n(\x,t^n) = \sum \limits_{\ell=0}^{\mathcal{M}-1} \psi_\ell(\x,t^n) \, \hat{\mathbf{w}}^{n}_{\ell,i} 
:= \psi_\ell(\x,t^n) \, \hat{\mathbf{w}}^{n}_{\ell,i}, \qquad \x \in P_i^n, \qquad \mathcal{M}=\mathcal{L}(M,d),
\label{eqn.wh}
\end{equation}
where we simply employ the same basis functions $\psi_l(\x,t^n)=\phi_l(\x,t^n)$ for the reconstruction according to~\eqref{eq.Dubiner_phi_spatial}, but with $ 0 \leq \ell \leq \mathcal{M}-1$ rather than
 $ 0 \leq \ell \leq \mathcal{N}-1$, see also~\cite{Dumbser2008}. 

With this notation, our method falls within the more general class of $P_NP_M$ schemes introduced in~\cite{Dumbser2008} for fixed unstructured simplex meshes in two and three space dimensions. 
In~\cite{Dumbser2008,ADERNSE,luo1,luo2} a new family of hybrid, reconstructed  discontinuous Galerkin methods was proposed, in which a Hermite-type reconstruction of degree $M \geq N$ is performed on cell data represented by polynomials of degree $N$. 
In this paper, however, we restrict ourselves to the two most common situations:
(i) $N=0$, with arbitrary high order reconstruction of degree $M > N$, which indeed corresponds to a FV scheme of order $M+1$, and (ii) $N=M$, which corresponds to a DG scheme of accuracy $N+1$.
Within the general $P_NP_M$ formalism one can simultaneously deal with arbitrary high order FV and DG schemes inside a unified framework, with only very few differences between the two schemes. 

For the sake of uniform notation, when $N>0$ and hence $M=N$, we trivially impose that 
the reconstruction polynomial is given by the DG polynomial, i.e.  
$\mathbf{w}_h^n(\x,t^n)=\mathbf{u}_h^n(\x,t^n)$, which automatically implies
that in the case $N=M$ the reconstruction operator is simply the identity.

\subsection{CWENO reconstruction}
\label{ssec.CWENO}
For finite volume schemes ($N=0$) the reconstruction procedure allows us to compute a high order non-oscillatory polynomial representation $\w_h^n(\x,t^n)$ of the solution $\Q(\x,t^n)$ for each Voronoi polygon $P_i^n$, starting from the values of $\mathbf{u}_h^n(\x,t^n)$ in $P_i^n$ and its neighbors. It should be employed in the case $N=0, M>0$. As already stated above, the total number of unknown degrees of freedom $\w_h^n(\x,t^n)$ is $\mathcal{M} = \mathcal{L}(M,d)$, with $M$ denoting the degree of the reconstruction polynomial $\w_h$.

In order to achieve high accuracy, a large stencil centered in $P_i^n$ is required, but this choice produces oscillations close to discontinuities, the well-known Gibbs phenomenon.
Indeed, for linear reconstruction operators, the requirements of high order of accuracy and non-oscillatory behavior are
in contrast with each other, due to the well-known Godunov theorem~\cite{godunov}.  
In order to fulfill also the requirement of non-oscillatory behavior, a nonlinear reconstruction operator has to be adopted. In this paper we rely on the CWENO  reconstruction strategy first introduced in~\cite{LPR:99,LPR:00,LPR:02}, and which can be cast in the general framework described in~\cite{cravero2018cweno}. 
Here, we closely follow the work outlined in~\cite{dumbser2017central} for unstructured triangular and tetrahedral meshes. For the sake of completeness, we report here the entire algorithm: the differences with respect to~\cite{dumbser2017central} are highlighted in the last paragraph of this section.

The reconstruction starts from the computation of a so-called \textit{central polynomial} $\Popt$ of degree $M$. 
In order to define $\Popt$ in a robust manner, following~\cite{dumbser2017central,barthlsq,kaeserjcp,SCR:CWENOquadtree}, we consider a stencil
$\mathcal{S}_i^0$ which is filled with  a total number of $n_e= f \cdot \mathcal{M} = f \cdot \mathcal{L}(M,d)$ elements, containing cell $P_i^n$ and its neighbors
\begin{equation}
\mathcal{S}_i^0 = \bigcup \limits_{k=1}^{n_e} P^n_{i_k},
\label{stencil}
\end{equation}
with the safety factor $f \ge 1.5$.
Stencil $S_i^0$ includes the current Voronoi polygon $P_i^n$, the first layer of Voronoi neighbors (node neighbors of $P_i^n$) denoted by $\mathcal{V}(P_i^n)$, and is filled by recursively adding neighbors of elements that have been already selected, until the desired number $n_e$ is reached. 
The polynomial $\Popt(\x,t^n)$ is then defined by imposing that its average on each cell $P^n_{i_k}$ matches the known cell average  
$\hat{\mathbf{u}}^{n}_{0,i_k}$. 
Since $n_e>\mathcal{M}$, this of course leads to an overdetermined linear system, which is solved using a constrained least-squares technique (CLSQ)~\cite{Dumbser2007693} as
\begin{equation}
\label{CWENO:Popt}
\begin{aligned}
\Popt(\x,t^n) = \underset{{\mathbf{p}\in\mathcal{P}_i}}{\argmin}  
\sum_{P_{i_k}^n \in \mathcal{S}_i^0} 
\left(\hat{\mathbf{u}}^{n}_{0,i_k} - \frac{1}{|P_{i_k}^n|} \int_{P^n_{i_k}} \mathbf{p}(\x,t^n) d\x \right)^2, \
\text{ with} \quad 
\mathcal{P}_i = \left\{\mathbf{p} \in \mathbb{P}_M: \frac{1}{|P_{i}^n|} \int_{P^n_{i}} \mathbf{p}(\x,t^n) d\x = \hat{\mathbf{u}}^{n}_{0,i} \right\},
\end{aligned}
\end{equation}
where $\mathbb{P}_M$ is the set of all polynomials of degree at most $M$. 
In other words, the polynomial $\Popt$ has exactly the cell average $\hat{\mathbf{u}}^{n}_{0,i}$ on the polygon $P_i^n$ and matches all the other cell averages of the remaining stencil elements in the least-square sense. 
The polynomial $\Popt$ is expressed in terms of the basis functions~\eqref{eq.Dubiner_phi_spatial} of degree $M$, hence
\begin{equation}
\label{eqn.recpolydef} 
\Popt(\x,t^n) = \sum \limits_{\ell=0}^{\mathcal{M}-1} \psi_\ell(\x,t^n) \hat{\mathbf{p}}^{n}_{\ell,i} ,   
\end{equation}
and the integrals appearing in~\eqref{CWENO:Popt} are computed in each Voronoi polygon $P^n_{i_k}$ 
by summing the contribution of each of its sub-triangles $T~\in~\mathcal{T}(P^n_{i_k})$. On the sub-triangles we employ $(M+1)^2$ quadrature points defined by the conical product of the one-dimensional Gauss-Jacobi formula, see~\cite{stroud}. 

\begin{figure}[!b]
	\centering
	\includegraphics[width=0.33\linewidth]{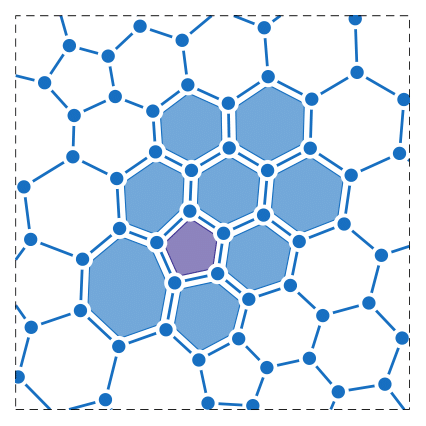}%
	\includegraphics[width=0.33\linewidth]{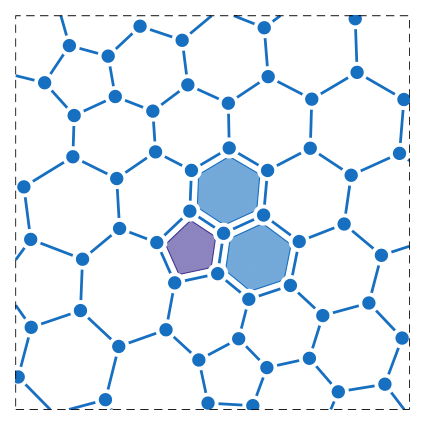}%
	\includegraphics[width=0.33\linewidth]{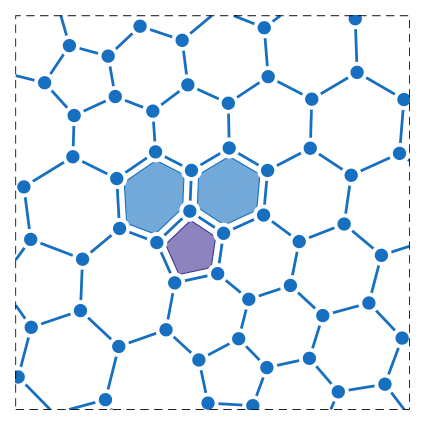}\\
	\includegraphics[width=0.33\linewidth]{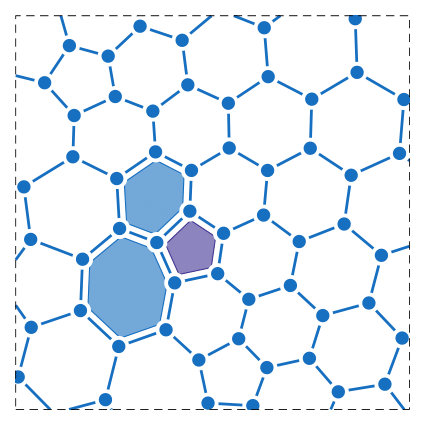}%
	\includegraphics[width=0.33\linewidth]{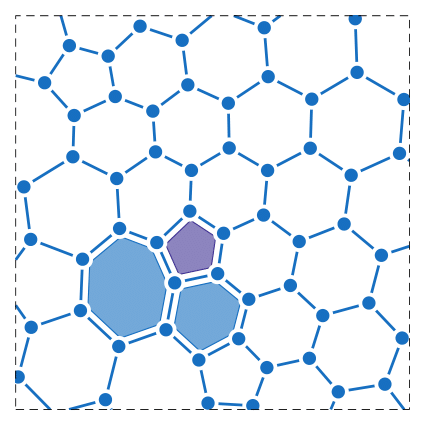}%
	\includegraphics[width=0.33\linewidth]{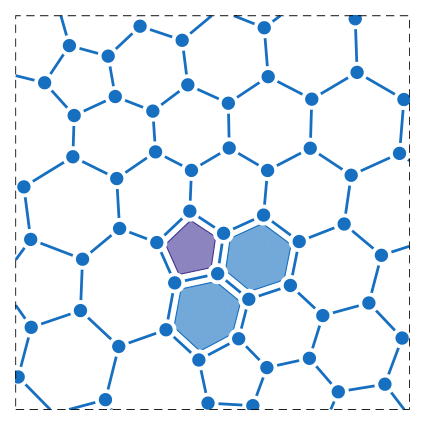}%
	\caption{Stencils for the CWENO reconstruction of order three ($M=2$) with $f=1.5$ for a pentagonal element $P_i^n$. Top-left: central stencil made of the element itself $P_i^n$ (in violet) and $n_e-1 = 8$ of its neighbors (in blue). In the other panels we report the $N_{V_i}^n =5$ sectorial stencils containing the element itself and two consecutive neighbors belonging to~$\mathcal{V}(P_i^n)$.}
	\label{fig.cweno_stencils}
\end{figure}

To make the reconstruction operator nonlinear, which is required in the presence of shock waves, the CWENO algorithm makes use of \textit{other polynomials} of lower degree. 
Given a Voronoi polygon $P_i^n$ with $N_{V_i}^n$ Voronoi neighbors $\mathcal{V}(P_i^n)$, 
we construct $N_{V_i}^n$ interpolating polynomials of degree $M^s=1$ referred to as \textit{sectorial polynomials}. 
More precisely, we consider $N_{V_i}^n$ stencils $S_i^s$ with $s\in[1,N_{V_i}^n]$, each of them containing exactly $\hat{n}_e=\mathcal{L}(M^s,d) = (d+1)$ cells. 
Each $S_i^s$ includes always the central cell $P^n_i$ and two consecutive neighbors belonging to~$\mathcal{V}(P_i^n)$.
An example of stencils $S_i^0$ and $S_i^s$ for a polygon with $N_{V_i}^n=5$ and $M=2$ is reported in Figure~\ref{fig.cweno_stencils}.

For each stencil $S_i^s$ we compute a linear polynomial $\mathbf{P}_s(\x,t^n)$ by solving the reconstruction systems 
\begin{equation}
\label{CWENO:Ps}
\mathbf{P}_s(\x,t^n) \in\mathbb{P}_1 
\text{  s.t. } \forall P^n_{i_k} \in S_i^s: 
\frac{1}{|P_{i_k}^n|} \int_{P^n_{i_k}} \mathbf{P}_s(\x,t^n) \, d\x = \hat{\mathbf{u}}^{n}_{0,i_k},
\end{equation}
which are not overdetermined and therefore have a unique solution for non-degenerate locations of the Voronoi barycenters. 
Following the general framework introduced in~\cite{cravero2018cweno}, we select a set of positive coefficients $\lambda_0,\ldots,\lambda_{N_p}$ such that
\begin{equation}
 \sum_{s=0}^{N_{V_i}^n}\lambda_s=1
 \label{eqn.sumCWENO}
\end{equation}
and we define a new polynomial 
\begin{equation}
\label{CWENO:P0}
\mathbf{P}_0(\x,t^n) = \frac{1}{\lambda_0}\left(\Popt(\x,t^n) - \sum_{s=1}^{N_p} \lambda_s \mathbf{P}_s(\x,t^n) \right) \in\mathbb{P}_M,
\end{equation}
so that the linear combination of the polynomials $\mathbf{P}_0,\ldots,\mathbf{P}_{N_{V_i}^n}$ with the coefficients $\lambda_0,\ldots,\lambda_{N_p}$ is equal to $\Popt$ and conservation is ensured. Specifically, we consider the linear weights used in~\cite{ADERFV}, namely $\lambda_0=10^5$ for $\mathcal{S}_i^0$ and $\lambda_s=1$ for the sectorial stencils. These weights are later normalized in order to sum to unity, according to the requirement~\eqref{eqn.sumCWENO}.
Finally, the sectorial polynomials $\mathbf{P}_s$ with $s\in[1,N_{V_i}^n]$ are nonlinearly hybridized with $\mathbf{P}_0$, as it is done also in other WENO schemes ~\cite{shu_efficient_weno,hu1999weighted,BalsaraAO}. We thus obtain $\w_{h}(\x,t^n)$ in $P_i^n$ as
\begin{equation}
\label{CWENO:Prec}
\w_{h}(\x,t^n) = \sum_{s=0}^{N_p} \omega_s \mathbf{P}_s(\x,t^n), \qquad \x \in P_i^n, 
\end{equation}
where the normalized {\em nonlinear weights} $\boldsymbol{\omega}_s$ are given by
\begin{equation}
\label{eqn.weights}
\boldsymbol{\omega}_s = \frac{\tilde{\boldsymbol{\omega}}_s}{ {\left({\sum \limits_{m=0}^{N_{V_i}^n} \tilde{\boldsymbol{\omega}}_m}\right)} }, 
\qquad \textnormal{ with } \qquad 
\tilde{\boldsymbol{\omega}}_s = \frac{\lambda_s}{\left(\boldsymbol{\sigma}_s + \epsilon \right)^r}. 
\end{equation} 
In the above expression the non-normalized weights $\tilde{\boldsymbol{\omega}}_s$ depend on the linear weights $\lambda_s$ and the oscillation indicators $\boldsymbol{\sigma}_s$ 
with the parameters $\epsilon=10^{-14}$ and $r=4$ chosen according to~\cite{Dumbser2007693}.
Note that in smooth areas, $\omega_s\simeq\lambda_s$ and then $\w_{h_i}\simeq\Popt$, so that we recover optimal accuracy. On the other hand, close to a discontinuity, $\mathbf{P}_0$ and some of the low degree polynomials $\mathbf{P}_s$ would be oscillatory and have high oscillation indicators, leading to $\omega_s\simeq0$ and in these cases only lower order non-oscillatory data are employed in $\w_{h_i}$, guaranteeing the non-oscillatory property of the reconstruction. The oscillation indicators $\boldsymbol{\sigma}_s$ appearing in~\eqref{eqn.weights} are simply 
given by 
\begin{equation}
\boldsymbol{\sigma}_s = \sum_{l} \left(\hat{\mathbf{p}}^{n,s}_{l,i}\right)^2.
\label{eqn.OI}
\end{equation}  

The CWENO procedure adopted in this work is similar to the one presented in~\cite{dumbser2017central} and it has been adapted to Voronoi polygons and their connectivity. The needed modifications concern the computation of integrals in~\eqref{CWENO:Popt}, the number of \textit{sectorial polynomials}, and the fact that basis functions are rescaled Taylor monomials referred to the physical element and not to the reference element, hence yielding a different and very simple evaluation of the oscillation indicators~\eqref{eqn.OI}. 

\subsection{Evolution of the computational domain}
\label{sec.MeshEvolution}

At this point we have a high order spatial representation of the solution  $\Q(\x, t^n)$ at the current time $t^n$ given by the polynomial $\w_h^n=\w_h(\x,t^n)$ of degree $M$. We recall that if $N=M>0$ then $\w_h^n = \u_h^n$; if instead $N=0$ then $\w_h^n$ is obtained through the reconstruction procedure described in the previous Section~\ref{ssec.CWENO}.

By evaluating $\w_h^n$ at the generator points $\xcn$, i.e. $\w_h^n(\xcn, t^n)$ with~\eqref{eqn.wh}, we recover the local fluid velocity $\mbf{v}(\xcn)$, that can be used to compute the new coordinates of the generator points simply as
\be
\xcinp = \xcin + \Delta t \, \mbf{v}(\xcin).
\label{eqn.xcnew} 
\ee
Note that in our ALE formalism, the mesh can be moved with \textit{any} velocity, hence it is not necessary to always integrate the above relation~\eqref{eqn.xcnew} with high order of accuracy. 
Moreover, for the sake of simplicity, all along this manuscript we do not move boundary elements.

The Delaunay triangulation connecting the new coordinates of the generator points $\xcnp$ is now recomputed, as well as the  corresponding updated Voronoi tessellation.
Note that the \textit{only connection} between the tessellations at time $t^n$ and $t^{n+1}$ is the number $N_P$ of generator points (i.e. of Voronoi polygons)  and their \textit{global numbering}. Instead, the shape of each polygon is allowed to change, i.e. $N_{V_i}^n \ne N_{V_i}^{n+1}$, and consequently also the connectivities, i.e. for example $\mathcal{V}(P_i^n) \ne \mathcal{V}(P_i^{n+1})$.

This change of the grid topology is actually the \textit{strength} of the present algorithm, since it allows \RIIIcolor{us} to maintain a high mesh quality without distorted elements, as depicted in Figures~\ref{fig.ShuDG4} and~\ref{fig.ShuDG4_standardALE}, where we show a comparison between the results obtained by allowing topology changes and by imposing a fixed connectivity, respectively. However, more care is needed in order to update the solution from time $t^n$ to $t^{n+1}$. 
In particular, to obtain a \textit{high order direct} ALE scheme we need a complete knowledge of the space--time structure between the two time levels, i.e. we need to construct the so called space--time control volumes and their space--time connectivity.
We would like to \RIIIcolor{emphasize} that up to Finite Volume schemes of order $2$, one could avoid the procedure that we are going to introduce (see~\cite{Springel, pakmor2015improving}), but starting from order $3$ it is essential. 

\RShashkovcolor{
\subsubsection{High order integration of the trajectories of the generator points}
\label{ssec.HighOrderTraj}

Due to the ALE framework of the present work, the mesh can in principle be moved with an \textit{arbitrary} velocity, and there is not a specific necessity of moving the grid in a fully Lagrangian fashion.
Nevertheless,~\eqref{eqn.xcnew} can also be replaced by a high order Taylor method~\cite{boscheri2013semi, boscheri2019high, tavellihigh}, leading to a high order approximation of the Lagrangian trajectories of the generators points. 
The use of this technique for example improves mesh quality in vortex flow, as clearly shown 
in Section~\ref{test.ShuVortex}, and also improves the overall Lagrangian behavior of the algorithm. 

In what follows we detail the high order approach used for the integration of the flow trajectories.

The Taylor expansion of the new position $\xcinp$ of a generator point at time $t^{n+1}$ with respect to its position at time $t^n$ can be written as
\be
\xcinp = \xcin + \Delta t \frac{d \x}{dt}  + \frac{\Delta t^2}{2} \frac{d^2\x}{dt^2} 
+ \frac{\Delta t^3}{6} \frac{d^3\x}{dt^3}  + \frac{\Delta t^4}{24} \frac{d^4\x}{dt^4} + \mathcal{O}(5),
\label{eqn.xcnew_high} 
\ee
which achieves \textit{fourth} order of accuracy in time.
Now, the high order time derivatives in~\eqref{eqn.xcnew_high} are replaced by high order spatial derivatives, via the \textit{Cauchy-Kovalevskaya} procedure, using repeatedly the trajectory equation 
\be
\frac{d \x}{dt} = \mbf{v} (\mbf{x}(t)),
\ee 
and assuming a stationary velocity field (i.e. $\de{\mbf{v}}{t} = \0$), hence
\be \label{eq:trajectory12}
\dfrac{d\mbf{x}}{dt}     = \mbf{v} = v_i,\qquad
\dfrac{d^2\mbf{x}}{dt^2} = \frac{d}{dt}\left(\dfrac{d\mbf{x}}{dt}\right) = \de{\mbf{v}}{\mbf{x}}\,\de{\mbf{x}}{t} = (\nabla\mbf{v})\,\mbf{v} = \frac{\partial v_i}{\partial x_j}\,{v_j}.
\ee
The chain rule, as written in~\eqref{eq:trajectory12}, can be applied iteratively to obtain the third derivative of the position
\be
\dfrac{d^3\mbf{x}}{dt^3} = \frac{d}{dt}\left(\dfrac{d^2\mbf{x}}{dt^2}\right) = 
\nabla(\nabla\mbf{v})\,\mbf{v}\,\mbf{v} + (\nabla\mbf{v})\,(\nabla\mbf{v})\,\mbf{v} = 
    \frac{\partial v_i}{\partial{x_j}\,\partial{x_k}}\,v_j\,v_k + \frac{\partial{v_i}}{\partial{x_j}}\,\frac{\partial{v_j}}{\partial{x_k}}\,v_k,
\ee
and similarly, the fourth derivative reads
\be
\dfrac{d^4\mbf{x}}{dt^4} &= 
\nabla\left(\nabla\left(\nabla\mbf{v}\right)\right)\,\mbf{v}\,\mbf{v}\,\mbf{v} + 
\nabla\left(\nabla\mbf{v}\right)\,\left(\nabla\mbf{v}\right)\,\mbf{v}\,\mbf{v} + 
2\,\nabla\left(\nabla\mbf{v}\right)\,\mbf{v}\,\left(\nabla\mbf{v}\right)\,\mbf{v} + 
\left(\nabla\mbf{v}\right)\,\nabla\left(\nabla\mbf{v}\right)\,\mbf{v}\,\mbf{v} + 
\left(\nabla\mbf{v}\right)\,\left(\nabla\mbf{v}\right)\,\left(\nabla\mbf{v}\right)\,\mbf{v} = \\[1pt]
& = \frac{\partial v_i}{\partial x_j\,\partial x_k\, \partial x_l}\,v_j\,v_k\,v_l + 
\frac{\partial v_i}{\partial x_j\,\partial x_k}\,\frac{\partial v_k}{\partial x_l}\,v_l\,v_j + 
2\,\frac{\partial v_i}{\partial x_j\,\partial x_k}\,v_k\,\frac{\partial{v_j}}{\partial x_l}\,v_l + 
\frac{\partial v_i}{\partial x_j}\,\frac{\partial v_j}{\partial x_k\,\partial x_l}\,v_l\,v_k + 
\frac{\partial v_i}{\partial x_j}\,\frac{\partial v_j}{\partial x_k}\,\frac{\partial v_k}{\partial x_l}\,v_l.
\ee
Finally, the partial derivatives of $\v$ are recovered from the local fluid velocities $\u$ 
through the high order polynomials $\w_h$~\eqref{eqn.wh} which represent the conserved variables $\mathbf{Q}$ inside each cell with high order of accuracy.
Since  $\w_h$ is given via \textit{modal} basis functions, the coefficients $\hat{\mathbf{w}}^{n}_{\ell,i}$ 
already represent the values of the partial derivatives with respect to $\x$ of the conserved variables, if a 
sufficiently high order accurate $P_NP_M$ method is employed.
Then, the chain rule should be applied in order to recover the partial derivatives of the primitive variable $\u$
from those of the conserved variables $\rho\u$ and $\rho$.

}

\RVcolor{
\subsection{Mesh optimization}
\label{ssec.smoothing}
The ALE framework also allows to apply some mesh optimization techniques, since the mesh velocity is not constrained to follow 
the local fluid velocity \textit{exactly}. Furthermore, an additional level of liberty in the choice of
type of smoothing scheme is introduced by the possibility of changing the grid connectivity between consecutive time levels.
In this work, the mesh optimization methods are implemented by 
slightly modifying, at each time step, the motion of the Voronoi generator points 
(that is, the vertexes of the dual Delaunay triangulation).
This aims at improving the overall robustness of the method, as well as reducing numerical errors
and spurious mesh imprinting. 

In general, the target polygonal elements will have a locally uniform edge length ({i.e.}
no small edges for a given element, smoothly graded mesh size) and will not be excessively stretched
in one direction only, so that anisotropies due to differences in numerical diffusion are minimized. More importantly,  
this increases the robustness of the matrix computations involved in the polynomial data 
reconstruction and in the fully discrete update formulae. Also, we must note that these objectives
shall be pursued in conjunction with the interest of preserving an accurate mesh motion that follows the 
local flow field, maintaining the Lagrangian character of the numerical method as far as possible. This means that one must 
allow a certain degree of anisotropy in the mesh, which might be desirable
to resolve flow discontinuities or strong gradients.

The mesh regularization procedure begins by computing 
all the new positions $\xcinp$ for the generator points of the Voronoi grid and then recovering,
for each generator, the position $\xcistar$ that is prescribed by a simple smoothing technique applied to the
candidate positions $\xcinp$. 
We say that $\xcistar$ is a location for the generator that is optimal 
in the sense of mesh quality, as opposed to optimal
in following the flow of the fluid, which would be the role taken by $\xcinp$.
The candidate position $\xcinp$ is subsequently replaced by a corrected value $\xcinpstar$ that is given by the weighted average 
$\xcinpstar = (1 - \mu)\,\xcinp + \mu\,\xcistar$, with $\mu$ being a blending factor that yields the balance between 
the amount of mesh motion due to fluid flow with the one due to smoothing.

Concerning the determination of $\xcistar$, we decided to 
simply compute it from the application to $\xcinp$ of one iteration of a Lloyd-type algorithm; 
that is, after updating the Delaunay triangulation
of the generator points taking into account their new candidate positions $\xcinp$, 
we evaluate the quality-optimal 
position $\xcistar$ for generator $\mbf{c}_i$ as
\begin{equation}
     \xcistar = {\left(\displaystyle\sum_{a_k\,\in\,\mathcal{A}_i}\omega_{a_k}\right)^{-1}}\,
     {\displaystyle\sum_{a_k\, \in\, \mathcal{A}_{i}}\frac{1}{2}\left(\xcijnp + 
     \xciknp\right)\,\omega_{a_k}}.
\end{equation}
We define $\mathcal{A}_i$ to indicate the set of Delaunay triangles $a_k$ that share $\xcinp$ as a vertex, while 
we denote with $\xcijnp$ and $\xciknp$ the other two vertexes of the triangle $a_k$, that is, the two
that do not coincide with $\xcinp$.
The choice of weights yields different smoothing methods, 
and in this work we mainly employ $\omega_{a_k} = \norm{\xcijnp - \xciknp}$
to obtain an algorithm that is reminiscent of Lloyd smoothing \cite{lloyd1982least}, as this would prescribe that each generator shall
be moved to the centroid of the polygonal chain obtained by connecting all vertexes $\xcijnp$ and $\xciknp$ of all 
Delaunay triangles in $\mathcal{A}_i$. This choice tends to eliminate small edges just like the algorithms forwarded in \cite{bubblesmoothing, pliantsmoothing, distmesh}.
Alternatively, we can set $\omega_{a_k} = 1$ and obtain Laplacian smoothing \cite{laplaciansmoothing,laplaciansmoothingbis}, that is, 
the generic generator $\xcinp$ is moved to the center of mass of the system of point masses defined by the vertexes 
of the above described polygonal chain. Laplacian smoothing yields nicely rounded cells and tends to preserve the grading of the mesh.

Once a set of quality-aware node positions $\xcistar$ has been determined, the algorithm must 
choose how to compromise between such positions and those prescribed by the fluid motion.
Instead of simply fixing the value of $\mu$ as a simulation parameter, we chose to let $\mu$ vary 
with time by recomputing it as a function of the solution data and of the current grid configuration, 
as well as by accounting for the specific explicit time step restriction in use.
Specifically, we compute
the relaxation parameter $\mu$ as
\begin{equation} \label{eq:relaxationparameter}
    \mu = \min\left(1,\ \sqrt{\frac{{U}_\ast\,\Delta t}{{\Delta s}}\,\mathcal{F}}\right), 
\end{equation}
with ${U}_\ast$ being a rough scaling estimate for the fluid velocity, computed at each 
time step as the maximum velocity encountered for all generator points, $\Delta t$ the time step size, 
and $\Delta s$ an indicator for the mesh spacing, 
given by the minimum value of 
the ratio between the area and the perimeter of all Voronoi polygons, in analogy to how the time step
duration is determined in \eqref{eq:timestep}. The underlying idea is that we want to balance, 
during each time step, the spatial scaling 
of fluid flow, with a characteristic length representative of the mesh motion due to pure
smoothing in the smallest cells of the domain, which we implicitly assume to be the most delicate.
In this way, we have replaced the blending factor 
$\mu$ with another non-dimensional \textit{smoothing parameter} $\mathcal{F}$, that fixes the
strength of smoothing in the small cells that are those that might otherwise compromise the stability of
the computation. The square root is arbitrarily introduced in order to reduce the sensitivity 
of Equation~\eqref{eq:relaxationparameter} to sudden variations in the flow speed $U_\ast$.

Note that, although in a very approximate form, the formula~\eqref{eq:relaxationparameter} 
scales with the square root of a characteristic Mach number, at least when $U_\ast$ is 
negligible with respect to sound speed or vice versa; 
one can verify this by substituting \eqref{eq:timestep} in \eqref{eq:relaxationparameter}
and noting that it simplifies into an expression
that includes the degree of the polynomial data $N$, the CFL coefficient, 
and an approximate Mach number. 
Further investigations on more complex scaling expressions that correct 
for such residual dependencies, as well as space-dependent formulae, 
are left for future work.

The results presented in this paper are obtained by moving the generators with 
the local fluid velocity and by applying one of the smoothing techniques described here,
with different values of the smoothing parameter $\mathcal{F}$.
}

\subsection{Space--time connectivity}
\label{ssec.SpaceTimeConnection}

\begin{figure}[!bp]
	\centering{
		\subfloat[]{\includegraphics[width=0.33\linewidth]{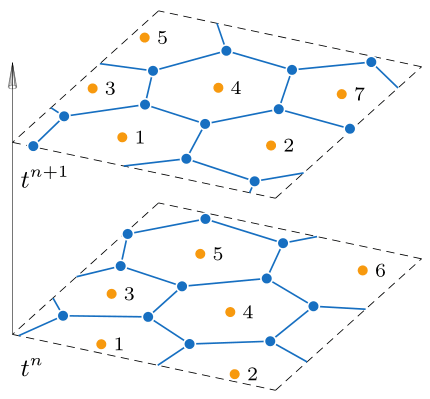}}%
		\subfloat[]{\includegraphics[width=0.33\linewidth]{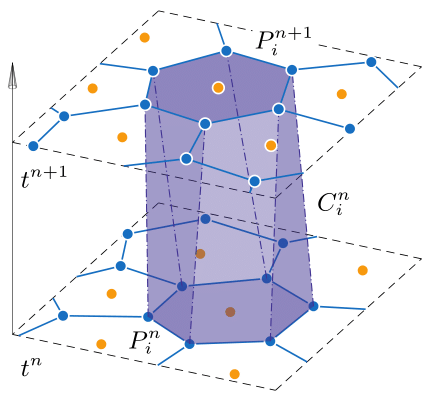}}%
		\subfloat[]{\includegraphics[width=0.33\linewidth]{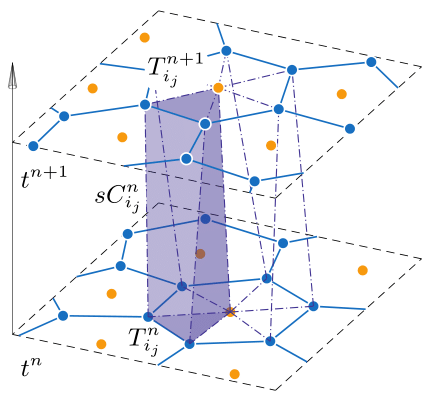}}%
	}
	\caption{Space time connectivity \textit{without} topology changes. (a) The tessellation at time $t^n$ and time $t^{n+1}$. (b) $P_i^n$ is connected with $P_i^{n+1}$ to construct the space--time control volume $C_i^n$. (c) The sub-triangle $T_{i_j}^n$ is connected with $T_{i_j}^{n+1}$ to construct the sub--space--time control volume $sC_{i_j}^n$. }
	\label{fig.ControlVolume_noTopChange}
\end{figure}

For the sake of clarity, let us first consider the simple case in which no topology changes have occurred between $t^n$ and $t^{n+1}$, i.e. $N_{V_i}^n = N_{V_i}^{n+1}$ and $\mathcal{V}(P_i^n) = \mathcal{V}(P_i^{n+1})$, as illustrated in Figure~\ref{fig.ControlVolume_noTopChange}.
Here, the space--time control volume $C_i^n$ is easily obtained by connecting each node of the polygon $P_i^n$ via \textit{straight} line segments with the \textit{corresponding} node of $P_i^{n+1}$. 
Moreover, each sub--triangle $T_{i_j}^n \in \mathcal{T}(P_i^n)$ is connected with the corresponding $T_{i_j}^{n+1} \in \mathcal{T}(P_i^{n+1})$ obtaining a \textit{sub}--space--time control volume, denoted by $sC_{i_j}^n$ in the following, which has the form of an oblique prism in space--time, with triangular faces on the bottom ($t^n$) and the top ($t^{n+1}$). 

We underline that each space--time element $C_i^n$ is given by a volume that is closed by the polygon $P_i^n$ at time $t^n$, the polygon $P_i^{n+1}$ at $t^{n+1}$ and by 
the lateral space-time faces $\partial C_{i_j}^n\  j=1,\dots,N_{V_i}^{n,st}$ which are quadrilaterals in space--time and represent the time evolution of the edges $e_{i_j}^n \in \mathcal{E}(P_i^n)$. Here, $N_{V_i}^{n,st} = N_{V_i}^n = N_{V_i}^{n+1}$ denotes the number of space--time neighbors of $C_i^n$.
The total surface of $C_i^n$ is denoted with $\partial C_i^n$ 
\be 
\partial C_i^n =  \bigcup_{j=1}^{N_{V_i}^{n,st}} \partial C_{i_j}^n \cup P_i^n \cup P_i^{n+1}.
\label{eqn.dC}
\ee

\paragraph{Technical details 1}
We recall that the node numbering (i.e. the numbering of the blue points in Figure~\ref{fig.ControlVolume_noTopChange}) could be in principle different at the two time levels so the \textit{correspondence} between the nodes at time level $t^n$ and $t^{n+1}$ is not obvious. Nevertheless, it can be recovered from the numbering of the Voronoi neighbors $\mathcal{V}(P_i^{n/n+1})$ that on the contrary remains the same.
Therefore, we loop over $\mathcal{V}(P_i^{n/n+1})$, we find the edges $e_{i_j}^{n/n+1}$ shared between $\mathcal{V}(P_{i_j}^{n/n+1})$ and $P_i^{n/n+1}$, and we put in correspondence their end points, so that the space--time control volume $C_i^n$ can be defined.
Besides, the surface obtained by connecting the end points of $e_{i_j}^{n}$ and $e_{i_j}^{n+1}$ is noted as $\partial C_{i_j}^n$, see Figure~\ref{sfig.quad_surface}. {\scriptsize{$\blacksquare$}}

\medskip

\begin{figure}[!p]
	\centering
	\subfloat[]{\includegraphics[width=0.5\linewidth]{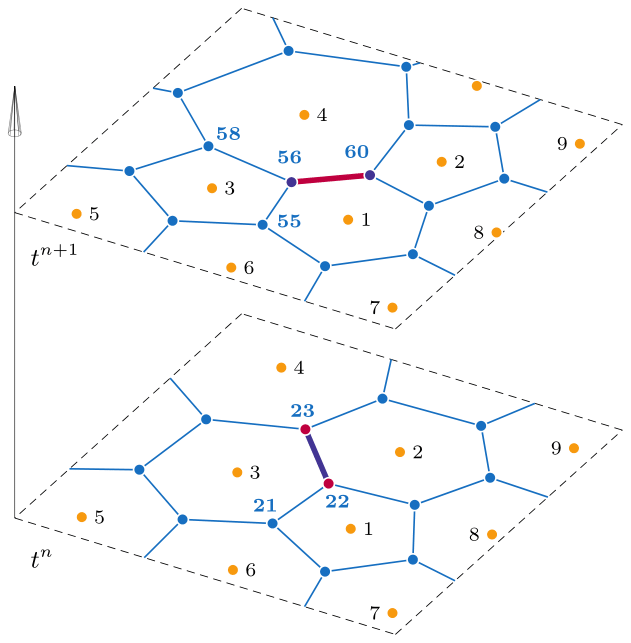}\label{sfig.Crazy_a}}%
	\subfloat[]{\includegraphics[width=0.5\linewidth]{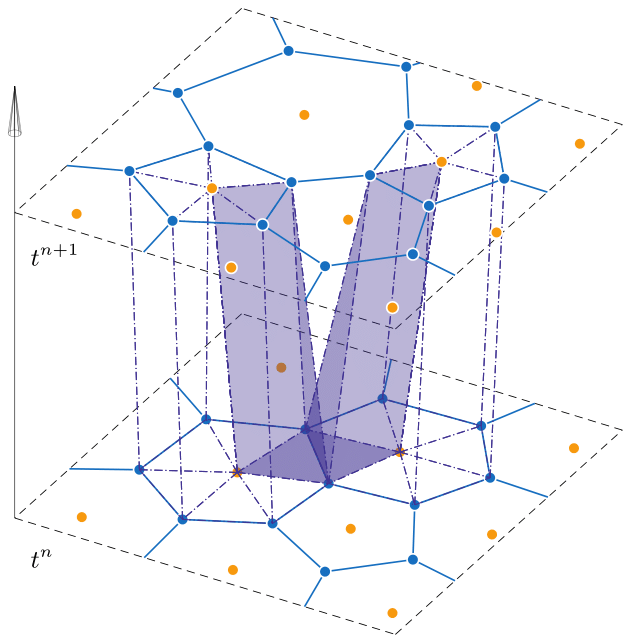}\label{sfig.Crazy_b}}\\%
	\subfloat[]{\includegraphics[width=0.5\linewidth]{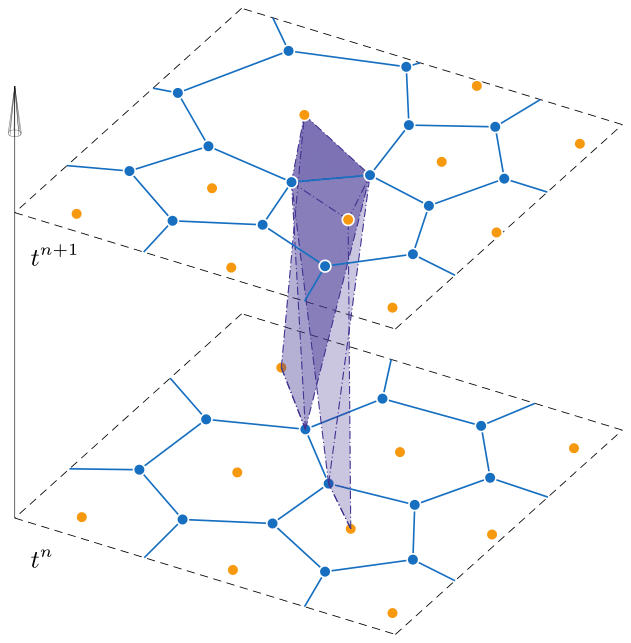}\label{sfig.Crazy_c}}%
	\subfloat[]{\includegraphics[width=0.5\linewidth]{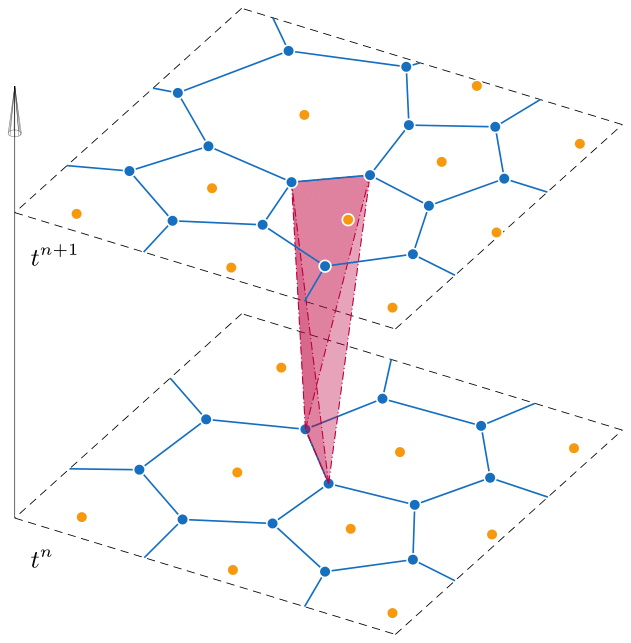}\label{sfig.Crazy_d}}%
	\caption{Space time connectivity \textit{with} topology changes, degenerate sub--space--time control volumes and sliver element. Panel (a): at time $t^n$ the polygons $P_2^n$ and $P_3^n$ are neighbors and share the highlighted edge, instead at time $t^{n+1}$ they do not touch each other; the opposite situation occurs for polygons $P_1^n$ and $P_4^n$. This change of topology causes the appearance of degenerate elements of different types. The first type is given by degenerate sub--space--time control volumes colored in violet in Panels (b) and (c). The second type of degenerate elements are called \textit{space--time sliver elements}, an example is colored in magenta in Panel (d).
		The sub--space--time control volumes of Panels (b) and (c) are triangular prisms with one of their faces collapsed to just a line: they do not pose particular problems because they are part of a standard control volume, so everything is naturally well defined on them (basis functions, quadrature points, values of $\u_h^n, \w_h^n, \q_h^n$). 
		On the contrary, the sliver element in panel (d) is a completely new control volume which does neither exist at time $t^n$, nor at time $t^{n+1}$, since it coincides with an edge of the tessellation and, as such, has zero areas in space. However, it has a \textit{non-negligible volume} in space--time. The difficulties associated to this kind of element are due to the fact that $\w_h$ is not clearly defined for it at time $t^n$ and that contributions across it should not be lost at time $t^{n+1}$ in order to guarantee conservation.  } 
	\label{fig.ControlVolume_TopChange_1sliver}
\end{figure}

Let us now consider $P_i^n$ and $P_i^{n+1}$ in the case $N_{V_i}^n \ne N_{V_i}^{n+1}$. Now, the space--time connection between them induces the appearance of degenerate elements of two types: (i) degenerate \textit{sub}--space--time control volumes $sC_{i_j}^n$, where either their top or bottom faces are degenerate triangles that are collapsed just to a line, see Figures~\ref{sfig.Crazy_b}-\ref{sfig.Crazy_c}; (ii) and also sliver space--time elements, see Figure~\ref{sfig.Crazy_d}.
Technical details on their construction (intended for the reader interested in reproducing the algorithm) are reported in the following paragraph. 
The main characteristics of this kind of elements are described in next Section~\ref{ssec.Sliver}. 

\paragraph{Technical details 2}
First, we order $\mathcal{V}(P_i^n)$ and $\mathcal{V}(P_i^{n+1})$ starting from the \textit{first} common neighbor (evidence that this choice does \textit{not} affect the results are shown in Table~\ref{tab.IndipendenceOfNeighsNumb}).
Then, we merge the two set of neighbors to compute $\mathcal{V}(C_i^n)$ which, in this case, does not coincide neither with $\mathcal{V}(P_i^n)$ nor with $\mathcal{V}(P_i^{n+1})$. $\mathcal{V}(C_i^n)$ contains all the polygons of $\mathcal{V}(P_i^n)$ and $\mathcal{V}(P_i^{n+1})$ counted once (i.e. without multiple entries) and counterclockwise ordered respecting the order of both $\mathcal{V}(P_i^n)$ and $\mathcal{V}(P_i^{n+1})$. It represents the set of $N_{V_i}^{n,st}$ space--time neighbors of $C_i^n$.

Next, we have to find the node connections in order to build $C_i^n$, which are not obviously determined and are recovered from $\mathcal{V}(C_i^n)$. We loop on $P_{i_j}\in \mathcal{V}(C_i^n)$ \RIcolor{(this loop assures that we account for \textit{all} the nodes of $C_i^n$, since by considering \textit{all} the neighbors we also consider \textit{all} the edges of both $P_i^{n}$ and $P_i^{n+1}$)} and we proceed as follows:
\begin{itemize}
	\item[I.]
	If $P_{i_j}$ belongs both to $\mathcal{V}(P_i^n)$ and to $\mathcal{V}(P_i^{n+1})$, the node connection procedure falls into the previous one, and a standard $sC_{i_j}^n$ and $\partial C_{i_j}^n$ can be recovered by connecting the end points of the edges shared between $P_i^n - P_{i_j}^n$ and $P_i^{n+1} - P_{i_j}^{n+1}$.
	Referring to $P_3^n$ depicted in Figure~\ref{fig.ControlVolume_TopChange_1sliver}, we could fix as first common neighbor $P_1^n$ because $P_1^n \in \mathcal{V}(P_3^n)$ and $P_1^n \in \mathcal{V}(P_3^{n+1})$: nodes $21-55$ and $22-56$ can be easily connected.
	\item[II.]  
	If $P_{i_j} \in \mathcal{V}(P_i^n)$ but $P_{i_j} \notin \mathcal{V}(P_i^{n+1})$, then the end points of the edge shared between $P_i^n - P_{i_j}^n$ will be connected to a unique node at time $t^{n+1}$, namely the \textit{top} node which is \textit{common} to $P_{i_{j-1}}$ and $P_{i_{j+1}}$ at time $t^{n+1}$.
	Referring to Figure~\ref{fig.ControlVolume_TopChange_1sliver}, both nodes $22$ and $23$ will be connected with node $56$.
	In this case , $\partial C_{i_j}^n$ is degenerate: it does not have a rectangular shape but a triangular one. Also $sC_{i_j}^n$ is degenerate because its top face is just given by a line connecting the barycenter of $P_i^{n+1}$ with the \textit{common top} node (node $56$ in Figure~\ref{fig.ControlVolume_TopChange_1sliver}).
	\item[III.] 
	If $P_{i_j} \in \mathcal{V}(P_i^{n+1})$ but $P_{i_j} \notin \mathcal{V}(P_i^{n})$, then the end points of the edge shared between $P_i^{n+1} - P_{i_j}^{n+1}$ will be connected to a unique node \RIIIcolor{at} time $t^{n}$, namely the \textit{bottom} node which is \textit{common} to $P_{i_{j-1}}$ and $P_{i_{j+1}}$ at time $t^n$.	
	Referring to $P_4^n$ shown in Figure~\ref{fig.ControlVolume_TopChange_1sliver}, both nodes $56$ and $60$ will be connected with node $23$.
	As in the previous case, $\partial C_{i_j}^n$ has a degenerate triangular shape and also $sC_{i_j}^n$ is degenerate because its bottom face is just given by a line connecting the barycenter of $P_i^{n}$ with the \textit{common bottom} node (node $23$ in Figure~\ref{fig.ControlVolume_TopChange_1sliver}). {\scriptsize{$\blacksquare$}}
\end{itemize}

\medskip

Note that when a change of topology occurs in a Voronoi polygon, the same happens to three of its neighbors and a total of four degenerate \textit{sub}--space--time control volumes will be originated, two of type (II) and two of type (III), refer to Figures~\ref{sfig.Crazy_b}-\ref{sfig.Crazy_c}.
Moreover, a void is left between them: to fill it and recover a fully conservative discretization, we insert a new element called \textit{space--time sliver element}, depicted in Figure~\ref{sfig.Crazy_d}, whose bottom and top faces just coincide with an edge of the tessellation at time $t^{n}$ and $t^{n+1}$, respectively. 
We denote this kind of element with $S_i^n$, its total lateral surface with $\partial S_i^n$ and each of the four lateral faces with $\partial S_{i_j}^n, j=1, \dots, 4$.

\paragraph{Technical details 3} The nodes of a \textit{sliver element} are given by the end points of those edges that \textit{flip} between the two time steps and are ordered in such a way that the volume of $S_i^n$ is positive.
Let us consider case~(II) in which $P_{i_j} \in \mathcal{V}(P_i^n)$ but $P_{i_j} \notin \mathcal{V}(P_i^{n+1})$: the edge between $P_i^n - P_{i_j}^n$ is taken as bottom face for the sliver.
Then, we loop over the edges outgoing from the \textit{common top} node: two of them belong to $P_i^{n+1}$, the third one will be taken as top face of the sliver element. 
If that edge connects  $P_i^{n+1} \rightarrow P_{i_j}^{n+1}$ then \textit{one} sliver element is enough to fill the space--time hole left from the topology change.

\begin{figure}[!bp]
	\centering
	\subfloat[]{\includegraphics[width=0.5\linewidth]{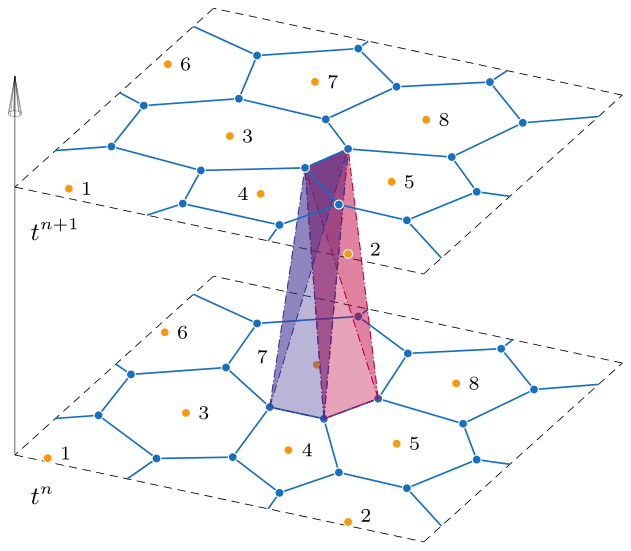}\label{sfig.consSliver_a}}%
	\subfloat[]{\includegraphics[width=0.5\linewidth]{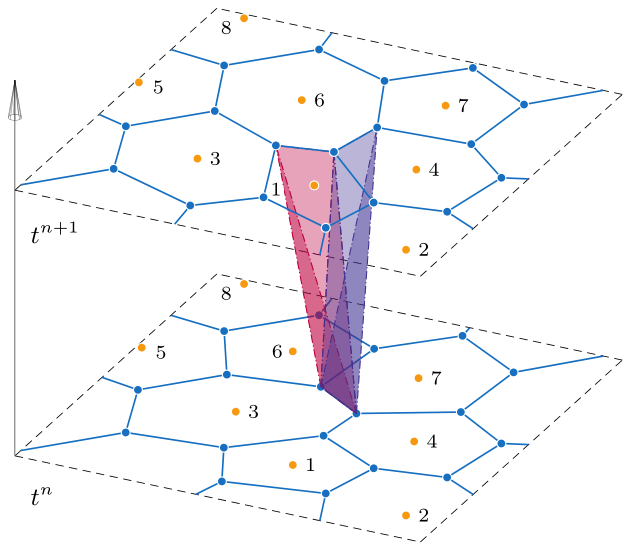}\label{sfig.consSliver_b}}\\%
	\subfloat[]{\includegraphics[width=0.5\linewidth]{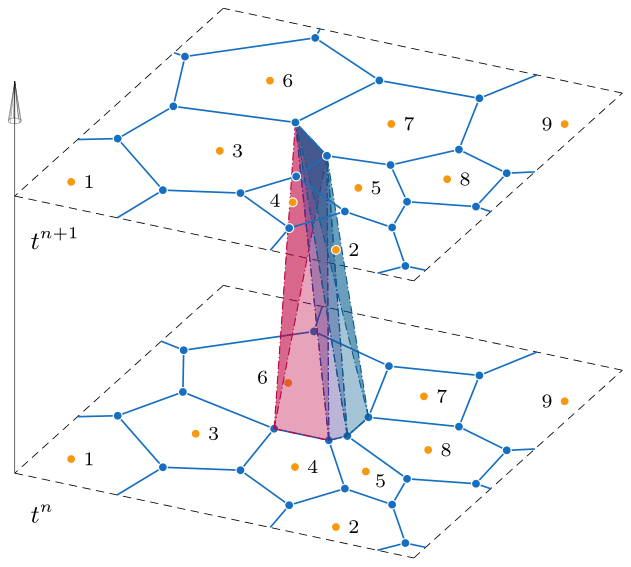}\label{sfig.consSliver_c}}%
	\subfloat[]{\includegraphics[width=0.5\linewidth]{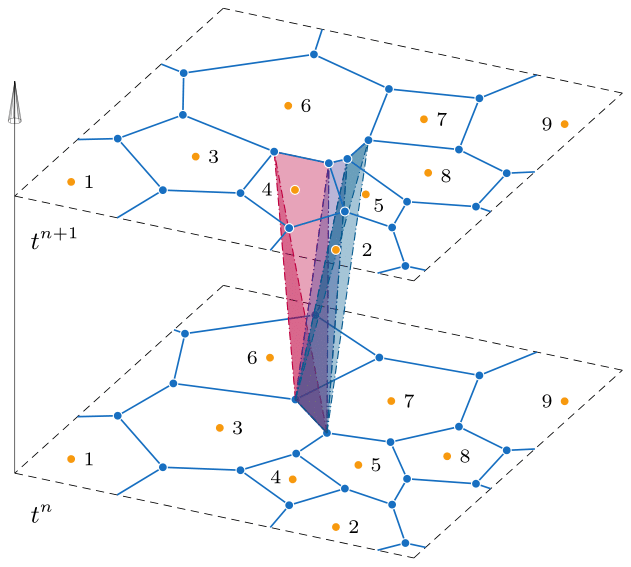}\label{sfig.consSliver_d}}%
	\caption{Consecutive space--time sliver elements. Refer for example to Panel (d): $P_3^n$ and $P_7^n$ are neighbors at time $t^n$ but this is no longer the case at time $t^{n+1}$ and moreover $P_4^{n+1}$, $P_5^{n+1}$, $P_6^{n+1}$ and $P_8^{n+1}$ are among them; this complex change of topology causes the appearance of $3$ space--time sliver elements. A similar situation with $3$ space--time sliver elements is depicted in Panel (c). In Panels (a) and (b) we show a change of topology with $2$ space--time sliver elements.}
	\label{fig.consecutiveSlivers}
\end{figure}

If this is not the case, as illustrated in Figures~\ref{sfig.consSliver_b}-\ref{sfig.consSliver_d}, \textit{more consecutive sliver elements} will be necessary to fill the space--time holes. These consecutive sliver elements have the bottom face in common, given by the edge between $P_i^n - P_{i_j}^n$, 
and the top faces given respectively by the edges composing the path connecting $P_i^{n+1} \rightarrow P_{i_j}^{n+1}$.
A similar procedure is employed for situations depicted in Figures~\ref{sfig.consSliver_a}-\ref{sfig.consSliver_c}, corresponding to case (III).
We allow a maximum of three consecutive sliver elements. {\scriptsize{$\blacksquare$}}

\medskip

Two \textit{problems} can arise while assembling the space--time connectivity: 
$\mathcal{V}(C_i^n)$ could be not sortable respecting both the order of $\mathcal{V}(P_i^n)$ and $\mathcal{V}(P_i^{n+1})$, or more than three sliver elements could be necessary to complete the connection path. In this case a MOOD~\cite{boscheri2015direct,ALEMOOD2} procedure described in Section~\ref{ssec.MOOD} will be adopted.

\subsection{Degenerate sub--space--time control volumes and sliver space--time elements}
\label{ssec.Sliver}

The change of topology induces the appearance of degenerate elements in the space--time connectivity. 

As is evident from Figures~\ref{sfig.Crazy_b}-\ref{sfig.Crazy_c}, some of the sub--space--time control volumes $sC_{i_j}^n$ of $C_i^n$, are triangular prisms with one of their top or bottom faces collapsed to just a line, and with the lateral space--time surface $\partial{C_{i_j}^n}$ being of triangular shape (instead of the standard quadrilateral shape).
They do not pose particular problems because they are part of a standard control volume, so everything is naturally well defined on them (basis functions, quadrature points, values of the numerical solution $\u_h^n$, of the reconstruction polynomials $\w_h^n$, \RIIcolor{and of the space--time predictor $\q_h^n$ defined below in~\eqref{eqn.qh}}). 

On the contrary, the space--time sliver element in Figure~\ref{sfig.Crazy_d} is a completely new control volume which does neither exist at time $t^n$, nor at time $t^{n+1}$, since it coincides with an edge of the tessellation at the old and at the new time levels, and, as such, has zero area in space at $t^n$ and $t^{n+1}$. 
However, it has a \textit{non-negligible volume} in space--time. The difficulties related to this kind of elements are due to the fact that $\w_h$ is not clearly defined for them at time $t^n$ and that contributions across them should not be lost at time $t^{n+1}$, in order to ensure conservation. Space--time sliver elements always have four neighbors, 
namely the two Voronoi polygons that share their degenerate bottom face (edge) and the two Voronoi polygons that share their degenerate top face (edge).

Note that the computation of numerical fluxes across degenerate triangular space--time faces has already been treated in~\cite{gaburro2016direct}. In the same paper a proof of concept was given, that situations like those shown in Figures~\ref{sfig.Crazy_b}-\ref{sfig.Crazy_c} could be handled up to second order of accuracy. Instead, the treatment of sliver elements is a completely new topic.

\section{Numerical method II: high order fully-discrete direct ALE FV-DG scheme}
\label{sec.NumMethod_ALE_FV-DGscheme}

The governing equations~\eqref{eq.generalform} are now solved with the aid of a high order fully-discrete \textit{one-step predictor-corrector} ADER FV-DG method obtained by generalizing the scheme first presented in~\cite{Dumbser2008} to our moving Voronoi meshes with topology change. 
ADER finite volume schemes go back to the pioneering work of Toro and Titarev~\cite{toro3,toro4,titarevtoro,schwartzkopff,Toro:2006a} on approximate solvers for the generalized Riemann problem (GPR). ADER schemes have been originally developed in the Eulerian framework on fixed grids~\cite{toro3,toro4,titarevtoro,schwartzkopff,Toro:2006a,DumbserEnauxToro,busto2018projection} and have subsequently also been extended to moving meshes in the ALE context~\cite{boscheri2017high,Lagrange3D,ALEMQF,LagrangeISO}.

We recall that high order of accuracy in \textit{space} is provided by the piecewise polynomial data representation $\w_h^n$, which for $N=M>0$ coincides with the DG polynomial, i.e. $\w_h^n = \u_h^n$, while, in the Finite Volume case ($N=0$), $\w_h^n$ is obtained through the reconstruction procedure described in Section~\ref{ssec.CWENO}.
In any case, $\w_h^n$ only depends on the mesh configuration at time $t^n$, so that an eventual degeneracy of the space--time geometry does not affect this first step.

Then, the \textit{predictor} step consists in a {\it local} solution of the governing PDE~\eqref{eq.generalform} \textit{in the small}, see~\cite{eno}, inside each space-time element $C_i^n$, thus including the sliver elements $S_i^n$. 
It is called \textit{local} because it is obtained by only considering cell $C_i^n$ with initial data $\w_h^n$ on $P_i^n$, the governing equations~\eqref{eq.generalform} and the geometry of $C_i^n$, without taking into account any interaction between $C_i^n$ and its neighbors.  
It provides, for each space--time control volume $C_i^n$, a polynomial data representation $\q_h^n$ (see below for the details) of high order both in space and time, which serves as a predictor solution, only valid inside $C_i^n$,
to be used for evaluating the numerical fluxes and sources when integrating the PDE 
in the final corrector step of the ADER scheme.  

Lastly, the \textit{corrector} step integrates the weak form of the PDE over the  space-time control volumes $C_i^n$, making use of the predictor solution 
$\q_h^n$, and returns $\u_h^{n+1}$ by taking care of the coupling with neighbors through the numerical flux computations across $\partial C_i^n$. 
It ensures high order of accuracy in space and time, provided the high order of accuracy of  $\q_h^n$. The scheme is by construction conservative since it takes into account all the flux contributions over $\partial C_i^n$, including those across the sliver elements (see Section~\ref{sssec.NumScheme_sliver_Flux}). 
Moreover, the method is stable if the time-step size $\Delta t$ satisfies an explicit CFL stability condition, which reads 
\begin{equation}
\Delta t < \textnormal{CFL} \left( 
\frac{ |P_i^n| }{ (2N+1) \, |\lambda_{\max,i}| \, \sum_{\partial P_{i_j}^n} |\ell_{i_j}| } 
\right), \qquad \forall P_i^n \in \Omega^n. 
\label{eq:timestep}
\end{equation}
In the above formula, $\ell_{i_j}$ is the length of the edge $j$ of $P_i^n$ and $|\lambda_{\max,i}|$ is the spectral radius of the Jacobian of the flux $\mathbf{F}$.
On unstructured meshes the CFL stability condition requires the inequality $\textnormal{CFL} < \frac{1}{d}$ to be satisfied, see~\cite{Dumbser2008}.

\subsection{High order in time: space--time predictor}
\label{ssec.predictor}

In what follows, a predictor of the solution is recovered, which is valid \textit{locally} inside $C_i^n$ and is given by high order piecewise space-time  polynomials $\q_h^n(\x,t)$ of degree $M$ that are expressed as 
\be
\q_h^n(\x, t) = \sum_{\ell=0}^{\mathcal{Q}-1} \theta_\ell (\x, t) \hat{\q}_\ell^n, \qquad (\x,t) \in C_i^n, \qquad \mathcal{Q} = \mathcal{L}(M,d+1).  
\label{eqn.qh}
\ee 
with $\theta_\ell(\x, t)$  being a \textit{modal space--time} basis of the polynomials of degree $M$ in $d+1$ dimensions ($d$ space dimensions plus time), which read
\begin{equation}
\label{eq.Dubiner_phi}
 \theta_\ell(x,y,t)|_{C_i^n} = \frac{(x - \xbixn)^{p_\ell}}{{p_\ell}! \, h_i^{p_\ell}} \, \frac{(y - \xbiyn)^{q_\ell}}{{q_\ell}! \, h_i^{q_\ell}}
\, \frac{(t - t^n)^{q_\ell}}{{q_\ell}! \, h_i^{q_\ell}}, 
 \qquad \ell = 0, \dots, \mathcal{L}(M,d+1), 
 \quad 0 \leq p_\ell + q_\ell + r_\ell \leq M. \\
\end{equation}  
\RIIIcolor{These basis functions $\theta$ are redefined at the beginning of each time step in function of the current position $\xbin$, thus they are directly linked to the current mesh configuration; however, contrarily to the test functions of Equation~\eqref{eq.Dubiner_phi_movingspatial} that are used in the corrector step (see the next section), there is no need to \textit{move} them during each time step, since they allow to represent information at the predictor step, which is only valid \textit{locally} inside each $C_i^n$.}

The predictor $\q_h^n$ is computed through an iterative procedure that looks for the polynomial satisfying a weak form of~\eqref{eq.generalform} obtained for any control volume $C_i^n$ as follows. We multiply the governing PDE~\eqref{eq.generalform} by a test function $\theta_k$,  integrate over $C_i^n$ and insert the discrete solution $\q_h^n$ instead of $\mathbf{Q}$, hence
obtaining 
\begin{equation}
\int_{C_i^n} \theta_k(\x,t) \de{\q_h^n}{t} \, d\x dt  + \int_{C_i^n} \theta_k(\x,t) \nabla \cdot \F(\q_h^n) \, d\x dt  = \int_{C_i^n} \theta_k(\x,t) \, \S(\q_h^n) \, d\x dt.
\label{eqn.PDEweak1}
\end{equation}
Differently from what has been proposed in~\cite{Dumbser2008,DumbserEnauxToro,Lagrange2D,Lagrange3D}, here we do \textit{not} integrate the first term in~\eqref{eqn.PDEweak1} by parts in time. Instead, we take into account potential jumps of $\q_h$ on the boundaries of $C_i^n$ in the sense of distributions, combined with upwinding of the fluxes in time. This approach is similar  to the path-conservative schemes proposed in~\cite{Pares2006,Castro2006,Castro2008}, but much simpler, since the test functions are only taken from within $C_i^n$ and there is no need to define a non-conservative product on $\partial C_i^n$. 
Therefore, the integral containing the time derivative in~\eqref{eqn.PDEweak1} is rewritten as 
\begin{equation}
\int_{C_i^n} \theta_k(\x,t) \de{\q_h^n}{t} \, d\x dt = \int_{C_i^n \backslash \partial C_i^n} \theta_k(\x,t) \de{\q_h^n}{t} \, d\x dt + \int_{\partial C_i^n} \theta_k(\x,t) \, \left( \q_h^{n,+} - \q_h^{n,-} \right) \, \tilde{\mbf{n}}_t^- \, dS,
\label{eqn.TimeInt}
\end{equation}
\RIIIcolor{where $C_i^n \backslash \partial C_i^n$ denotes the interior of $C_i^n$}.
Here, $\q_h^{n,-}$ and $\q_h^{n,+}$ denote the boundary-extrapolated inner and outer states across the jump on $\partial C_i^n$.
Furthermore, $\tilde{\mbf{n}}^-$ are only those outward pointing unit-normal vectors on  $\partial C_i^n$ that point \textit{back} in time and $\tilde{\mbf{n}}_t^-$ is their time component, i.e. $\tilde{\mbf{n}}_t^- = \min( 0, \, \tilde{\mbf{n}} \cdot (0,0,1) ) \leq  0$. 
Upwinding in time is therefore automatically guaranteed, since we only consider the contributions coming from the \textit{past}, according to the causality principle. 
In other words, only time fluxes that \textit{enter} the space--time control volume $C_i^n$ contribute to the jump term in~\eqref{eqn.TimeInt}, and they are easily  identified by checking the sign of the time component of the space--time 
normal vector $\tilde{\mbf{n}}$.  

\subsubsection{Space--time predictor on standard space--time elements} 
\label{sec.pred.standard} 

For standard elements, we apply the jump term only on the bottom surface $P_i^n$ of the space--time element $C_i^n$ under consideration, where it then simplifies to 
\begin{equation}
\left. \left( \q_h^{n,+} - \q_h^{n,-} \right) \, \tilde{\mbf{n}}_t^- \ \right|_{P_i^n}  = -\left(  \w_h^{n}(\x,t^n) - \q_h^{n}(\x,t^n) \right) =  \q_h^{n}(\x,t^n) - \w_h^{n}(\x,t^n), 
\label{eqn.pathint3} 
\end{equation}
with $\q_h^{n,+} = \w_h(\x,t^n)$ being simply given by the reconstruction polynomial at time $t^n$
and obviously $\tilde{\mbf{n}}^- = (0,0,-1)$ on $P_i^n$ and thus $\tilde{\mbf{n}}_t^- = -1$. In this case,~\eqref{eqn.TimeInt} reduces to 
\begin{equation}
\int_{C_i^n} \theta_k(\x,t) \de{\q_h^n}{t} \, d\x dt = \int_{C_i^n \backslash P_i^n} \theta_k(\x,t) \de{\q_h^n}{t} \, d\x dt + \int_{P_i^n} \theta_k(\x,t^n) \, \left( \q_h^{n}(\x,t^n) - \w_h(\x,t^n) \right) \, d\x
\label{eqn.TimeInt.std}
\end{equation}
for standard space--time elements. 
The reason for this choice is that in this manner, all space--time predictors of the standard elements are \textit{decoupled} from each other, since they only require the initial data $\w_h^n$ and no information from the neighbor elements. 
This will not be the case for sliver elements, for which we do not have any reconstruction polynomial available at $t^n$.
If we considered the jump terms also on lateral surfaces of standard space--time elements, the space--time predictors would no longer be independent of each other, since our mesh is \textit{moving} and there will be in general always a non--empty subset of $\partial C_i^n$ with $\tilde{\mathbf{n}}_t^- < 0$. This would require a proper ordering of the execution sequence of the space--time predictors on the standard elements, but this is something we want to avoid. 
%
With the following definitions
\begin{eqnarray}
&&\mathbf{K}_{1} = \int_{C_i^n \backslash P_i^n} \theta_k \frac{\partial \theta_\ell}{\partial t} \, d\x dt, \quad \mathbf{K}_{x} = \int_{C_i^n} \theta_k \frac{\partial \theta_\ell}{\partial x} \, d\x dt, \quad \mathbf{K}_{y} = \int_{C_i^n} \theta_k \frac{\partial \theta_\ell}{\partial y} \, d\x dt,  \nonumber \\
&&\mathbf{M} = \int_{C_i^n} \theta_k \theta_\ell \, d\x dt, \quad 
\mathbf{F}_{0} =\int_{P_i^n} \theta_k(\x,t^n) \psi_\ell(\x,t^n)  \, d\x, \quad  \mathbf{F}_{1} =\int_{P_i^n} \theta_k(\x,t^n) \theta_\ell(\x,t^n)  \, d\x,
\end{eqnarray}  
the weak form~\eqref{eqn.PDEweak1}-\eqref{eqn.TimeInt} can be compactly rewritten as
\begin{equation}
\left(\mathbf{K}_{1} + \mathbf{F}_{1} \right) \hat{\q}_i^n = \mathbf{F}_{0} \hat{\w}_i^n -\mathbf{K}_{x} \, \f(\hat{\q}_i^n) - \mathbf{K}_{y} \, \g(\hat{\q}_i^n) + \mathbf{M} \, \S(\hat{\q}_i^n),
\label{eqn.PDEpredictor}
\end{equation}
where $\hat{\q}_i^n$ and  $\hat{\w}_i^n$ contain all the expansion coefficients of $\hat{\q}_{\ell,i}^n$ in~\eqref{eqn.qh} and $\hat{\w}_{\ell,i}^n$ in~\eqref{eqn.wh}, respectively. The solution of~\eqref{eqn.PDEpredictor} can be found via a simple and fast converging fixed point iteration (a discrete Picard iteration), as detailed in~\cite{Dumbser2008, hidalgo2011ader}. Here, as initial guess we simply impose $\hat{\q}_{\ell,i}^n = \hat{\w}_{\ell,i}^n$ for the common spatial degrees of freedom (with $\ell \le \mathcal{M}$) and zero for the other ones. 
For linear homogeneous systems, the discrete Picard iteration
converges in a finite number of at most $M+1$ steps, since the involved iteration matrix is nilpotent,
see~\cite{jackson2017eigenvalues}. In the nonlinear case we allow a maximum of $10$ iterations if convergence is not reached before, being $M+1$ iterations enough for obtaining the correct order $M$ of convergence.

Notice again that in~\eqref{eqn.TimeInt.std}  and therefore in~\eqref{eqn.PDEpredictor} we have considered only one jump term, namely the contribution \textit{coming from the past} through the bottom face $P_i^n$ of $C_i^n$, where $\w_h^n=\w_h(\x,t^n)$ is \textit{known} and \textit{well defined}. This allows us to couple~\eqref{eqn.PDEweak1} with the initial condition  $\w_h(\x,t^n) |_{P_i^n}$ via~\eqref{eqn.TimeInt.std}. No other information (as neighbors values) is taken into account in this local phase. Indeed, neighbor data will be considered later in the corrector step (Section~\ref{ssec.corrector}).  

\begin{figure}[!bp]
	\centering
	\subfloat[]{\includegraphics[width=0.33\linewidth]{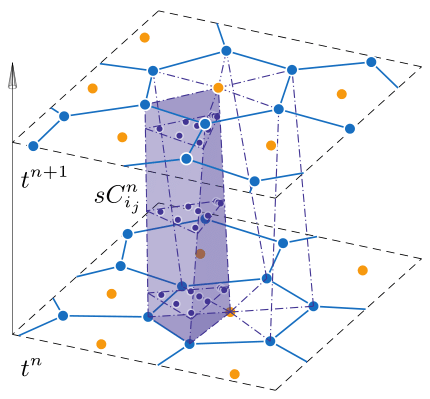}\label{sfig.quad_volume}}%
	\subfloat[]{\includegraphics[width=0.33\linewidth]{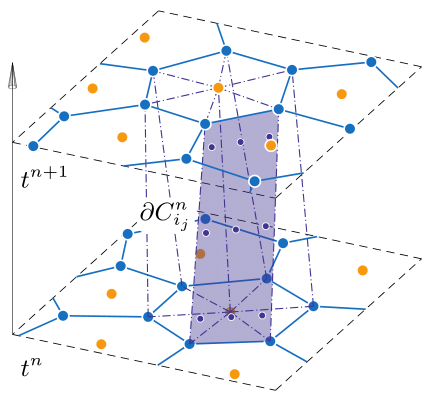}\label{sfig.quad_surface}}%
	\subfloat[]{\includegraphics[width=0.33\linewidth]{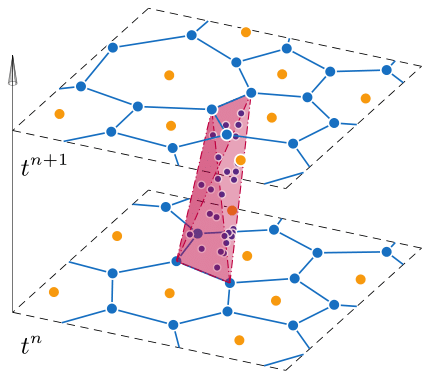}}
	\caption{Space--time quadrature points for third order methods, i.e. $M=2$. (a) Quadrature points for the volume integrals and the space--time predictor. (b) Quadrature points for the surface integrals, i.e. for flux computation. (c) Quadrature points for the volume integrals and the space--time predictor for a sliver element.}
	\label{fig.quadraturePoints}
\end{figure}

The integrals above are evaluated using multidimensional Gaussian quadrature rules of suitable order of accuracy, see~\cite{stroud} and Figure~\ref{fig.quadraturePoints} for details. In order to carry out the integration, we split the space-time volume $C_i^n$ into a set of sub--space-time volumes $sC_{i_j}^n$ of $C_i^n$, whose shape is an oblique triangular prism.
Note that for degenerate sub--space--time control volumes, as those of Figures~\ref{sfig.Crazy_b} and~\ref{sfig.Crazy_c}, the above quadrature formulae
remain well defined, hence the predictor procedure over them does not pose any problem and does not need any adaptation. 

We emphasize that we \textit{first} carry out the space--time predictor for all standard elements, which can be computed independently of each other, and only  subsequently process the remaining space--time sliver elements. The reason for this will become clear in the next section. 

\subsubsection{Space--time predictor on the space--time sliver elements}
\label{sssec.NumScheme_sliver_Predictor}  

The predictor procedure on space--time sliver elements, as those shown in Figures~\ref{sfig.Crazy_d} and~\ref{fig.consecutiveSlivers}, needs particular care.
The main problem connected with the space--time sliver elements is the fact that 
their bottom face is degenerate and consists only in a line segment, hence the 
spatial integral over $P_i^n$ vanishes, i.e. there is no possibility to introduce
the initial condition of the local Cauchy problem at time $t^n$ into the 
predictor for space--time sliver elements. 

Furthermore, the degenerate bottom faces are the edges of the Voronoi tesselation at $t^n$ and are thus at the interface between two adjacent elements, which have in  principle a discontinuous solution $\w_h^n$. Therefore, an initial value for a sliver element is in general not easy to define.  
Thus, in order to couple~\eqref{eqn.PDEweak1} with some known data from the past we have to slightly modify the algorithm detailed previously.

In particular, the upwinding in time approach is not only used for the surface $P_i^n$, as done in~\eqref{eqn.pathint3}, but we actually use the jump terms on the \textit{entire} part of the space--time surface $\partial C_{i}^n$ that closes a sliver control volume. As already stated in the previous section, the information needed to feed the predictor is allowed to come only from the past, i.e. only from those space--time neighbors $C_j^n$ whose common surface $\partial C_{ij}^n = C_i^n \cap C_j^n$ exhibits a \textit{negative} time component of the outward pointing space--time normal vector ($ \tilde{\mbf{n}}_t^- < 0$). In this way, we can introduce information from the past into the space--time sliver elements by considering also its neighbor elements, but respecting at the same time the causality principle in time, hence using again upwinding for the flux evaluation of the jump term in~\eqref{eqn.TimeInt}.   
As a consequence, the predictor solution $\q_{h}^{n}$ is again obtained by means of~\eqref{eqn.PDEweak1}, but treating the \textit{entire} space--time surface $\partial C_{i}^n$ with the upwind in time approach, hence leading to
\begin{equation}
\left(\mathbf{K}_{1}^* - \mathbf{F}_{1}^* \right) \hat{\q}_i^n = -\sum_j \mathbf{F}_{j}^* \, \hat{\q}_j^{n} - \mathbf{K}_{x}^* \, \f(\hat{\q}_i^n) - \mathbf{K}_{y}^* \, \g(\hat{\q}_i^n) + \mathbf{M}^* \, \S(\hat{\q}_i^n),
\label{eqn.PDEpredictor_sliver}
\end{equation}
where the following definitions for the sliver element hold
\begin{eqnarray}
&&\mathbf{K}_{1}^* = \int_{C_i^n \backslash \partial C_{i}^n} \theta_k \frac{\partial \theta_\ell}{\partial t} \, d\x dt, \quad \mathbf{K}_{x}^* = \int_{C_i^n} \theta_k \frac{\partial \theta_\ell}{\partial x} \, d\x dt, \quad \mathbf{K}_{y}^* = \int_{C_i^n} \theta_k \frac{\partial \theta_\ell}{\partial y} \, d\x dt, \nonumber \\
&& \mathbf{M}^* = \int_{S_i^n} \theta_k \theta_\ell \, d\x dt, \quad 
\mathbf{F}_{1}^* =\int_{\partial C_{i}^n} \theta_k \theta_\ell \, \tilde{\mbf{n}}_t^- \, dS, \quad 
\mathbf{F}_{j}^* =\int_{\partial C_{ij}^n} \theta_k \theta_\ell \, \tilde{\mbf{n}}_t^- \, dS.
\end{eqnarray} 

This is slightly different from what is done for standard elements in~\eqref{eqn.PDEpredictor}, where only the space--time surface at time $t^n$, i.e. $P_i^n$, is considered for introducing the initial condition $\w_h^n$.  
Here, the information from the past comes through the upwind fluxes contained in the term $\mathbf{F}_{j}^* \, \hat{\q}_j^{n}$ in~\eqref{eqn.PDEpredictor_sliver} and thus requires the knowledge of the predictor solution $\hat{\q}_j^{n}$ in the neighbor $C_j^n$. 
This is the reason why the predictor step must \textit{first} be performed over all the standard elements using~\eqref{eqn.PDEpredictor}, so that the predictor solution $\q_h^n$ is always available to feed the temporal fluxes with the quantities  $\hat{\q}_j^{n}$ that are needed for solving~\eqref{eqn.PDEpredictor_sliver} in the 
case of the space--time sliver elements. We underline again that a space--time sliver element has always four standard Voronoi elements as neighbors.  
This closes the description of the predictor step for the space--time sliver elements.

\subsection{Corrector step: direct ALE FV-DG scheme}
\label{ssec.corrector}

This section contains the core of our direct ALE FV-DG scheme used to solve~\eqref{eq.generalform} on regenerating moving meshes.

Following~\cite{Lagrange2D, Lagrange3D, LagrangeISO}, the PDE system~\eqref{eq.generalform} is rewritten in a space-time divergence form as
\begin{equation}
\tilde \nabla \cdot \tilde{\F} = \mathbf{S}, 
\label{eqn.st.pde}
\end{equation}
with $\tilde \nabla = \left( \partial_x, \, \partial_y, \, \partial_t \right)$ denoting the space-time divergence operator and $\tilde{\F} = \left( \mathbf{f}, \, \mathbf{g}, \, \Q \right)$ being the corresponding space-time flux tensor.
Then, we multiply~\eqref{eqn.st.pde} by a set of \textit{moving} spatial modal test functions $\phitilde_k(\x,t)$, 
which coincide with~\eqref{eq.Dubiner_phi_spatial} at $t=t^n$ and at $t=t^{n+1}$, 
i.e. $\phitilde_k(\x,t^n) = \phi_k(\x,t^n)$ and 
$\phitilde_k(\x,t^{n+1}) = \phi_k(\x,t^{n+1})$. The test functions are tied to the motion of the barycenter $\x_{\mathbf{b}_i}(t)$ and move together with $P_i(t)$ in such a way that at time $t = t^{n+1}$ they refer to the new barycenter $\xbinp$. Thus, the test functions explicitly read as follows:  
\begin{eqnarray} 
\label{eq.Dubiner_phi_movingspatial}
&& \phitilde_\ell(x,y,t)|_{C_i^{n}} = 
\frac{(x - x_{b_i}(t))^{p_\ell}}{{p_\ell}! \, h_i^{p_\ell}} \, \frac{(y - y_{b_i}(t))^{q_\ell}}{{q_\ell}! \, h_i^{q_\ell}}, \quad 
 \text{ with } \ \x_{\mbf{b}_i} (t) = \frac{t-t^n}{\Delta t} \xbin + \left(1-\frac{t-t^n}{\Delta t}\right) \xbinp, \\ 
&& \ell = 0, \dots, \mathcal{N}, \quad  0 \leq p+q \leq N. \nonumber 
\end{eqnarray} 
These \textit{moving modal basis functions} are \textit{essential} for the approach presented in this paper. They \textit{naturally} allow for topology changes, without the need of any remapping steps, which we want to avoid in a direct ALE formulation. 

Next, integration over the closed space-time control volume $C_i^n$ yields
\begin{equation}
\int_{C_i^n} \phitilde_k \tilde \nabla \cdot \tilde{\F}(\Q) \, d\mathbf{x} dt = \int_{C_i^n} \phitilde_k \mathbf{S} (\Q)\, d\x dt.
\label{eqn.intPDE}
\end{equation}
Application of the Gauss theorem leads to the following weak form that is the basis of our fully-discrete ALE scheme 
\begin{equation}
\int_{\partial C_{i}^n} 
\phitilde_k  \tilde{\F}(\Q) \cdot \ \mathbf{\tilde n} \, dS   
- \int_{C_i^n} 
\tilde \nabla \phitilde_k \cdot \tilde{\F}(\Q) \, d\mathbf{x} dt = \int_{C_i^n} \phitilde_k\mathbf{S}(\Q) \, d\x dt,
\label{eqn.GaussPDE}
\end{equation}
where $\mathbf{\tilde n} = (\tilde n_x,\tilde n_y,\tilde n_t)$ denotes the outward pointing space-time unit normal vector on the space-time faces composing the boundary $\partial C_{i}^n$ of the space-time control volume. Moreover, the surface integral can be decomposed over the faces of $\partial C_{i}^n$ given by~\eqref{eqn.dC}. 

\subsubsection{Corrector step for standard space--time elements}
\label{sec.corrector.standard} 

We first describe the corrector step for standard space--time control volumes. 
After introducing the discrete solution $\u_h$, the space--time predictor $\q_h$ and a two-point numerical flux function on the element boundaries of the type 
\be 
\label{eqn.fluxFunctionGeneric}
\tilde{\F} (\Q) \cdot \ \mathbf{\tilde n} := \mathcal{F}(\q_h^{n,-},\q_h^{n,+}) \cdot \mathbf{\tilde n},
\ee 
into~\eqref{eqn.GaussPDE}, where $\q_h^{n,-}$ and $\q_h^{n,+}$ are the inner and outer boundary-extrapolated data respectively,
\RIIcolor{(i.e. the values assumed by the predictors of the two neighbor elements at a point on the shared space--time lateral surface)}, 
we obtain the final direct ALE scheme: 
\begin{equation}
\int \limits_{P_i^{n+1}} 
\phitilde_k  \u_h(\x,t^{n+1}) \, d\x   = 
\int  \limits_{P_{i}^n} 
\phitilde_k  \u_h(\x,t^{n})  \, d\x   
- \sum_{j=1}^{N_{V_i}^{n,st}} \int \limits_{\partial C_{ij}^n} 
\phitilde_k  \mathcal{F}(\q_h^{n,-},\q_h^{n,+}) \cdot \mathbf{\tilde n} \, dS  
+ \int \limits_{C_i^n} 
\tilde \nabla \phitilde_k \cdot \tilde{\F} (\q_h)\, d\mathbf{x} dt + \int \limits_{C_i^n} \phitilde_k \mathbf{S}(\q_h) \, d\x dt,
\label{eqn.GaussPDE2}
\end{equation}
where the unknown solution at the new time step $\u_h(\x,t^{n+1})$ can be computed  \textit{directly} from the solution at the previous time step $\u_h(\x,t^{n})$ 
through the integration of the fluxes and source terms over $C_i^n$, 
without needing any further remapping/remeshing steps.  

Our scheme is high order accurate in space and time because 
the predictor solution $\q_h^n$, which is given by piecewise space--time polynomials of degree $M$, is employed for a high order accurate space-time integration of all remaining terms in~\eqref{eqn.GaussPDE2}, namely the numerical 
surface flux integral on $\partial C_{i_j}^n$ and the volume integrals on $C_i^n$ for the fluxes and the source terms. 

The boundary fluxes are obtained by a Riemann solver, 
thus providing the coupling between neighbors, which was neglected in the predictor step. 
The ALE Jacobian matrix w.r.t. the normal direction in space reads
\begin{equation} 
\label{eq.ALEjacobianMatrix}
\mathbf{A}^{\!\! \mathbf{V}}_{\mathbf{n}}(\Q)=\left(\sqrt{\tilde n_x^2 + \tilde n_y^2}\right)\left[\frac{\partial \mathbf{F}}{\partial \Q} \cdot \mathbf{n}  - 
(\mathbf{V} \cdot \mathbf{n}) \,  \mathbf{I}\right], \qquad    
\mathbf{n} = \frac{(\tilde n_x, \tilde n_y)^T}{\sqrt{\tilde n_x^2 + \tilde n_y^2}},  
\end{equation} 
with $\mathbf{I}$ representing the identity matrix and $\mathbf{V} \cdot \mathbf{n}$ denoting the local normal mesh velocity. Furthermore, $\mathbf{n}$ is the spatial normalized normal vector, which is different from the space-time normal vector $\mathbf{\tilde{n}}$. We adopt either a simple and robust Rusanov-type~\cite{Rusanov:1961a} ALE scheme, 
\begin{equation}
\mathcal{F}(\q_h^{n,-},\q_h^{n,+}) \cdot \mathbf{\tilde n} =  
\frac{1}{2} \left( \tilde{\F}(\q_h^{n,+}) + \tilde{\F}(\q_h^{n,-})  \right) \cdot \mathbf{\tilde n}_{ij}  - 
\frac{1}{2} s_{\max} \left( \q_h^{n,+} - \q_h^{n,-} \right),  
\label{eq.rusanov} 
\end{equation} 
where $s_{\max}$ is the maximum eigenvalue of $\mathbf{A}^{\!\! \mathbf{V}}_{\mathbf{n}}(\q_h^{n,+})$ and $\mathbf{A}^{\!\! \mathbf{V}}_{\mathbf{n}}(\q_h^{n,-})$,
or a less dissipative Osher-type~\cite{osherandsolomon,OsherUniversal} ALE flux
\begin{equation}
\mathcal{F}(\q_h^{n,-},\q_h^{n,+}) \cdot \mathbf{\tilde n} = 
\frac{1}{2} \left( \tilde{\F}(\q_h^{n,+}) + \tilde{\F}(\q_h^{n,-})  \right) \cdot \mathbf{\tilde n}_{ij}  - 
\frac{1}{2} \left( \int_0^1 \left| \mathbf{A}^{\!\! \mathbf{V}}_{\mathbf{n}}(\boldsymbol{\Psi}(s)) \right| ds \right) \left( \q_h^{n,+} - \q_h^{n,-} \right),  
\label{eqn.osher} 
\end{equation} 
where we choose to connect the left and the right state across the discontinuity using a simple straight--line segment path  
\begin{equation}
\boldsymbol{\Psi}(s) = \q_h^{n,-} + s \left( \q_h^{n,+} - \q_h^{n,-} \right), \qquad 0 \leq s \leq 1.  
\label{eqn.path} 
\end{equation} 
The absolute value of $\mathbf{A}^{\!\! \mathbf{V}}_{\mathbf{n}}$ is evaluated as usual as 
$ \mathbf{R} |\boldsymbol{\Lambda}| \mathbf{R}^{-1}$,  
where $\mathbf{R}$, $\mathbf{R}^{-1}$ and $\boldsymbol{\Lambda}$ denote, respectively, the right eigenvector matrix, its inverse and the diagonal matrix of the eigenvalues of $\mathbf{A}^{\!\! \mathbf{V}}_{\mathbf{n}}$.

\bigskip
Finally, using the definitions~\eqref{eqn.uh} and~\eqref{eqn.wh}, our arbitrary high order one-step direct ALE FV-DG scheme becomes 
\begin{equation}
\left( \!
\int_{P_i^{n+1}} \!\!\phitilde_k \phi_\ell \, d\x 
\right) 
\hat{\mathbf{u}}^{n+1}_{\ell} 
\! = \left( \!
\int_{P_i^{n}} \!\!\phitilde_k \psi_\ell \, d\x 
\right) 
\hat{\mathbf{w}}^{n}_{\ell} 
- \sum_{j=1}^{N_{V_i}^{n,st}} \int_{\partial C_{i_j}^n} \!\!\!
\phitilde_k \mathcal{F}(\q_h^{n,-},\q_h^{n,+}) \cdot \mathbf{\tilde n} \, dS   
+ \!\int_{C_i^n} \!\!\! \tilde \nabla 
\phitilde_k \cdot \tilde{\F}(\q_h^n) \, d\mathbf{x} dt + \! \int_{C_i^n} \!\! \phitilde_k \mathbf{S}(\q_h^n) d\x dt.
\label{eqn.ALE-DG}
\end{equation}
The volume integrals in the above expression~\eqref{eqn.ALE-DG} can be easily computed directly on the physical space-time element $C_i^n$ by summing up the contributions on each sub-volume $sC_{i_j}^n$ and employing Gaussian quadrature rules of sufficient precision, see~\cite{stroud}. 
The lateral space--time surfaces of $\partial C_{ij}^n$ instead are parameterized using a set of bilinear basis functions~\cite{Lagrange2D}, that is 
\be
\partial C_{ij}^n = \mathbf{\tilde{x}} \left( \chi,\tau \right) = 
\sum\limits_{k=1}^{4}{\beta_k(\chi,\tau) \, \mathbf{\tilde{X}}_{ij,k}^n },	
\qquad 0 \leq \chi \leq 1,  \quad	0 \leq \tau \leq 1, 										 
\label{eq.SurfParBeta}
\ee
where the $\mathbf{\tilde{X}}_{ij,k}^n$ represent the physical space--time coordinates of the four vertexes of $\partial C_{i_j}^n$, and the functions $\beta_k(\chi,\tau)$ are defined as follows 
\begin{equation}
\beta_1(\chi,\tau) = (1-\chi)(1-\tau), \qquad 
\beta_2(\chi,\tau) = \chi(1-\tau), \qquad 
\beta_3(\chi,\tau) = \chi\tau, \qquad  
\beta_4(\chi,\tau) = (1-\chi)\tau.
\label{BetaBaseFunc}
\end{equation}
The mapping in time is given by the transformation 
\be
t = t_n + \tau \, \Delta t, \qquad  \tau = \frac{t - t^n}{\Delta t}.
\label{eq.timeTransf}
\ee 
In this way, every $\partial C_{i_j}^n$ (even if degenerate, i.e. with a triangular shape) can be mapped to a reference square $[0,1]\times[0,1]$ and surface integrals can be computed.

We close this section remarking that the integration of the governing PDE over the space-time volume $C_i^n$ automatically satisfies the geometric conservation law (GCL) for all test functions $\phitilde_k$. This simply follows from Gauss theorem applied to \textit{closed} space--time control volumes and we refer to~\cite{Lagrange3D} for a complete proof. 
\RIIIcolor{The satisfaction of the GCL property up to machine precision has been numerically 
verified in \textit{each} simulation presented in this paper and a series of test cases aimed at demonstrating its validity is presented in Section~\ref{test.GCLTest}}

\subsubsection{Corrector step on sliver elements}
\label{sssec.NumScheme_sliver_Flux}  

Let us now consider the numerical scheme given by~\eqref{eqn.ALE-DG} in the case of a sliver element $C_i^n = S_i^n$: 
\begin{equation}
0_\ell \, \hat{\mathbf{u}}^{n+1}_{\ell} = 
0_\ell \hat{\mathbf{w}}^{n}_{\ell} 
- \sum_{j=1}^{4} \int_{\partial S_{i_j}^n} \!\!\!
\phitilde_k \mathcal{F}(\q_h^{n,-},\q_h^{n,+}) \cdot \mathbf{\tilde n} \, dS   
+ \int_{S_i^n} \!\!\! \tilde \nabla 
\phitilde_k \cdot \tilde{\F}(\q_h^n) \, d\mathbf{x} dt + \int_{S_i^n} \!\! \phitilde_k \mathbf{S}(\q_h^n) d\x dt, 
\label{eqn.ALE-DG_sliver}
\end{equation}
Since for sliver elements $|P_i^n|=|P_i^{n+1}|=0$, the first two terms vanish. 
However, since the method is explicit and $\q_h^n$ only depends on information coming from the past, the remaining terms in~\eqref{eqn.ALE-DG_sliver} are in general not equal to zero, i.e. 
\begin{equation}
- \sum_{j=1}^{4} \int_{\partial S_{i_j}^n} \!\!\!
\phitilde_k \mathcal{F}(\q_h^{n,-},\q_h^{n,+}) \cdot \mathbf{\tilde n} \, dS    
+ \int_{S_i^n} \!\!\! \tilde \nabla 
\phitilde_k \cdot \tilde{\F}(\q_h^n) \, d\mathbf{x} dt + \int_{S_i^n} \!\! \phitilde_k \mathbf{S}(\q_h^n) d\x dt \ne  \0. \label{eqn.ALE-DG_lostQuantity}
\end{equation}
We underline that computing these quantities does not pose any problem, since $\q_h^n$ on $S_i^n$ is well defined (refer to Section~\ref{sssec.NumScheme_sliver_Predictor}), and the shape of a space--time sliver element is that of a tetrahedron in space--time, hence allowing standard quadrature rules to be used for integral evaluations.

The problem here arises from the fact that, using~\eqref{eqn.ALE-DG_sliver}, the non-null quantity~\eqref{eqn.ALE-DG_lostQuantity} will be lost at time $t^{n+1}$ because it plays a role only in the evolution of $S_i^n$, which exists between $t^n$ and $t^{n+1}$, but is null at $t^{n+1}$.
In order to be conservative, we must avoid losing any contribution from the sliver elements. We therefore couple the weak formulation on $S_i^n$ 
with the weak form of \textit{one} of its standard space--time neighbors. Here, we always choose the one with the biggest space--time volume, referred to as $C_\mathrm{big}$. The choice of the biggest volume is not mandatory, it only represents our way to uniquely fix the choice of a particular neighbor of the sliver element.
The test function $\phitilde_k$ of~\eqref{eqn.ALE-DG_sliver} is then referred to the barycenter of $C_\mathrm{big}$. Conservation is guaranteed by adding the contribution~\eqref{eqn.ALE-DG_lostQuantity} of the sliver element $S_i^n$ to the neighbor $C_\mathrm{big}$, hence
\be
\left( \int_{P_\mathrm{big}^{n+1}} \!\phitilde_k \phi_\ell \, d\x \right) \hat{\mathbf{u}}^{n+1}_{\ell}   = \, & 
\left( \int_{P_\mathrm{big}^{n}} \!\!\!\phitilde_k \psi_\ell \, d\x \right) \hat{\mathbf{w}}^{n}_{\ell}   -
\sum_{j=1}^{N_{C_\mathrm{big}}^{n,st}} \int_{\partial C_{\mathrm{big}_j}^n} \!\!\!
\phitilde_k \mathcal{F}(\q_h^{n,-},\q_h^{n,+}) \cdot \mathbf{\tilde n} \, dS   
+ \! \int_{C_\mathrm{big}^n} \!\!\!\! \tilde \nabla 
\phitilde_k \cdot \tilde{\F}(\q_h^n) \, d\mathbf{x} dt + \!\int_{C_\mathrm{big}^n} \!\!\! \phitilde_k \mathbf{S}(\q_h^n) d\x dt \\[10pt]
& + 
\sum_{j=1}^{4} \int_{\partial S_{i_j}^n} \!
\phitilde_k \mathcal{F}(\q_h^{n,-},\q_h^{n,+}) \cdot \mathbf{\tilde n} \, dS    
+ \int_{S_i^n} \!\! \tilde \nabla 
\phitilde_k \cdot \tilde{\F}(\q_h^n) \, d\mathbf{x} dt +  \int_{S_i^n} \! \phitilde_k \mathbf{S}(\q_h^n) d\x dt. \\
\ee

We would like to remark that sliver elements \textit{only} exist in between two consecutive time levels and are degenerate both at $t^n$ and $t^{n+1}$, hence they introduce some complexity in the algorithm. In particular, i) the fact that they coincide with an edge at time $t^n$ makes it difficult to fix a valid initial condition in the predictor step necessary for the high order of accuracy in time, and ii) the fact that they coincide with an edge at time $t^{n+1}$ could prevent  conservation in an explicit scheme. Nevertheless, with the strategy outlined in Sections~\ref{sssec.NumScheme_sliver_Predictor} and ~\ref{sssec.NumScheme_sliver_Flux}, no space-time contributions are lost while advancing the numerical solution in time, i.e. our proposed ADER ALE FV-DG schemes are fully conservative and keep their formal high order of accuracy even in the presence of space--time sliver elements.

Furthermore, notice that the presence of degenerate elements
is \textit{strictly unavoidable} in order to connect meshes in space \textit{and} time that include topology changes. They are also needed to collect enough geometrical  information for ensuring high order of accuracy in a direct ALE framework. 
For comparison purposes, let us consider the work presented in~\cite{re2017interpolation}, 
where the authors, in order to connect meshes with topology changes (within a different framework w.r.t. this work), have introduced some \textit{pyramidal} degenerate elements instead of our sliver elements. 
The strategy proposed in the aforementioned reference is indeed interesting and could in principle be applied also to the framework of our \textit{explicit} high order direct ALE schemes. However, besides the same complexities described for our sliver elements, an additional difficulty would arise, since a degeneracy would occur at the midpoint of the time step.

\subsection{A posteriori sub--cell finite volume limiter}
\label{ssec.limiter}

Up to now, the presented $P_NP_M$ scheme is high order accurate in space and time and, formally, 
the differences between the FV case ($N=0$) and the DG case ($N=M$) are only due to the procedure 
for achieving high order of accuracy in space, which is obtained through a CWENO reconstruction in the FV case and is instead automatic for DG. 
But there is actually one major difference, because the CWENO operator provides a nonlinear stabilization of the FV scheme, while the DG scheme presented so far is unlimited and, as such, it is affected by the so-called Gibbs phenomenon, i.e. oscillations are likely to appear in presence of shock waves or other discontinuities, which typically occur while solving nonlinear hyperbolic systems. These oscillations can be explained by the Godunov theorem~\cite{godunov}, because the presented high order DG scheme is linear in the sense of Godunov. 

As a consequence, a limiting technique is required. Our strategy is based on the MOOD approach~\cite{CDL1,CDL2,CDL3}, which has already been successfully applied in the framework of ADER finite volume schemes~\cite{ADERMOOD,ALEMOOD1,ALEMOOD2}. Specifically, the numerical solution is checked \textit{a posteriori} for nonphysical values and spurious oscillations and, instead of applying a limiter to the already computed solution, the solution is {\em locally recomputed} with a more robust scheme in the so-called \textit{troubled cells}. Troubled elements are those that do not pass the admissibility detection criteria, given by both physical and numerical requirements which mark the numerical solution as \textit{acceptable} or \textit{not acceptable}. If the solution in a cell is discarded, it is recomputed relying on a first order finite volume method applied to a fine sub-grid generated within each troubled cell. A second order TVD scheme has been used as limiter in~\cite{DGLimiter3,ALEDG,SonntagDG}, while higher order ADER-WENO subcell finite volume limiters are presented in~\cite{DGLimiter1,DGLimiter2,DGCWENO,Rannabauer,DeLaRosaMunzDGMHD}.

\medskip

We refer to the aforementioned references for an exhaustive description of the \textit{a posteriori} finite volume subcell limiter. Here, for the sake of clarity, we briefly recall the main concepts and we underline the differences introduced for dealing with moving Voronoi elements and topology changes.

First, using the notation adopted in~\cite{ALEDG}, the numerical solution computed so far is assumed to be a \textit{candidate} solution and is denoted with $\mathbf{u}^{n+1,*}_{h}(\x,t^{n+1})$. Then, we define a sub-triangulation of $P_i^n$ made of a set of non-overlapping so called \textit{small sub-triangles}. Consequently, each control volume $C_i^n$ is split into sub-triangular prisms, called \textit{small sub-volumes}, as follows.
\begin{itemize}
	\item For $N=1$ we consider a total number of small sub-triangles $\mathcal{S}_i$ which is equal to $N_{C_i}^n$, i.e. $\mathcal{S}_i=N_{C_i}^n$. The small sub-triangles are given by $T_{i_j}^n$ and the associated small sub-volumes are $sC_{i_j}^n$, as defined in Section~\ref{ssec.SpaceTimeConnection}.	
	\item If a topology change happens with $N=1$, i.e. $\mathcal{V}(P_i^n) \ne \mathcal{V}(P_i^{n+1})$, degenerate small sub-triangles/sub-volumes are considered as well, thus including also sub-triangles which can be given by a line.  
	\item For $N\ge2$ we further subdivide each $T_{i_j}^n$ into $N^2$ small sub-triangles,
	which are defined through the sub-nodes provided by standard nodes of classical high order conforming finite elements on triangular meshes. In this way, a total number of $\mathcal{S}_i=N_{C_i}^n \cdot N^2$ small sub-triangles is taken into account. The splitting of $sC_{i_j}^n$ is consequently defined.
	\item Even in the case $N\ge2$, degenerate sub-triangles/sub-volumes are counted if a topology change happens, i.e. $\mathcal{V}(P_i^n) \ne \mathcal{V}(P_i^{n+1})$. 
    This results in small sub-triangles which may be given by a portion of a line.
\end{itemize} 
We denote each small sub-triangle of $P_i^n$ with $s_{i,\alpha}^n$, where $\alpha \in [1, \mathcal{S}_i]$. 
Next, we define the corresponding subcell average of the numerical solution at time $t^n$  
\begin{equation}
\mathbf{v}_{i,\alpha}^n(\x,t^n) = \frac{1}{|s_{i,\alpha}^n|} \int_{s_{i,\alpha}^n} \mathbf{u}_{h}^n(\x,t^n) \, d\x = \frac{1}{|s_{i,\alpha}^n|} \int_{s_{i,\alpha}^n} \phi_\ell(\x) \, d\x \, \hat{\mathbf{u}}^{n}_{l}:=\mathcal{P}(\mathbf{u}_h^n) \qquad \forall \alpha \in [1,\mathcal{S}_i],
\label{eqn.vh}
\end{equation}
where $|s_{i,\alpha}^n|$ denotes the volume of subcell $s_{i,\alpha}^n$ of element $P_i^n$ and the definition $\mathcal{P}(\mathbf{u}_h)$ is the $L_2$ projection operator. We fix also the \textit{candidate} subcell average of the numerical solution at time $t^{n+1}$ as $\mathbf{v}_{i,\alpha}^{n+1,*}(\x,t^{n+1}) = \mathcal{P}(\mathbf{u}_h^{n+1,*})$.

Now, we mark the troubled cells. The candidate solution $\mathbf{v}^{n+1,*}_h(\x,t^{n+1})$ is checked against a set of detection criteria.
Here we follow the criteria described in~\cite{ALEDG}, \RIIcolor{however also other and more elaborate choices could be considered, see for example the recent work of Guermond \textit{et al.}~\cite{guermond2018second} on invariant domain preserving schemes for hyperbolic systems.}

Thus, our first criterion is the requirement that the computed solution is physically acceptable, i.e. belongs to the phase space of the conservation law being solved. For instance, if the compressible Euler equations for gas dynamics are considered, density and pressure should be positive and in practice we require that they are greater than a prescribed tolerance $\epsilon=10^{-12}$. 
Then, a relaxed discrete maximum principle (DMP) is applied, hence we verify
\begin{equation}
\min_{m \in \mathcal{V}(C_i^n)} \left( \, \min_{\beta\in[1,\mathcal{S}_m]} (\mathbf{v}_{m,\beta\,}^n) \right) 
- \delta \leq \ \mathbf{v}^{n+1,*}_{i,\alpha} \ \leq 
\max_{m \in \mathcal{V}(C_i^n)} \left( \, \max_{\beta\in[1,\mathcal{S}_m]} (\mathbf{v}_{m,\beta\,}^n) \right) 
+ \delta \qquad \forall \alpha \in [1,\mathcal{S}_i],
\label{eqn.RDMP}
\end{equation}
where $\delta$ is a parameter which, according to~\cite{ALEDG, DGLimiter1, DGLimiter2}, reads
\begin{equation}
\delta = \max \, \Biggl( \delta_0\, , \ \epsilon \cdot 
\biggl[
\max_{m \in \mathcal{V}(C_i^n)} \biggl(\, \max_{\beta\in[1,\mathcal{S}_m]} (\mathbf{v}_{m,\beta\,}^n) \biggr) -
\min_{m \in \mathcal{V}(C_i^n)} \biggl(\, \min_{\beta\in[1,\mathcal{S}_m]} (\mathbf{v}_{m,\beta\,}^n) \biggr) 
\biggr]
\Biggr),
\label{eqn.deltaRDMP}
\end{equation}
with $\delta_0=10^{-4}$ and $\epsilon=10^{-3}$. 

If a cell fulfills the detection criteria in \textit{all} its subcells, then the cell is marked as \textit{good}, otherwise the cell is \textit{troubled}. We emphasize that this step is performed independently in each element and thus the projection $\mathbf{v}^*_h(\x,t^{n+1})$ does not need to be retained after the cell has been assigned its marker.

The following step consists in re-computing the solution \textit{only} in the troubled cells with a first order FV scheme, applied in each small sub-triangle/sub-volume, that evolves the cell averages $\mathbf{v}_{i,\alpha}^n$ in order to obtain $\mathbf{v}_{i,\alpha}^{n+1}$. 

We do not report the details on the very well-known first order ALE-FV scheme, but we add some remarks on flux computation at the space--time lateral surfaces of each $s_i^n$. 
i) The same numerical flux function, i.e.~\eqref{eq.rusanov} or~\eqref{eqn.osher}, used in the rest of the scheme is adopted here as well.
ii) The employed quadrature rule is a simple mid-point rule that makes use of the space--time barycenters $g_i^n$ of the space--time lateral faces of the sub-volume. 
iii) The normal vectors are also computed at $g_i^n$. 
iv) Referring to~\eqref{eqn.fluxFunctionGeneric}, when computing the flux between the sub-volume $\alpha$ of $C_i^n$ and the neighboring sub-volume $\beta$ (of $C_i^n$ or of any other $C_{i_j}^n$), boundary data are simply given by $\q_h^{n,-}=\v_{i,\alpha}^n$ and $\q_h^{n,+}=\v_{i/ i_j,\,\beta}^n\,$.
v) If instead the neighbor is \textit{not} a troubled Voronoi element $C_{i_j}^n$ (which thus has not been sub-triangulated), then $\q_h^{n,-}=\v_{i,\alpha}^n$ and $\q_h^{n,+} = \q_h^n|_{C_{i_j}^n}(\g_{i_j}^n)$.

A first order finite volume scheme always provides a valid solution, hence $\mathbf{v}_{i,\alpha}^{n+1}$ is acceptable.
Moreover, since the FV scheme is not directly applied to the Voronoi element but to each of its sub-triangles, the sub-mesh resolution does not completely spoil the solution of the DG scheme. Nevertheless, the method does not maintain the formal order of accuracy of the $P_NP_M$ scheme, but it is only used and activated across shock waves and strong discontinuities.
Note also that for a troubled cell the mesh motion is not recomputed because it has been fixed using only information coming from space at time $t^n$, which are, as such, not affected by any problem.

Finally, the DG polynomial for the Voronoi cell $P_i^{n+1}$ is recovered from the robust and stable solution on the sub-grid level $\mathbf{v}_{i,\alpha}^{n+1}$ by applying the reconstruction operator $\mathcal{R}(\mathbf{v}_{i,\alpha}^{n+1}(\x,t^n))$, that is
\begin{equation}
\int_{S_{i,\alpha}^n} \mathbf{u}_{h}^{n+1}(\x,t^{n+1}) \, d\x = \int_{S_{i,\alpha}^n} \mathbf{v}_{i,\alpha}^{n+1}(\x,t^n) \, d\x :=\mathcal{R}(\mathbf{v}_{i,\alpha}^{n+1}(\x,t^n))  \qquad \forall \alpha \in [1,\mathcal{S}_i].
\label{eqn.intRec}
\end{equation}
The reconstruction is imposed to be \textit{conservative} on the main cell $P_i^{\RIIcolor{n+1}}$, hence yielding the additional linear constraint
\begin{equation}
\int_{P_i^{\RIIcolor{n+1}}} \mathbf{u}_{h}(\x,t^{n+1}) \, d\x = \int_{P_i^{\RIIcolor{n+1}}} \mathbf{v}_{h}(\x,t^{n+1}) \, d\x.
\label{eqn.LSQ}
\end{equation}
Moreover, the projection operator $\mathcal{P}$ in~\eqref{eqn.vh} and the reconstruction operator $\mathcal{R}$ in~\eqref{eqn.intRec} satisfy the property $\mathcal{P} \cdot \mathcal{R}=\mathcal{I}$, with $\mathcal{I}$ being the identity operator. \RIIcolor{The reconstruction operator \eqref{eqn.intRec}-\eqref{eqn.LSQ} might still lead to an oscillatory solution, since it is based on a linear unlimited least squares technique. If this is the case, the cell $P_i^{n+1}$ will be detected as troubled while performing the time marching at the next time level $t^{n+2}$, therefore the finite volume subcell limiter will be used again in that cell. However, the projection operator $\mathcal{P}$ \eqref{eqn.vh} applied to the oscillatory reconstructed solution in cell $P_i^{n+1}$ does not produce an acceptable set of subcell averages $\mathbf{v}_{i,\alpha}$, so that one could not guarantee that the solution coming from the subcell finite volume scheme is valid at the next time level $t^{n+2}$. In order to overcome this issue, the subcell averages $\mathbf{v}_{i,\alpha}^{n+1}$ are always kept in memory till the cell $P_i$ is again marked as valid. In other words, if a cell is detected to be troubled for the second time step in a row, then the starting subcell averages are not obtained via the projection operator, but they are given by the solution of the subcell finite volume scheme at the previous time step.}

If a cell $C_i^n$ is acceptable but has at least one troubled neighbor cell $C_{i_j}^n$ in its $\mathcal{V}(C_i^n)$, 
then we cannot accept its candidate solution $\mathbf{u}^{n+1,*}_{h}(\x,t^{n+1})$ because the scheme would become nonconservative.
Indeed, at the common space--time lateral surface $\partial C_{i_j}^n$, 
the flux computed from $C_i^n$ would be obtained through the DG scheme (i.e. high order predictor and high order corrector), while the one coming from the troubled neighbor $C_{i_j}^n$ would be updated using the first order FV scheme.
Thus, the DG solution in these cells is recomputed in a \textit{mixed way}: the volume integral and the surface integrals on good faces are kept, while the numerical flux across the troubled faces is always provided by the first order limiter.


\paragraph{Neighborhood of a sliver element}

\hfill \break 

At the subcell level, the difficulties associated with degenerate small sub-volumes are the same stated at the 
end of Section~\ref{sssec.NumScheme_sliver_Flux} for degenerate big elements: how to impose
an initial condition for cells with zero area at $t^n$ and how not to lose 
any contribution computed through elements with zero area at $t^{n+1}$. In order to activate and apply the limiter, the following strategy is proposed. 

Firstly, the sliver elements are \textit{not} sub-triangulated. If one neighbor of a sliver $S_i^n$ is troubled, we mark as troubled also the remaining three neighbors. Among the four neighbors of $S_i^n$, we select the one with the biggest volume which we call $C_\mathrm{big}^n$. 

Next, we need to provide the values $\q_h^{+,-}$ when computing the fluxes~\eqref{eqn.fluxFunctionGeneric}.
\begin{itemize} 
\item For a degenerate $s_{i,\alpha}^n$ with zero area at $t^n$ we take the value obtained by evaluating $\u_h^{n}$ at the mid point of $s_{i,\alpha}^n|_{t^n}$ (this value is well defined because $s_{i,\alpha}^n \subset P_i^n$ and so  $\u_h^{n}$ is continuous).
\item For a sliver element $S_i^n$ we take the value obtained by evaluating $\u_h^{n}$ of $C_\mathrm{big}$ at the mid point of $S_i^n|_{t^n}$; this arbitrary choice is justified by the fact that here we simply employ a first order method, which is stable even if the sliver elements are neglected (see~\cite{Springel}).
\end{itemize}

Finally, we need to redistribute the fluxes computed across the degenerate elements when they disappear at $t^{n+1}$.
\begin{itemize} 
	\item For a degenerate $s_{i,\alpha}^n$ with zero area at $t^{n+1}$ we assign the sum of the fluxes computed through its space--time lateral surfaces to the closest $s_{i,\beta}^n$ that is not degenerate at $t^{n+1}$ (the concept of \textit{closest} is uniquely fixed through a specific numbering of the sub-volumes). 
	\item For a sliver element $S_i^n$ we assign its fluxes to $C_\mathrm{big}^n$. 
\end{itemize}

Besides, we remark that the space--time geometry definition in itself does not pose any problem: indeed, the configuration of big elements has already been fixed in Section~\ref{ssec.SpaceTimeConnection} and the subdivision has been deduced just above.
Therefore, quadrature formulae, normal vectors and bilinear mapping are always well defined.

\subsection{MOOD approach to verify the consistency of the space-time connectivity}
\label{ssec.MOOD}

As already stated at the end of Section~\ref{ssec.SpaceTimeConnection}, it may not always be possible to connect two consecutive meshes in a consistent way if the associated topology changes are too strong.
However, these situations are immediately detected at the beginning of the new time step, when the space--time connectivity is built. 
Indeed, if i) the set $\mathcal{V}(C_i^n)$ cannot be ordered consistently with both the order of $\mathcal{V}(P_i^n)$ and $\mathcal{V}(P_i^{n+1})$, or if ii) more than three sliver elements are necessary to complete a path between elements which are neighbors at one time level but not at the previous or at the next one, or if iii) \RIIcolor{the path involving the minimum number of slivers is not unique}, then the algorithm
detects the problem. To overcome it, the current time step is simply restarted with a smaller time step size $\Delta t$ (reduced by a factor of $2$ for example). Eventually, more restarts are needed, until the connection between the two meshes is coherent.

Since the mesh generation and the connectivity construction are not expensive, the performances of the algorithm are not negatively influenced by this additional MOOD-type procedure (which applies \textit{before} the evolution in time). 
However, future work will consider the possibility of remeshing only \textit{locally}, in the neighborhood of a connectivity problem without reducing the time step size or more sophisticated mesh optimization algorithms. 
We \RIIIcolor{emphasize} that such problems are encountered very rarely, \RIVcolor{see Tables~\ref{tab.GCLTest_CPUtimes}, \ref{tab.explosion_percentage}, \ref{tab.RT_sliver_percentage} and the notes of Section~\ref{test.MHDVortex} for some statistics}.

\section{Numerical results}
\label{sec.Results}

The numerical results presented in this section will show the following properties of our new algorithms.
\begin{enumerate}
	\item[\textbf{i.}] Our method 
    has been implemented as general purpose code, 
      in the sense that any kind of 
        hyperbolic system cast in the form~\eqref{eq.generalform} can be easily studied: 
        for this reason we test it on several models, namely the standard Euler equations 
        of gas dynamics (Section~\ref{ssec.EulerEq}), the Euler equations with gravity source 
        term (Section~\ref{ssec.EulerEq+S}) and the ideal magnetohydrodynamics (MHD) system (Section~\ref{ssec.MHDEq}).
	\smallskip
	\item[\textbf{ii.}] Next, we show the capability of our scheme in maintaining a high quality mesh 
	for very long times, 
	even in the case of strong shear flows and vortices, 
	thanks to its high \textit{robustness and adaptability} to complex flow patterns, 
	see Sections~\ref{test.ShuVortex} and~\ref{test.MHDVortex}. In particular, we show that 
    we can preserve the accuracy of the high order trajectory integration of generator points, 
    as well as sharply fit strong shocks in Section~\ref{test.Sedov}.
	\smallskip
	\item[\textbf{iii.}] Then, we compute numerically the \textit{order of convergence} of both Finite Volume and Discontinuous Galerkin schemes for two different test problems, see  Tables~\ref{tab.orderOfconvergenceFV_shu},~\ref{tab.orderOfconvergenceDG_shu},~\ref{tab.orderOfconvergenceFV_MHD} and~\ref{tab.orderOfconvergenceDG_MHD}.
	\item[\textbf{iv.}] 
	Moreover, for \textit{all} the presented test cases we have numerically verified that mass and volume conservation is respected up to \textit{machine precision} at \textit{any} time step, and that the same holds true for the GCL condition on \textit{each} element, \RIIIcolor{even when topology changes occur.}
	\RIIIcolor{In addition, to provide even more evidence on the fact that the GCL condition is satisfied by construction, we refer to the set of test cases of Section~\ref{test.GCLTest}, where constant states are preserved up to machine precision for long times over moving meshes where topology changes regularly occur}.
	\smallskip
	\item[\textbf{v.}] Finally, we study some more complicated test problems (see Sections~\ref{test.ExpPb},~\ref{test.Sedov},~\ref{test.triplepoint} and~\ref{test.MHDRotor}) to show the \textit{robustness} of our method, concerning both the mesh quality in presence of arbitrary and strong velocity fields as well as the consistency/stability of our high order schemes. In particular, we test the \textit{a posteriori} sub--cell finite volume \textit{limiter} used to stabilize the DG scheme that indeed avoids undesirable oscillations by activating only where needed (see Figures~\ref{fig.ExpPb_P2P2_DG} and~\ref{fig:sedovp1p180x800-0}).	
\end{enumerate}

The great variety of the presented tests is intended to show both the wide range of applicability of the proposed high order ALE scheme on moving Voronoi meshes with topology changes and its level of novelty with respect to the state of the art.

\subsection{Euler equations of gasdynamics}
\label{ssec.EulerEq}

A well-known example of a hyperbolic system of the form~\eqref{eq.generalform} is given by the homogeneous Euler equations of compressible gas dynamics with  
\begin{equation}
\label{eulerTerms}
\Q = \left( \begin{array}{c} \rho   \\ \rho u  \\ \rho v  \\ \rho E \end{array} \right), \quad
\mathbf{F} = \left( \begin{array}{ccc}  \rho u       & \rho v        \\ 
\rho u^2 + p & \rho u v          \\
\rho u v     & \rho v^2 + p      \\ 
u(\rho E + p) & v(\rho E + p)   
\end{array} \right), \qquad \S = 0.  
\end{equation}
The vector of conserved variables $\Q$ involves the fluid density $\rho$, the momentum density vector $\rho \v=(\rho u, \rho v)$ and the total energy density $\rho E$. 
The fluid pressure $p$ is related to conservative quantities $\Q$ using the equation of state for an ideal 
gas    
\begin{equation}
\label{eqn.eos} 
p = (\gamma-1)\left(\rho E - \frac{1}{2} \rho \mathbf{v}^2 \right), 
\end{equation}
where $\gamma$ is the ratio of specific heats so that the speed of sound takes the form $c=\sqrt{\frac{\gamma p}{\rho}}$. Where not otherwise specified we employ the Rusanov-type ALE flux~\eqref{eq.rusanov} as numerical flux function and we move the generator points using the local fluid velocity obtained from $\w_h^n$ (see Section~\ref{sec.MeshEvolution}). Furthermore, we set $\gamma=1.4$.

\subsubsection{Isentropic vortex}
\label{test.ShuVortex}

To verify the order of convergence of the proposed ALE FV-DG scheme we consider 
a smooth isentropic vortex flow according to~\cite{HuShuVortex1999}. The initial computational domain is the square $\Omega=[0;10]\times[0;10]$ with wall  boundary conditions set everywhere. 
The initial condition is given by some perturbations $\delta$ that are superimposed onto a homogeneous background field $\Q_0=(\rho,u,v,p)=(1,0,0,1)$, assuming that the entropy perturbation is zero, i.e. $\delta S= 0$. The perturbations for density and pressure are
\begin{equation}
\label{rhopressDelta}
\delta \rho = (1+\delta T)^{\frac{1}{\gamma-1}}-1, \quad \delta p = (1+\delta T)^{\frac{\gamma}{\gamma-1}}-1, 
\end{equation}
with the temperature fluctuation $\delta T = -\frac{(\gamma-1)\epsilon^2}{8\gamma\pi^2}e^{1-r^2}$ and the vortex strength is $\epsilon=5$.
The velocity field is affected by the following perturbations
\begin{equation}
\label{ShuVortDelta}
\left(\begin{array}{c} \delta u \\ \delta v  \end{array}\right) = \frac{\epsilon}{2\pi}e^{\frac{1-r^2}{2}} \left(\begin{array}{c} -(y-5) \\ \phantom{-}(x-5)  \end{array}\right).
\end{equation}
This is a stationary equilibrium of the system so the exact solution coincides with the initial condition at any time.

\begin{figure*}[!p]
\centering
\includegraphics[width=0.20\linewidth]{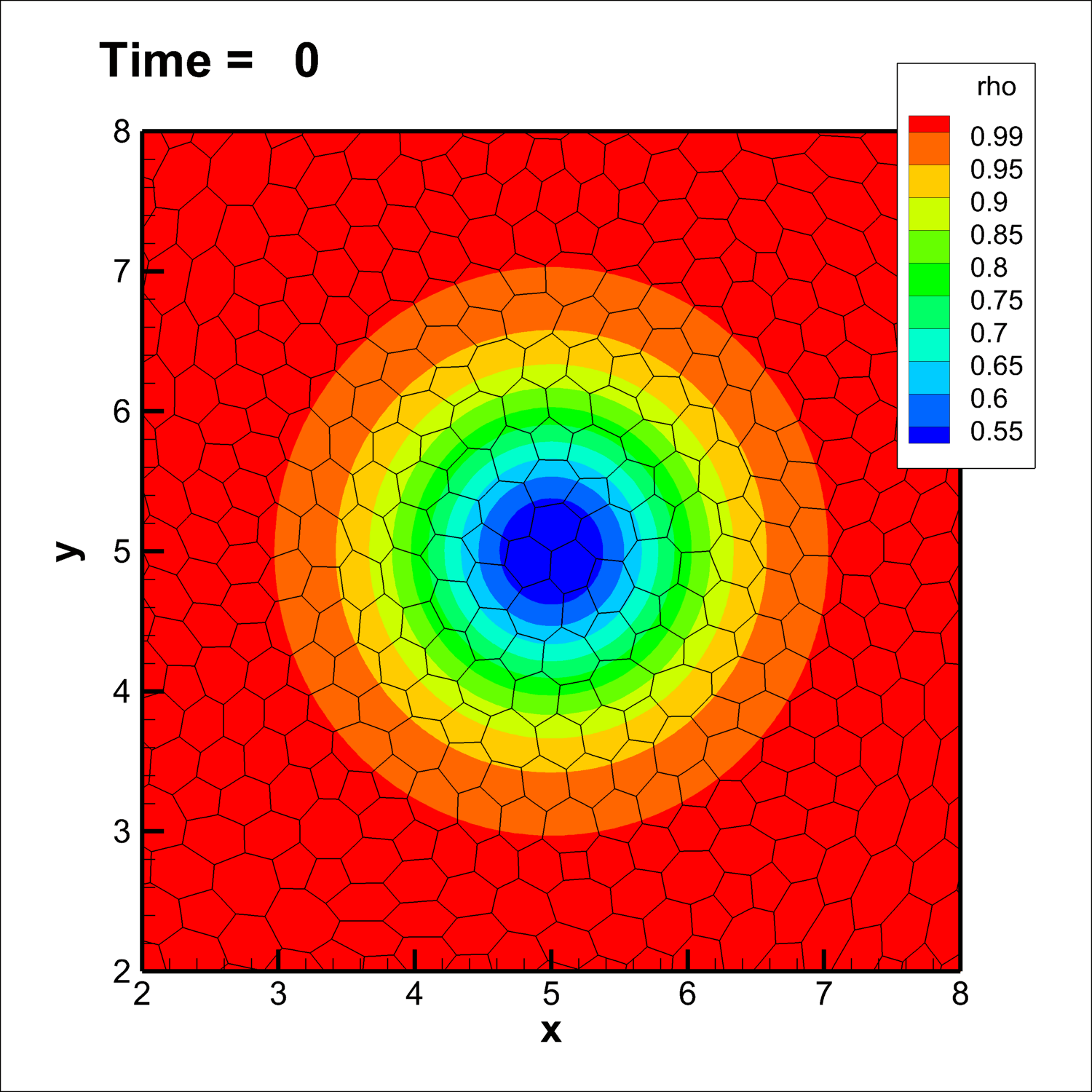}%
\includegraphics[width=0.20\linewidth]{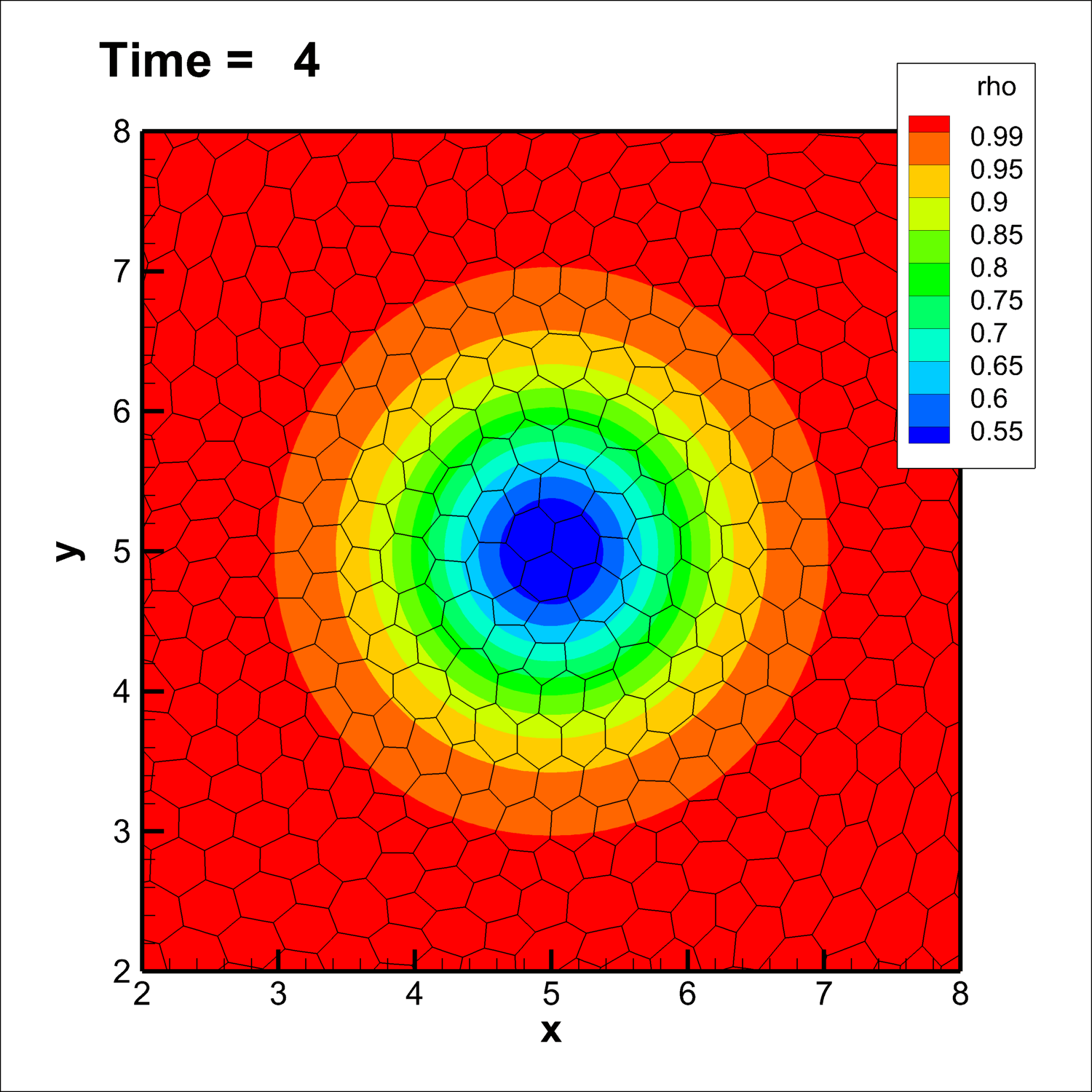}%
\includegraphics[width=0.20\linewidth]{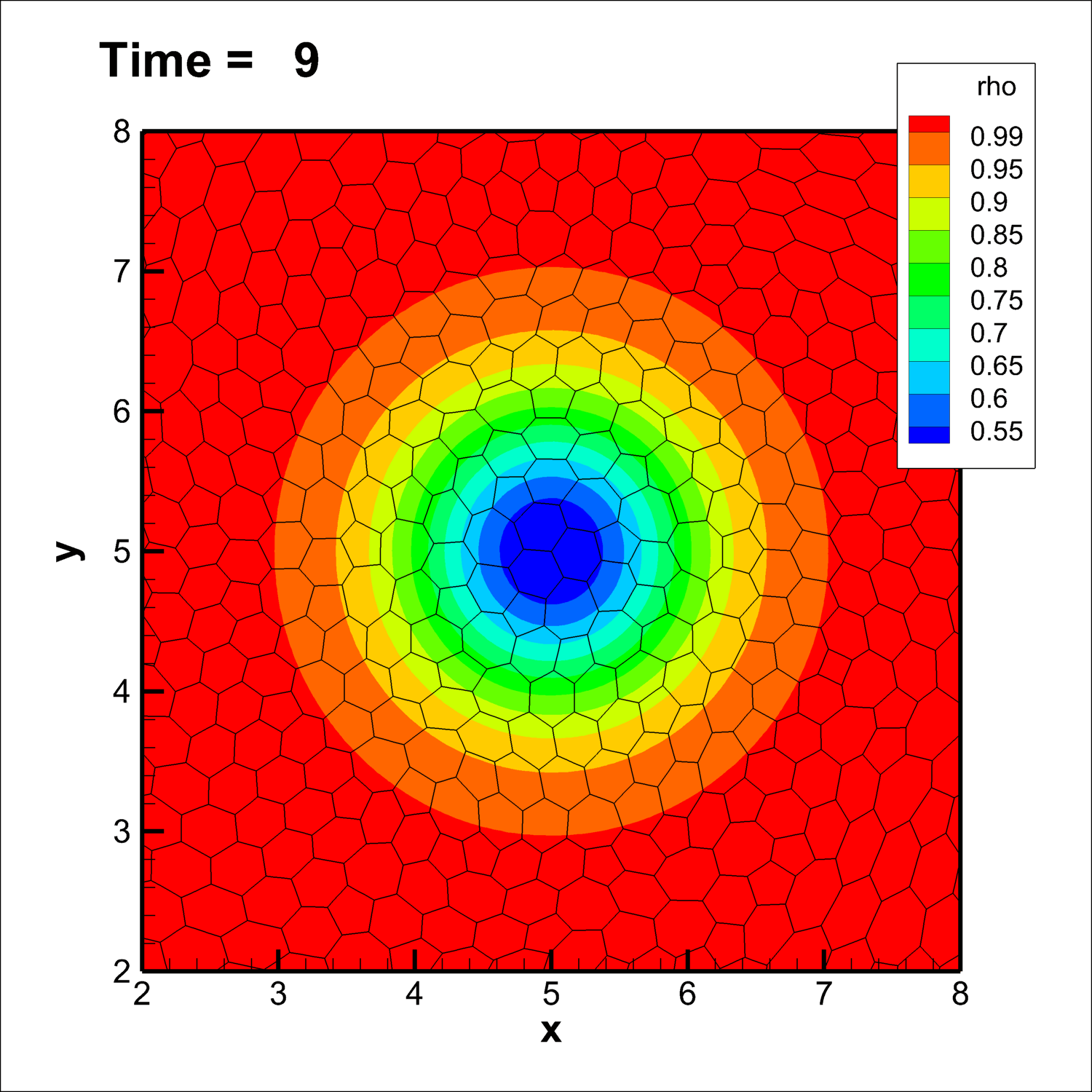}%
\includegraphics[width=0.20\linewidth]{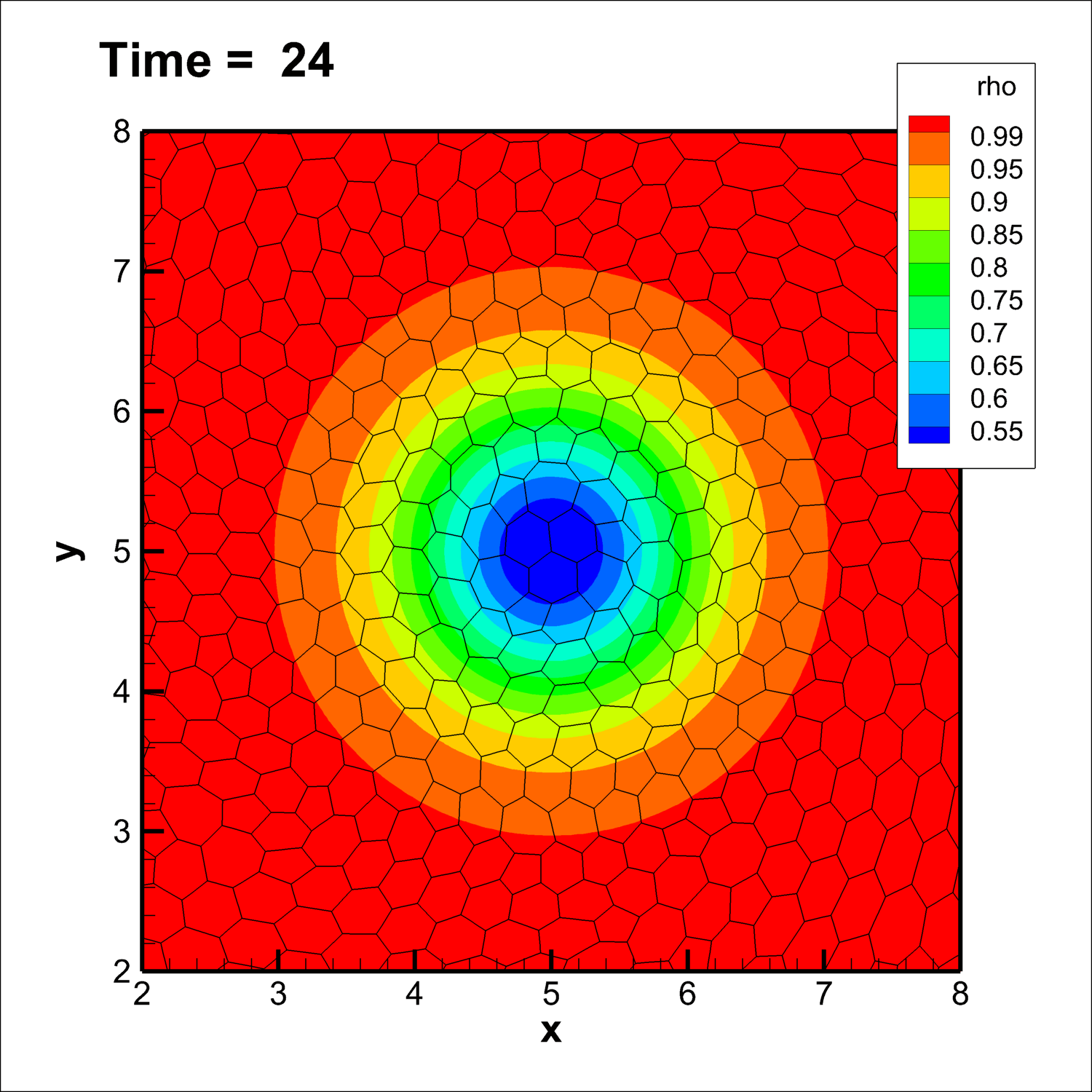}%
\includegraphics[width=0.20\linewidth]{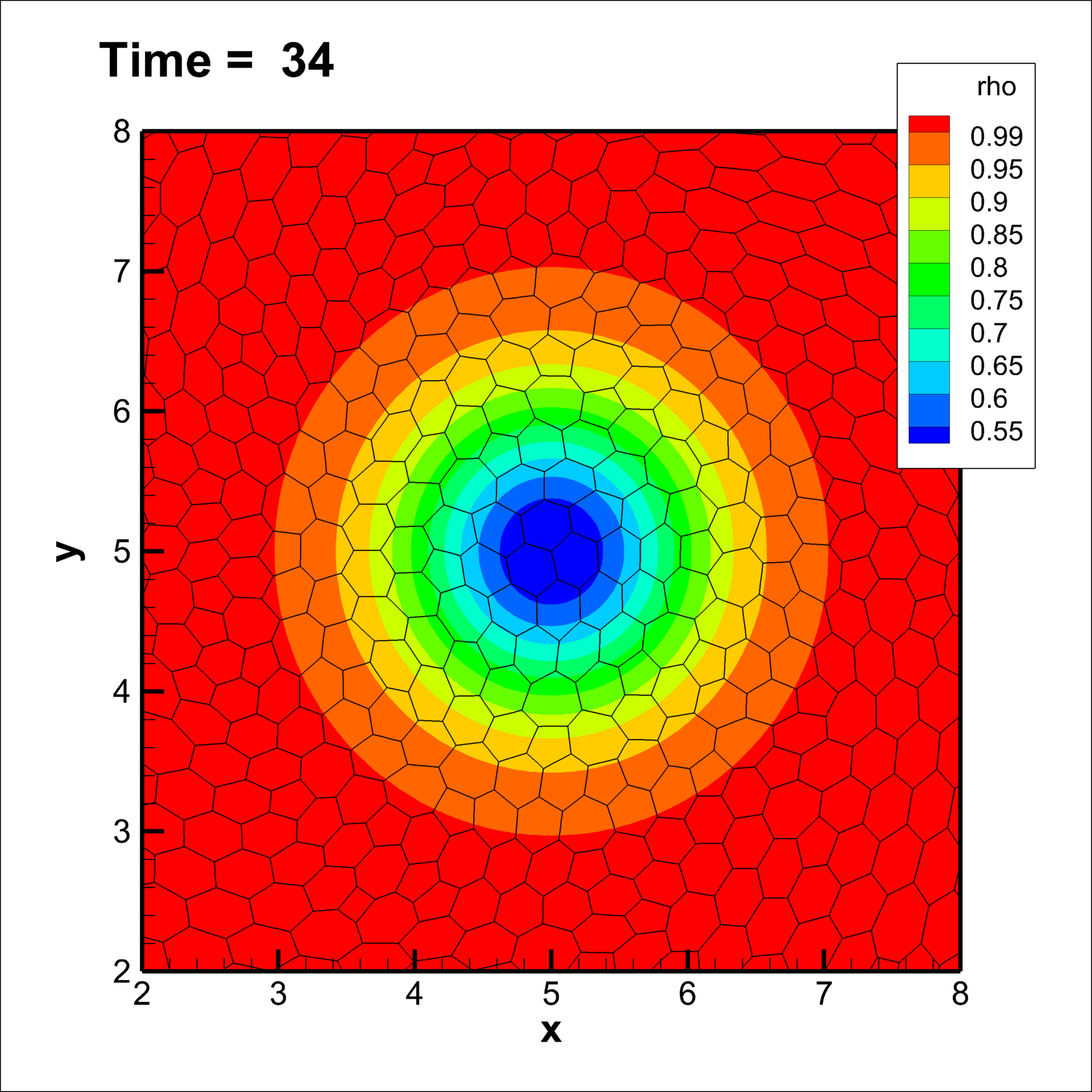}\\[-1.55pt]
\includegraphics[width=0.20\linewidth]{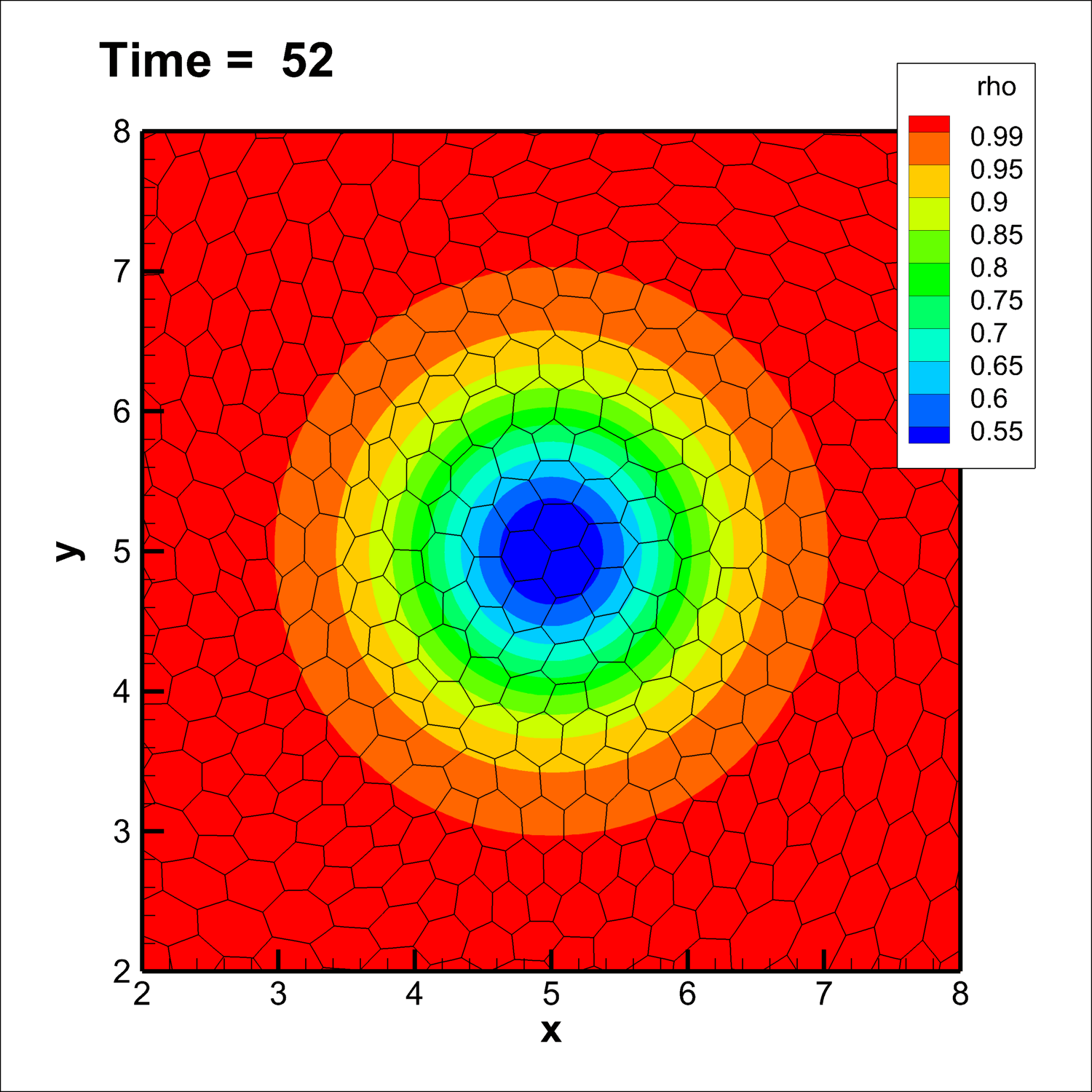}%
\includegraphics[width=0.20\linewidth]{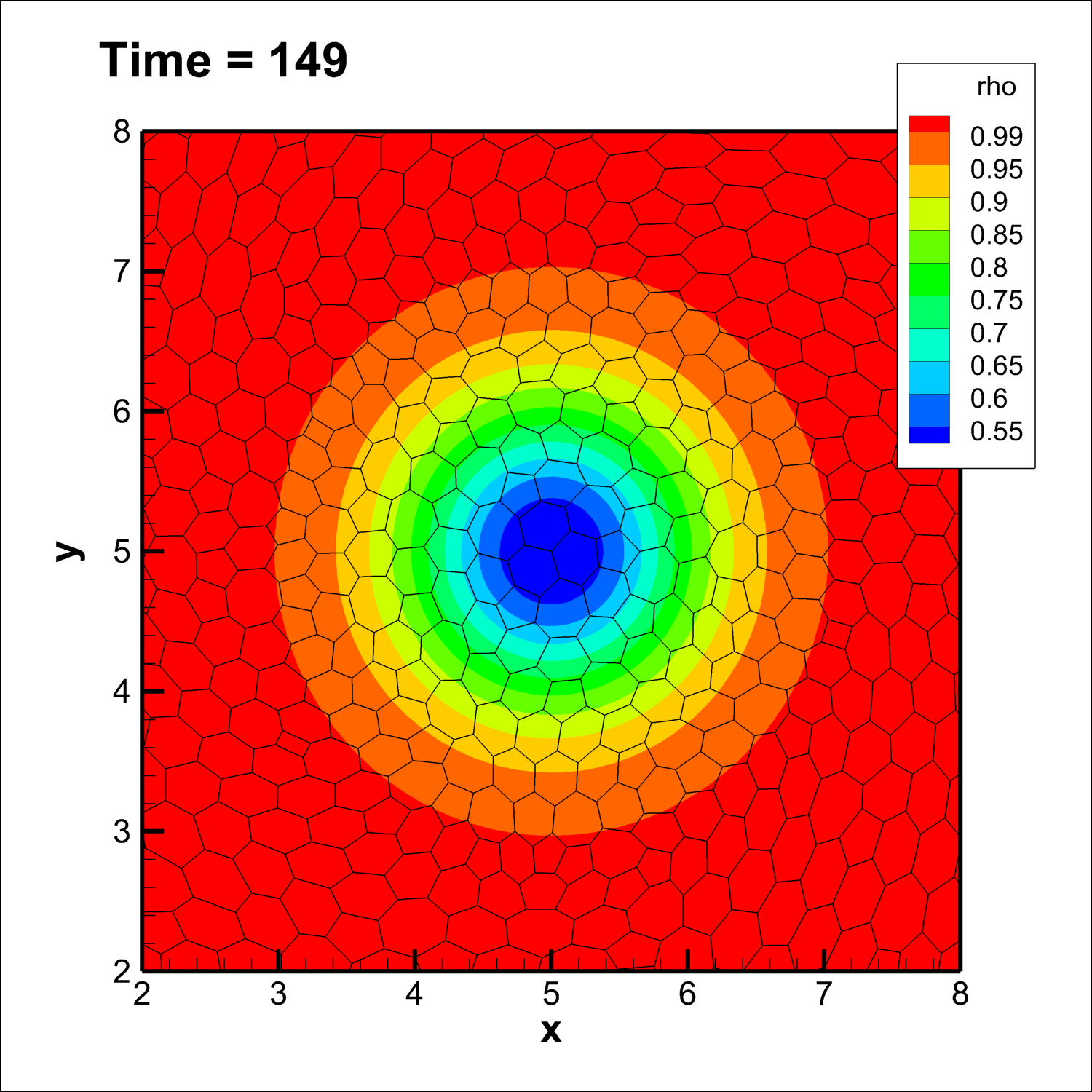}%
\includegraphics[width=0.20\linewidth]{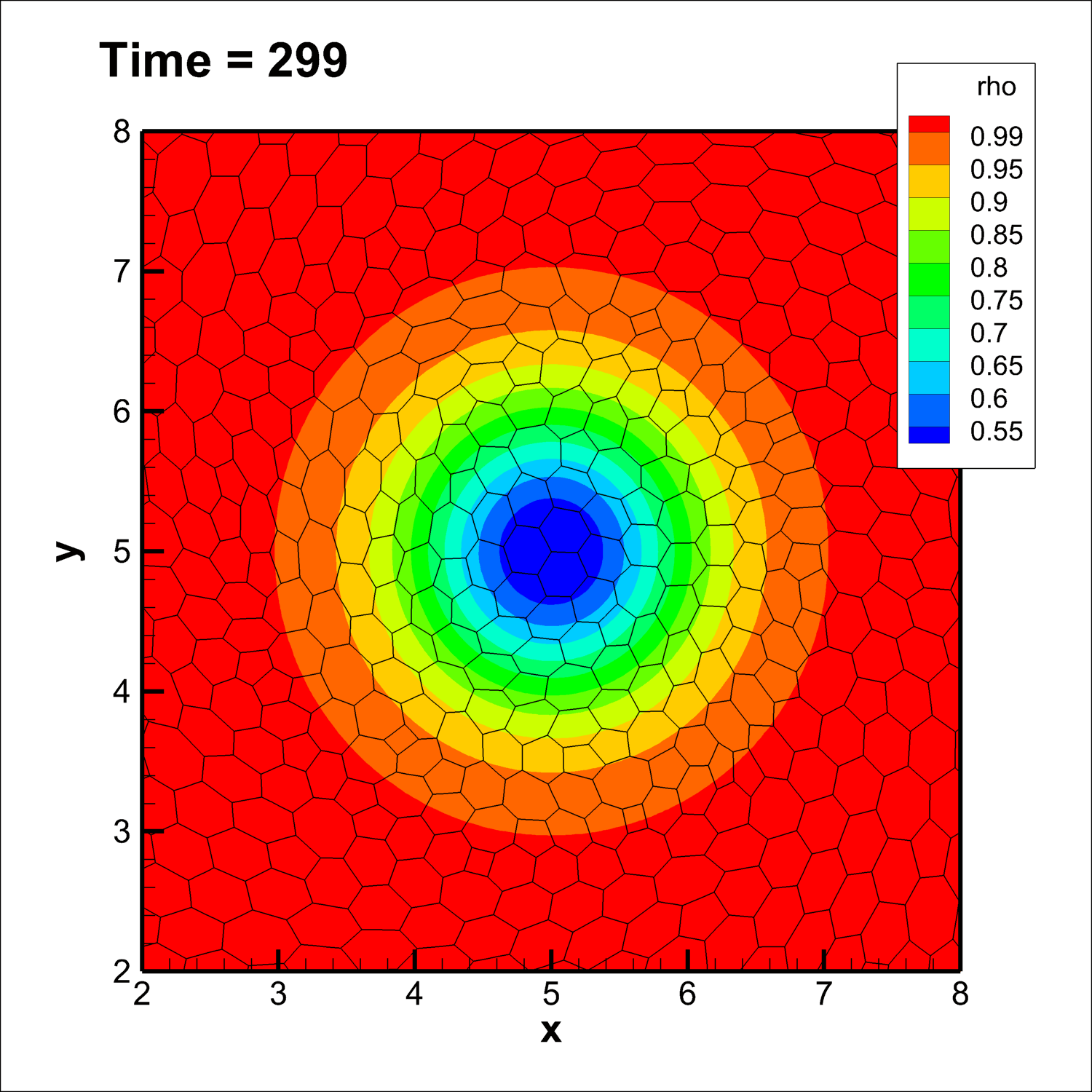}%
\includegraphics[width=0.20\linewidth]{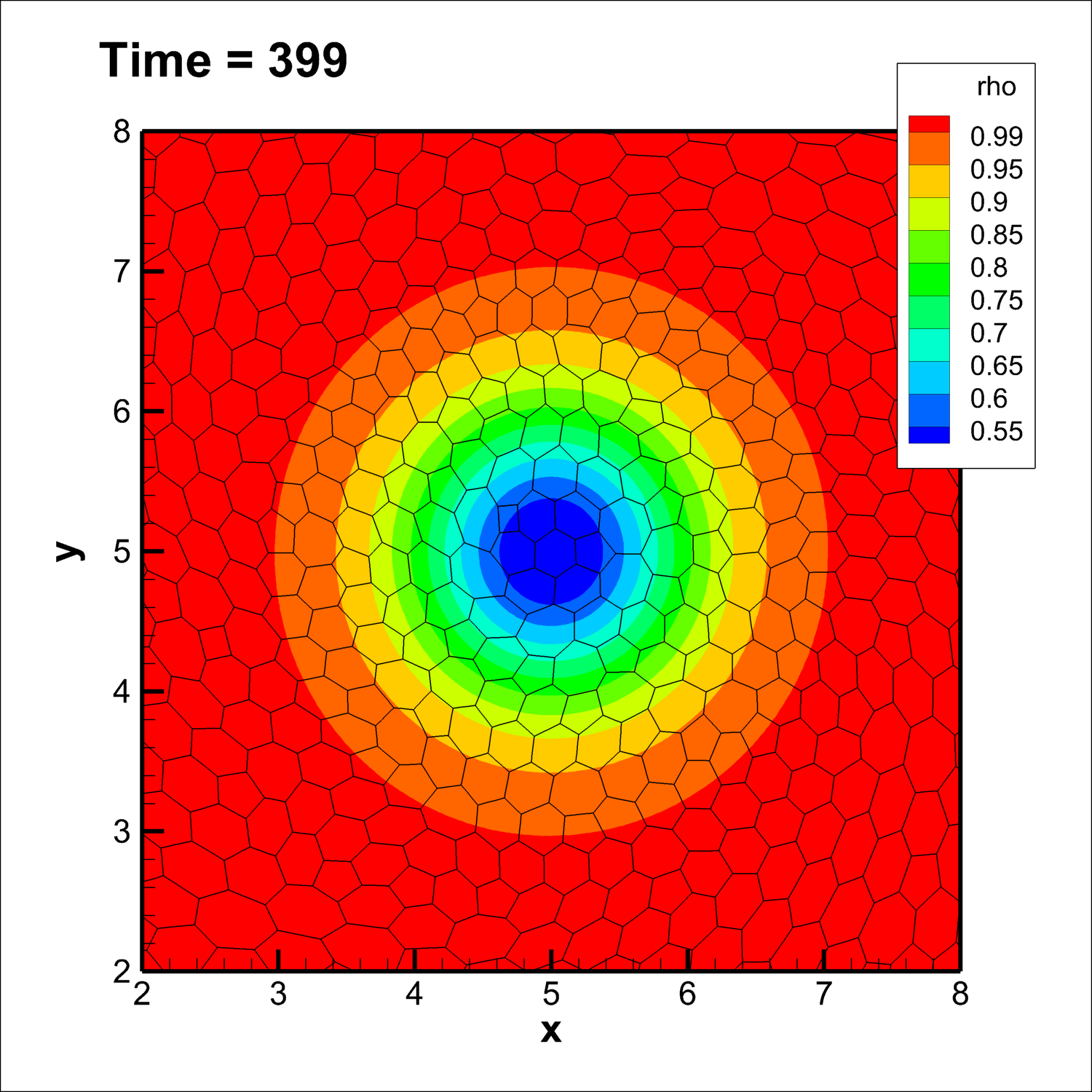}%
\includegraphics[width=0.20\linewidth]{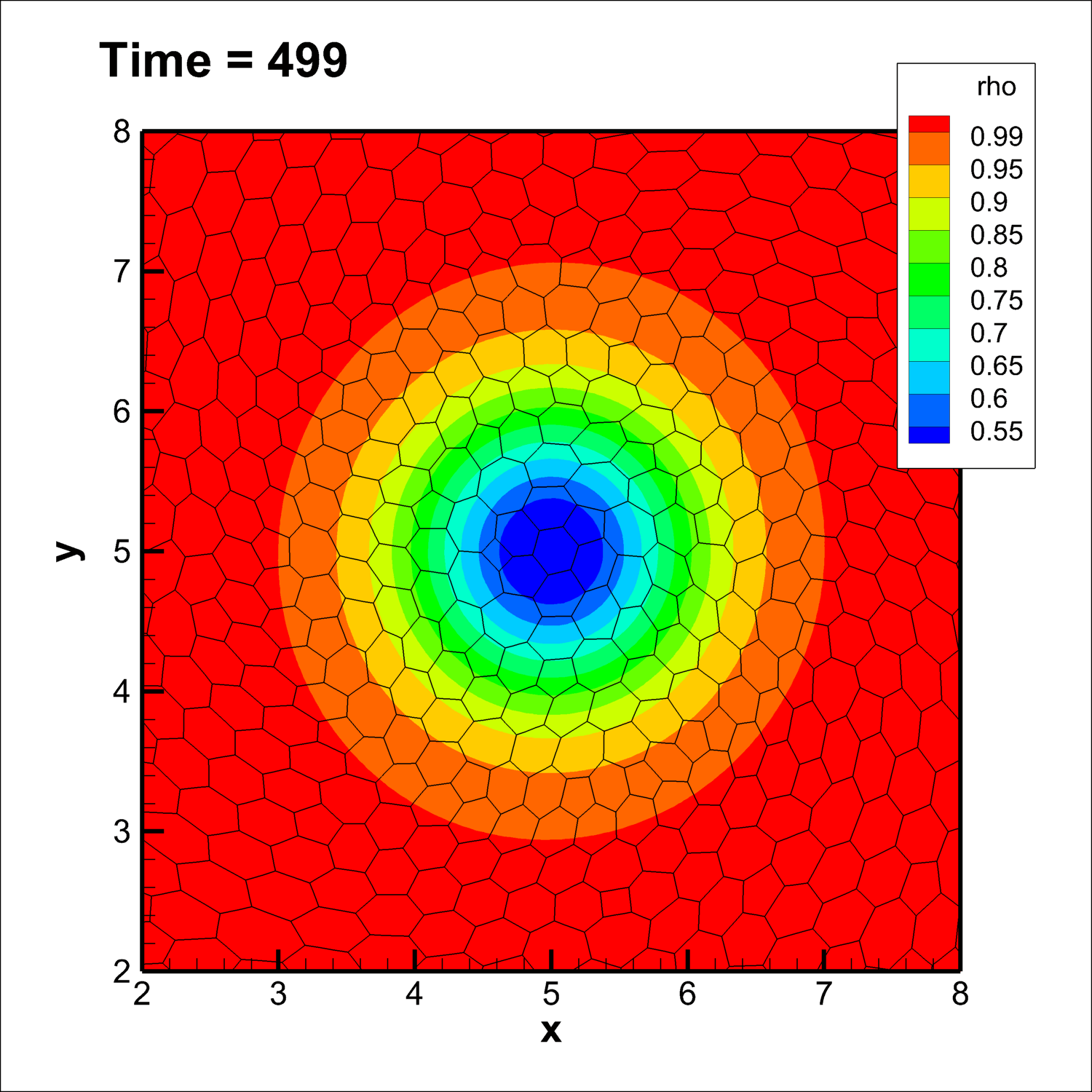}\\[-1.55pt]
\includegraphics[width=0.20\linewidth]{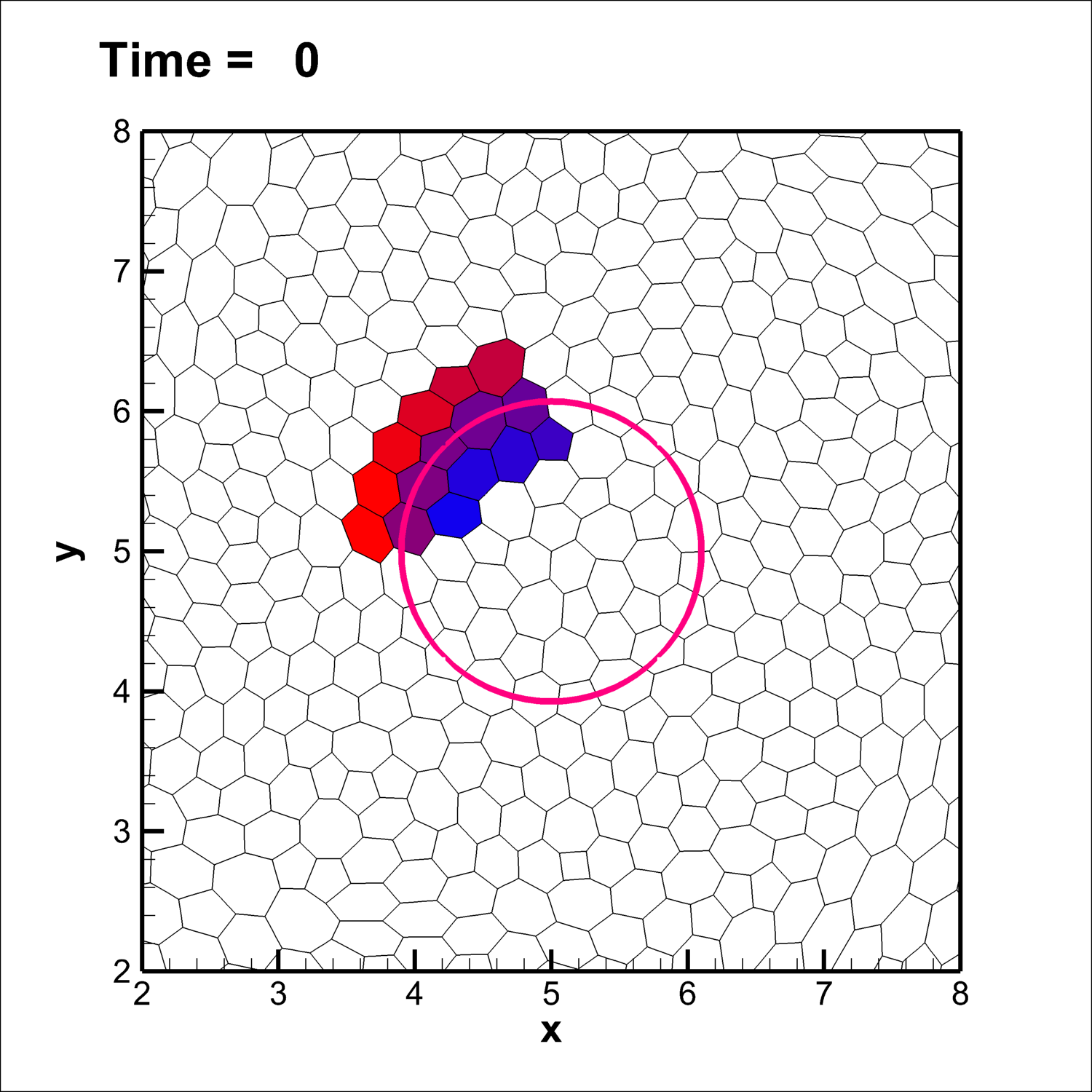}%
\includegraphics[width=0.20\linewidth]{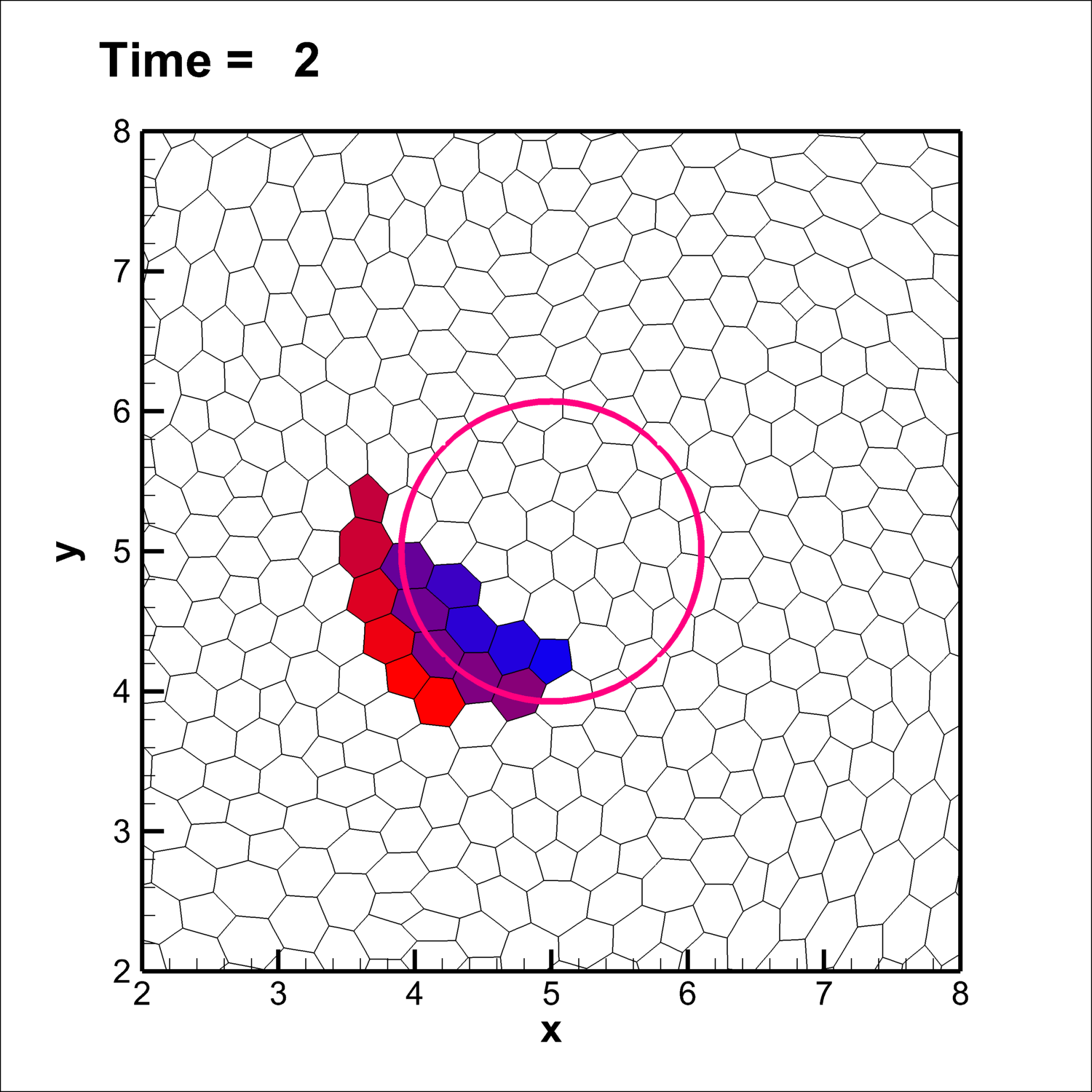}%
\includegraphics[width=0.20\linewidth]{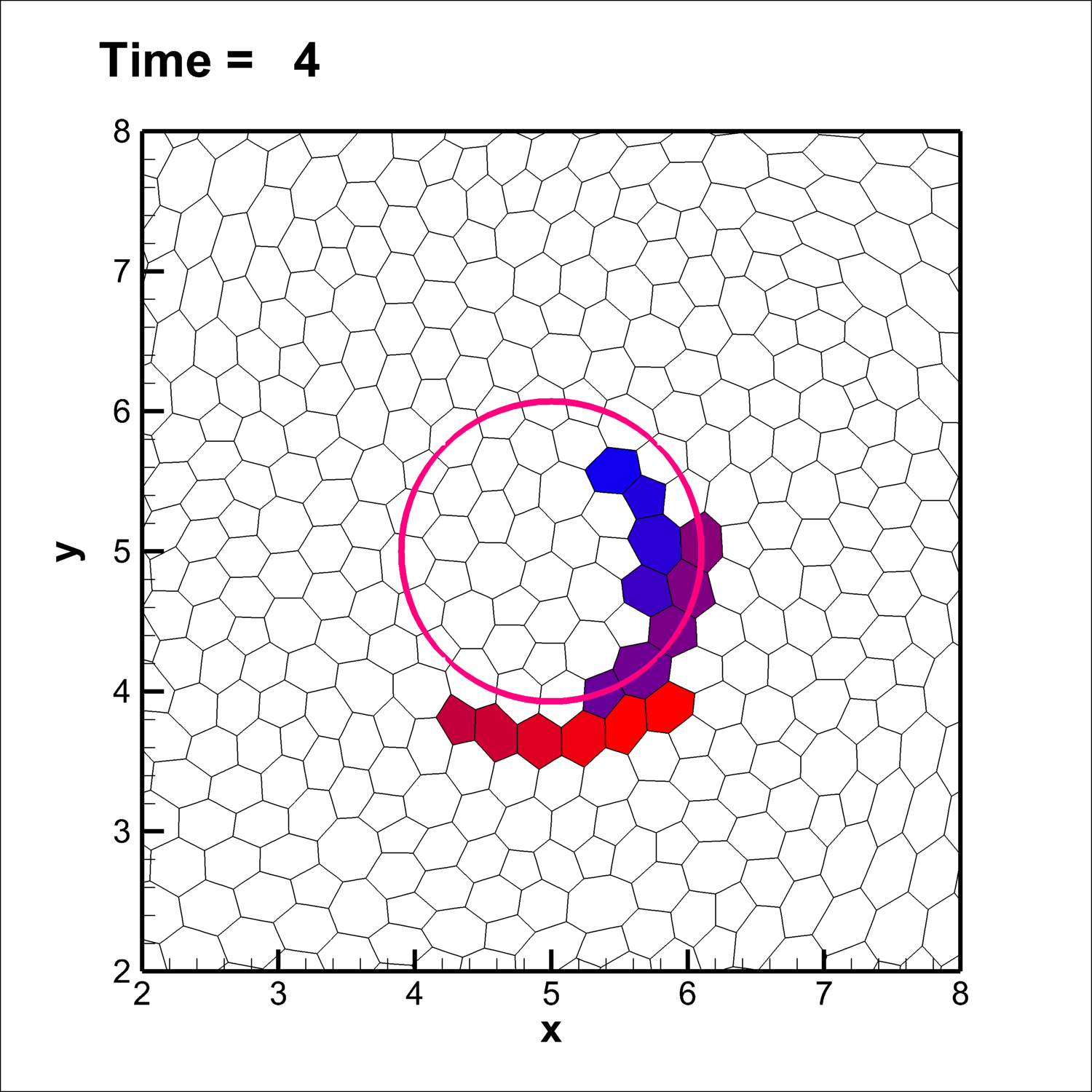}%
\includegraphics[width=0.20\linewidth]{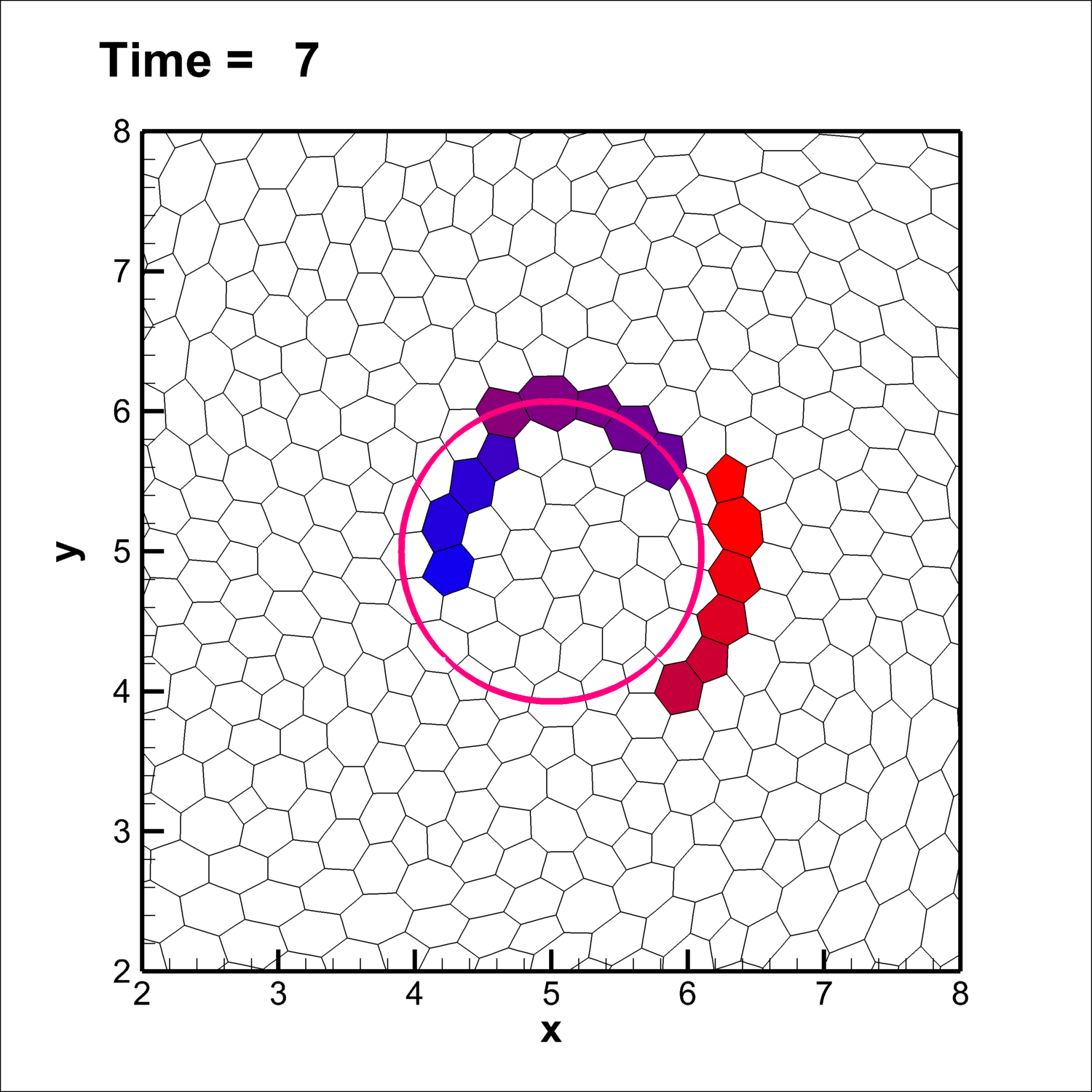}%
\includegraphics[width=0.20\linewidth]{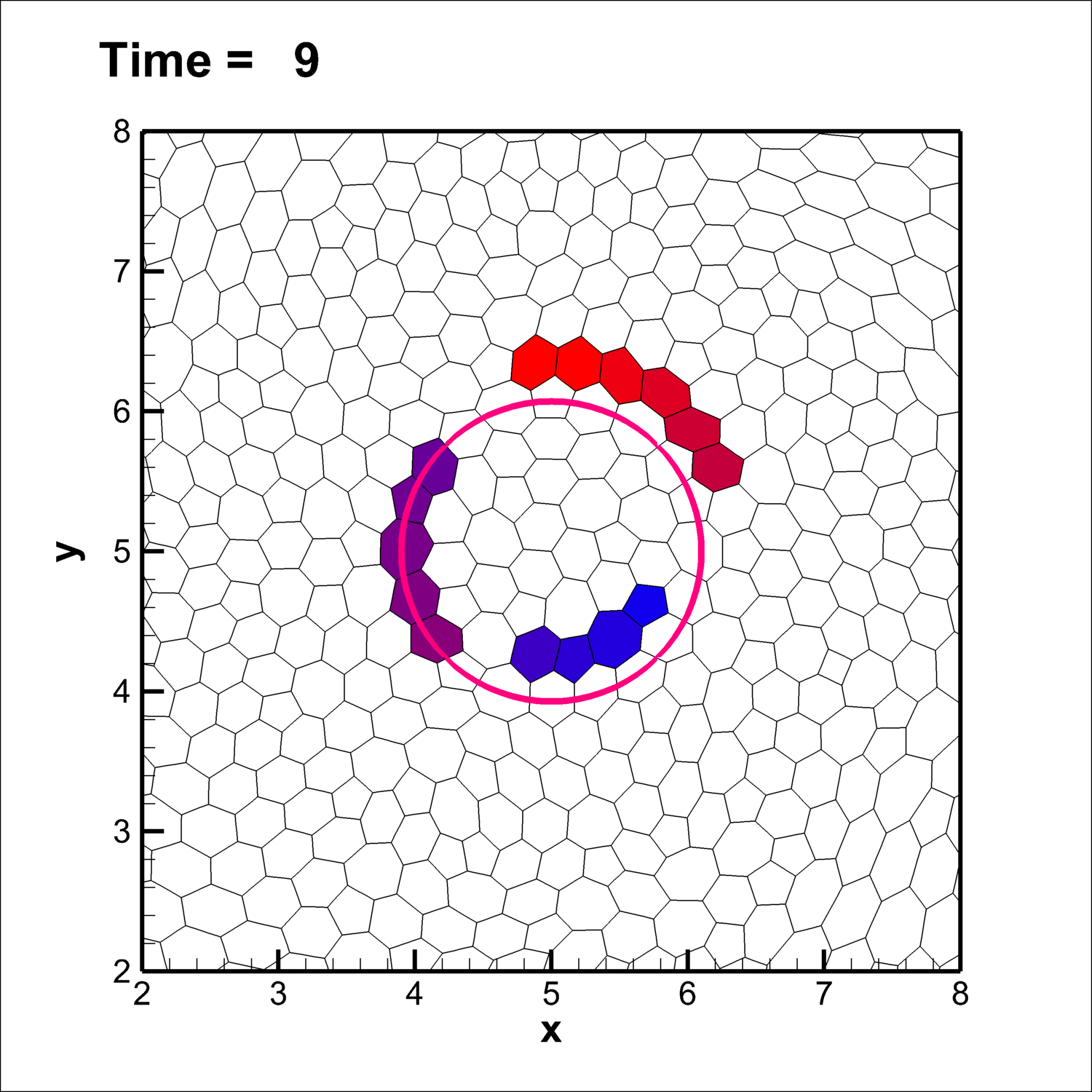}\\[-1.55pt]
\includegraphics[width=0.20\linewidth]{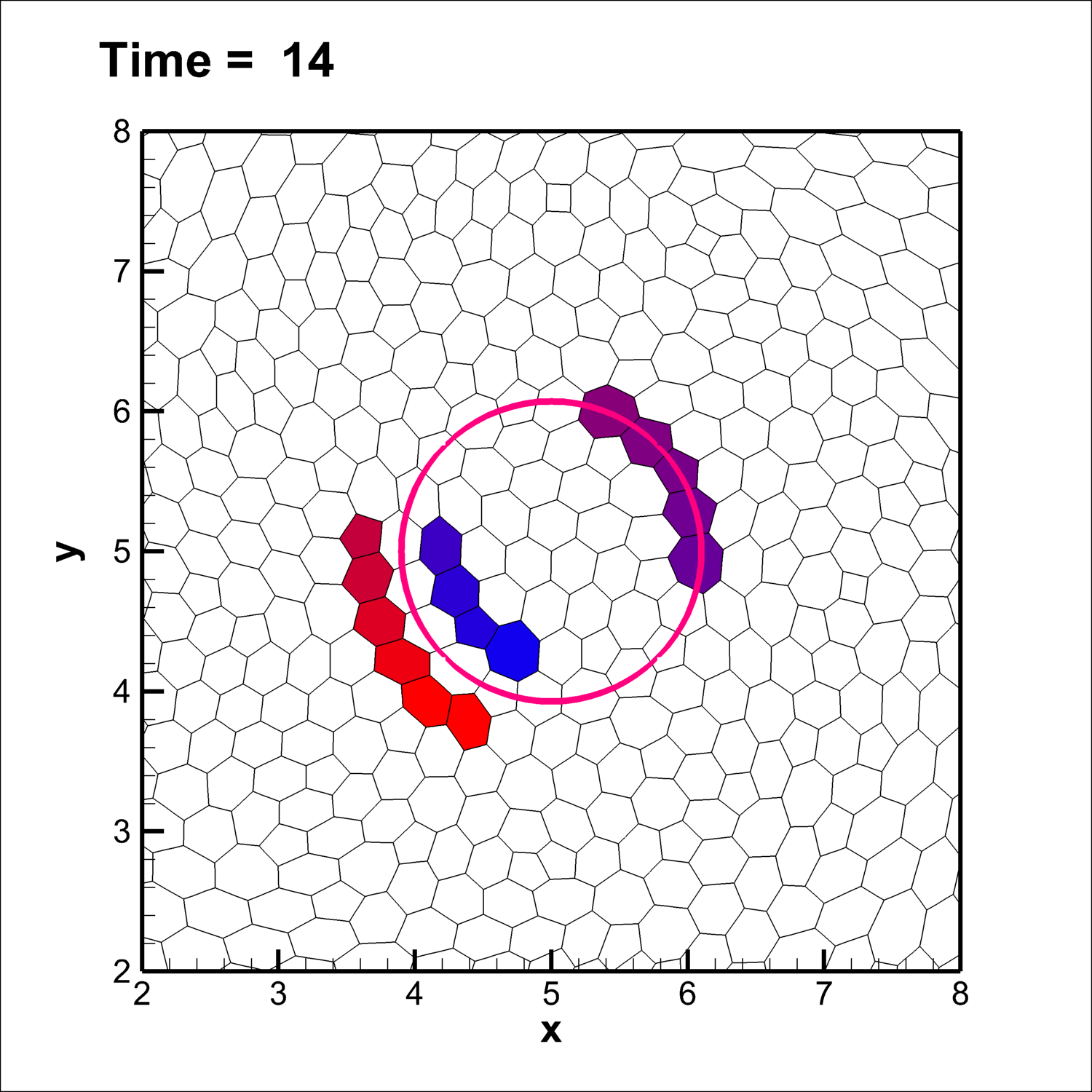}%
\includegraphics[width=0.20\linewidth]{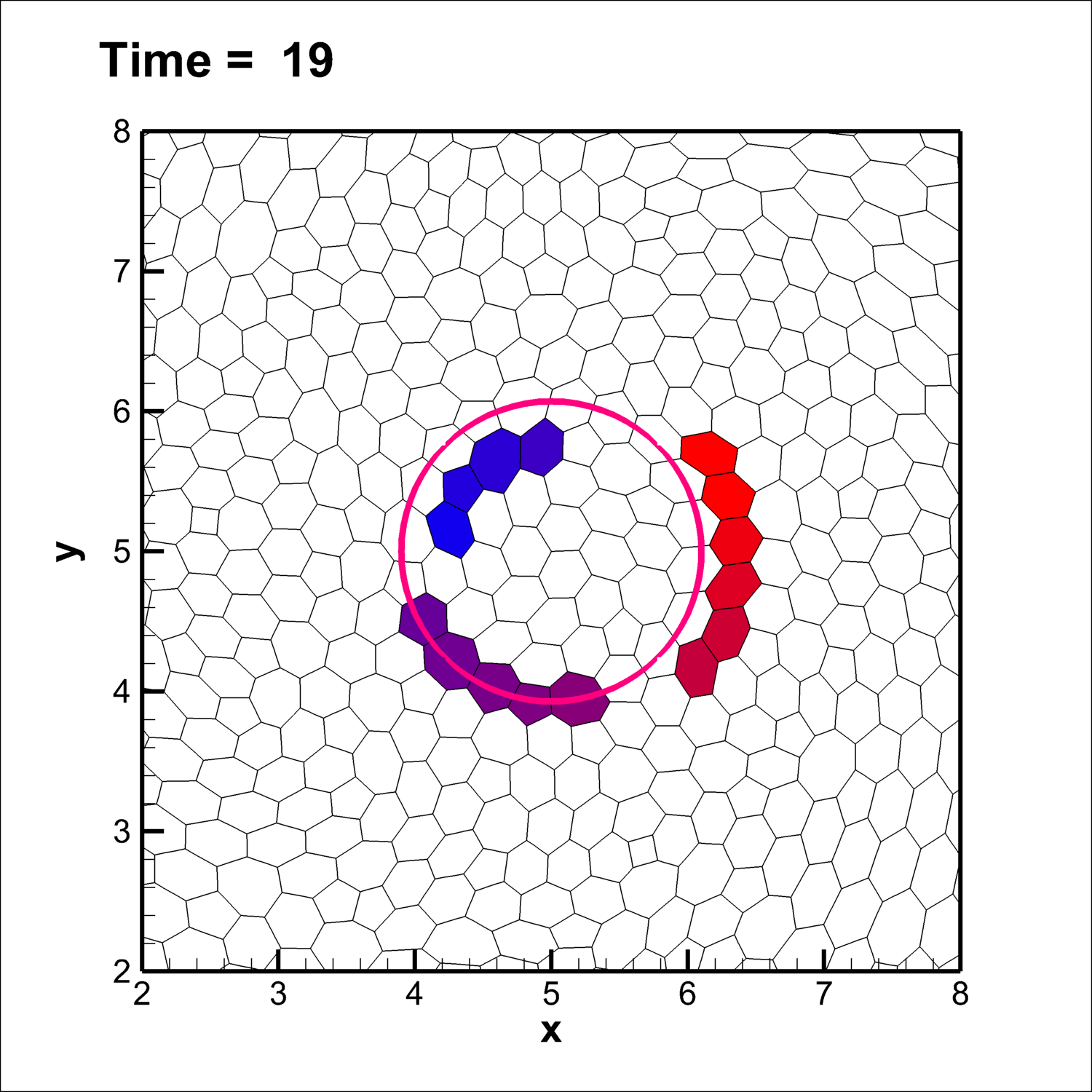}%
\includegraphics[width=0.20\linewidth]{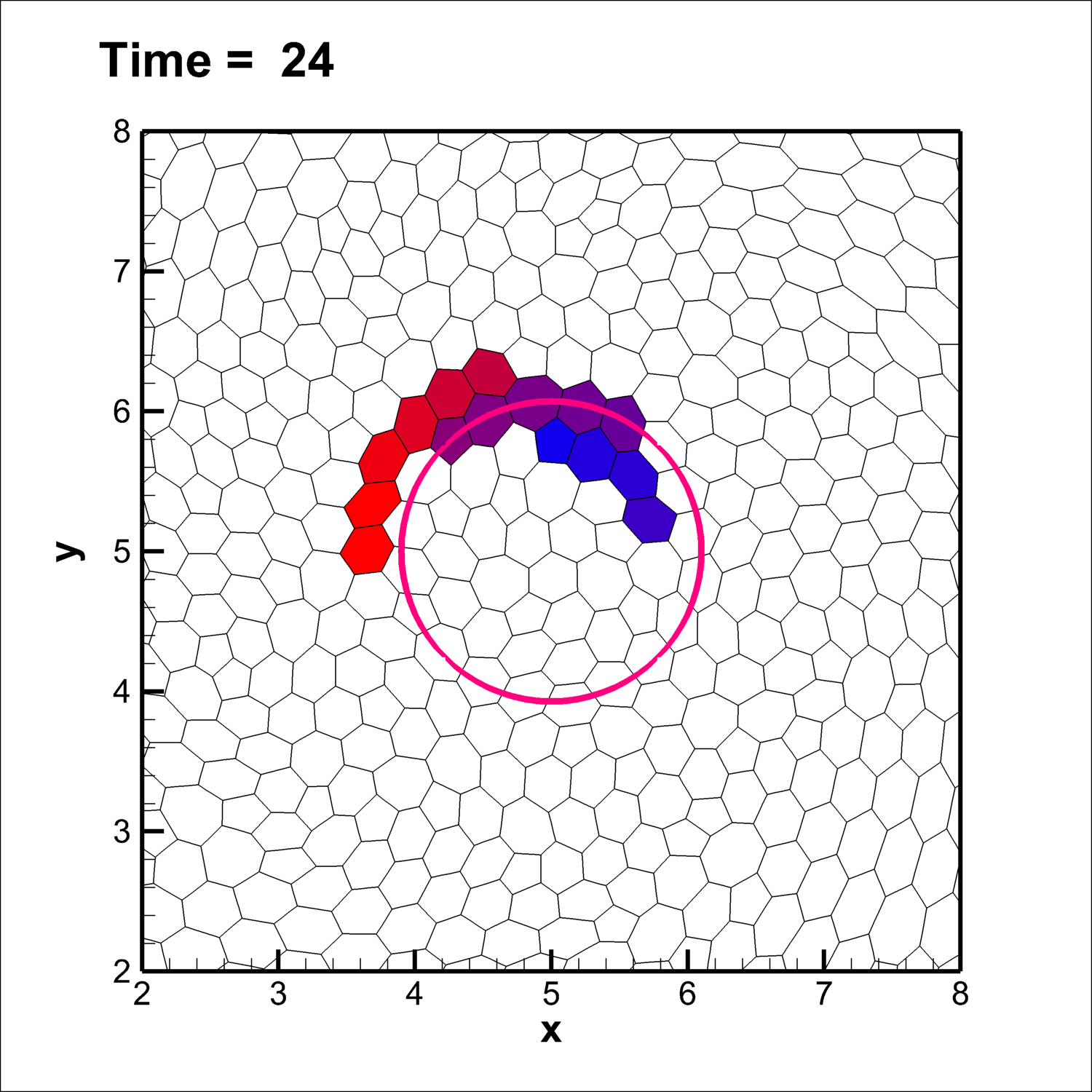}%
\includegraphics[width=0.20\linewidth]{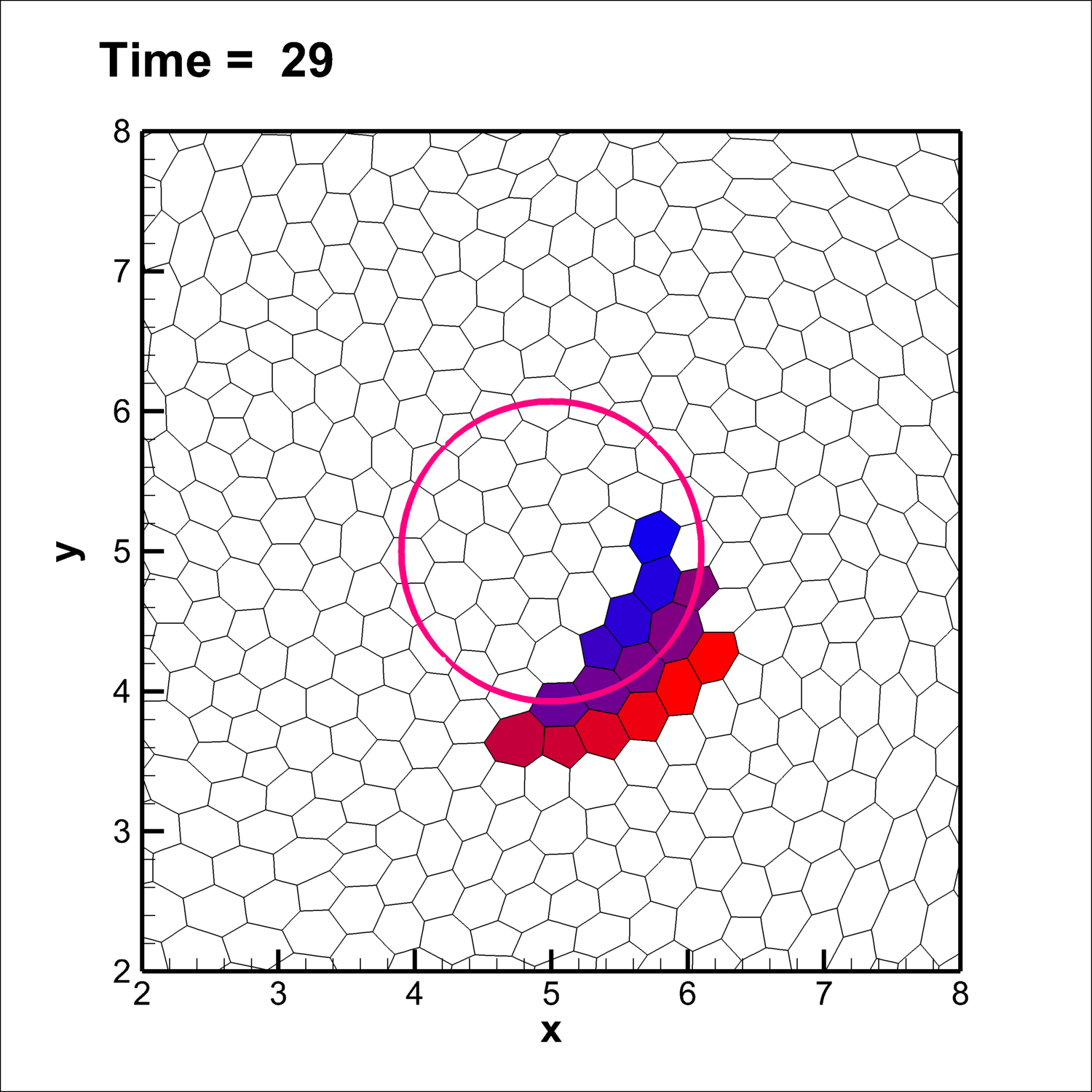}%
\includegraphics[width=0.20\linewidth]{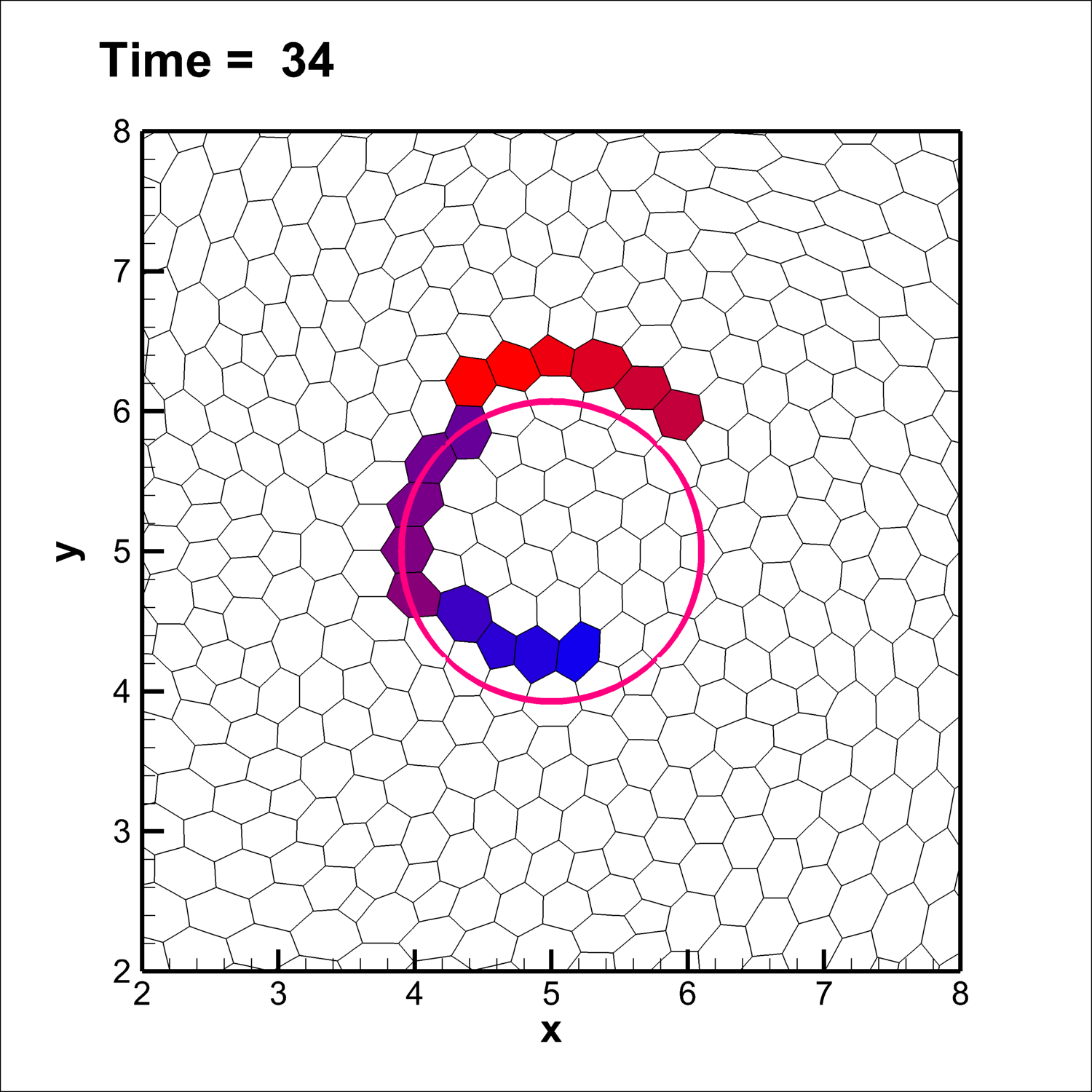}\\[-1.55pt]
\includegraphics[width=0.20\linewidth]{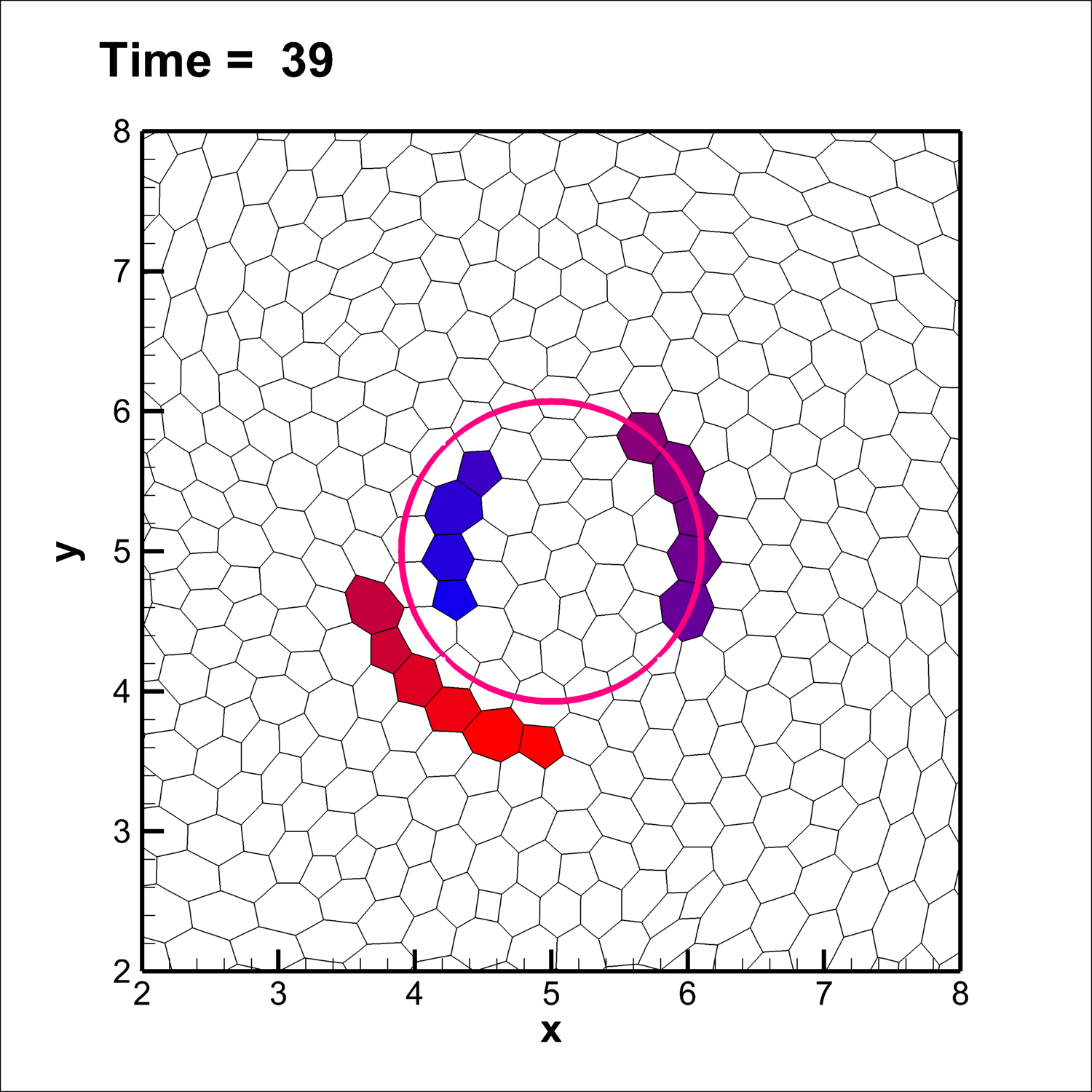}%
\includegraphics[width=0.20\linewidth]{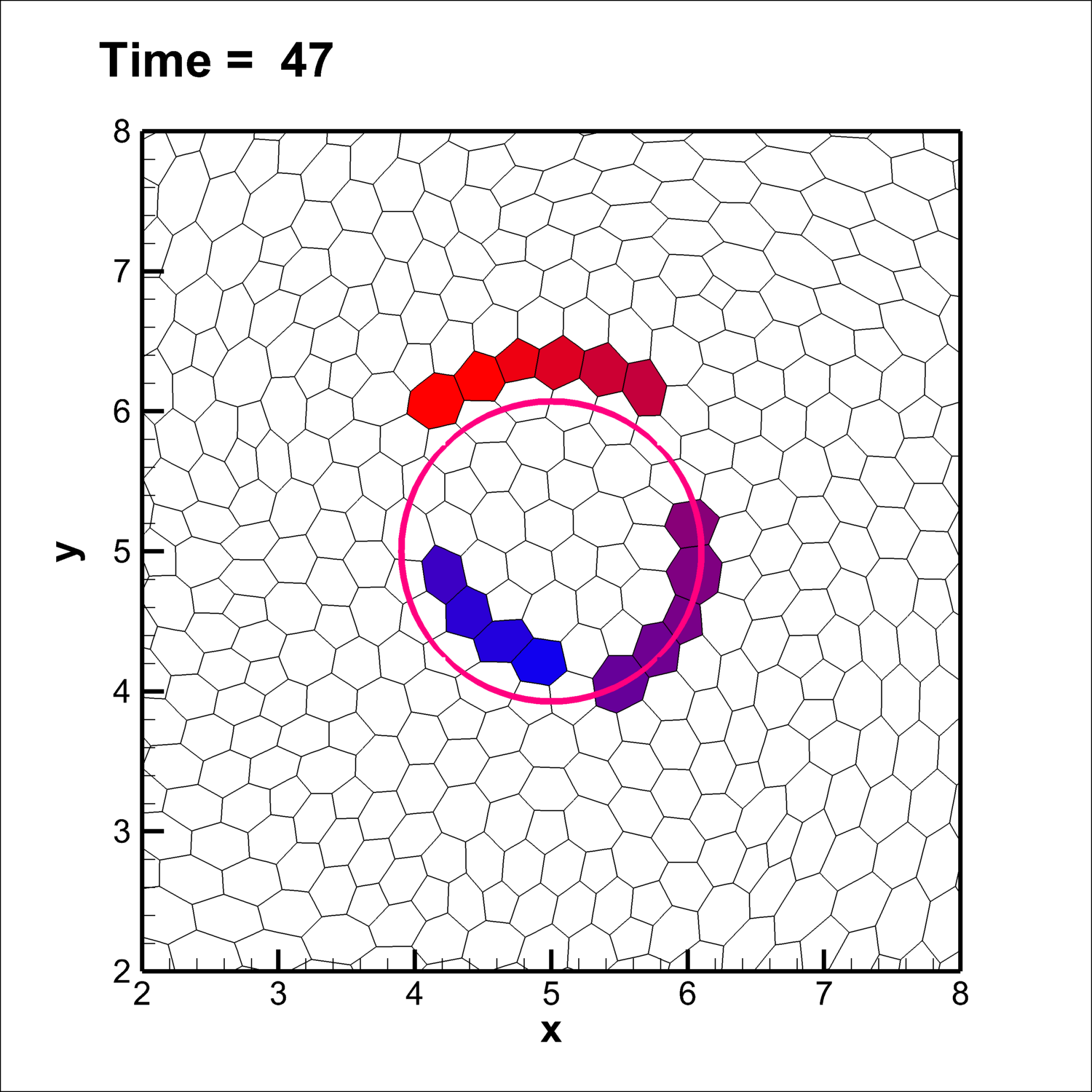}%
\includegraphics[width=0.20\linewidth]{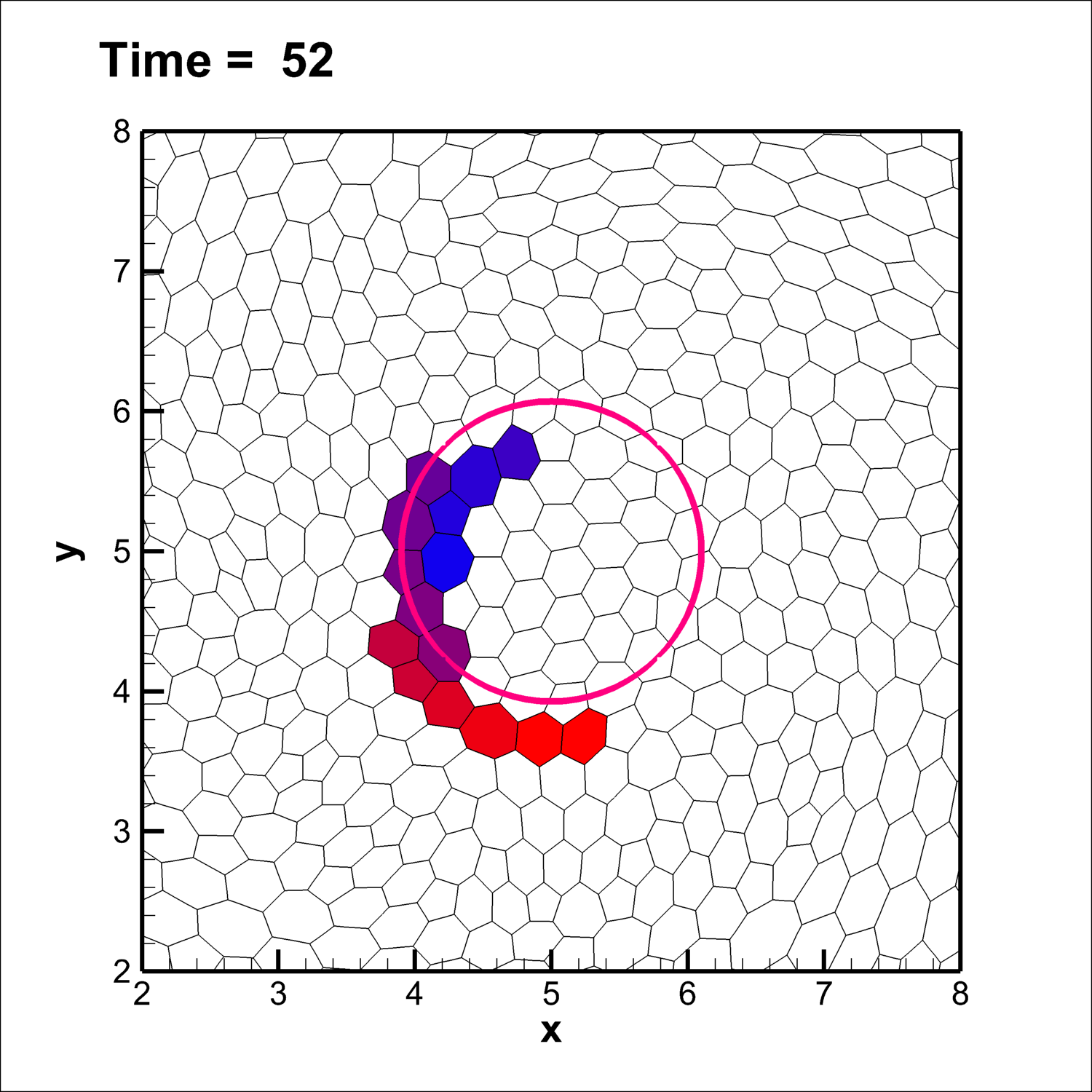}%
\includegraphics[width=0.20\linewidth]{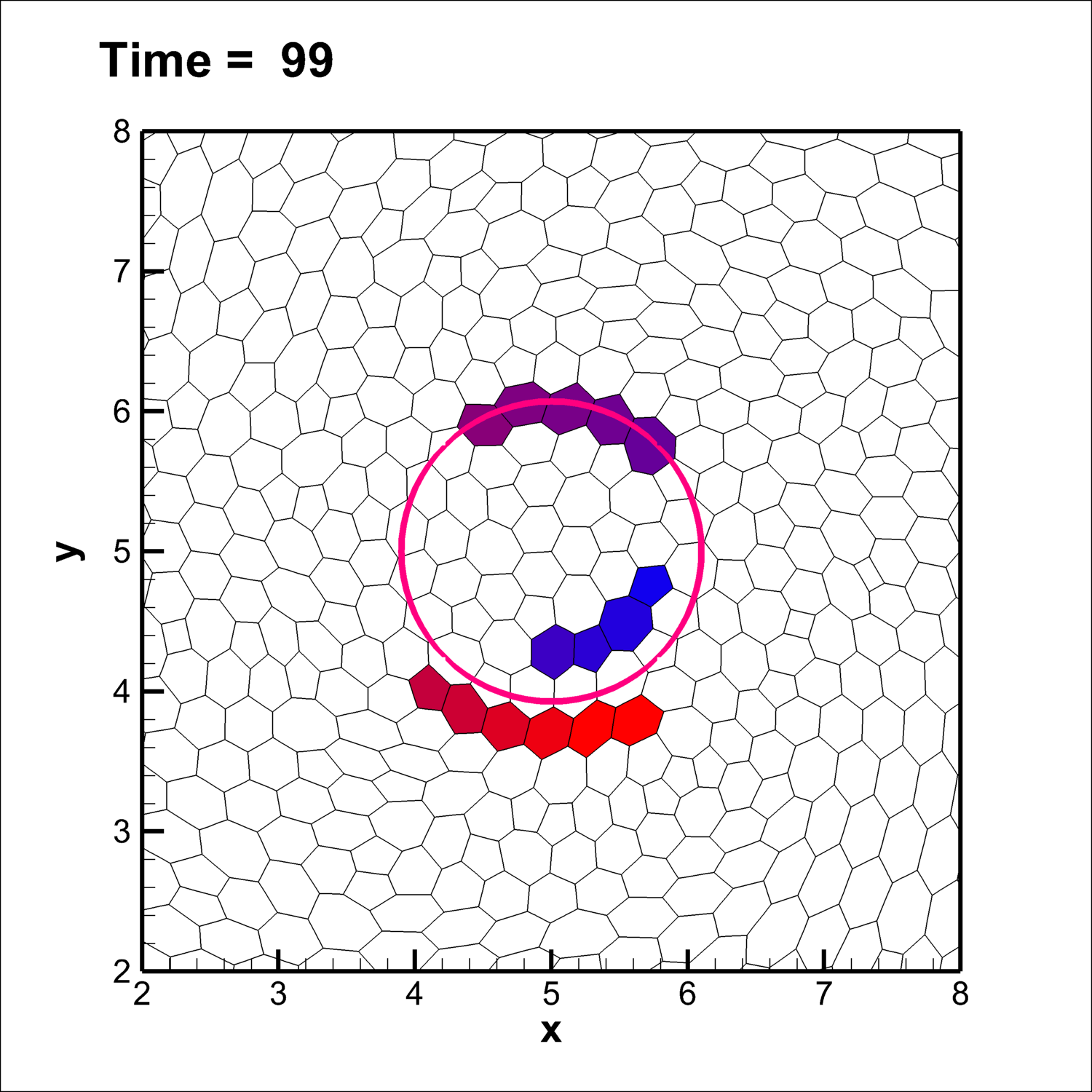}%
\includegraphics[width=0.20\linewidth]{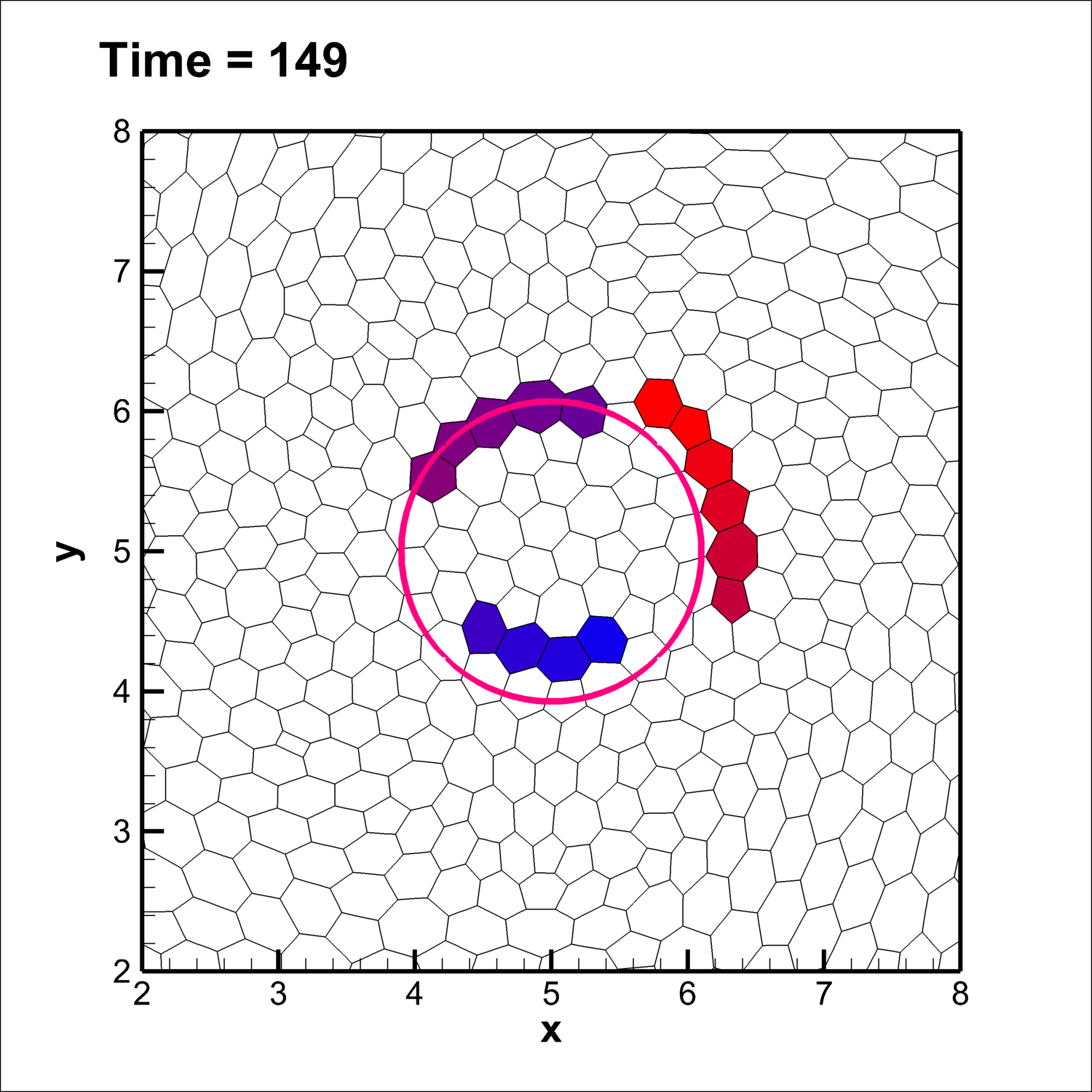}\\[-1.55pt]
\includegraphics[width=0.20\linewidth]{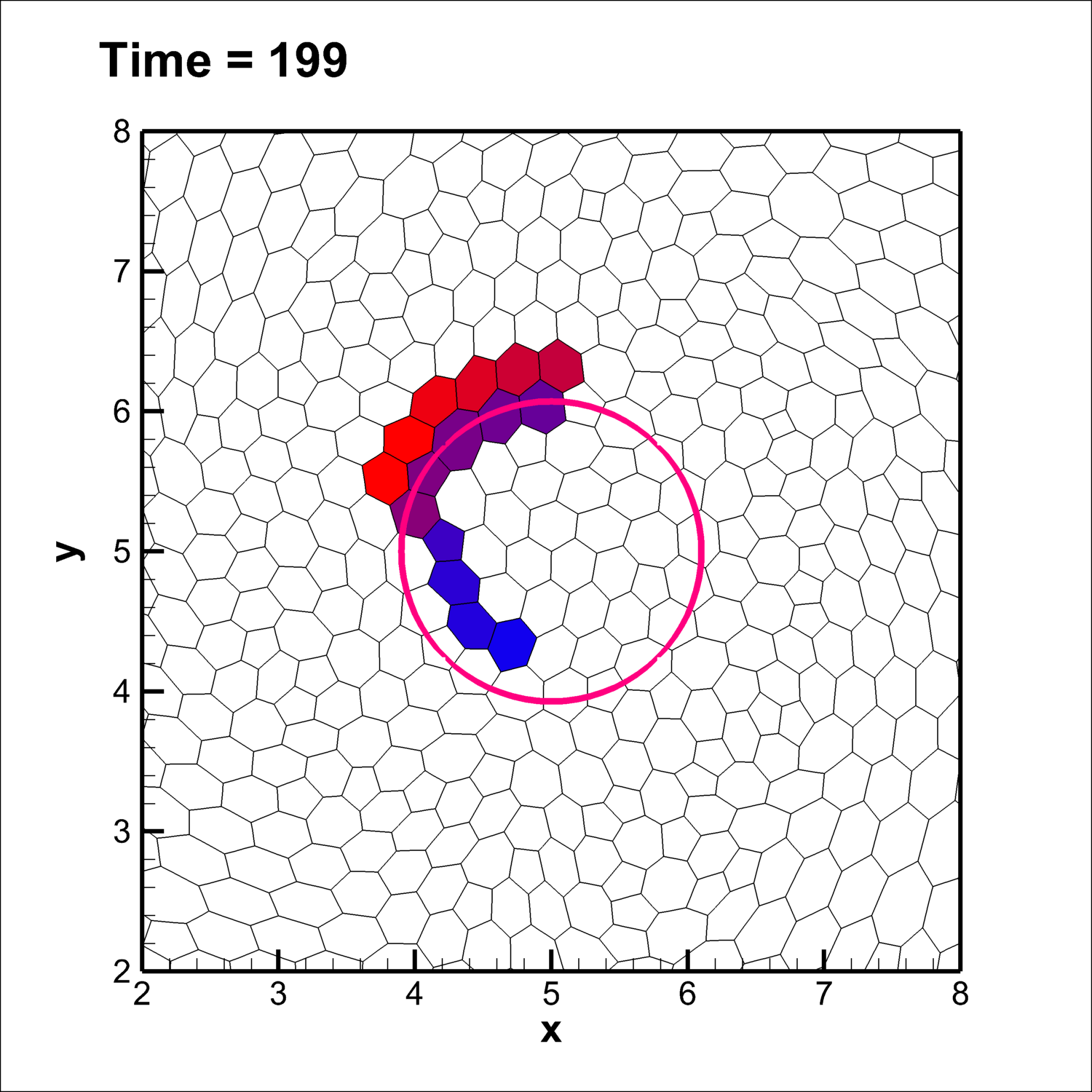}%
\includegraphics[width=0.20\linewidth]{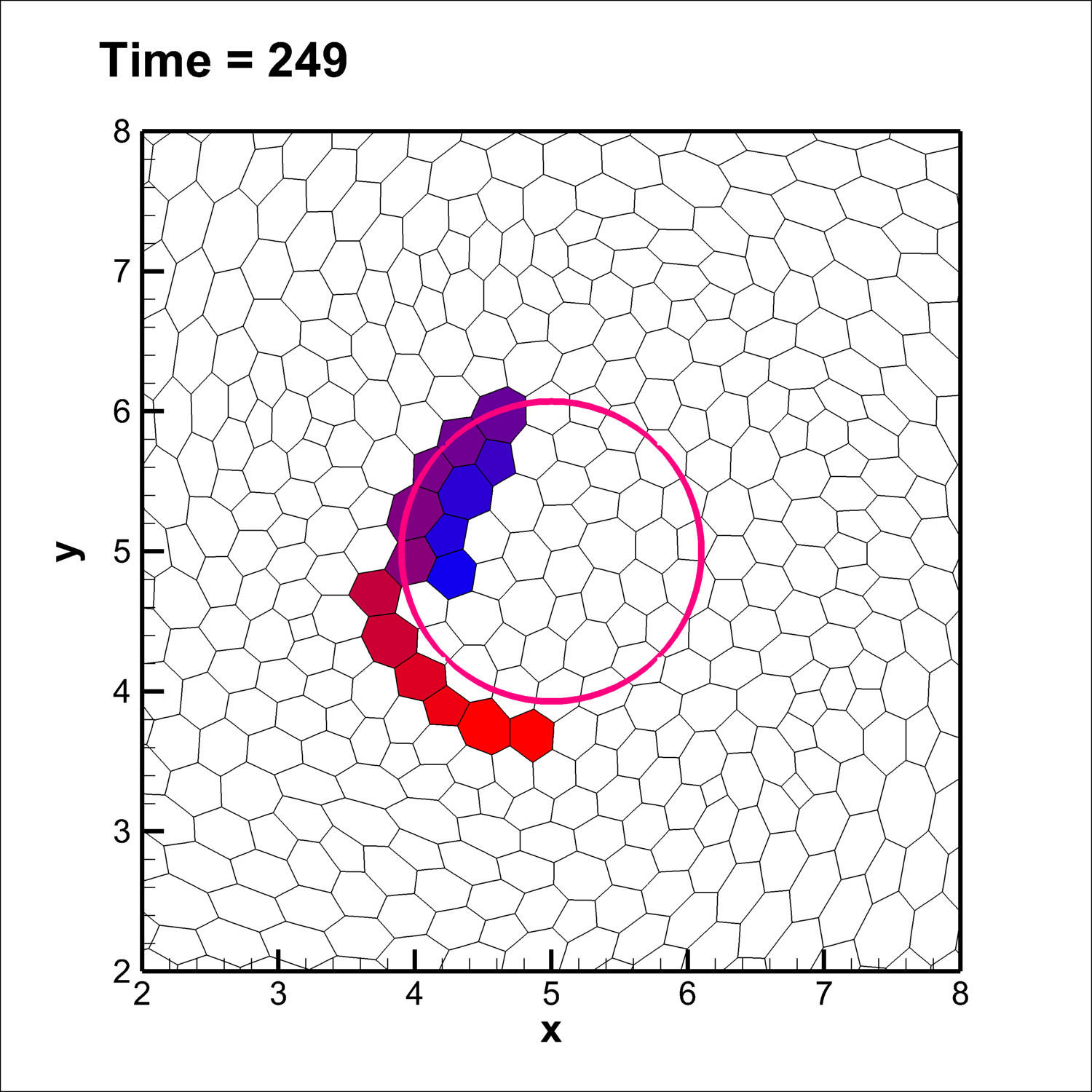}%
\includegraphics[width=0.20\linewidth]{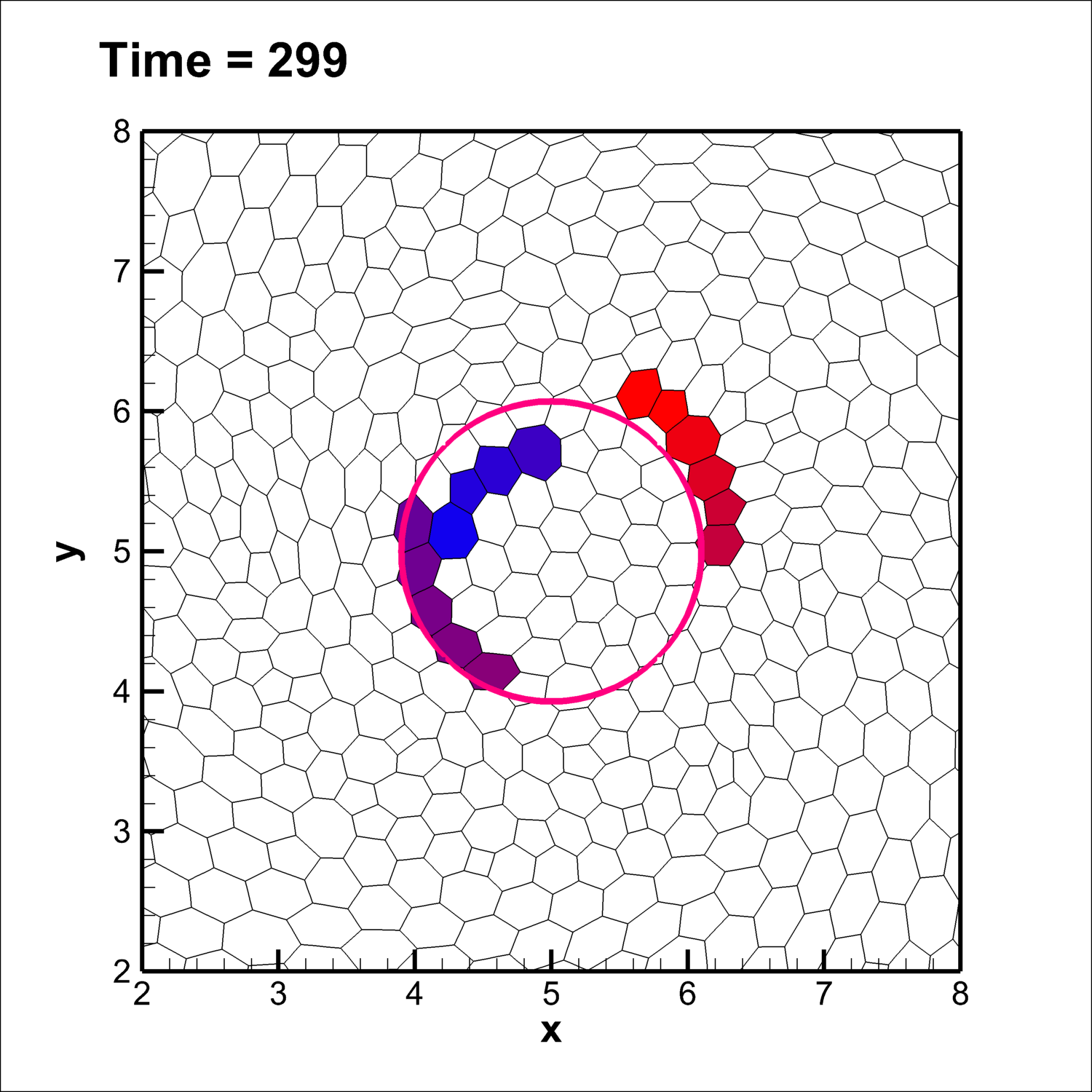}%
\includegraphics[width=0.20\linewidth]{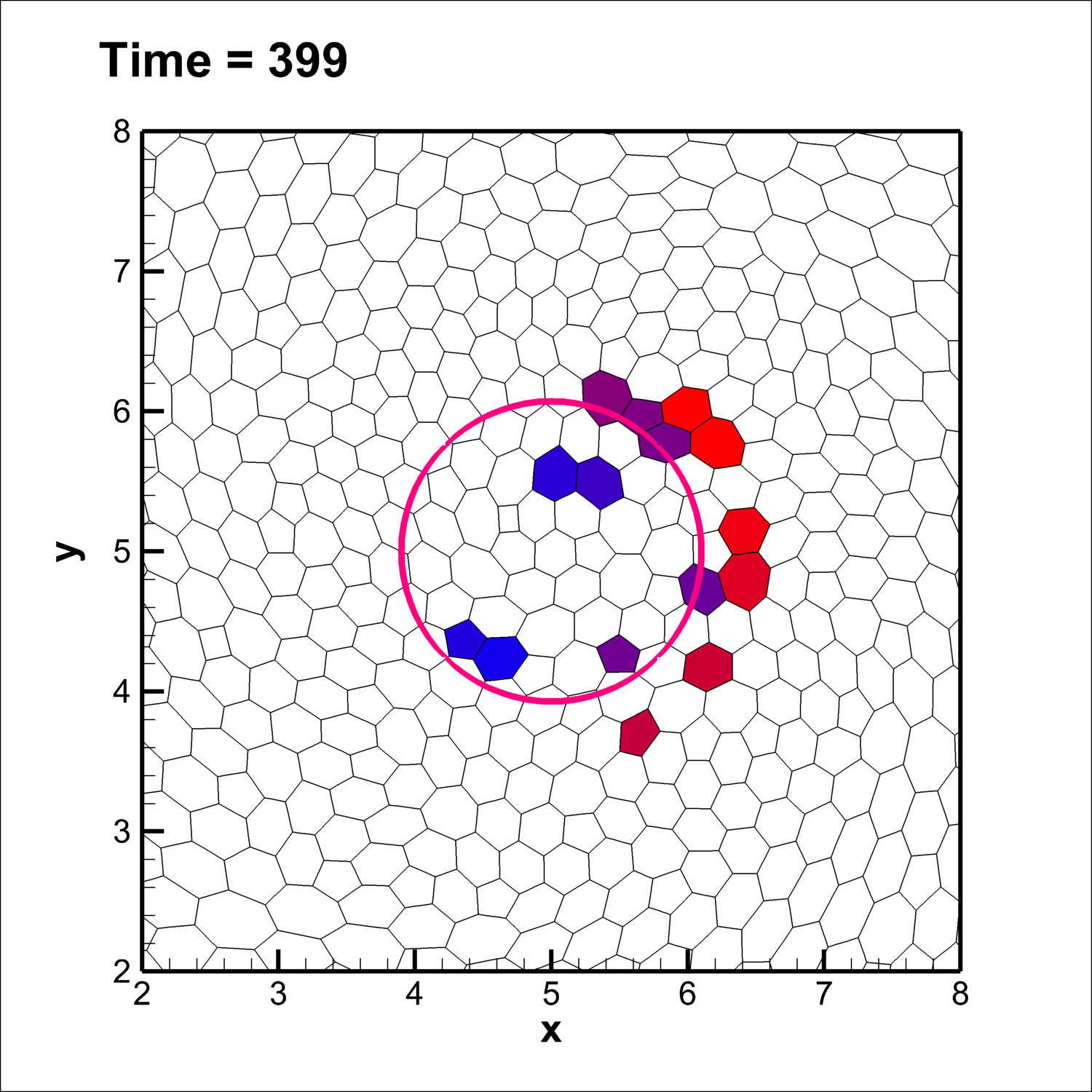}%
\includegraphics[width=0.20\linewidth]{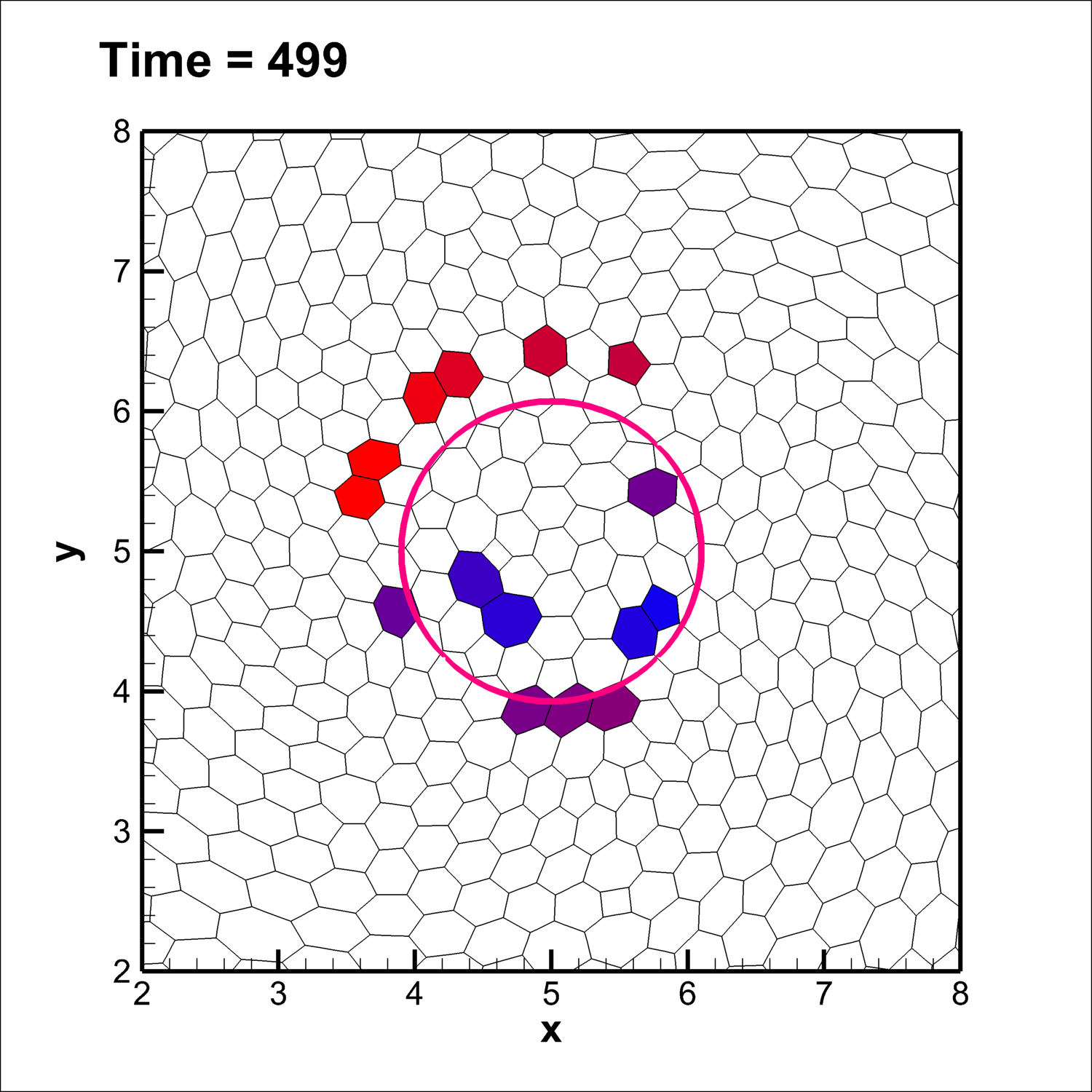}%
\vspace{-4pt}
\caption{Stationary rotating vortex solved with our fourth order $P_3P_3$ ALE-DG scheme on a moving Voronoi mesh of $957$ 
elements with dynamical change of connectivity and with the  generators trajectories computed with fourth order of accuracy 
and Lloyd-like smoothing algorithm with $\mathcal{F} = 10^{-4}$.    
        Density contours (top) and the position of a bunch of highlighted elements (bottom) are provided at different times. 
        The dynamic change of connectivity makes it possible to \textit{substantially improve the mesh quality} w.r.t. standard conforming ALE schemes without
          topology change, for which mesh tangling would lead to a stop of the simulation. In this simulation, after $500$s a point located at $r=1.05$ (pink circle) completes $56$ revolutions. Moreover, it does not collapse in the center and instead maintains a perfect circular trajectory (see also next Figure.)}
    \label{fig.ShuDG4}
\end{figure*}

\paragraph*{Convergence} 
Tables~\ref{tab.orderOfconvergenceFV_shu} and~\ref{tab.orderOfconvergenceDG_shu} report the convergence rates from second up to fifth order of accuracy for the vortex test problem run on a sequence of successively refined meshes \RIIcolor{(finer meshes are obtained by simply increasing the number of generators)}.
For each element, its characteristic size $h_i^n$ at time $t^n$ is given by the diameter of the circumcircle and we denote with $h(\Omega(t_f))$ the average of $h_i^n$ at the final time of the simulation $t_f=0.5$. Thus, $h(\Omega(t_f))$ represents the characteristic mesh size of our mesh.
The optimal order of accuracy is achieved both in space and time for the FV schemes as well as for the DG schemes. 
We would like to underline that this is not trivial for moving Voronoi meshes, 
because the changing characteristic mesh sizes could affect the convergence results: the mesh is \textit{not} stationary at all and
no mesh smoothing is applied for this test case. 

\begin{figure*}[!bp]
\includegraphics[width=0.495\linewidth]{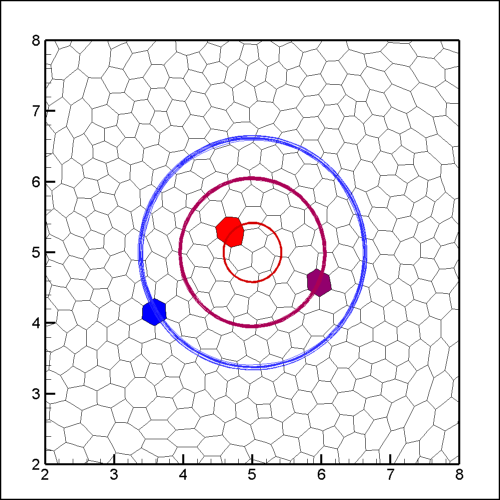}%
\includegraphics[width=0.495\linewidth]{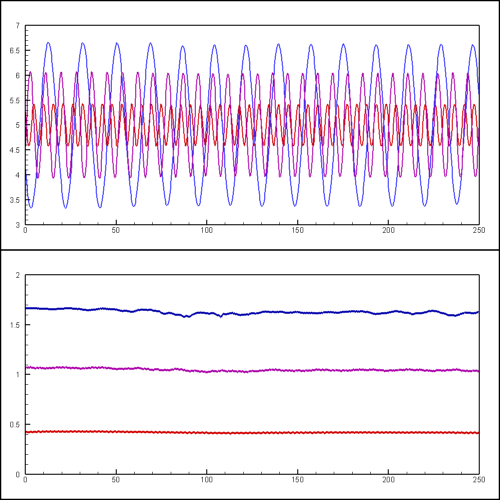}
\caption{Stationary rotating vortex solved with our fourth order $P_3P_3$ ALE-DG scheme on a moving Voronoi mesh of $957$ 
elements with dynamical change of connectivity and with the  generators trajectories computed with fourth order of accuracy, 
refer to Equation~\eqref{eqn.xcnew_high}. The test uses the Lloyd-like smoothing algorithm with $\mathcal{F} = 10^{-4}$.
Left: We depict the trajectories (in Cartesian coordinates) of the generators of $3$ Voronoi elements 
(those highlighted respectively in blue, violet and red) from time $t=0$ up to time $t=250$. During this 
time interval the red Voronoi element makes $30$ complete rotations around the origin. Right: we depict the $y$
 coordinates of the $3$ generators (top) and their radial coordinates (bottom). We would like to emphasize that the
  trajectories are circular (their radius is indeed almost constant) for a very long evolution time.}
\label{fig.Shu_traj}
\end{figure*}

\RIIcolor{We remark that for similar spatial resolution and for the same order of accuracy,
 the numerical errors of the DG scheme are always less than those of the FV method. 
 This is natural and expected since a DG scheme involves many more DOFs w.r.t. a FV algorithm and, in this sense it is more accurate on the same mesh.} 

\begin{table*}[!p]
    \caption{Isentropic vortex. Numerical convergence results for the finite volume algorithm on moving meshes with topology changes. The error norms refer to the variable $\rho$ at time $t=0.5$ in $L_1$ norm.} 
    \label{tab.orderOfconvergenceFV_shu}
    \centering  
        \begin{tabular}{ccc|ccc|ccc|ccc} 
            \hline 
            \multicolumn{3}{c}{$P_0P_1 \rightarrow \mathcal{O}2$} & \multicolumn{3}{c}{$P_0P_2\rightarrow \mathcal{O}3$}    & \multicolumn{3}{c}{$P_0P_3\rightarrow \mathcal{O}4$} & \multicolumn{3}{c}{$P_0P_4\rightarrow \mathcal{O}5$}    \\ 
            \hline
            \!\!$h(\Omega(t_f))$ \!\!&\!\! $\epsilon(\rho)_{L_1}$ \!\!&\!\!\!\! $\mathcal{O}(L_1)$ \!\!\!&  $h(\Omega(t_f))$ \!\!&\!\! $\epsilon(\rho)_{L_1}$ \!\!&\!\!\!\! $\mathcal{O}(L_1)$ \!\!\!& $h(\Omega(t_f))$ \!\!&\!\! $\epsilon(\rho)_{L_1}$ \!\!&\!\! $\mathcal{O}(L_1)$  \!\!&\!\! $h(\Omega(t_f))$ \!\!&\!\! $\epsilon(\rho)_{L_1}$ \!\!&\!\!\!\! $\mathcal{O}(L_1)$ \!\!\!\\
            \hline
            3.8E-01 \!\!&\!\! 3.1E-02  \!\!&\!\! -   \!& 3.8E-01  \!\!&\!\! 2.9E-02 \!\!&\!\! -    \!& 1.9E-01 \!\!&\!\!   1.6E-03 \!\!&\!\!   -  \!& 4.7E-01 \!\!&\!\! 4.0E-02 \!\!& -   \\
            2.0E-01 \!\!&\!\! 6.2E-03  \!\!&\!\! 2.4 \!& 1.9E-01  \!\!&\!\! 4.6E-03 \!\!&\!\! 2.8  \!& 1.3E-01 \!\!&\!\!   4.1E-04 \!\!&\!\!  3.4 \!& 3.8E-01 \!\!&\!\! 1.4E-02 \!\!& 4.8 \\
            1.3E-01 \!\!&\!\! 2.4E-03  \!\!&\!\! 2.4 \!& 1.3E-01  \!\!&\!\! 1.4E-03 \!\!&\!\! 2.9  \!& 9.9E-02 \!\!&\!\!   1.4E-04 \!\!&\!\!  3.8 \!& 1.3E-01 \!\!&\!\! 2.5E-04 \!\!& 3.8 \\
            9.9E-02 \!\!&\!\! 1.3E-03  \!\!&\!\! 2.3 \!& 9.9E-02  \!\!&\!\! 6.1E-04 \!\!&\!\! 3.0  \!& 7.9E-02 \!\!&\!\!   6.0E-05 \!\!&\!\!  3.9 \!& 9.9E-02 \!\!&\!\! 6.7E-05 \!\!& 4.6 \\
            8.0E-02 \!\!&\!\! 7.8E-04  \!\!&\!\! 2.2 \!& 7.9E-02  \!\!&\!\! 3.1E-04 \!\!&\!\! 2.0  \!& 6.7E-03 \!\!&\!\!   3.0E-05 \!\!&\!\!  3.8 \!& 7.9E-02 \!\!&\!\! 2.4E-05 \!\!& 4.7 \\
            \hline 
        \end{tabular}       
\end{table*}

\begin{table*}[!p]
    \caption{Isentropic vortex. Numerical convergence results for the discontinuous Galerkin  algorithm on moving meshes with topology changes. The error norms refer to the variable $\rho$ at time $t=0.5$ in $L_1$ norm.} 
    \label{tab.orderOfconvergenceDG_shu}
    \centering
            \begin{tabular}{ccc|ccc|ccc|ccc} 
            \hline 
            \multicolumn{3}{c}{$P_1P_1\rightarrow \mathcal{O}2$} & \multicolumn{3}{c}{$P_2P_2\rightarrow \mathcal{O}3$} & \multicolumn{3}{c}{$P_3P_3\rightarrow \mathcal{O}4$}      & \multicolumn{3}{c}{$P_4P_4\rightarrow \mathcal{O}5$}     \\ 
            \hline
            \!\!$h(\Omega(t_f))$ \!\!&\!\! $\epsilon(\rho)_{L_1}$ \!\!&\!\!\!\! $\mathcal{O}(L_1)$ \!\!\!&  $h(\Omega(t_f))$ \!\!&\!\! $\epsilon(\rho)_{L_1}$ \!\!&\!\!\!\! $\mathcal{O}(L_1)$ \!\!\!& $h(\Omega(t_f))$ \!\!&\!\! $\epsilon(\rho)_{L_1}$ \!\!&\!\! $\mathcal{O}(L_1)$  \!\!&\!\! $h(\Omega(t_f))$ \!\!&\!\! $\epsilon(\rho)_{L_1}$ \!\!&\!\!\!\! $\mathcal{O}(L_1)$ \!\!\!\\
            \hline
            7.5E-01 \!\!&\!\! 6.3E-03  \!\!&\!\! -   \!&  7.5E-01 \!\!&\!\! 1.4E-02 \!\!&\!\! -    \!& 6.1E-01 \!\!&\!\!   1.4E-03 \!\!&\!\!   -   \!& 1.4E-00 \!\!&\!\! 1.1E-02 \!\!& -   \\
            6.1E-01 \!\!&\!\! 4.2E-04  \!\!&\!\! 1.9 \!&  6.1E-01 \!\!&\!\! 7.2E-03 \!\!&\!\! 3.4  \!& 5.2E-01 \!\!&\!\!   7.4E-04 \!\!&\!\!  3.7  \!& 1.0E-00 \!\!&\!\! 2.0E-03 \!\!& 5.9 \\
            3.2E-01 \!\!&\!\! 9.9E-04  \!\!&\!\! 2.2 \!&  3.2E-01 \!\!&\!\! 9.3E-04 \!\!&\!\! 3.2  \!& 4.7E-01 \!\!&\!\!   4.1E-04 \!\!&\!\!  5.9  \!& 9.8E-01 \!\!&\!\! 1.6E-03 \!\!& 4.7 \\
            2.2E-01 \!\!&\!\! 4.4E-04  \!\!&\!\! 2.0 \!&  2.2E-01 \!\!&\!\! 2.8E-04 \!\!&\!\! 3.0  \!& 3.2E-01 \!\!&\!\!   7.7E-05 \!\!&\!\!  4.4  \!& 8.9E-01 \!\!&\!\! 9.0E-04 \!\!& 5.9 \\
            1.6E-01 \!\!&\!\! 2.5E-05  \!\!&\!\! 2.0 \!&  1.6E-01 \!\!&\!\! 1.2E-04 \!\!&\!\! 3.0  \!& 2.2E-01 \!\!&\!\!   1.6E-05 \!\!&\!\!  4.0  \!& 8.5E-01 \!\!&\!\! 7.0E-04 \!\!& 5.1 \\
            \hline 
        \end{tabular}       
\end{table*}

\begin{table*}[!tp]
    \caption{Isentropic vortex. Numerical convergence results for the third order $P_2P_2$ discontinuous Galerkin algorithm on moving meshes with topology changes. The error norms refer to the variable $\rho$ at time $t=0.5$ in $L_2$ norm. The three groups of results refer to three different ways of ordering the space--time neighbors of each element. The fact that the errors are exactly the same up to machine precision proves that the algorithm is independent of the neighbor ordering used in the construction of the space--time elements.} 
    \label{tab.IndipendenceOfNeighsNumb}
    \centering
        \begin{tabular}{ccc} 
            \multicolumn{3}{c}{ordering from $1^{\text{st}}$ common neighbor} \\
            \hline
            $h(\Omega(t_f))$ & $\epsilon(\rho)_{L_2}$ & $\!\!\mathcal{O}(L_2)\!\!\!\!$\\ 
            \hline
            0.319411631217116 & 9.2414523328907E-04 &   -   \\
            0.242212163540348 & 3.9353901580992E-04 &  3.1  \\
            0.194949032600822 & 2.0616099552666E-04 &  3.0  \\
            0.163155447483668 & 1.1964571728528E-04 &  3.1  \\
            0.122985013713313 & 5.1270456290057E-05 &  3.0  \\      
            \hline \\[1pt]
            \multicolumn{3}{c}{ordering from $2^{\text{nd}}$ common neighbor}\\
            \hline
            $h(\Omega(t_f))$ & $\epsilon(\rho)_{L_2}$ & $\!\!\mathcal{O}(L_2)\!\!\!\!$\\ 
            \hline          
            0.319411631217114 & 9.2414523328982E-04  &  -    \\
            0.242212163540348 & 3.9353901581037E-04  &  3.1  \\
            0.194949032600822 & 2.0616099552752E-04  &  3.0  \\
            0.163155447483668 & 1.1964571728459E-04  &  3.1  \\
            0.122985013713313 & 5.1270456288495E-05  &  3.0  \\
            \hline \\[1pt]
            \multicolumn{3}{c}{ordering from $3^{\text{rd}}$ common neighbor}\\
            \hline
            $h(\Omega(t_f))$ & $\epsilon(\rho)_{L_2}$ & $\!\!\mathcal{O}(L_2)\!\!\!\!$\\ 
            \hline          
            0.319411631217116 &  9.2414523328907E-04& -     \\
            0.242212163540348 &  3.9353901580992E-04& 3.1   \\
            0.194949032600822 &  2.0616099552666E-04& 3.0   \\
            0.163155447483668 &  1.1964571728400E-04& 3.1   \\
            0.122985013713313 &  5.1270456291299E-05& 3.0   \\
            \hline 
        \end{tabular}       
\end{table*}

\begin{figure*}[!bp]
    \centering
    \includegraphics[width=0.2\linewidth]{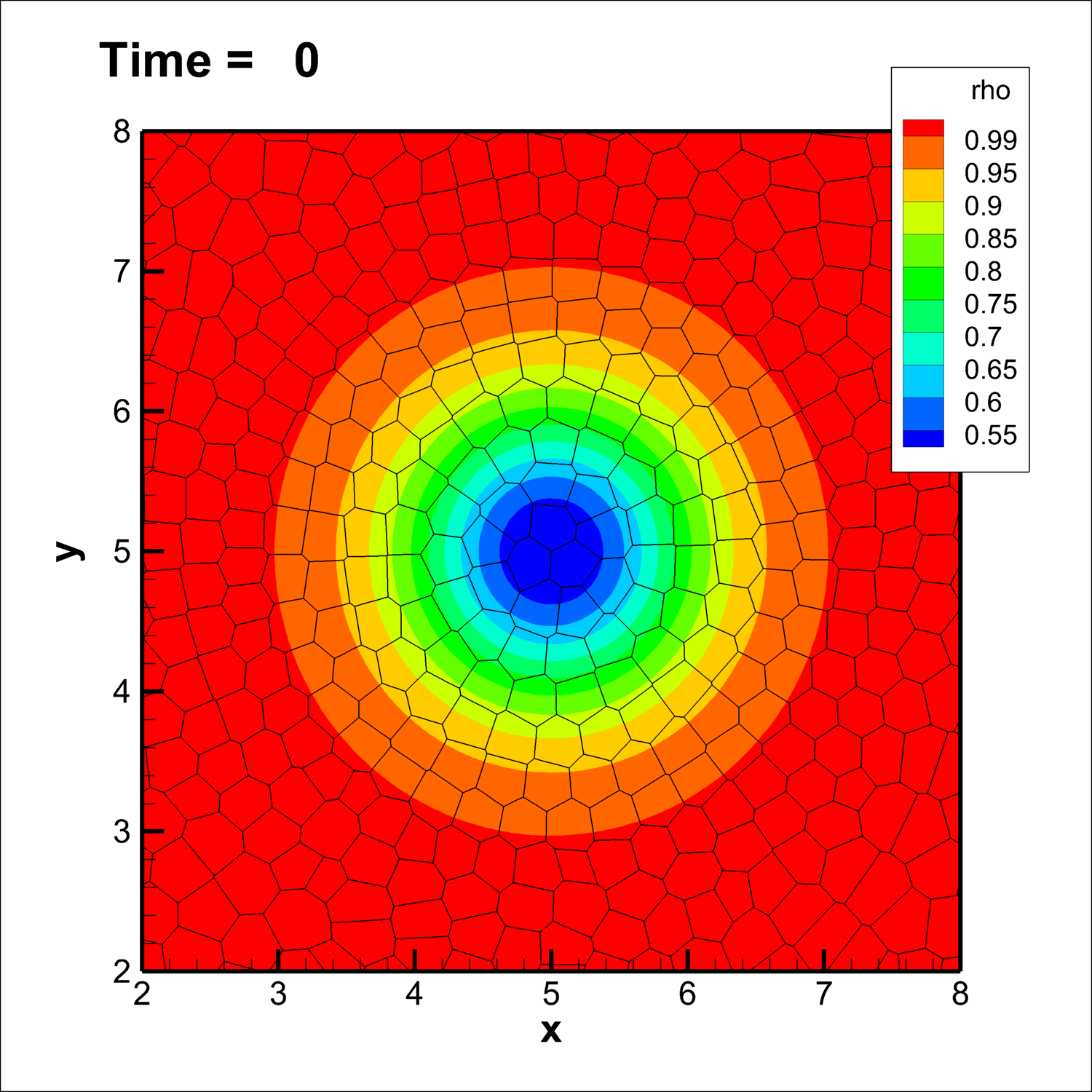}%
    \includegraphics[width=0.2\linewidth]{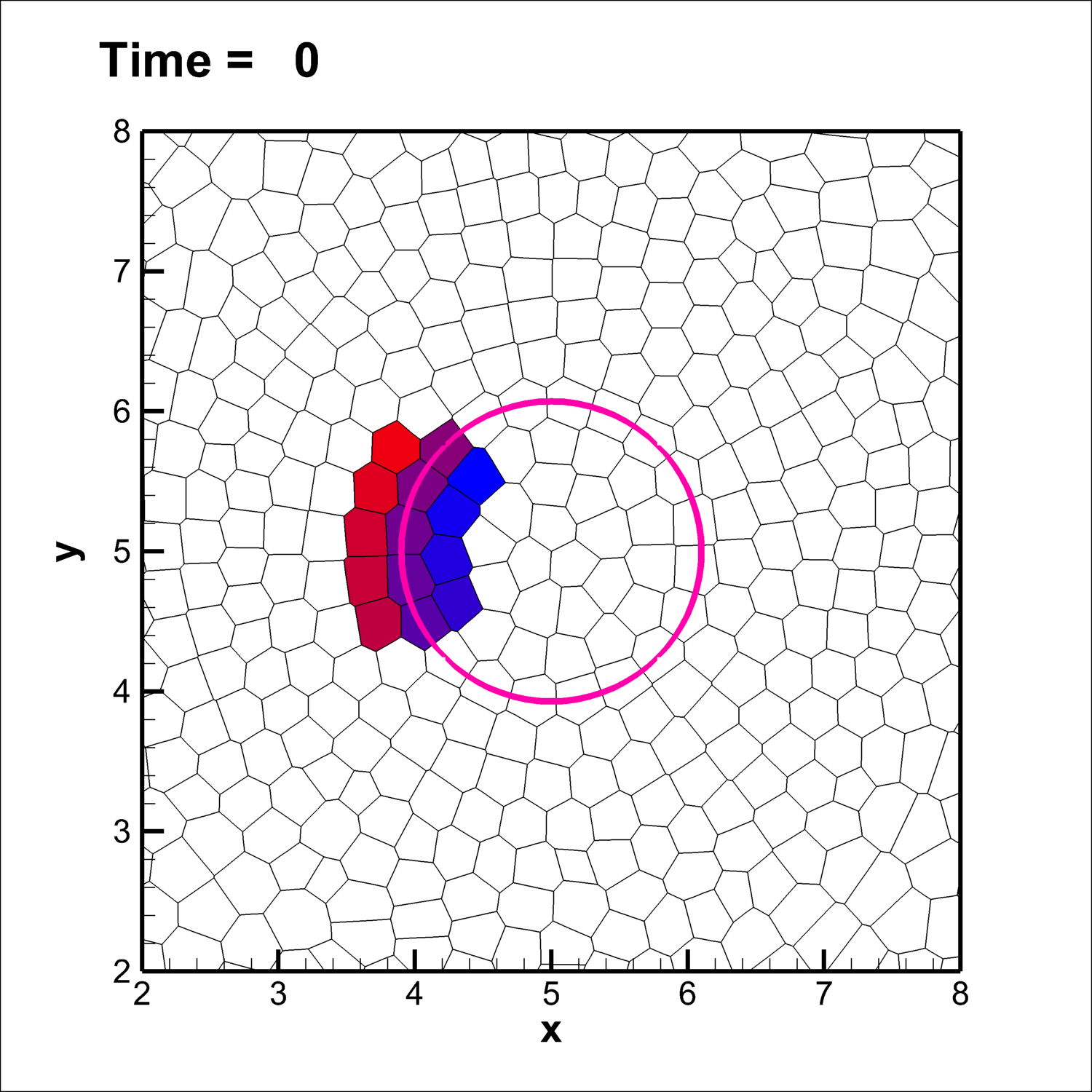}%
    \includegraphics[width=0.2\linewidth]{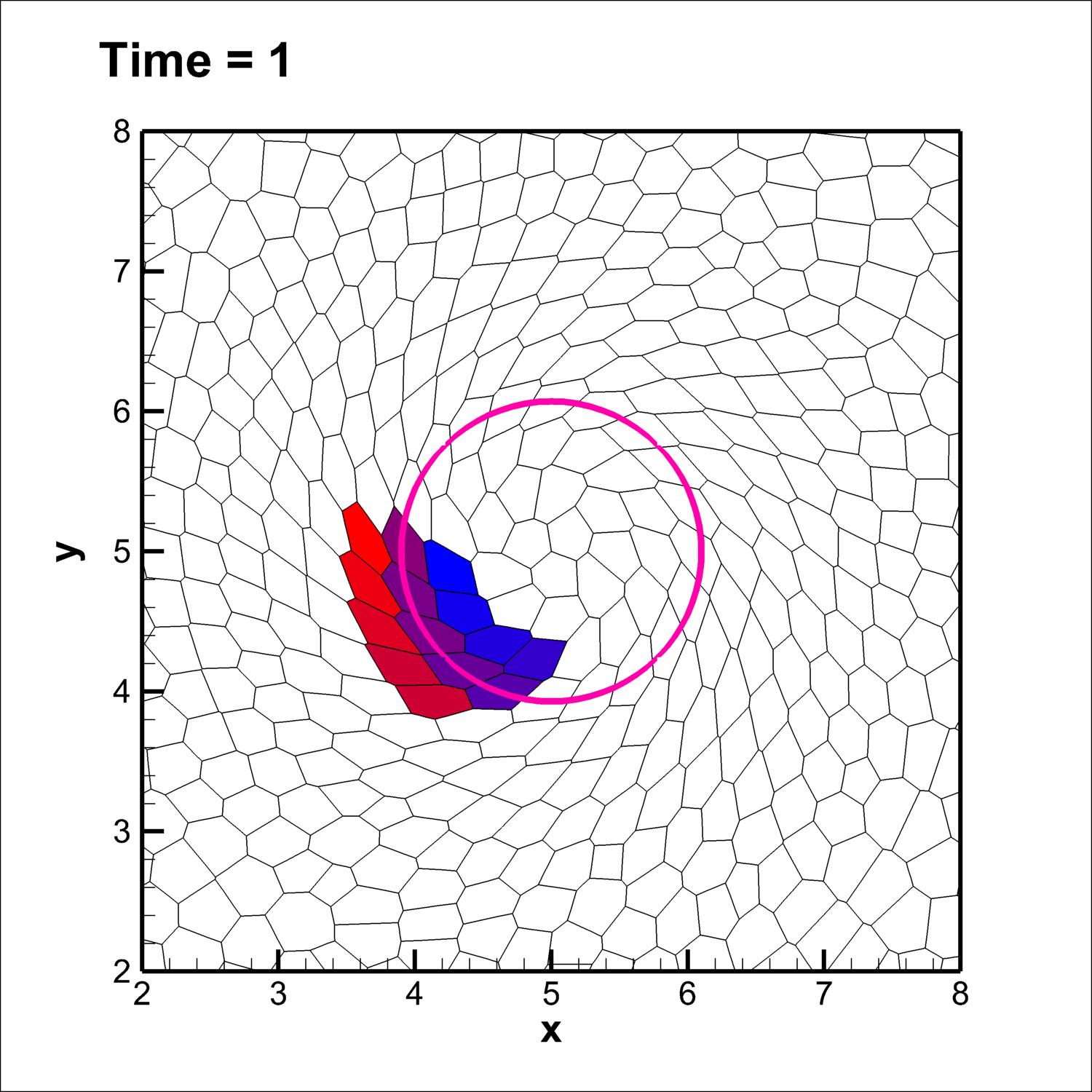}%
    \includegraphics[width=0.2\linewidth]{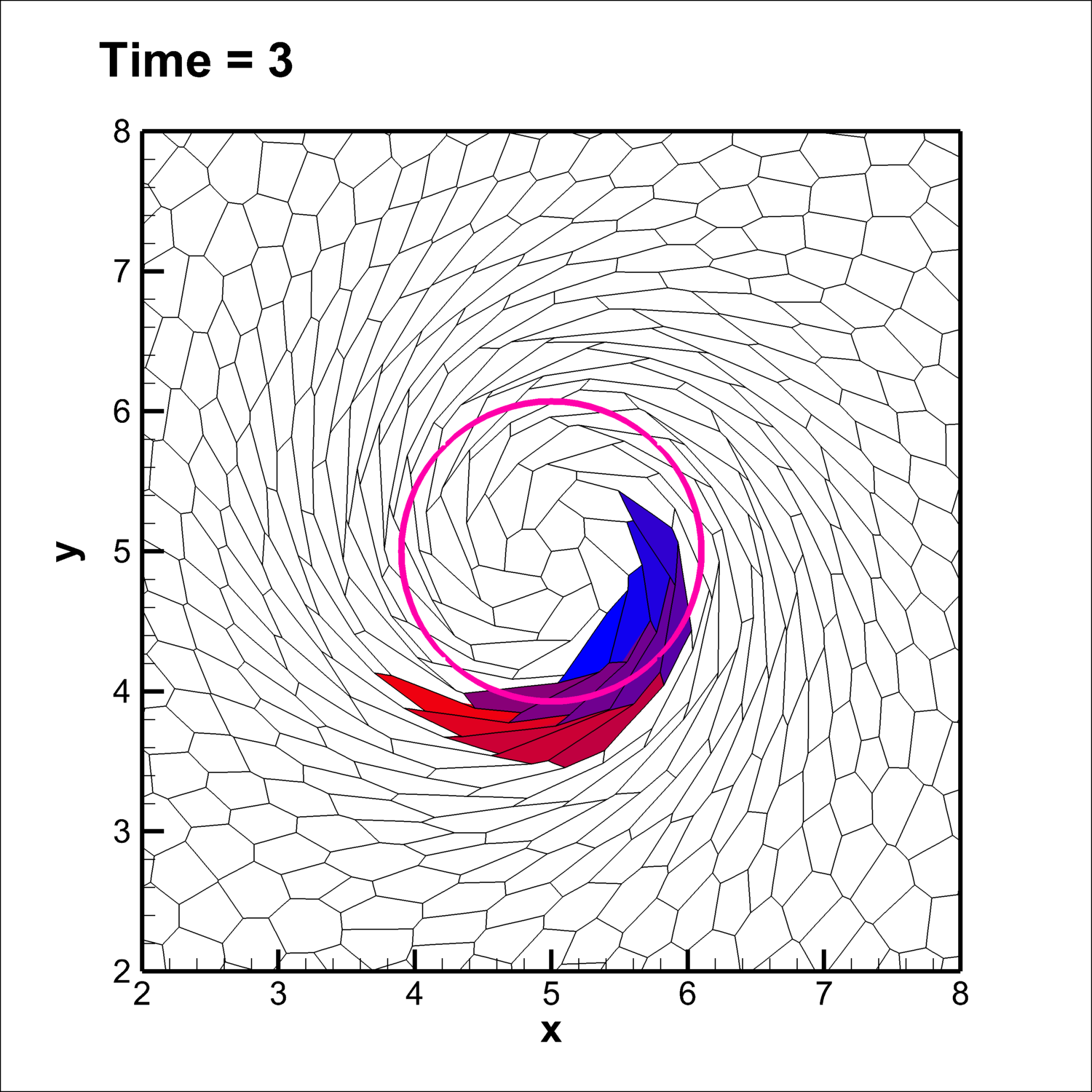}%
    \includegraphics[width=0.2\linewidth]{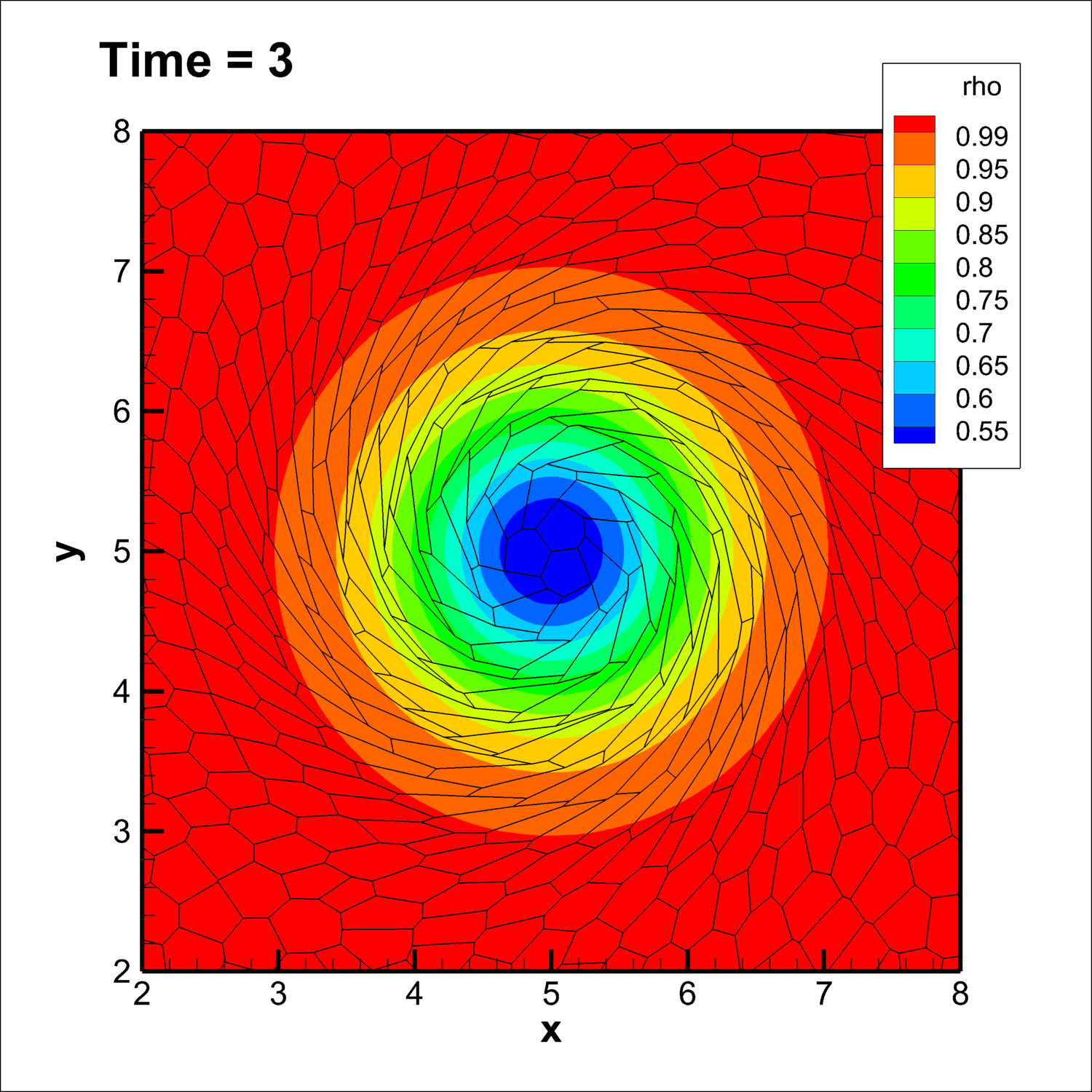}%
    \caption{Stationary rotating vortex solved with a fourth order $P_3P_3$ ALE-DG scheme
        on a moving Voronoi mesh of $957$ elements with \textit{fixed} connectivity. Density contours (left and right) and position of a
         bunch of highlighted elements (middle) are provided at different times. 
         The mesh quality rapidly deteriorates and the simulation ultimately stops at $t\simeq3.5431$ 
         due to tangling elements.}
    \label{fig.ShuDG4_standardALE}
\end{figure*}

\paragraph*{Quality} 
In Figure~\ref{fig.ShuDG4} we plot the density contours and the two-dimensional mesh configuration at various output 
times obtained with our fourth order ALE-DG scheme. We would like to attract the attention to the long duration of the 
simulation and on the high quality of the density profile obtained even after very long simulation times. 
The correct density profile and a high quality mesh are conserved for times that are two orders of magnitude larger 
with respect to standard conforming ALE schemes, where mesh tangling would occur and stop the simulation 
much earlier (see Figure~\ref{fig.ShuDG4_standardALE}). 

The second part of Figure~\ref{fig.ShuDG4} shows the position of a bunch of highlighted elements at different times: 
this makes clear how strong the differential rotation of the mesh elements is. It also highlights the importance of allowing topology changes in the computational grid, 
which needs to provide enough flexibility in order to preserve a high quality mesh over long simulation times. 
Indeed, if the preservation of the connectivity had been imposed, the elements would have been quite distorted 
after only rather short times (see Figure~\ref{fig.ShuDG4_standardALE}).

Finally, we would like to emphasize that generator trajectories are almost perfectly circular even for very long evolution
times (indeed the elements are not collapsing into the center), refer to Figure~\ref{fig.Shu_traj}. 
This is achieved thanks to the computation of the generator trajectories through a fourth order accurate Taylor method. 
In this test, we used the Lloyd-like smoothing method described in Section~\ref{ssec.smoothing}, with $\mathcal{F} = 10^{-4}$.
 
\paragraph*{Independence of the neighbor numbering}  
To show that our algorithm is also completely independent of the space--time neighbor numbering chosen when 
connecting the old mesh to the new one (see Section~\ref{ssec.SpaceTimeConnection}), we have carried out the following test.
In the framework of a third order $P_2P_2$ DG scheme we have simulated the isentropic vortex up to a final 
time of $t=0.5$ on a series of meshes, namely composed by $961$, $1681$, $2601$, $3721$ and $6561$ Voronoi elements 
moving with the exact velocity computed at the generator point of each element.
Then, we have run the algorithm for each mesh configuration by ordering the space--time neighbors in three different 
ways, namely starting first with the first common neighbor, next with the second common neighbor and last with the 
third common neighbor (if existing, otherwise we have used the first one again).

Table~\ref{tab.IndipendenceOfNeighsNumb} shows that not only the order of the algorithm does not depend on the 
neighbor numbering, but also that the final errors are the same up to machine precision.

\begin{figure}[!bp]
    \centering
    \includegraphics[width=0.49\linewidth]{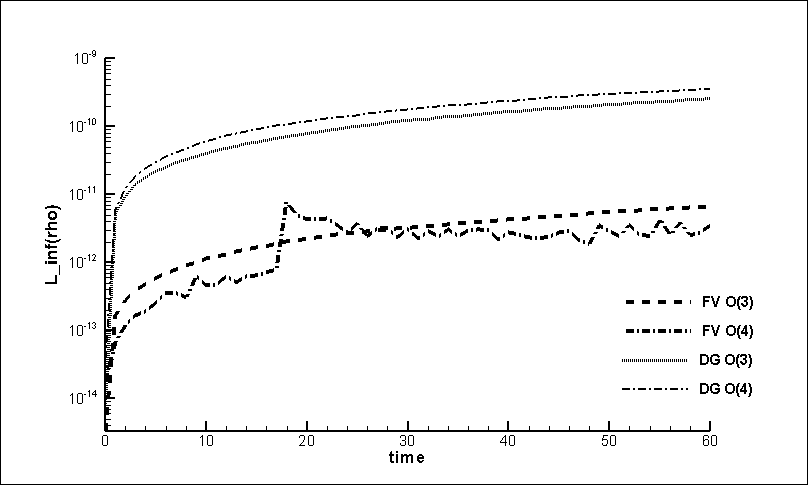}
    \includegraphics[width=0.49\linewidth]{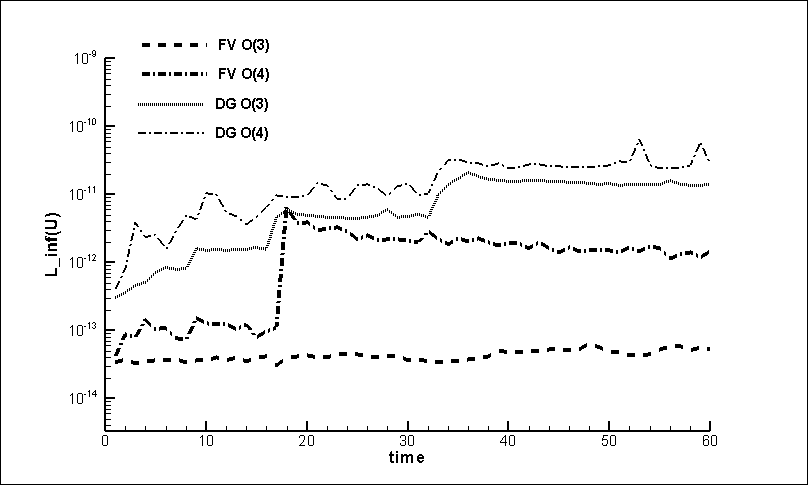}
    \caption{\RIIIcolor{GCL verification problem. In this Figure we report the error in $L-$infinity norm 
    of the error between our numerical results and the exact constant solution for the density $\rho$ 
    (left) and the velocity $U=\sqrt{u^2+v^2}$ (right) for four different methods, namely a third and 
    fourth order Finite Volume scheme and a third and fourth order Discontinuous Galerkin scheme. 
    We underline that the final time of $t_f=60$ is reached after \textit{thousands} 
    of time steps and involves the generations of \textit{thousands} of sliver elements, as reported in 
    Table~\ref{tab.GCLTest_CPUtimes}.    }}
    \label{fig.gclrhoconstant}
\end{figure}

\begin{figure}[!bp]
    \centering
    \includegraphics[width=0.20\linewidth]{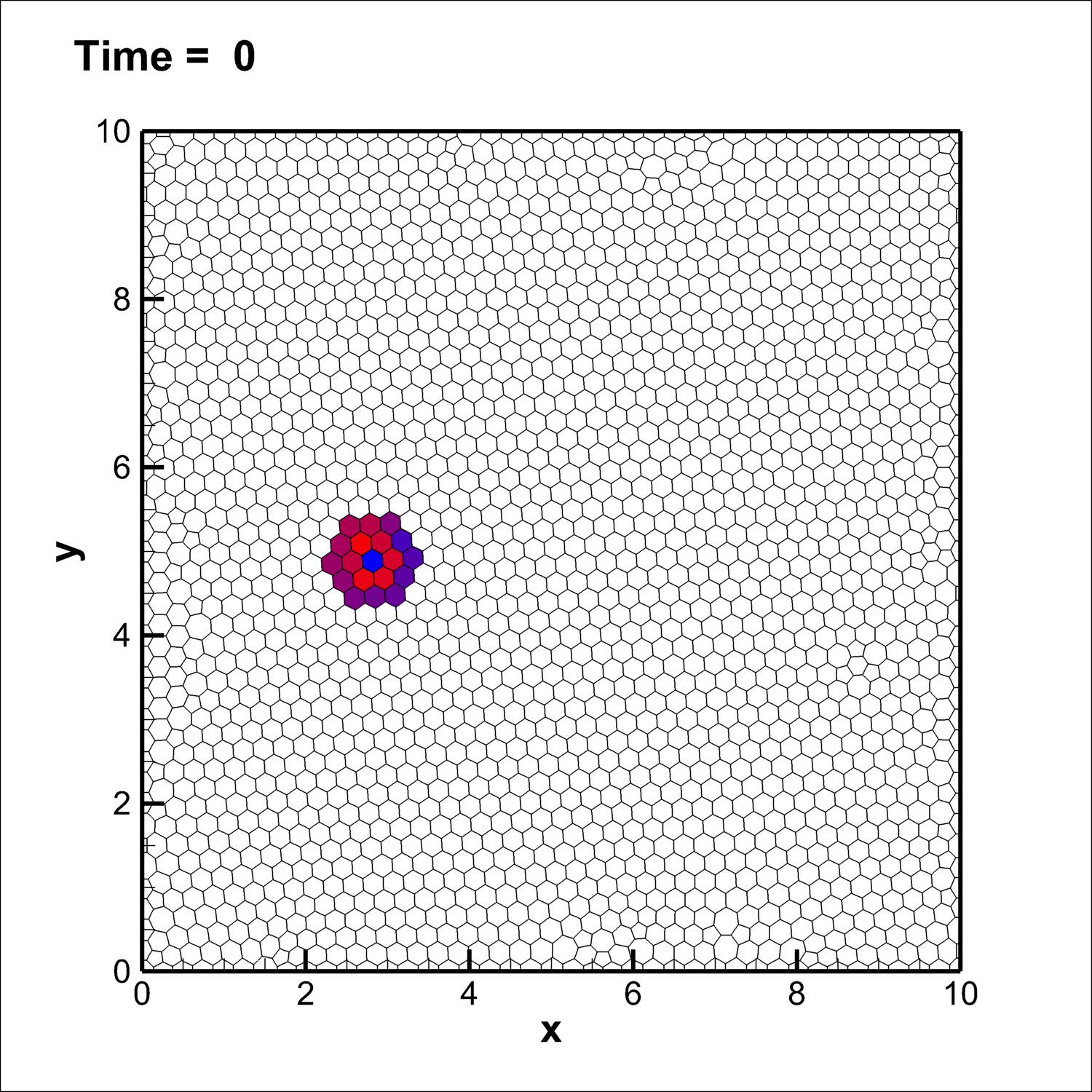}%
    \includegraphics[width=0.20\linewidth]{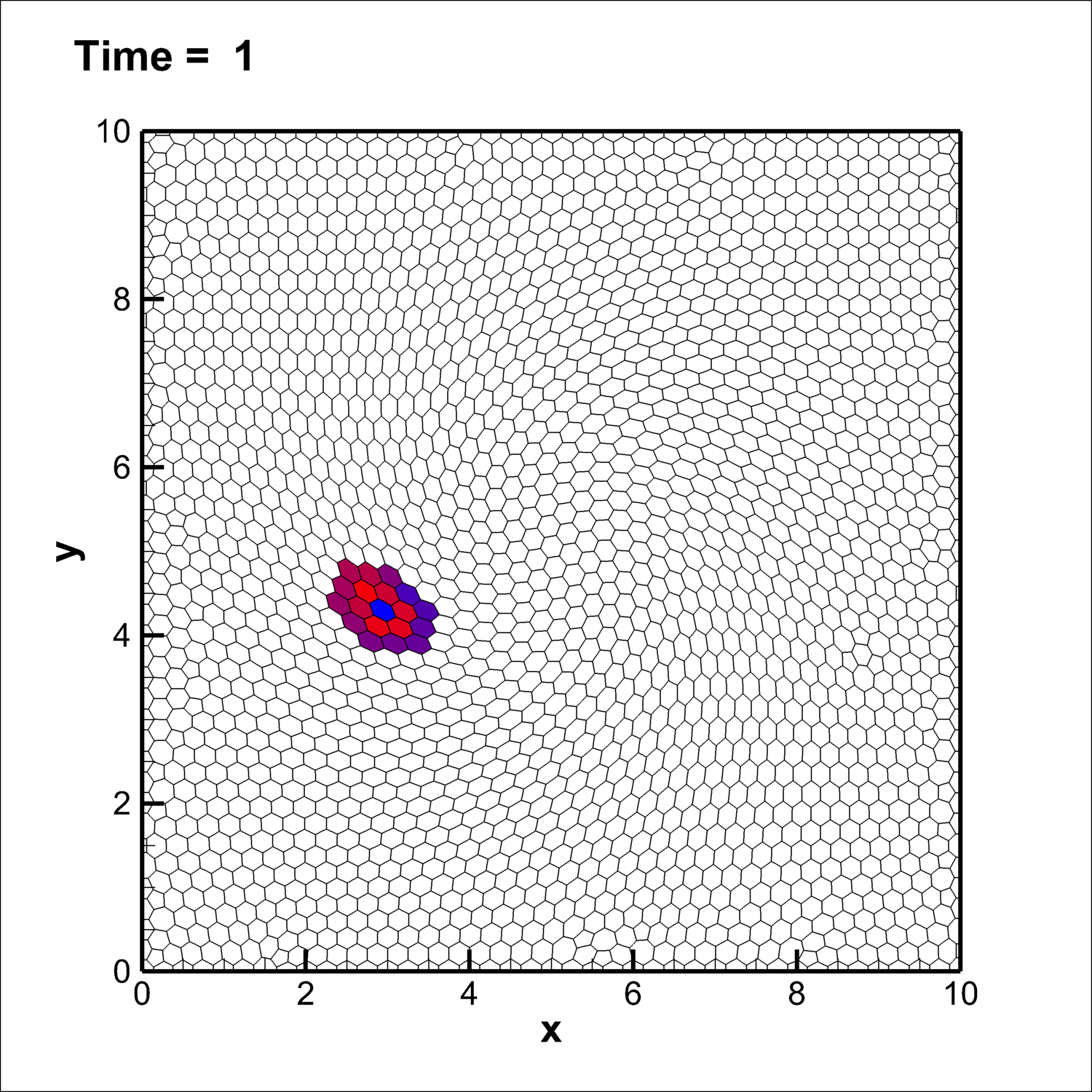}%
    \includegraphics[width=0.20\linewidth]{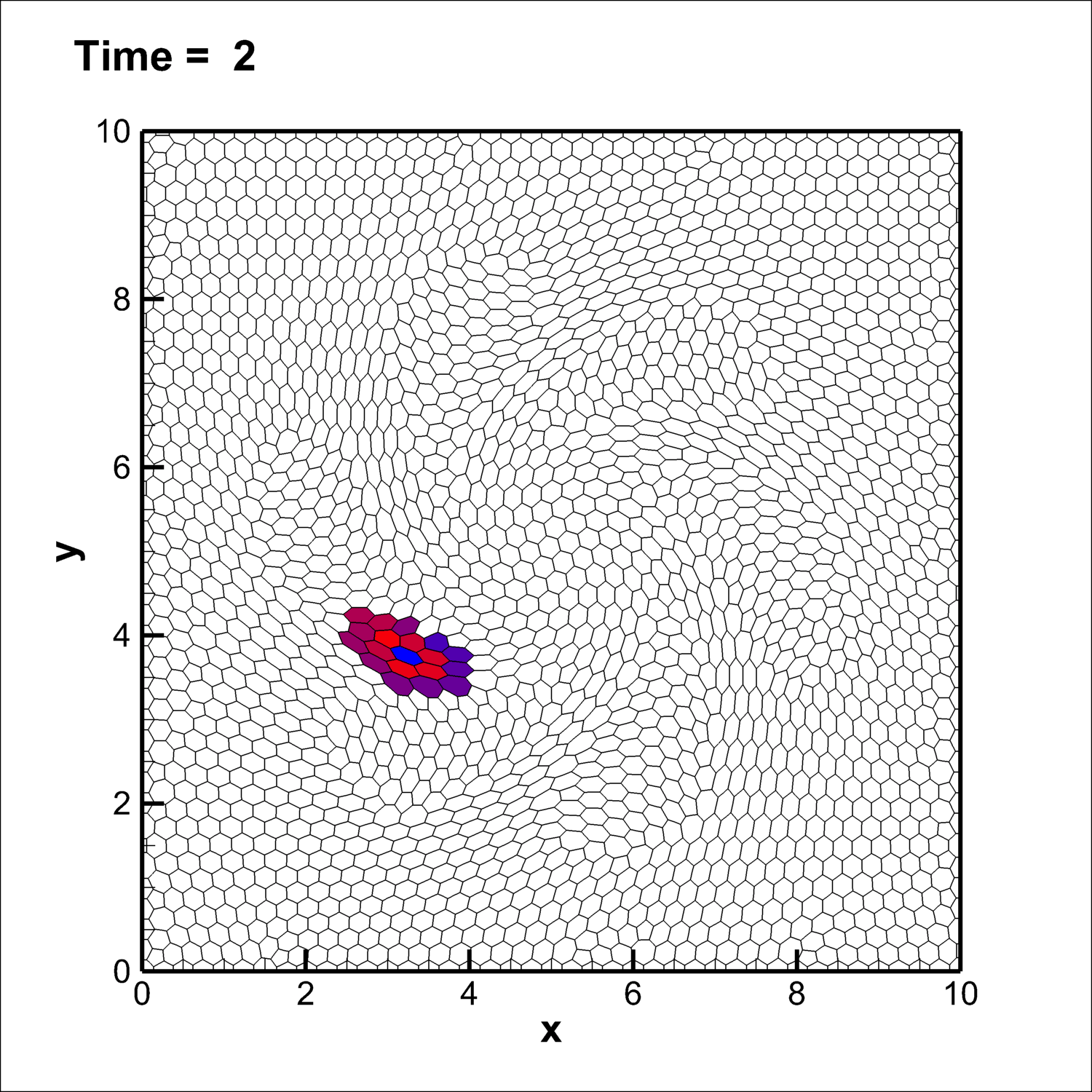}%
    \includegraphics[width=0.20\linewidth]{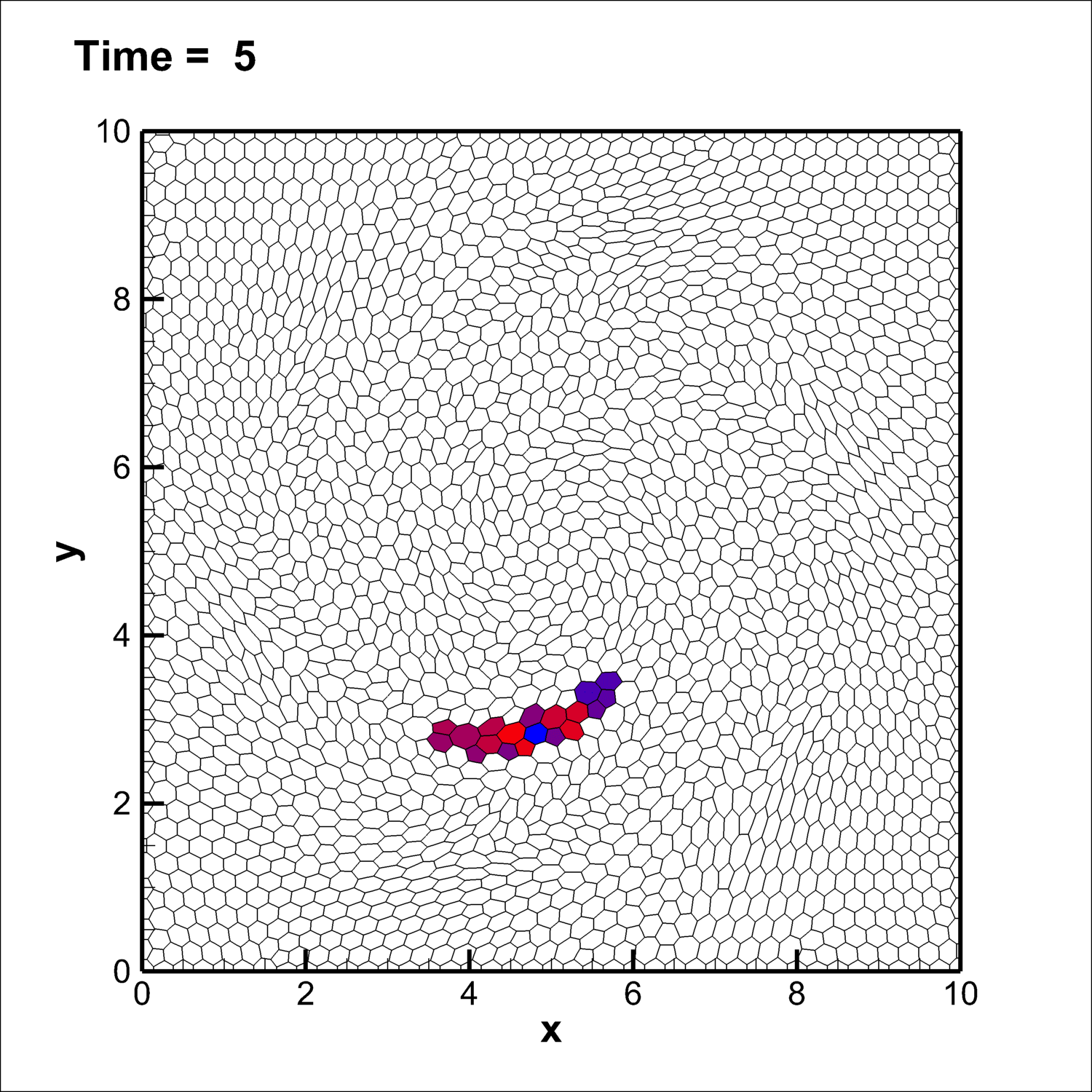}%
    \includegraphics[width=0.20\linewidth]{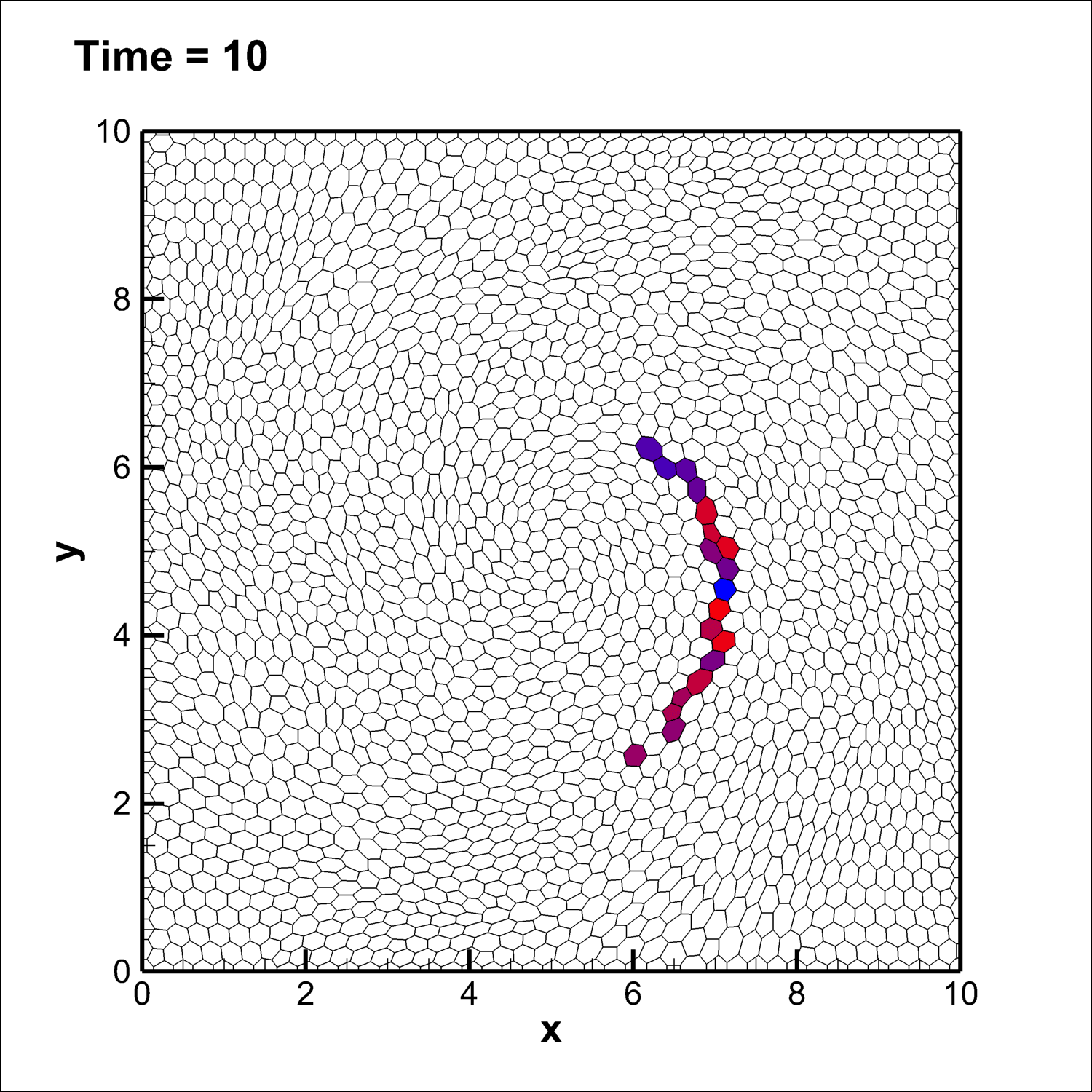}\\[-1pt]
    \includegraphics[width=0.20\linewidth]{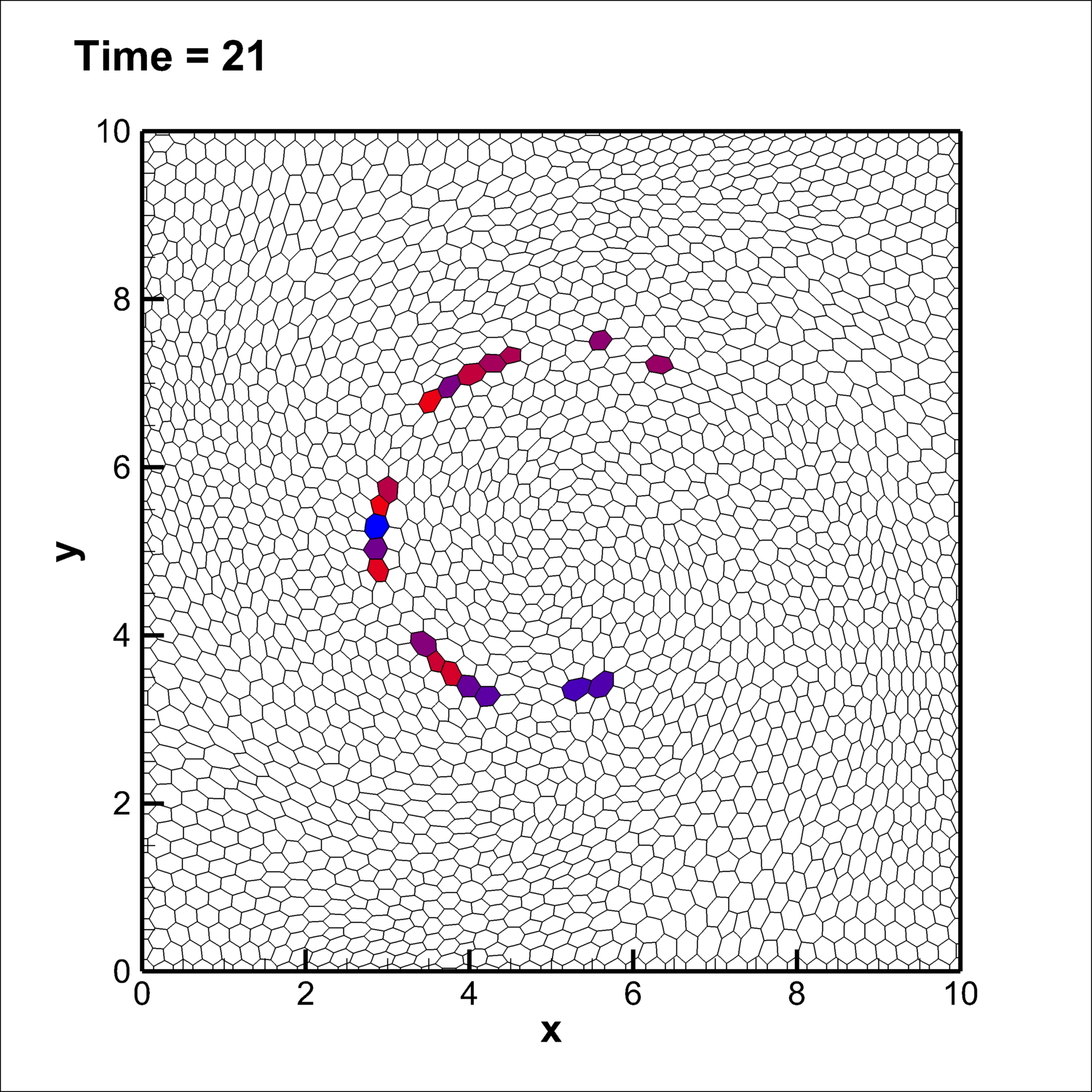}%
    \includegraphics[width=0.20\linewidth]{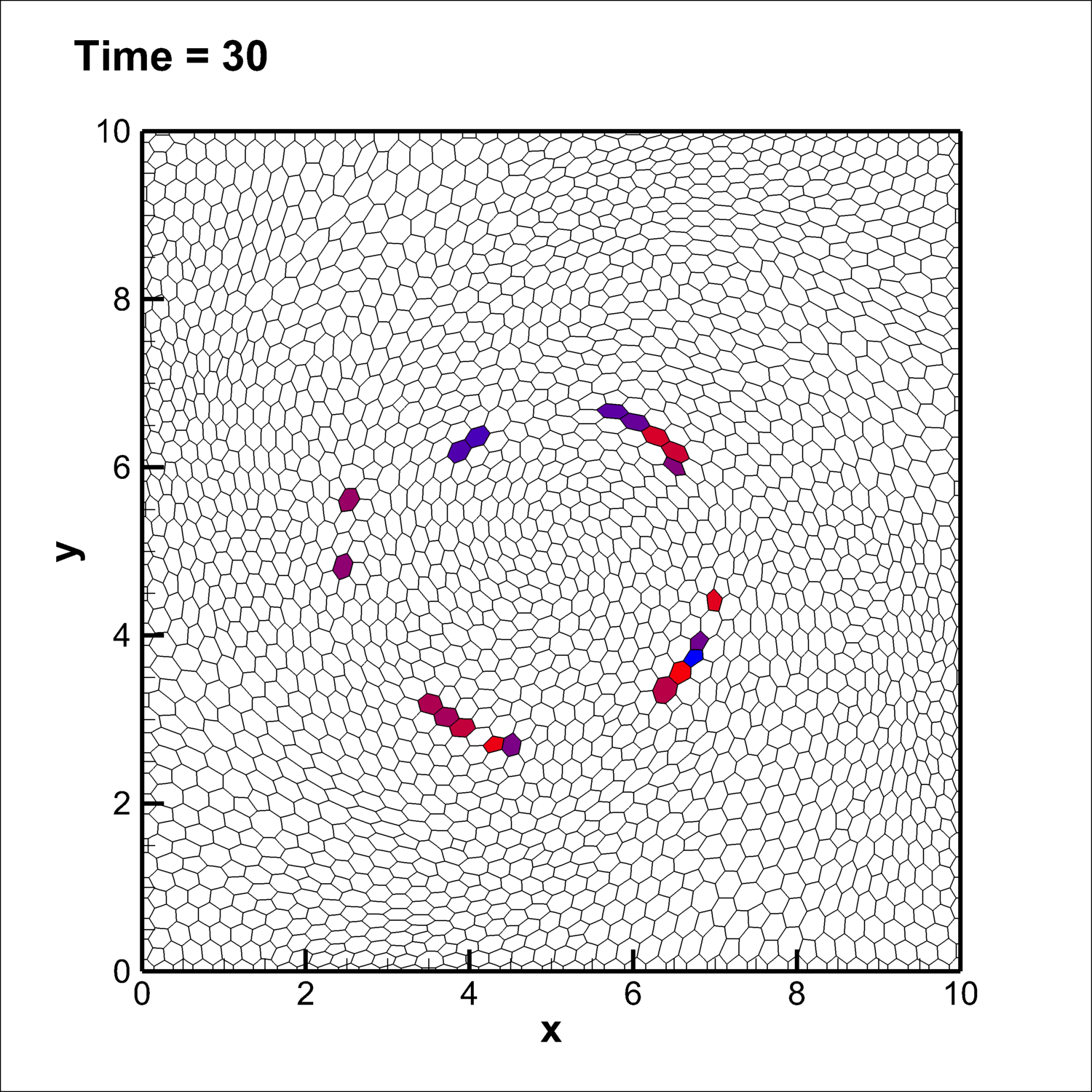}%
    \includegraphics[width=0.20\linewidth]{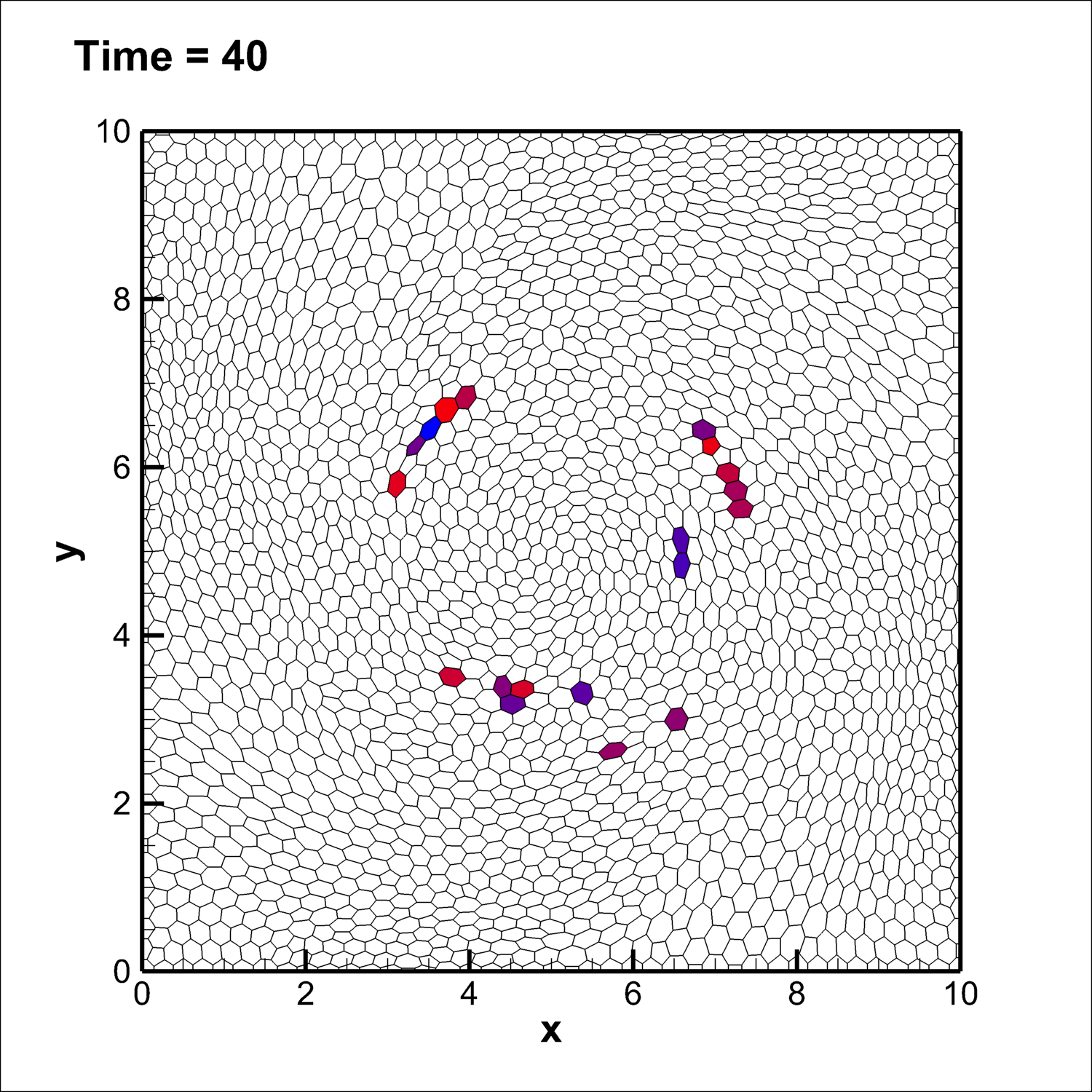}%
    \includegraphics[width=0.20\linewidth]{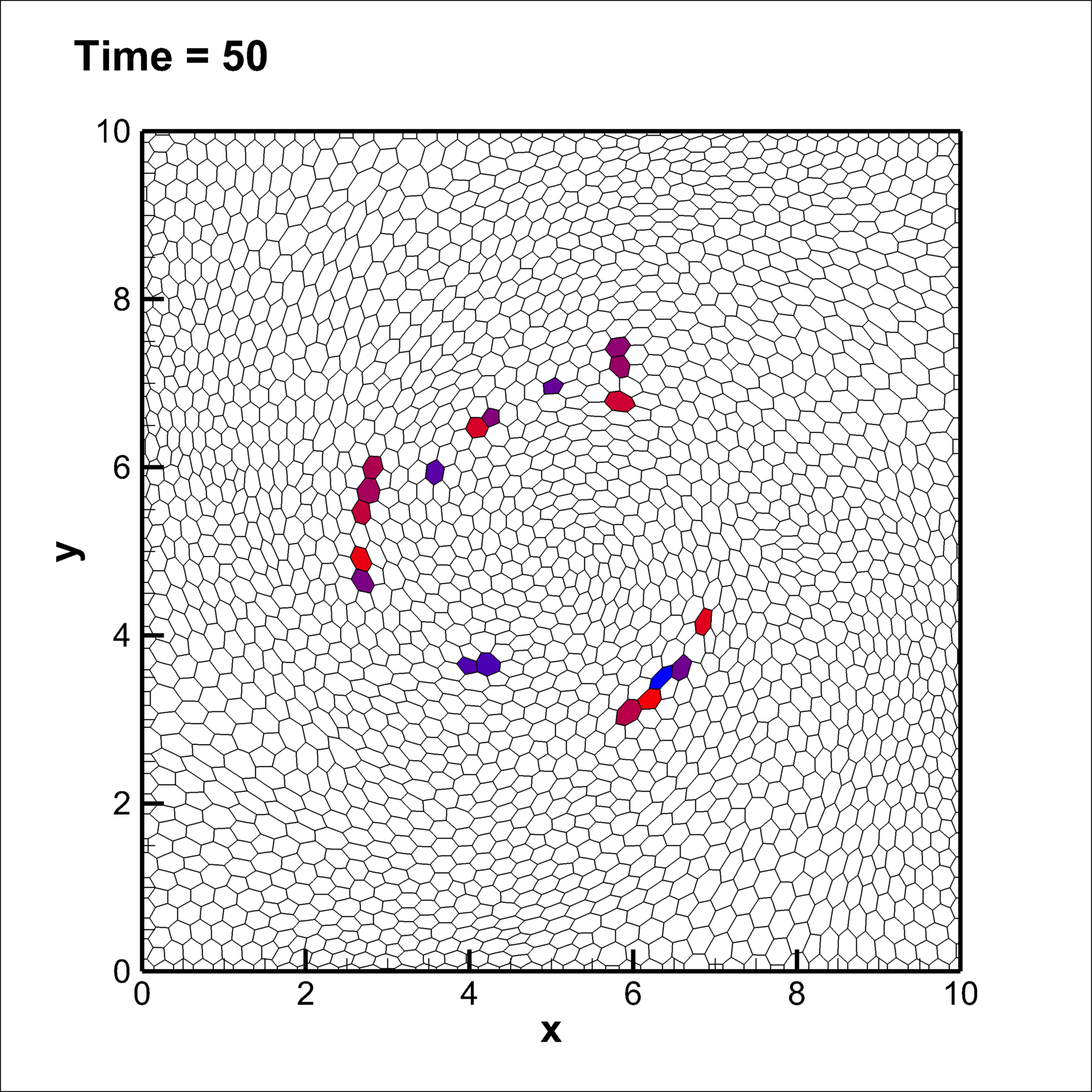}%
    \includegraphics[width=0.20\linewidth]{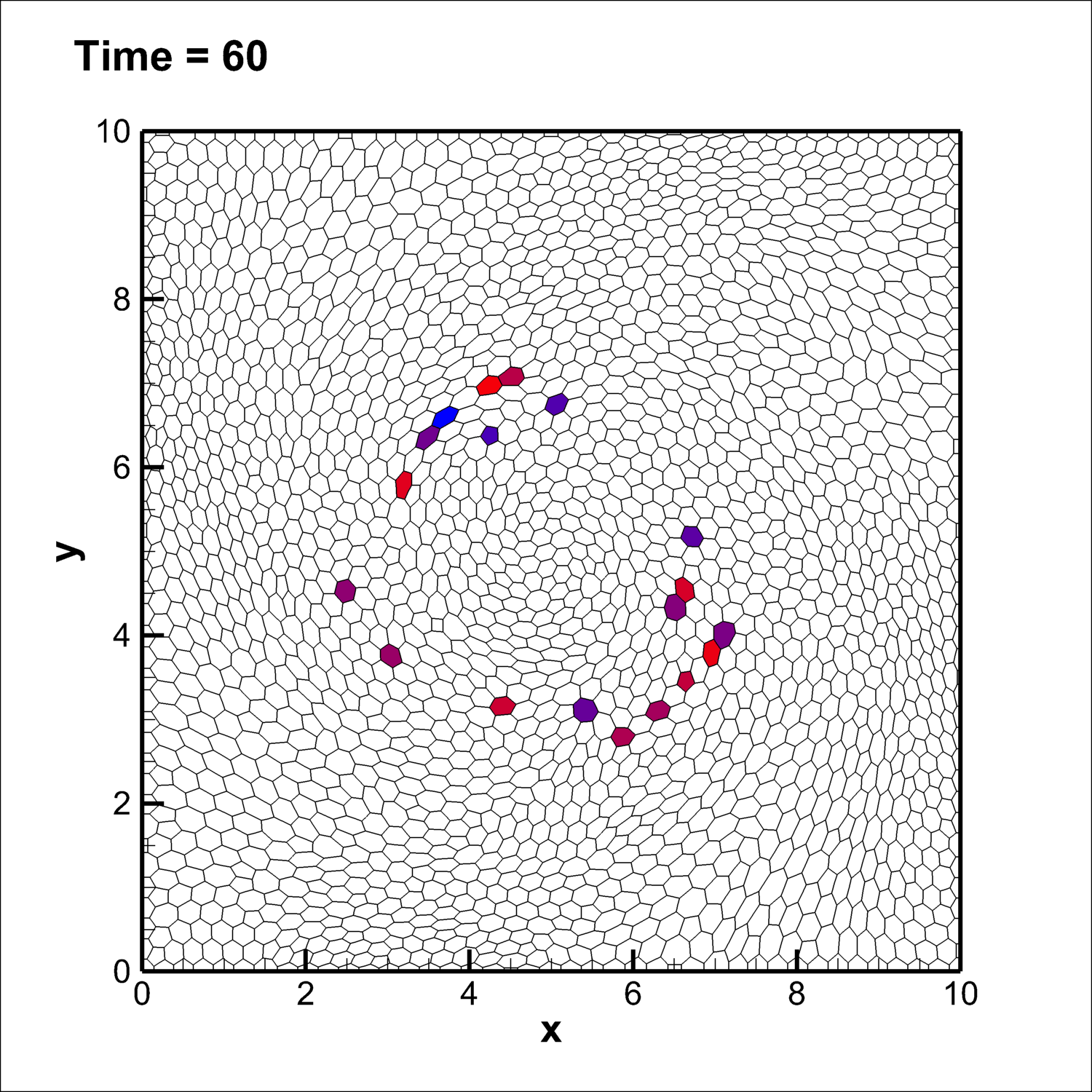}%
    \caption{\RIIIcolor{GCL verification problem. In this Figure we report the position of a bunch of elements initially neighbors subjected to
     the velocity field~\eqref{eq.VelVorticalSinu} at different times. 
     The Voronoi tessellation is regenerated at any of the $33019$ time steps of the employed third order $P_2P_2$ DG schemes. 
     The test uses the Lloyd-like smoothing algorithm with $\mathcal{F} = 10^{-3}$.
      In total $19392$ sliver elements have been treated during this simulation (and the MOOD procedure 
      of Section~\ref{ssec.MOOD} has never been activated).}}
    \label{fig.GCLTest_bunchOfElements}
\end{figure}

\RIIIcolor{
\subsubsection{Numerical verification of the GCL property}
\label{test.GCLTest}

In order to verify that our schemes satisfy the GCL property up to machine precision we consider the following standard test case, see for example~\cite{ma2015geometric}. 
 
The initial condition for this test is given by a \textit{constant state}, namely  $\Q_0=(\rho,u,v,p)=(1,0,0,1)$ that should remain constant even with moving meshes if the GCL condition is satisfied. The initial computational domain is the square $\Omega=[0;10]\times[0;10]$ covered with a mesh of $1979$ Voronoi elements and wall boundary conditions are set everywhere.

Then, in agreement with the ALE framework, the mesh is moved with a completely arbitrary velocity field (thus not with the fluid velocity which would be zero in this test case).
In particular, we have chosen a vortical velocity field varying in a sinusoidal fashion, given in the form
\be
\v(x,y) = \left(
-\sin\left( \frac{2\pi}{\ell} \left(y-y_0\right) \right) \cos\left( \frac{\pi}{\ell} \left(x - x_0\right) \right) \exp(-kr),\ 
\cos\left( \frac{\pi}{\ell} \left(y-y_0\right) \right) \sin\left( \frac{2\pi}{\ell} \left(x - x_0\right) \right) \exp(-kr) \right),
\label{eq.VelVorticalSinu}
\ee 
where $\x_0 =(x_0, y_0) =5, \ell =10, k=0.1$ and $r=\sqrt{(x-x_0)^2+(y-y_0)^2}$.

In Figure~\ref{fig.gclrhoconstant}  we show the error between our numerical results for density and velocity, with respect to the exact constant solution, 
obtained with our FV and DG schemes of order three and four, namely the $P_0P_2$, $P_0P_3$,  $P_2P_2$ and  $P_3P_3$ schemes.
We emphasize that these results are obtained over long times ($60$ time units, and three complete revolutions of a reference point), thus after \textit{thousands} of time steps.
Furthermore, topology changes regularly occur during the mesh motion: 
the number of sliver elements generated per simulation over the total number of time steps is reported in Table~\ref{tab.GCLTest_CPUtimes}, 
and a bunch of initially adjacent elements is plotted in Figure~\ref{fig.GCLTest_bunchOfElements} at different times, 
in order to make clear that elements really change their topology and connectivity during the simulation.

This test case proves numerically that the GCL property is satisfied by our scheme even when topology changes occur and sliver elements appear, 
which is indeed a property that our schemes satisfy \textit{by construction}, since PDE integration is always performed over closed space--time control volumes, see \cite{Lagrange3D} for a formal proof.
}

\begin{table*}[!tp]
    \caption{\RIVcolor{Scheme statistics. In this Table we report some statistics 
    on the computational cost of reaching the final time $t_f=60$ with four different schemes. 
    In the second to fourth column we show the the number of total time steps needed to reach the final time, 
    the number of sliver elements that have been originated during the simulations due to 
    the occurred topology changes and the number of time steps that have been repeated through the MOOD loop described in Section~\ref{ssec.MOOD}.
            In the fifth column we report the percentage of computational time employed by mesh 
            regeneration and space-time connectivity generation, in the sixth column the 
            percentage of time needed for the predictor-corrector step on standard elements, and 
            in the last column the percentage of time spent on sliver elements. It is evident that 
            the cost due to mesh rearrangement and sliver computations is minimal.
            (For what concerns FV schemes, the time for spatial reconstruction is not included in the 
            third column, in order to facilitate the comparison between the cost on standard elements 
            and sliver elements, for which spatial reconstruction is not performed.) 
    }} 
    \label{tab.GCLTest_CPUtimes}  
    \centering
        \begin{tabular}{c|ccc|ccc} 
            \hline
            Method & time steps & slivers & restarts & Mesh $\%$  & $P_NP_M$ standard $\%$  & $P_NP_M$ sliver $\%$  \\
            \hline
            FV $\mathcal{O}(3)$ & 5524  & 21930 &  5 & 0.31   & 56.03   &  0.01  \\
            FV $\mathcal{O}(4)$ & 4896  & 19819 &  4 & 0.31   & 56.00   &  0.01  \\
            DG $\mathcal{O}(3)$ & 33019 & 19392 &  0 & 1.29   & 91.07   &  0.004 \\
            DG $\mathcal{O}(4)$ & 55496 & 18995 &  2 & 0.25   & 96.31   &  0.001 \\
            \hline  
        \end{tabular}       
\end{table*}

\RIVcolor{Finally, in Table~\ref{tab.GCLTest_CPUtimes} we report the percentage of computational times employed by 
	i) all the the procedures necessary to rebuild a new mesh at each time step, 
    namely the mesh regeneration of Section~\ref{sec.MeshEvolution}, the mesh optimization of Section~\ref{ssec.smoothing}, 
    and the construction of the space time connectivity of Sections~\ref{ssec.SpaceTimeConnection} and~\ref{ssec.Sliver};
	ii) the predictor-corrector step of Section~\ref{sec.NumMethod_ALE_FV-DGscheme} performed on standard elements;
	iii) the predictor-corrector step of Section~\ref{sec.NumMethod_ALE_FV-DGscheme} performed on sliver elements.
		We remark that step i) actually consists in a complete regeneration of a new configuration without exploiting 
        the previous one, which would decrease the computational cost of the geometric part of the code. 
	Nevertheless, one can see that the computational cost of the geometric part of our scheme is minimal and 
    does not affect the final cost of the entire algorithm.		
}

\subsubsection{Explosion problem}
\label{test.ExpPb}

\begin{table*}[!p]
    \caption{\RIVcolor{Explosion problem. In this Table we report the number of total time steps needed to reach the final time $t_f=0.25$ with our $P_0P_3$ and $P_2P_2$ schemes, the number of sliver elements and the number of time steps that have been repeated through the MOOD loop described in Section~\ref{ssec.MOOD}. 
            Moreover, we report the percentage of computational time employed by mesh regeneration and space time connectivity generation, by the predictor-corrector step on standard elements and on sliver elements. The other part of the computational time is mostly spent on the reconstruction procedure for what concerns the FV scheme, and on the limiter for what concerns the DG scheme.}} 
    \label{tab.explosion_percentage}  
    \centering  
        \begin{tabular}{c|ccc|ccc} 
            \hline
            Method & time steps & slivers & restarts & Mesh $\%$  & $P_NP_M$ standard $\%$  & $P_NP_M$ sliver $\%$   \\
            \hline
            FV $\mathcal{O}(4)$ & 150   & 1785  & 0  & 0.16     & 79.69      & 0.003                 \\
            DG $\mathcal{O}(3)$ & 883   & 2238  & 3  & 0.70     & 64.59      & 0.002                 \\
            \hline  
        \end{tabular}       
\end{table*}

\begin{figure}[!p]
    \centering
    \includegraphics[width=0.45\linewidth]{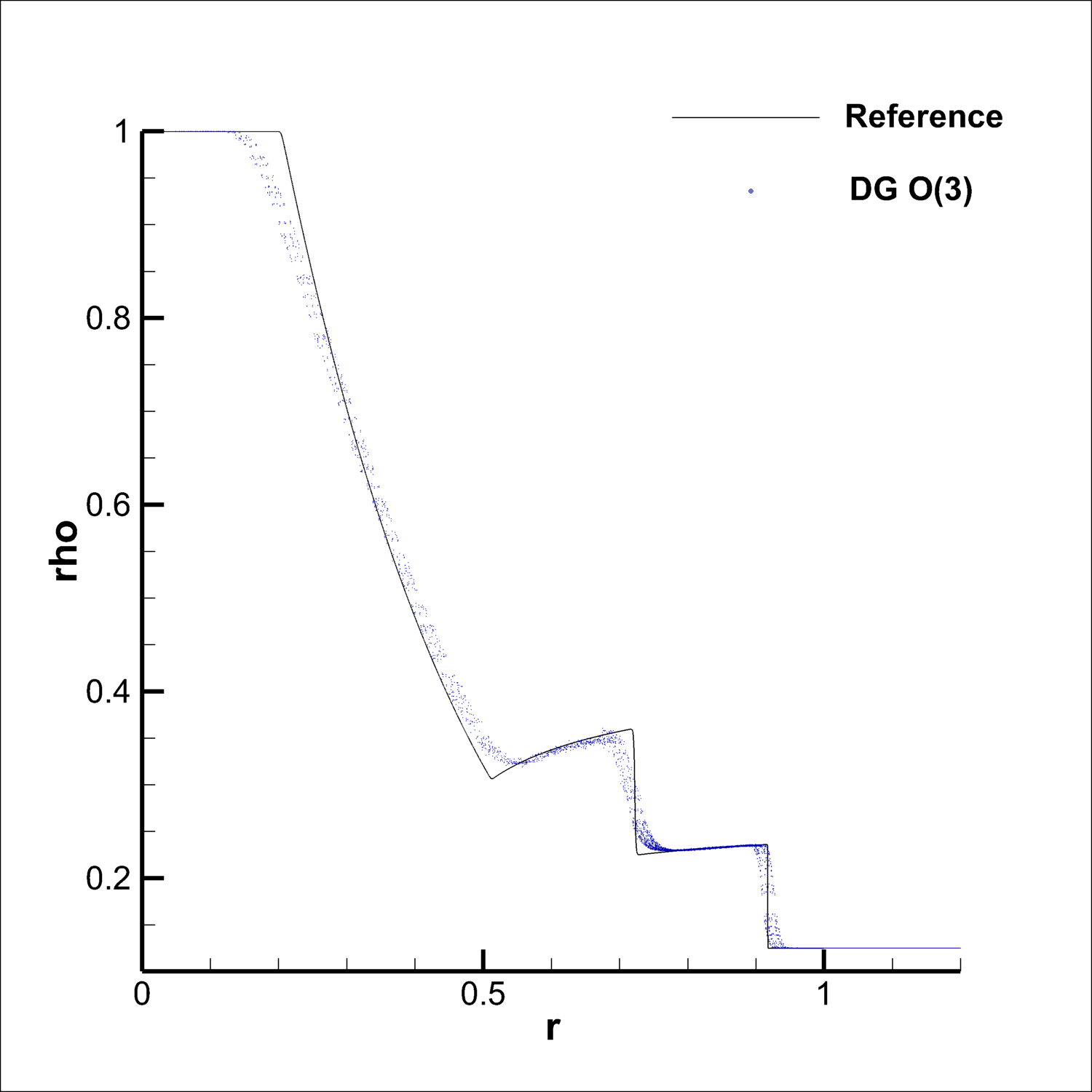} \quad
    \includegraphics[width=0.45\linewidth]{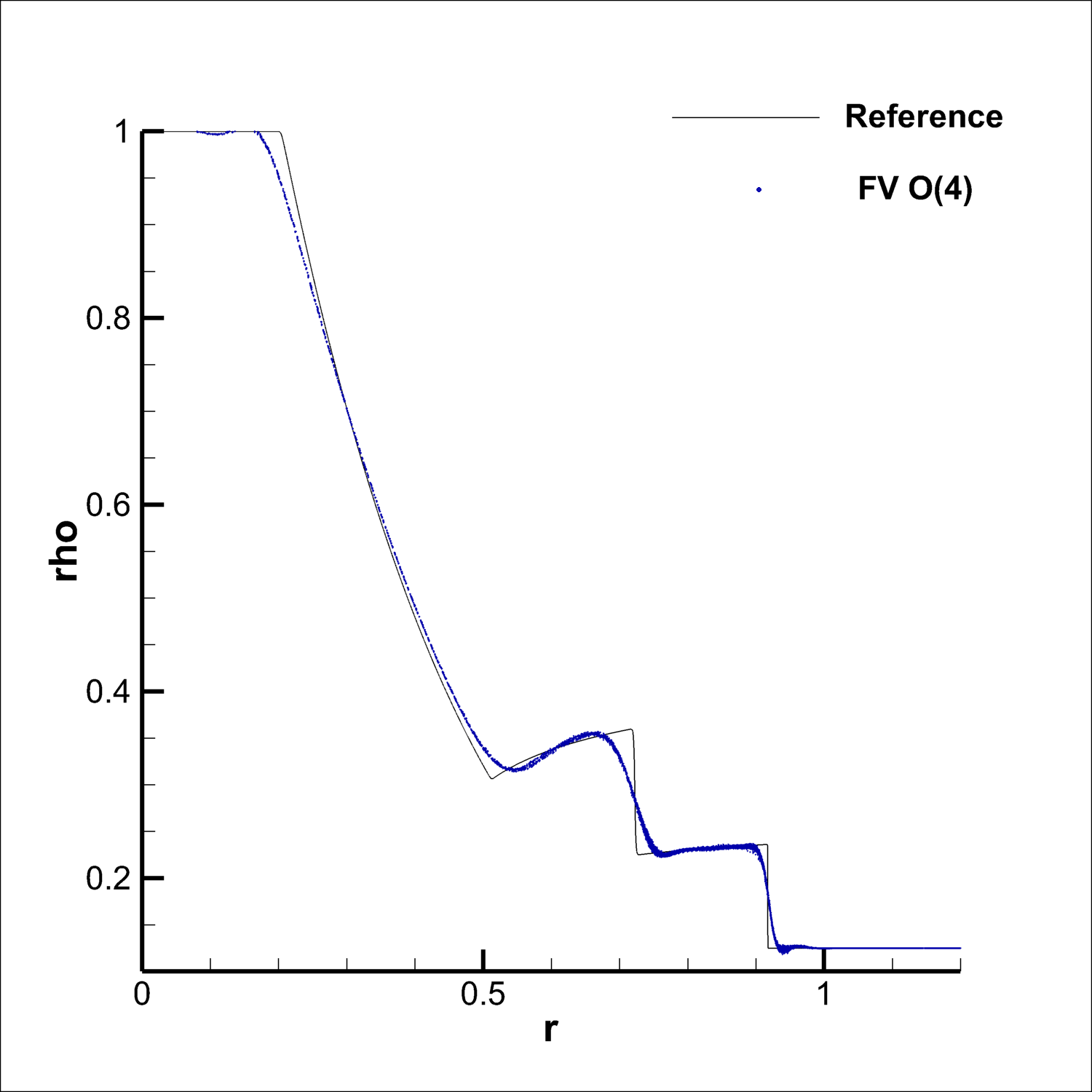}
    \caption{\RIIIcolor{Explosion problem: we compare the scatter plot of our numerical results (blue dots) 
    with the reference solutions (line) at time $t_f=0.25$. Left: results obtained with our $P_2P_2$ DG scheme on a moving Voronoi mesh of $9805$ elements 
    \RVIcolor{(corresponding to $58830$ total DOFs)}. 
    Right: results obtained with our $P_0P_3$ FV scheme on a moving Voronoi mesh of $19856$ elements \RVIcolor{(corresponding 
    to $19856$ total spatial DOFs, treated with a fourth order method)}. The represented values (squares) are obtained 
    from a cut of our numerical solutions along $y=0$. Both tests use the Laplacian smoothing algorithm with $\mathcal{F} = 10^{-3}$.}}
    \label{fig:ExpPb_P2P2_100x100_RhoScatterVsExact}
\end{figure}

\begin{figure}[!p]
    \centering
    \includegraphics[width=0.365\linewidth]{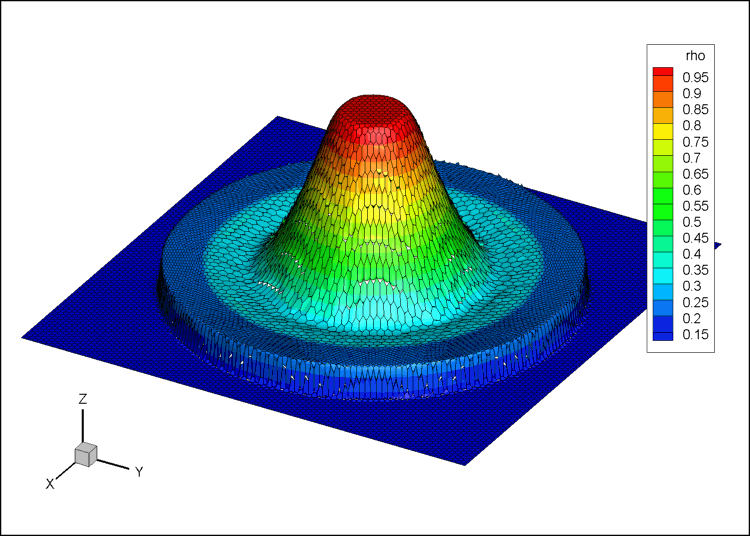}
    \includegraphics[width=0.365\linewidth]{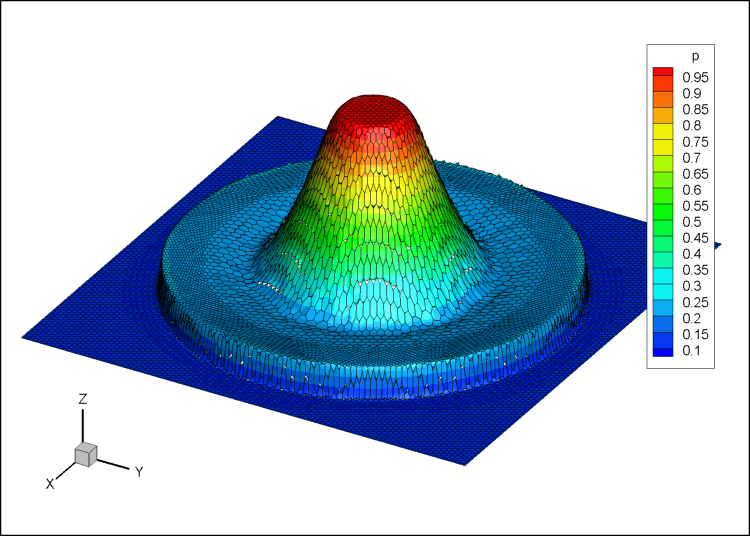}
    \includegraphics[width=0.26\linewidth]{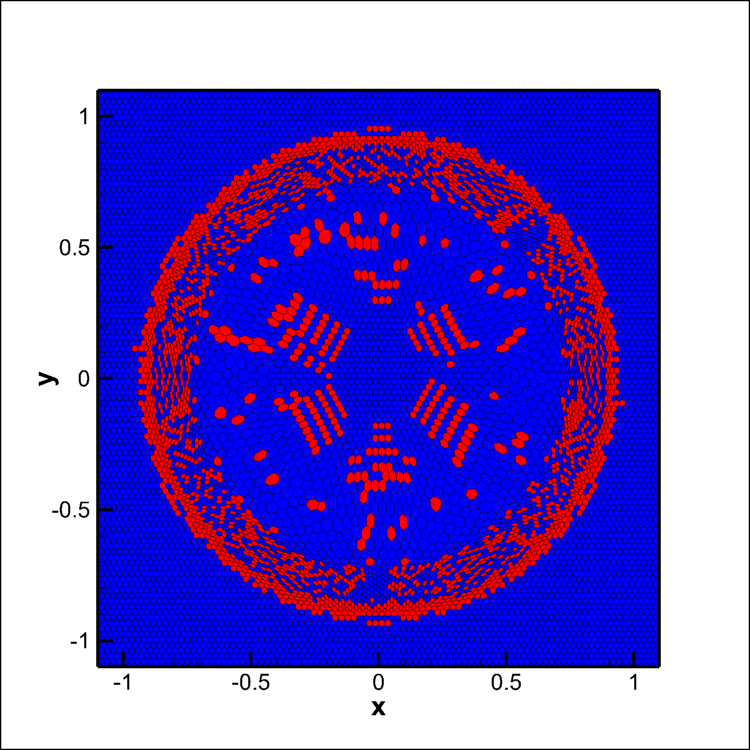}
    \caption{Explosion problem solved with our $P_2P_2$ method on a moving Voronoi mesh of $9805$ elements. 
    The test uses the Laplacian smoothing algorithm with $\mathcal{F} = 10^{-3}$. We show the density profile (left), the pressure profile (center) and we color in red the cells on which the limiter is active (right) at time $t=0.25$. 
}
\label{fig.ExpPb_P2P2_DG}
\end{figure}

The explosion problems can be seen as a multidimensional extension of the classical Sod test case.
Here, we consider as computational domain a square of dimension $[-1.1;1.1]\times[-1.1;1.1]$, and the initial condition is composed of two  different states, separated by a discontinuity at radius $r_d=0.5$ 
\be
\begin{cases} 
\rho_L = 1, \quad \mathbf{u}_L = 0, \quad p_L = 1, \quad  &\left\| \x \right\| \leq r_d \\
\rho_R = 0.125, \quad \mathbf{u}_R = 0, \quad  p_R = 0.1, \quad & \left\| \x \right\| > r_d .
\end{cases}
\ee
The final time is chosen to be $t_f=0.25$, so that the  shock  wave does not cross the external boundary of the domain, where a wall boundary condition is set.  
We run this problem with two different configurations.
\begin{description}
	\item[(a)]
	In the first case we use a third order $P_2P_2$ DG scheme on a mesh of $9805$ Voronoi elements \RVIcolor{(corresponding to $58830$ total DOFs)}.  
	The results are depicted in Figure~\ref{fig.ExpPb_P2P2_DG}. In particular, one can notice that the limiter activates in proximity of the shock waves where it is indeed essential, and only on a handful of other elements. 
	\smallskip
	\item[(b)]
	Then, we test our FV algorithm by employing a fourth order $P_0P_3$ scheme on a finer mesh of $19856$ Voronoi elements \RVIcolor{(corresponding to $19856$ total DOFs)}. 
\end{description} 
In both  cases, we can observe a good agreement between the numerical results and
the reference solution. The non perfect symmetry is justified by the non symmetric initial meshes. 
Statistics on generation of sliver elements, MOOD restarts, and computational times are reported in Table~\ref{tab.explosion_percentage}.

As in~\cite{Lagrange3D, ToroBook}, a reference solution can be obtained by making use of the rotational symmetry of the problem and by solving a reduced one-dimensional system with geometric source terms using a classical second order TVD scheme on a very fine one-dimensional mesh. The comparison between our numerical solutions and the reference solution is given with 
a scatter plot in Figure~\ref{fig:ExpPb_P2P2_100x100_RhoScatterVsExact}.
\RVIcolor{In order to obtain a similar resolution, the FV scheme needs one order more of accuracy w.r.t. the DG scheme and a finer mesh as well. 
We remark that for this comparison we have not exactly matched the total number of DOFs of the FV and the DG simulations because, 
even if the DOFs of a FV scheme are only one per cell, the FV algorithm involves a spatial reconstruction procedure
 (of the necessary order of accuracy) which increases the resolution (and the cost) of the scheme.}

We would like to underline that this test problem involves three different waves, therefore it allows 
each ingredient of our scheme to be properly checked. Indeed, we have
\begin{itemize}
	\item
	one cylindrical shock wave that is running towards the external boundary: our scheme does not exhibit spurious oscillations thanks to the CWENO reconstruction, in the case (b), and to the \textit{a posteriori} sub--cell finite volume limiter, in case (a);
	\item
	a rarefaction fan traveling in the opposite direction, which is well captured thanks to the high order of accuracy;
	\item
	an outward-moving contact wave in between, which is well resolved thanks to the Lagrangian framework of our scheme, 
    in which the mesh moves together with the fluid flow.
\end{itemize}

\subsubsection{Sedov problem}
\label{test.Sedov}

\begin{figure}[!p]
	\centering
	\includegraphics[width=0.32\linewidth]{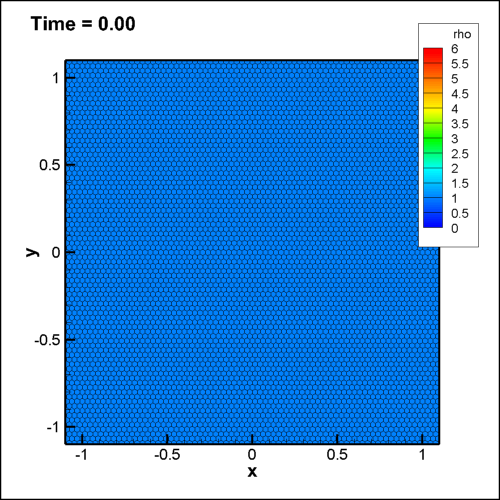}%
	\includegraphics[width=0.32\linewidth]{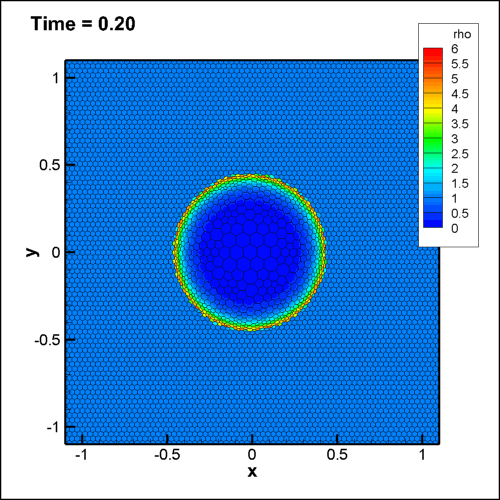}%
	\includegraphics[width=0.32\linewidth]{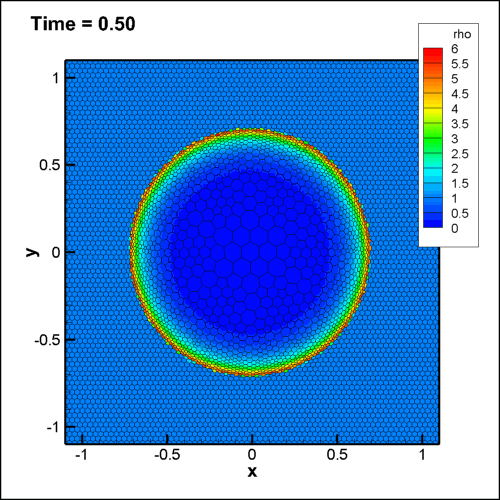}\\[-1pt]
	\includegraphics[width=0.32\linewidth]{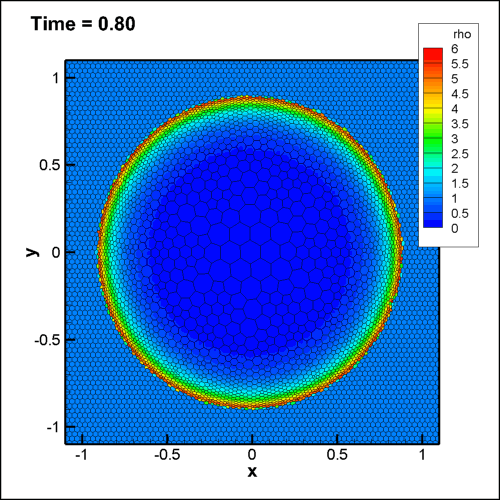}%
	\includegraphics[width=0.32\linewidth]{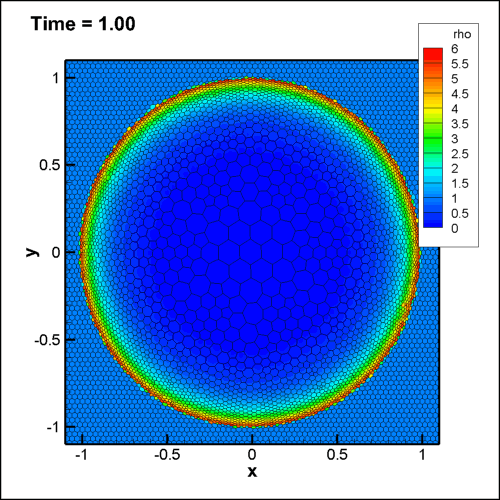}%
	\includegraphics[width=0.32\linewidth]{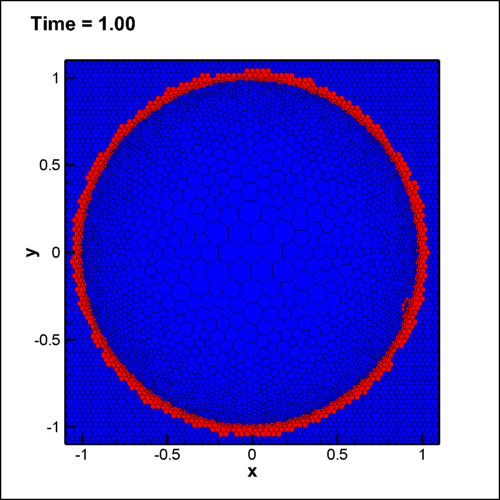}%
	\caption{Sedov problem solved with our $P_2P_2$ scheme on a moving Voronoi mesh of $7234$ elements. 
    The test uses the Laplacian smoothing algorithm with $\mathcal{F} = 10^{-2}$. We depict the density profile and the mesh configuration at times $t=0,0.2, 0.5, 0.8, 1$ and in the last images we show in red the cells on which the limiter is activated.}
	\label{fig:sedovp1p180x800-0}
\end{figure}

\begin{figure}[!p]
	\centering
	\includegraphics[width=0.45\linewidth]{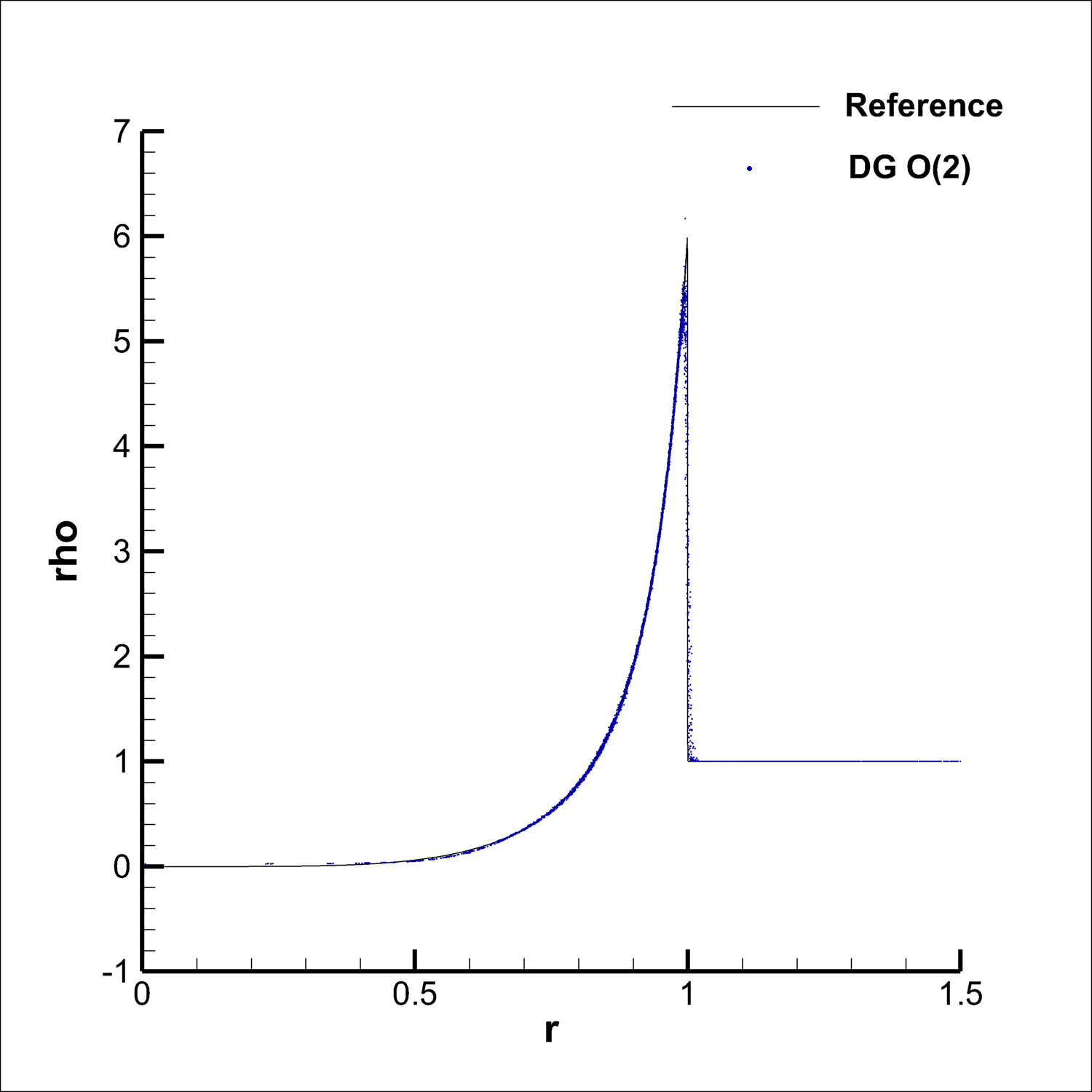} \quad
	\includegraphics[width=0.45\linewidth]{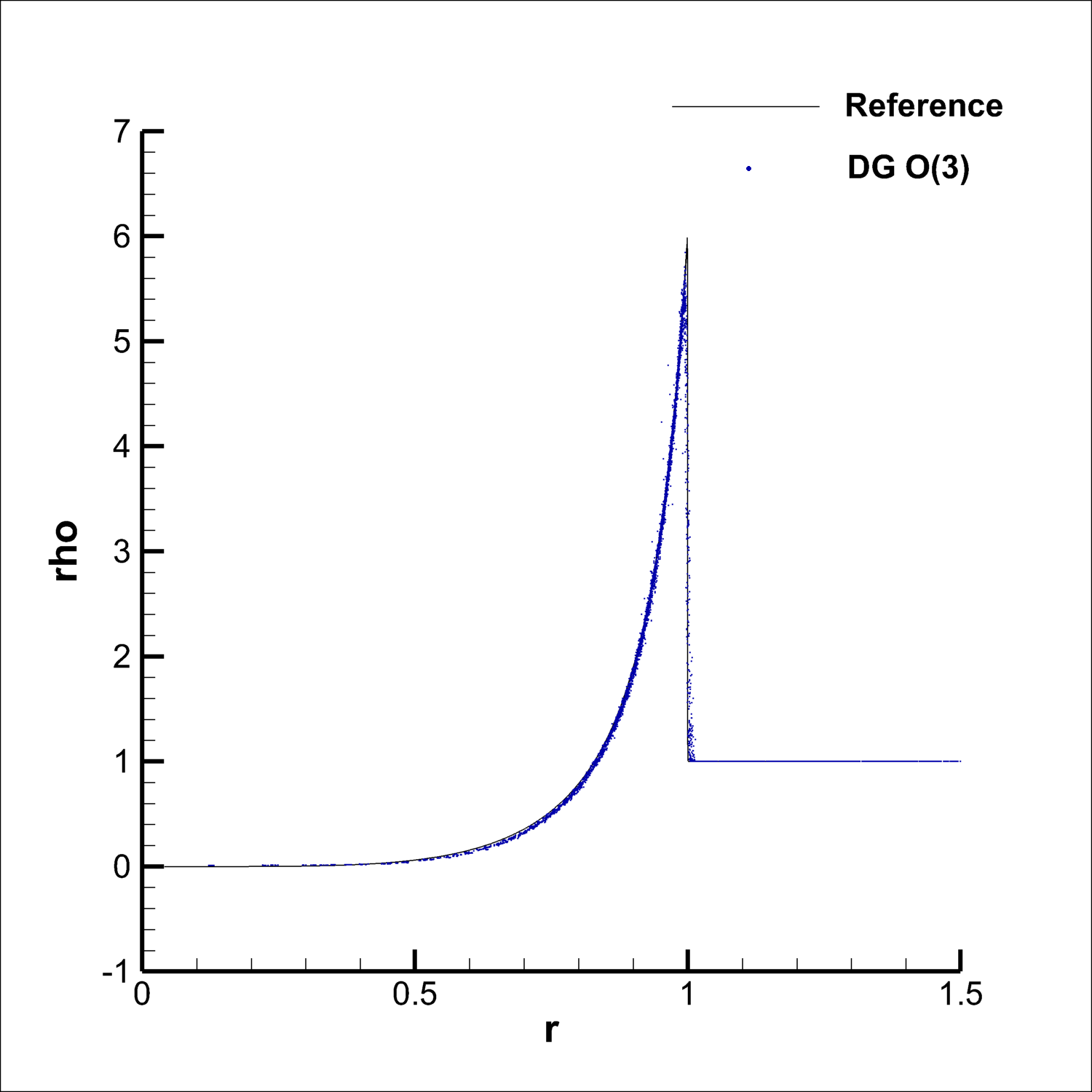}
		\caption{\RIIIcolor{Scatter plot of density values (blue dots) for the Sedov problem compared with the reference solution (line). 
			Left: the initial setting is with outer pressure equal to $p_0=10^{-6}$ and it is solved with our DG scheme of order $2$.
			Right: the initial setting is with outer pressure equal to $p_0=10^{-10}$ and it is solved with our DG scheme of order $3$.  
            In both cases, the test uses the Laplacian smoothing algorithm with $\mathcal{F} = 10^{-2}$.
			In the two cases the mesh moves with the fluid flow and topology changes occur.}}
	\label{fig:sedovp1p11e-6scatter}
\end{figure}

This test problem is widespread in the literature~\cite{LoubereSedov3D} and it describes the evolution of a blast wave that is generated at the origin $\mathbf{O}=(x,y)=(0,0)$ of the computational domain $\Omega(0)=[0;1.2]\times[0;1.2]$. An exact solution based on self-similarity arguments is available from~\cite{Sedov} and the fluid is assumed to be an ideal gas with $\gamma=1.4$, which is initially at rest and assigned with a uniform density $\rho_0=1$. The initial pressure is $p_0$ everywhere (with $p_0=10^{-6}$ or $p_0=10^{-10}$) except in the cell $V_{or}$ containing the origin $\mathbf{O}$ where it is given by
\begin{equation}
p_{or} = (\gamma-1)\rho_0 \frac{E_{tot}}{|V_{or}|}, \quad \textnormal{ with } \
E_{tot} =  0.979264,
\label{eqn.p0.sedov}
\end{equation}
being $E_{tot}$ the total energy concentrated at $\x=\mathbf{0}$.
We solve this numerical test with a second order $P_1P_1$ DG scheme for $p_0=10^{-6}$ and a third order DG scheme for $p_0=10^{-10}$; we employ a mesh of $7234$ Voronoi elements.
The density profiles are shown in Figure~\ref{fig:sedovp1p180x800-0} for various output times $t=0, 0.2, 0.5, 0.8, 1.0$. The obtained results are in good agreement with the literature. 
Moreover, one can refer to Figure~\ref{fig:sedovp1p11e-6scatter} for a comparison between our numerical solution (\RIIIcolor{scatter plot}) and the exact solution: the position of the shock wave and the density peak are perfectly captured. We remark that this is quite a challenging benchmark because of the low pressure and the strong shock. 

Finally, we refer to the last panel of Figure~\ref{fig:sedovp1p180x800-0} for the behavior of our \textit{a posteriori} sub--cell finite volume limiter, which activates only and exactly where the shock wave is located.

\RIIIcolor{
\subsubsection{Single material triple point problem}
\label{test.triplepoint}

A typical test in the ALE community (especially for comparing ReALE approaches) is given by the so called \textit{triple point} problem, first introduced in~\cite{ReALE2010} and widespread in the literature, see for example~\cite{ReALE2011,ShashkovMultiMat3,Kucharik2014,boscheri2015direct,ReALE2015}. 
Here in particular we follow the setting presented in~\cite{ReALE2015}, which has been slightly modified due to the fact that in this work we do not allow for moving boundaries and thus we have enlarged our computational domain, which is taken to be  $ [-1,8]\times[-1, 4]$; moreover, generators located outside $ [-0.5,7.5]\times[-0.5, 3.5]$ are kept fixed and slip wall boundary conditions are imposed everywhere.
 
To impose the initial condition we divide our domain in three regions, namely 
$R_1 = [-1, 1]\times[-1, 4]$, $R_2 = [1, 8]\times[-1, 1.5]$  and $R_3 = [1, 8]\times[1.5, 4]$, and we set
\be
(\rho, u, v, p)(\x) = 
\begin{cases} 
	(1,   0, 0, 1       ),  \ \text{ if } \  \x \in R_1 \\
	(1,   0, 0, 0.125   ),  \ \text{ if } \  \x \in R_2 \\
	(0.1, 0, 0, 0.125   ),  \ \text{ if } \  \x \in R_3, \\	
\end{cases}
\ee
with $\gamma = 1.4$. 

Due to the difference of density and pressure, two shocks propagate in the top and bottom domains with different speeds. This creates a shear flow along the initial horizontal contact discontinuity. 
On moving meshes capturing accurately the vorticity is the difficult part of such a simulation, involving the interaction of shocks, rarefactions, shear and contact waves. 

Since an exact solution does not exist, in order to visually see the convergence of our method, we present the results 
obtained with a third order DG scheme, namely a $P_2P_2$ scheme, on three different meshes ($M_1$, $M_2$, $M_3$) at time $t=3$ 
(see Figure~\ref{fig.triplepoint_3}) and $t=5$ (see Figure~\ref{fig.triplepoint_5}) 
(which are the same times at which results are given in Reference~\cite{ReALE2015}).

In particular, our coarsest mesh $M_1$ is made by a total number of $6140$ Voronoi elements ($2798$ of which in the usual reference domain, i.e. $[0,7]\times[0,3]$), the second mesh $M_2$ is made by a total number of $12706$ elements ($5823$ of which in the usual reference domain) and the finest mesh $M_3$ is made by $25158$ ($11597$ of which in the usual reference domain).

This test case demonstrates the robustness of our approach in dealing with complex mesh motion, its ability in following 
accurately the fluid flow, and its resolution due to the underlying high order DG method. 

\begin{figure}[!bp]
	\centering
	\includegraphics[width=1.0\linewidth]{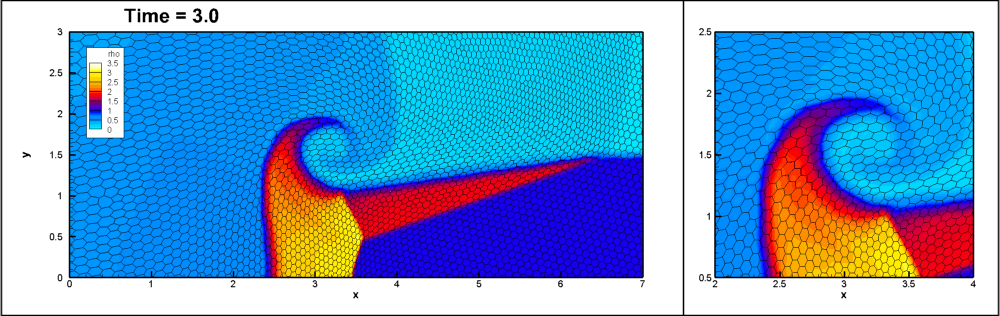}\\[-1pt] 
	\includegraphics[width=1.0\linewidth]{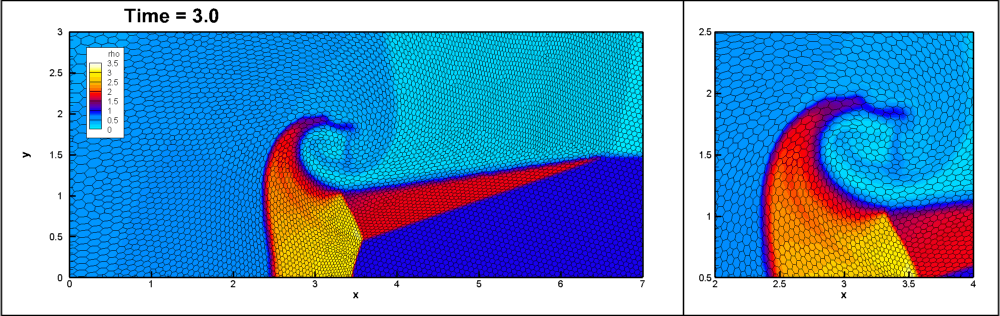}\\[-1pt] 
	\includegraphics[width=1.0\linewidth]{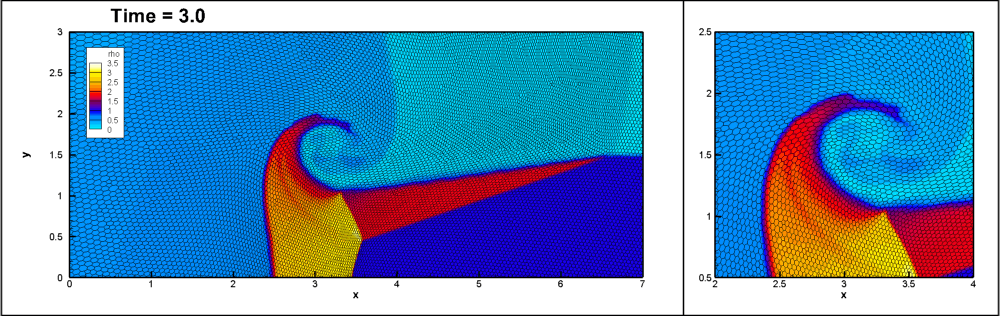}\\
	\caption{Triple point solved with our $P_2P_2$ third order DG scheme on three meshes. 
    The test uses the Lloyd-like smoothing algorithm with $\mathcal{F} = 10^{-1}$. Here we plot the density contour levels at time $t=3$ in the usual domain $[0,7]\times[0,3]$ (left) and a zoom in the region $[2,4]\times[0.5,2.5]$ (right) on three finer and finer meshes $M_1$ (top), $M_2$ (middle) and $M_3$ (bottom); refer to Section \ref{test.triplepoint} for the number of element of each mesh. }
	\label{fig.triplepoint_3}
\end{figure}

\begin{figure}[!bp]
	\centering
	\includegraphics[width=1.0\linewidth]{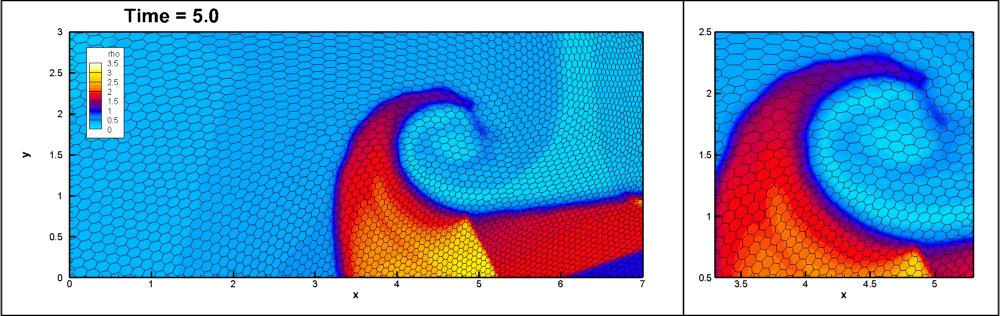}\\[-1pt] 
    \includegraphics[width=1.0\linewidth]{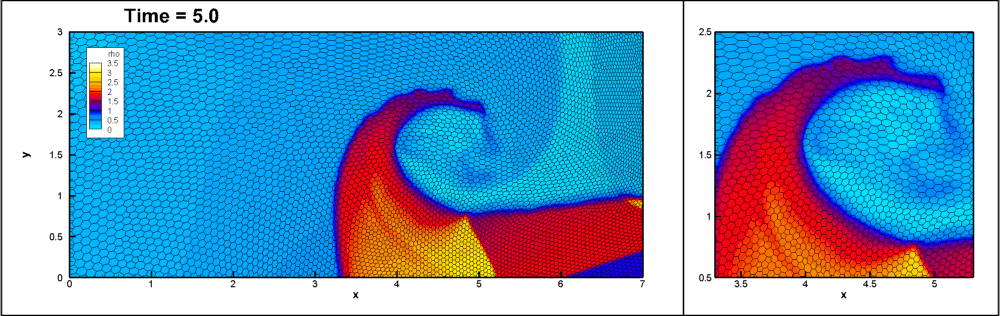}\\[-1pt] 
    \includegraphics[width=1.0\linewidth]{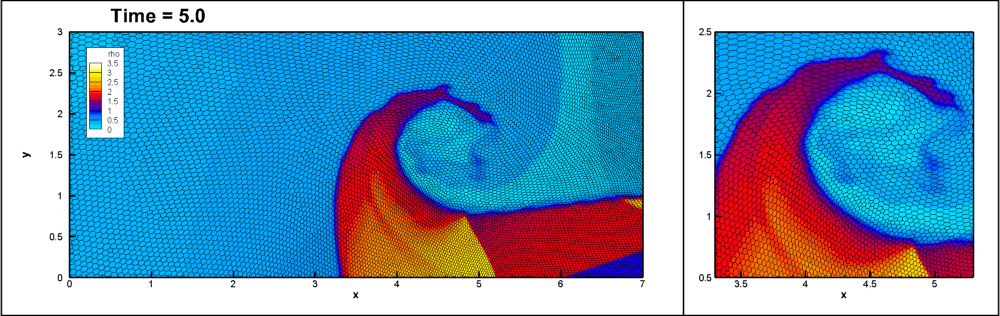}\\
		\caption{Triple point solved with our $P_2P_2$ third order DG scheme on three meshes. 
        The test uses the Lloyd-like smoothing algorithm with $\mathcal{F} = 10^{-1}$. Here we plot the density contour levels at time $t=5$ in the domain $\Omega_r=[0,7]\times[0,3]$ (left) and a zoom in the region $[3.3,5.3]\times[0.5,2.5]$ (right) on three finer and finer meshes $M_1$ (top), $M_2$ (middle) and $M_3$ (bottom); refer to Section \ref{test.triplepoint} for the number of element of each mesh. }
	\label{fig.triplepoint_5}
\end{figure}
}

\subsection{Euler equations with source term}
\label{ssec.EulerEq+S}

Next, we consider the Euler equations given in~\eqref{eulerTerms}, but with a gravity source term of the form 
\be
\S~=~\left( 0, \, 0, \,  -g\rho,  \, - g \rho v \right)^T.
\ee 
The Euler equations with gravity are of interest not only in  hydrodynamics~\cite{BottaKlein,Kapelli2014,Klingenberg2015,desveaux2016well,KlingenbergPuppo}, but also computational astrophysics~\cite{Springel, springel2010moving, mignone2011pluto}. 

\paragraph{Rayleigh-Taylor instability} 

With this test case we study an important type of fluid instability that arises in stratified atmospheres in approximate hydrostatic equilibrium if a denser fluid lies above a lighter phase.  In such a Rayleigh-Taylor unstable state, energy can be gained if the lighter fluid rises in the gravitational field, triggering buoyancy-driven fluid motions.
We consider here a simple test where we excite only one single Rayleigh-Taylor mode.  

Our  setup  is  a  small  variation (due to the fact that we do not allow for moving boundaries in this work)  of  a  similar  test considered in~\cite{liska2003comparison} and in~\cite{Springel}. 
The computational domain is $ [-0.25, 0.75]\times [-0.25, 1.75]$, with generators kept fixed outside the domain $ [-0.1, 0.6]\times [-0.1, 1.6]$ and wall boundary conditions set everywhere. 
The imposed initial condition is given by  
\be
\begin{cases} 
	\rho_B = 2,  \quad p_B = P_0 + g(y-0.75)\rho_B, \quad  &y \leq 0.75\\
	\rho_T = 1, \quad  p_T = P_0 + g(y-0.75)\rho_T, \quad & y > 0.75 ,
\end{cases}
\ee
with $P_0=2.5$ and $g=-0.1$.
The initial velocities are zero everywhere, i.e. $\mathbf{u} = (u,v) = \0$, except for a small perturbation that is designed to excite one single mode for the Rayleigh-Taylor instability
\be
v(x,y) = \omega_0 \left(1-\cos(4\pi x)\right) \left(1-\cos(4\pi y/3) \right)
\quad \text{if } \, 0 \le x \le 0.5, 
\ee 
where $\omega_0 = 0.0025$.
Next, we smooth the initial discontinuity in the density (in such a way that the limiter for the DG scheme will not be necessary) according to ~\cite{tavelli2014high}
\be
\rho(\x) = \frac{1}{2}\left( \rho_B+\rho_T\right) + \frac{1}{2} \left (\rho_T - \rho_B \right)\text{erf}\left( \frac{y-0.75}{\epsilon}\right).
\ee 

\begin{figure}[!bp]
	\centering
	\subfloat[FV $\mathcal{O}(3)$, $M_2$]{\includegraphics[width=0.25\linewidth]{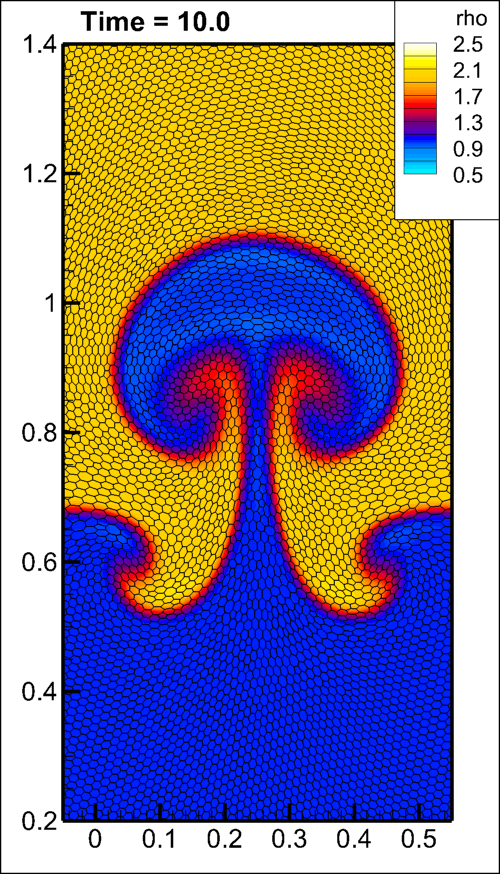}}%
	\subfloat[DG $\mathcal{O}(2)$, $M_1$]{\includegraphics[width=0.25\linewidth]{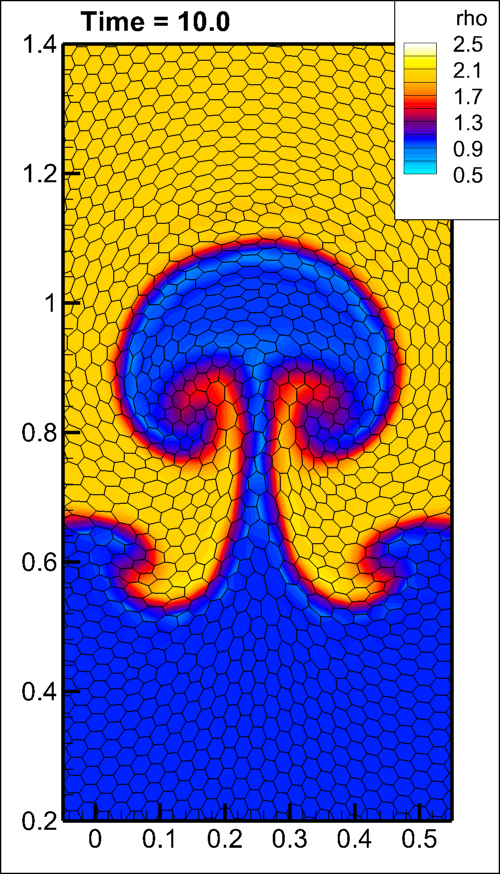}}%
	\subfloat[DG $\mathcal{O}(3)$, $M_1$]{\includegraphics[width=0.25\linewidth]{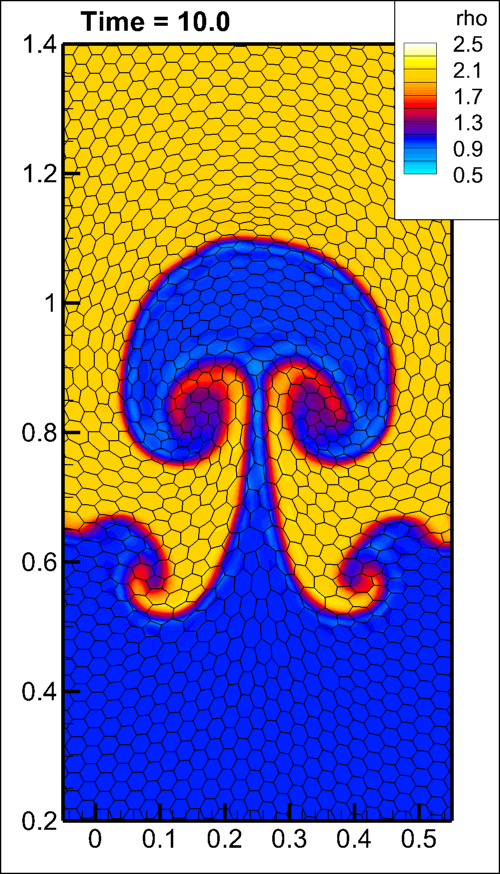}}%
	\subfloat[DG $\mathcal{O}(2)$, $M_2$]{\includegraphics[width=0.25\linewidth]{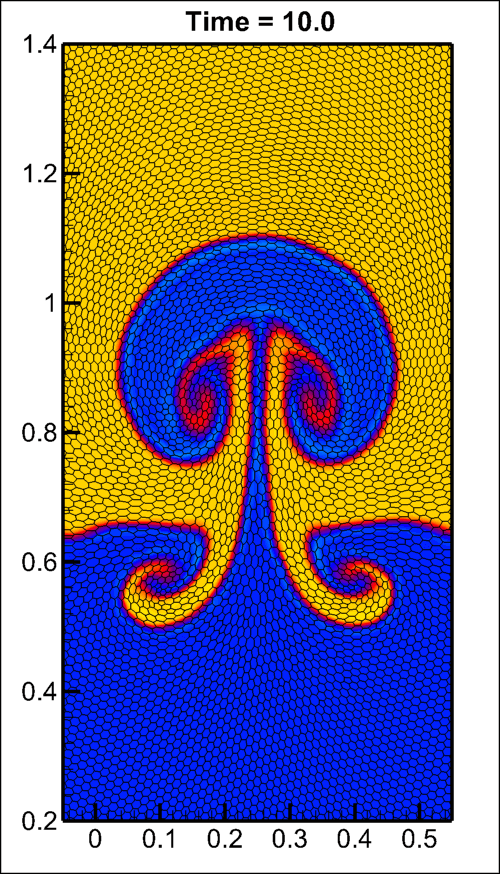}}\\
	\subfloat[FV $\mathcal{O}(3)$, $M_2$]{\includegraphics[width=0.25\linewidth]{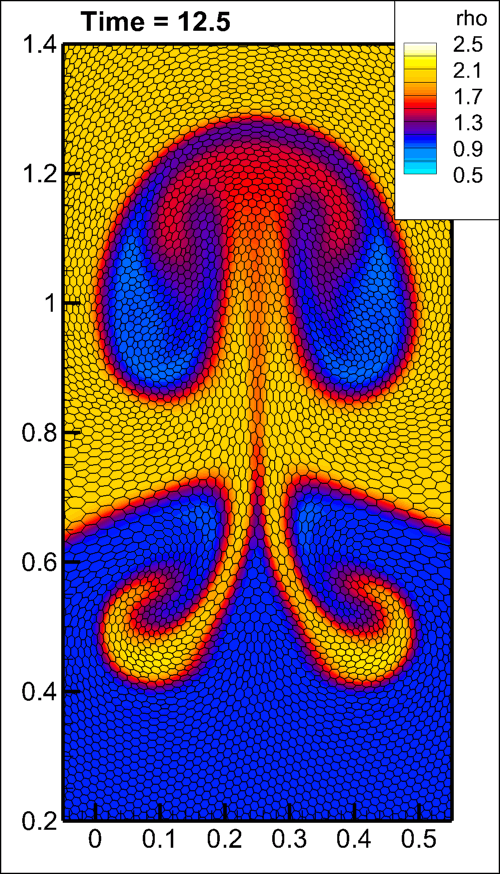}}%
	\subfloat[DG $\mathcal{O}(2)$, $M_1$]{\includegraphics[width=0.25\linewidth]{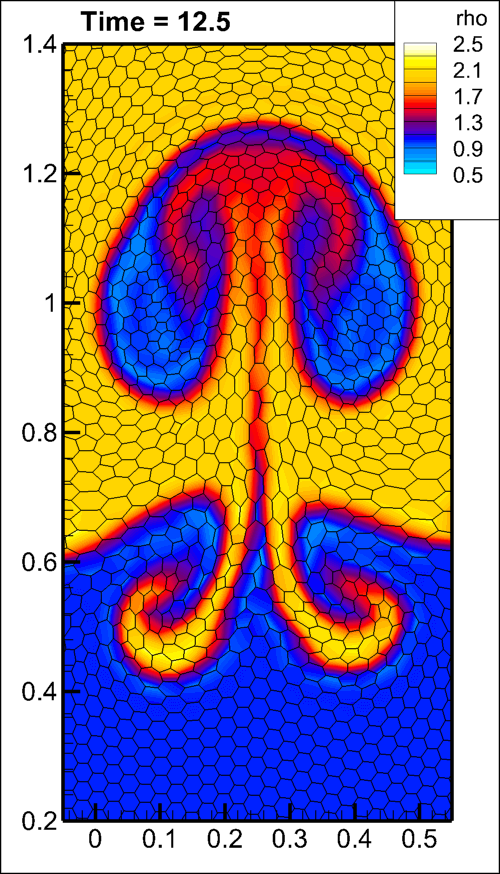}}%
	\subfloat[DG $\mathcal{O}(3)$, $M_1$]{\includegraphics[width=0.25\linewidth]{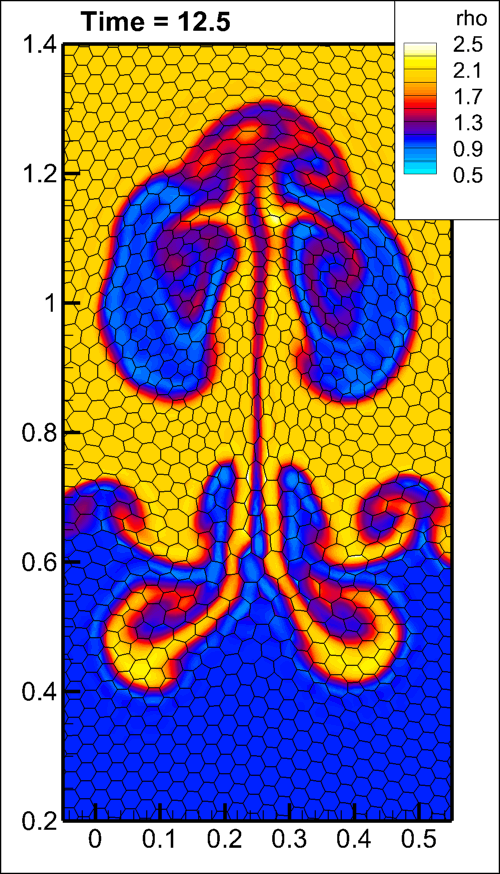}}%
	\subfloat[DG $\mathcal{O}(2)$, $M_2$]{\includegraphics[width=0.25\linewidth]{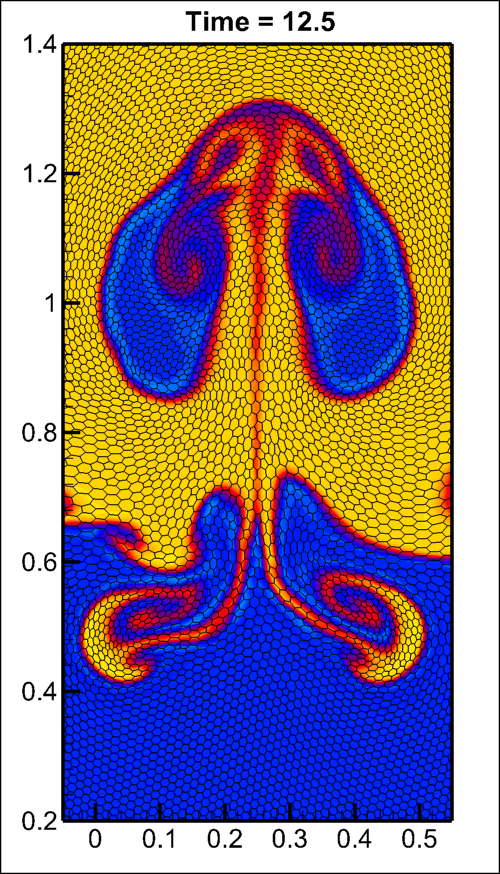}}%
	\caption{Rayleigh-Taylor instabilities at time $t=10$ (top) and time $t=12.5$ (bottom). The results in the panel are obtained by using two coarse meshes: 
		$M_1$ made of $3175$ elements and $M_2$ which is made of $13883$ elements and is $2$ times finer than $M_1$. The test uses the Lloyd-like smoothing algorithm with $\mathcal{F} = 10^{-3}$.
		We have employed our FV scheme of order $3$ (a,e) and our DG scheme of order $2$ (b,d,f,h) and $3$ (c,g), 
        \RIIcolor{and the total number of DOFs in the DG case corresponds to $9525$ (b,f), $19050$ (c,g), and $41649$ (d,h)}.
		We would like to \RIIIcolor{emphasize} that the use of a high order DG scheme makes secondary structures appear even on the 
		coarse mesh $M_1$ (g) which cannot be seen with standard third order FV schemes \RIIcolor{with a similar number of DOFs (a,e)}.}
	\label{fig.RTcomparison}
\end{figure}

We solve this problem deliberately on coarse meshes ($M_1$ made of $2\,706$ elements and $M_2$ made of $13\,340$ cells) and we compare the  resolution of the instabilities obtained with  our ALE FV-DG schemes with different order of accuracy, see Figure~\ref{fig.RTcomparison}.
Specifically, we compare third order FV scheme with second and third order DG schemes, i.e. $P_0P_1$, $P_1P_1$, $P_0P_2$, $P_2P_2$ and we employ the Osher-type ALE flux as approximate Riemann solver~\eqref{eqn.osher}; we note that secondary instability vortexes only appear within a high order DG method, being hidden by numerical dissipation in the other cases. 

\RIVcolor{Finally in Table~\ref{tab.RT_sliver_percentage},
 we report some statistics on the number of sliver elements created over the total number of time steps, 
 and on the percentage of computational time required both for the geometrical part of the code and for the $P_NP_M$ predictor-corrector algorithm.}

Comparing our results with those presented in~\cite{Springel}, we can remark the importance of coupling our 
new high order DG and FV  algorithms with the direct ALE framework with topology changes; indeed our results 
show a highly increased resolution with respect to~\cite{Springel} and the ability of capturing secondary instability structures on 
coarse meshes.  
 
 \RIVcolor{
\begin{table*}[!tp]
    \caption{\RIVcolor{Rayleigh-Taylor instabilities. In this Table we report the number of total time steps needed to reach the final time $t_f=12.5$ with our $P_0P_2$, $P_1P_1$ and $P_2P_2$ schemes, the number of sliver elements and the number of time steps that have been repeated through the MOOD loop described in Section~\ref{ssec.MOOD}. Moreover, we report the percentage of computational time employed by mesh regeneration and space time connectivity generation, by the predictor-corrector step on standard elements and on sliver elements. The other part of the computational time for what concerns the FV scheme is mostly spent on the reconstruction procedure. Instead the DG scheme was run without activating the limiter.}} 
    \label{tab.RT_sliver_percentage}  
    \centering
        \begin{tabular}{c|ccc|ccc} 
            \hline
            Method & time steps & slivers & restarts & Mesh $\%$  & $P_NP_M$ standard $\%$  & $P_NP_M$ sliver $\%$   \\
            \hline          
            FV $\mathcal{O}(3)$, $M_2$ & 15219   & 15956  & 0     & 1.00     & 74.19      & 9.7E-4             \\
            DG $\mathcal{O}(2)$, $M_1$ & 27951  & 2707   & 0   & 4.23     & 87.65      & 6.5E-4                \\
            DG $\mathcal{O}(3)$, $M_1$ & 54493  & 1297   & 0   & 0.65     & 95.64      & 9.0E-5                \\
            DG $\mathcal{O}(2)$, $M_2$ & 41919  & 17114   & 0   & 5.3     & 86.43      & 6.1E-4                \\
            \hline  
        \end{tabular}       
\end{table*}
}

\subsection{Ideal MHD equations}
\label{ssec.MHDEq}

We also consider the equations of ideal classical magnetohydrodynamics (MHD) that result in a more complicated system of hyperbolic conservation laws. The state vector $\Q$ and the flux tensor $\F$ for the MHD equations in the general form~\eqref{eq.generalform} are 
\begin{equation}
\label{MHDTerms}
\Q=\left( \begin{array}{c} \rho \\ \rho \mathbf{v} \\ \rho E \\ \B \\ \psi \end{array} \right), \quad \F(\Q) = 
\left( \begin{array}{c}  \rho \v  \\ 
\rho \v \otimes \v + p_t \mathbf{I} - \frac{1}{4 \pi} \B \otimes \B \\ 
\v (\rho E + p_t ) - \frac{1}{4 \pi} \B ( \v \cdot \B ) \\ 
\v \otimes \B - \B \otimes \v + \psi \mathbf{I} \\
c_h^2 \B  \end{array} \right).
\end{equation}
Here, $\B=(B_x,B_y,B_z)$ represents the magnetic field and $p_t=p+\frac{1}{8\pi}\mathbf{B}^2$ is the total pressure. The hydrodynamic pressure is given by the equation of state used to close the system, thus
\begin{equation}
p = \left(\gamma - 1 \right) \left(\rho E - \frac{1}{2}\mathbf{v}^2 - \frac{\mathbf{B}^2}{8\pi}\right).
\label{MHDeos}
\end{equation}
System~\eqref{MHDTerms} requires an additional constraint on the divergence of the magnetic field to be satisfied, that is
\begin{equation}
\nabla \cdot \mathbf{B} = 0.
\label{eq:divB}
\end{equation}  
Here,~\eqref{MHDTerms} includes one additional scalar PDE for the evolution of the variable $\psi$, which is needed to transport divergence errors outside the computational domain with an artificial divergence cleaning speed $c_h$, see~\cite{Dedneretal}. A more recent and more sophisticated methodology to fulfill this condition exactly on the discrete level also in the context of high order ADER WENO finite volume schemes on unstructured simplex meshes can be found in~\cite{BalsaraDivB2015}. A similar approach is adopted in~\cite{fambri2017space, boscheri2014high, boscheri2017efficient}.

\subsubsection{MHD vortex}
\label{test.MHDVortex}

\begin{table*}[!tp]
	\caption{MHD vortex. Numerical convergence results for the finite volume algorithm on moving meshes with topology changes. The error norms refer to the variable $\rho$ at time $t=1.0$ in $L_1$ norm.} 
    \centering
		\begin{tabular}{ccc|ccc|ccc|ccc} 
			\hline 
			\multicolumn{3}{c}{$P_0P_1 \rightarrow \mathcal{O}2$} & \multicolumn{3}{c}{$P_0P_2\rightarrow \mathcal{O}3$}	& \multicolumn{3}{c}{$P_0P_3\rightarrow \mathcal{O}4$} & \multicolumn{3}{c}{$P_0P_4\rightarrow \mathcal{O}5$}    \\ 
			\hline
		\!\!$h(\Omega(t_f))$ \!\!&\!\! $\epsilon(\rho)_{L_1}$ \!\!&\!\!\!\! $\mathcal{O}(L_1)$ \!\!\!&  $h(\Omega(t_f))$ \!\!&\!\! $\epsilon(\rho)_{L_1}$ \!\!&\!\!\!\! $\mathcal{O}(L_1)$ \!\!\!& $h(\Omega(t_f))$ \!\!&\!\! $\epsilon(\rho)_{L_1}$ \!\!&\!\! $\mathcal{O}(L_1)$  \!\!&\!\! $h(\Omega(t_f))$ \!\!&\!\! $\epsilon(\rho)_{L_1}$ \!\!&\!\!\!\! $\mathcal{O}(L_1)$ \!\!\!\\ 
			\hline
			4.6E-01 \!\!&\!\! 3.3E-02  \!\!&\!\! -   \!& 3.2E-01  \!\!&\!\! 1.0E-02 \!\!&\!\! -    \!& 4.7E-01 \!\!&\!\!   2.1E-02 \!\!&\!\!   -  \!& 6.0E-01 \!\!&\!\! 3.6E-0.2 \!\!& -   \\
			3.9E-01 \!\!&\!\! 1.6E-02  \!\!&\!\! 1.8 \!& 2.4E-01  \!\!&\!\! 5.5E-03 \!\!&\!\! 2.3  \!& 3.2E-01 \!\!&\!\!   6.0E-03 \!\!&\!\!  3.2 \!& 5.8E-01 \!\!&\!\! 3.0E-0.2 \!\!& 5.8 \\
			2.4E-01 \!\!&\!\! 8.9E-03  \!\!&\!\! 2.3 \!& 1.9E-01  \!\!&\!\! 2.7E-03 \!\!&\!\! 3.3  \!& 2.4E-01 \!\!&\!\!   2.0E-03 \!\!&\!\!  3.9 \!& 5.6E-01 \!\!&\!\! 2.7E-0.2 \!\!& 3.6 \\
			1.9E-01 \!\!&\!\! 5.3E-03  \!\!&\!\! 2.4 \!& 1.6E-01  \!\!&\!\! 1.5E-03 \!\!&\!\! 3.1  \!& 2.2E-01 \!\!&\!\!   1.3E-03 \!\!&\!\!  3.6 \!& 5.5E-01 \!\!&\!\! 2.3E-0.2 \!\!& 5.9 \\
			1.6E-01 \!\!&\!\! 3.4E-03  \!\!&\!\! 2.5 \!& 1.4E-01  \!\!&\!\! 1.0E-03 \!\!&\!\! 2.9  \!& 1.9E-01 \!\!&\!\!   8.1E-04 \!\!&\!\!  4.8 \!& 5.2E-01 \!\!&\!\! 1.8E-0.2 \!\!& 4.8 \\
			\hline  
		\end{tabular}		
	\label{tab.orderOfconvergenceFV_MHD}
\end{table*}

\begin{table*}[!tp]
	\caption{MHD vortex. Numerical convergence results for the discontinuous Galerkin algorithm on moving meshes with topology changes. The error norms refer to the variable $\rho$ at time $t=1.0$ in $L_1$ norm.} 
	\centering	
		\begin{tabular}{ccc|ccc|ccc|ccc} 
			\hline 
			\multicolumn{3}{c}{$P_1P_1\rightarrow \mathcal{O}2$} & \multicolumn{3}{c}{$P_2P_2\rightarrow \mathcal{O}3$}	& \multicolumn{3}{c}{$P_3P_3\rightarrow \mathcal{O}4$}   	& \multicolumn{3}{c}{$P_4P_4\rightarrow \mathcal{O}5$}     \\ 
			\hline
			\!\!$h(\Omega(t_f))$ \!\!&\!\! $\epsilon(\rho)_{L_1}$ \!\!&\!\!\!\! $\mathcal{O}(L_1)$ \!\!\!&  $h(\Omega(t_f))$ \!\!&\!\! $\epsilon(\rho)_{L_1}$ \!\!&\!\!\!\! $\mathcal{O}(L_1)$ \!\!\!& $h(\Omega(t_f))$ \!\!&\!\! $\epsilon(\rho)_{L_1}$ \!\!&\!\! $\mathcal{O}(L_1)$  \!\!&\!\! $h(\Omega(t_f))$ \!\!&\!\! $\epsilon(\rho)_{L_1}$ \!\!&\!\!\!\! $\mathcal{O}(L_1)$ \!\!\!\\ 
			\hline
			4.7E-01\!\! & \!\!8.5E-03 \!\! & \!\!-   \!&  6.1E-01 \!\!&\!\! 2.8E-03 \!\!&\!\! -    \!& 8.8E-01 \!\!&\!\!   1.1E-03 \!\!&\!\!   -   \!& 1.6E-00 \!\!&\!\! 6.9E-0.3 \!\!& -   \\
			3.2E-01\!\! & \!\!3.2E-04 \!\! & \!\!2.5 \!&  4.7E-01 \!\!&\!\! 1.3E-03 \!\!&\!\! 2.8  \!& 7.5E-01 \!\!&\!\!   6.2E-04 \!\!&\!\!  3.5  \!& 6.1E-01 \!\!&\!\! 1.3E-0.4 \!\!& 4.1 \\
			2.8E-01\!\! & \!\!2.1E-04 \!\! & \!\!2.9 \!&  3.8E-01 \!\!&\!\! 7.3E-04 \!\!&\!\! 2.7  \!& 6.1E-01 \!\!&\!\!   3.1E-04 \!\!&\!\!  3.4  \!& 5.2E-01 \!\!&\!\! 4.7E-0.5 \!\!& 5.8 \\
			2.4E-01\!\! & \!\!1.6E-04 \!\! & \!\!2.0 \!&  3.5E-01 \!\!&\!\! 5.6E-04 \!\!&\!\! 3.6  \!& 5.5E-01 \!\!&\!\!   1.9E-04 \!\!&\!\!  4.3  \!& 4.9E-01 \!\!&\!\! 3.1E-0.5 \!\!& 8.1 \\
			1.9E-01\!\! & \!\!9.7E-05 \!\! & \!\!2.4 \!&  3.2E-01 \!\!&\!\! 4.1E-04 \!\!&\!\! 3.0  \!& 3.2E-01 \!\!&\!\!   2.3E-05 \!\!&\!\!  3.9  \!& 4.7E-01 \!\!&\!\! 2.4E-0.5 \!\!& 5.3 \\
			\hline 
		\end{tabular}		
	\label{tab.orderOfconvergenceDG_MHD}
\end{table*}

For the numerical convergence studies, we solve the vortex test problem proposed by Balsara in~\cite{Balsara2004}.  
The computational domain is given by the box $\Omega=[0;10]\times[0;10]$ with wall boundary conditions imposed everywhere. 
The initial condition is given in terms of the vector of primitive variables $\mathbf{V} = ( \rho, u, v, w, p, B_x, B_y,  B_z, \Psi )^T$ as 
\begin{equation}
\mathbf{V}(\mathbf{x},0) =
( 1, \delta u, \delta v, 0, 1+\delta p, \delta B_x, \delta B_y,  0, 0 )^T, 
\end{equation}
with $\delta \mathbf{v} = (\delta u, \delta v, 0)^T$, $\ \delta \mathbf{B} = ( \delta B_x, \delta B_y, 0 )^T$ and 
\begin{equation}
\begin{aligned} 
\label{eqn.mhd3d.ic1}
&\delta \mathbf{v} = \frac{\kappa}{2\pi} e^{ q(1-r^2)} \mathbf{e}_z \times \mathbf{r}  \\ 
&\delta \mathbf{B} = \frac{\mu}{2\pi}    e^{ q(1-r^2)} \mathbf{e}_z \times \mathbf{r},  \\ 
&\delta p = \frac{1}{64 q \pi^3} \left( \mu^2 (1 - 2 q r^2) - 4 \kappa^2 \pi \right) e^{2q(1-r^2)}. 
\end{aligned} 
\end{equation}
We have $\mathbf{e}_z = (0,0,1)$, $\mathbf{r} = (x-5,y-5,0)$ and $r = \left\| \mathbf{r} \right\| = \sqrt{ (x-5)^2 + (y-5)^2 }$. The divergence cleaning speed 
is chosen as $c_h=3$. The other parameters are $q=\halb$, $\kappa=1$ and $\mu=\sqrt{4 \pi}$, according to~\cite{Balsara2004}. 


\begin{figure*}[!bp]
    \centering
    \includegraphics[width=0.249\linewidth]{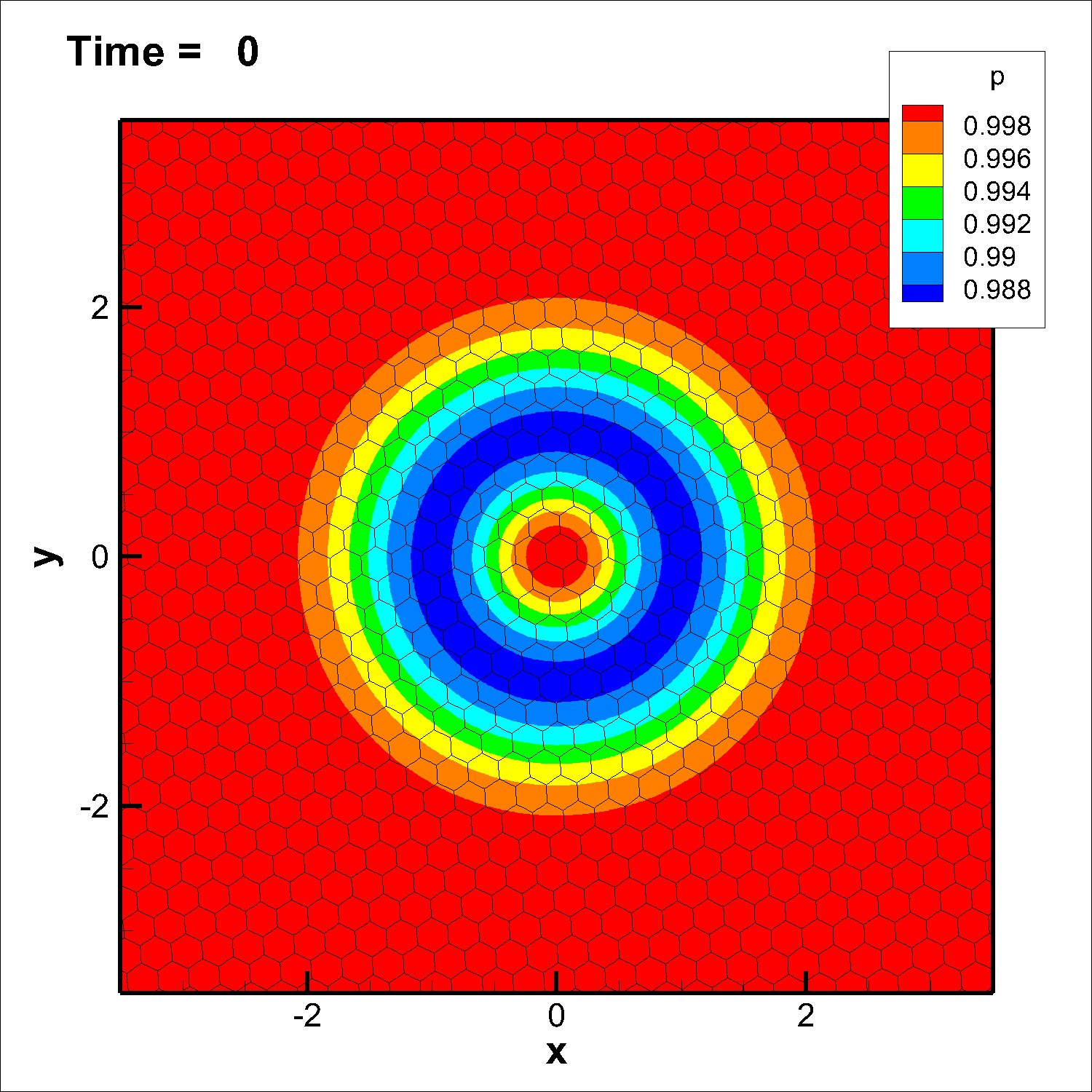}%
    \includegraphics[width=0.249\linewidth]{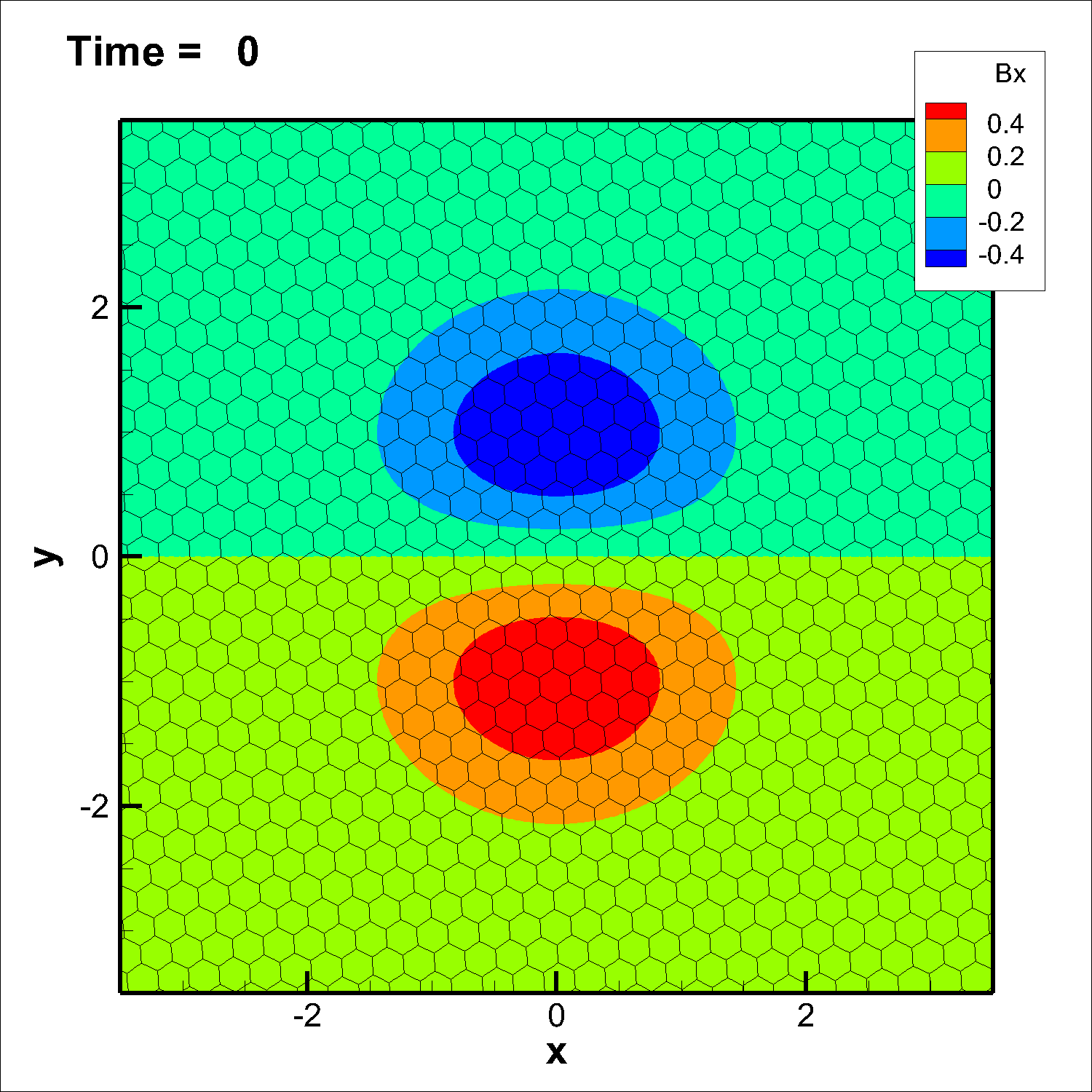}%
    \includegraphics[width=0.249\linewidth]{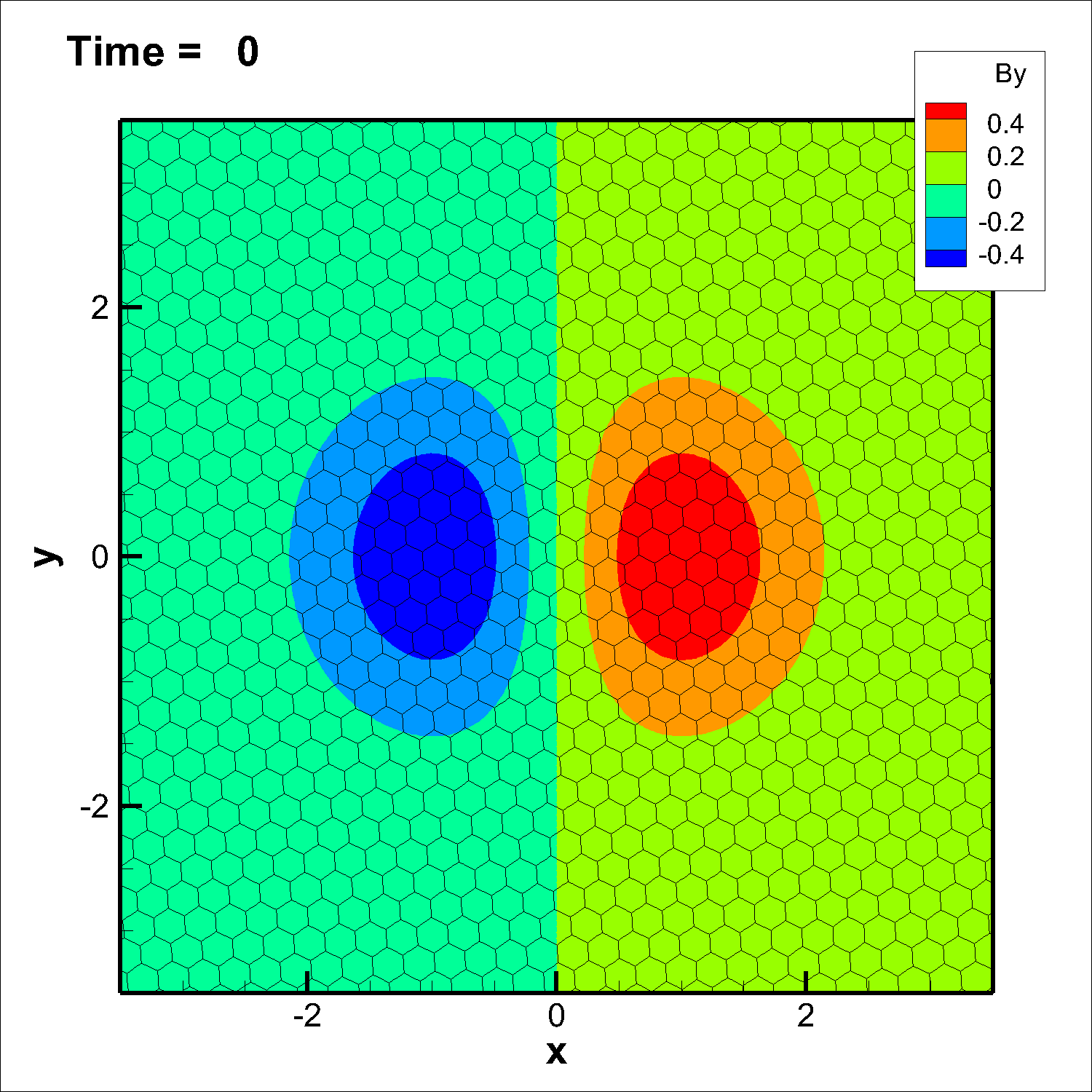}%
    \includegraphics[width=0.249\linewidth]{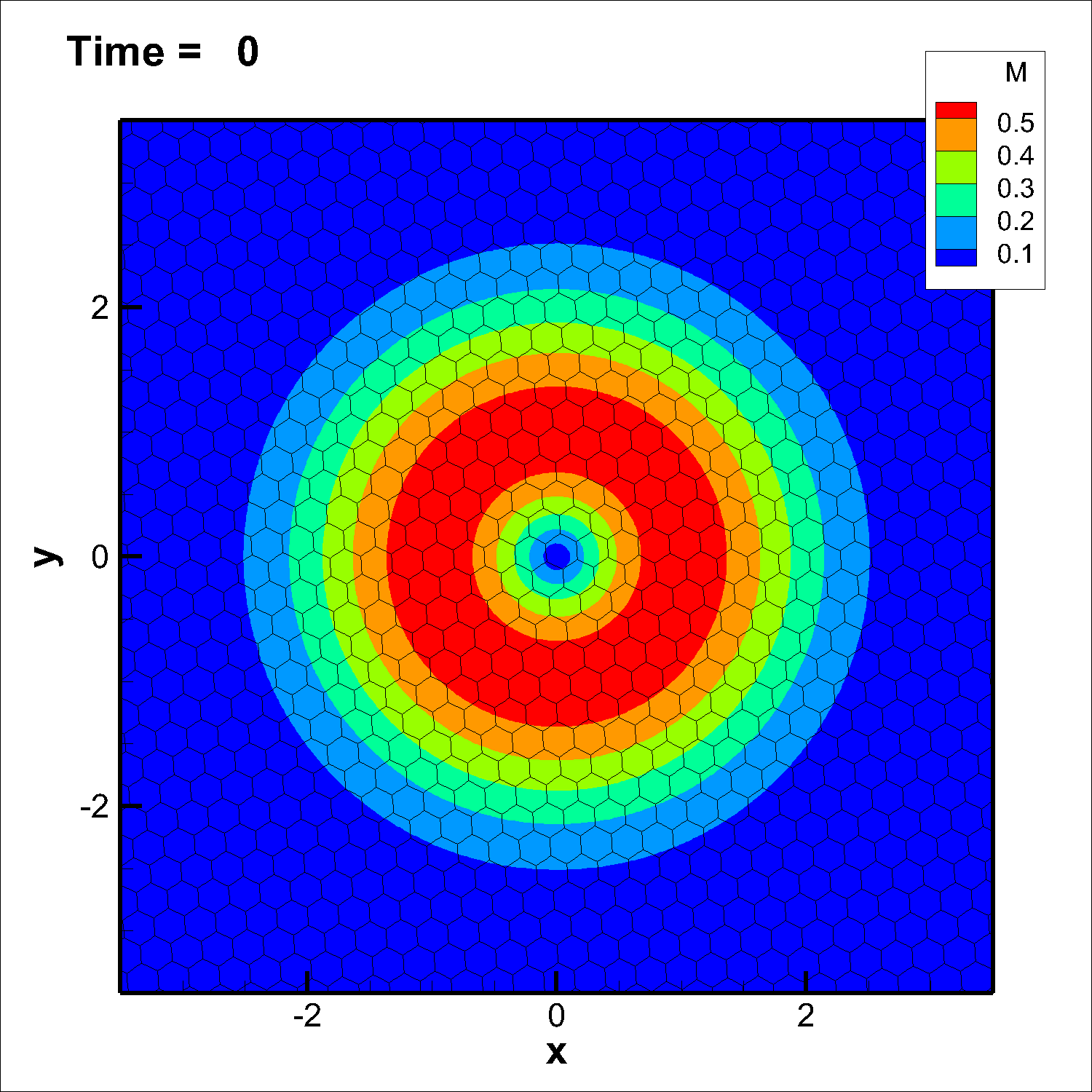}\\[-1pt]
    \includegraphics[width=0.249\linewidth]{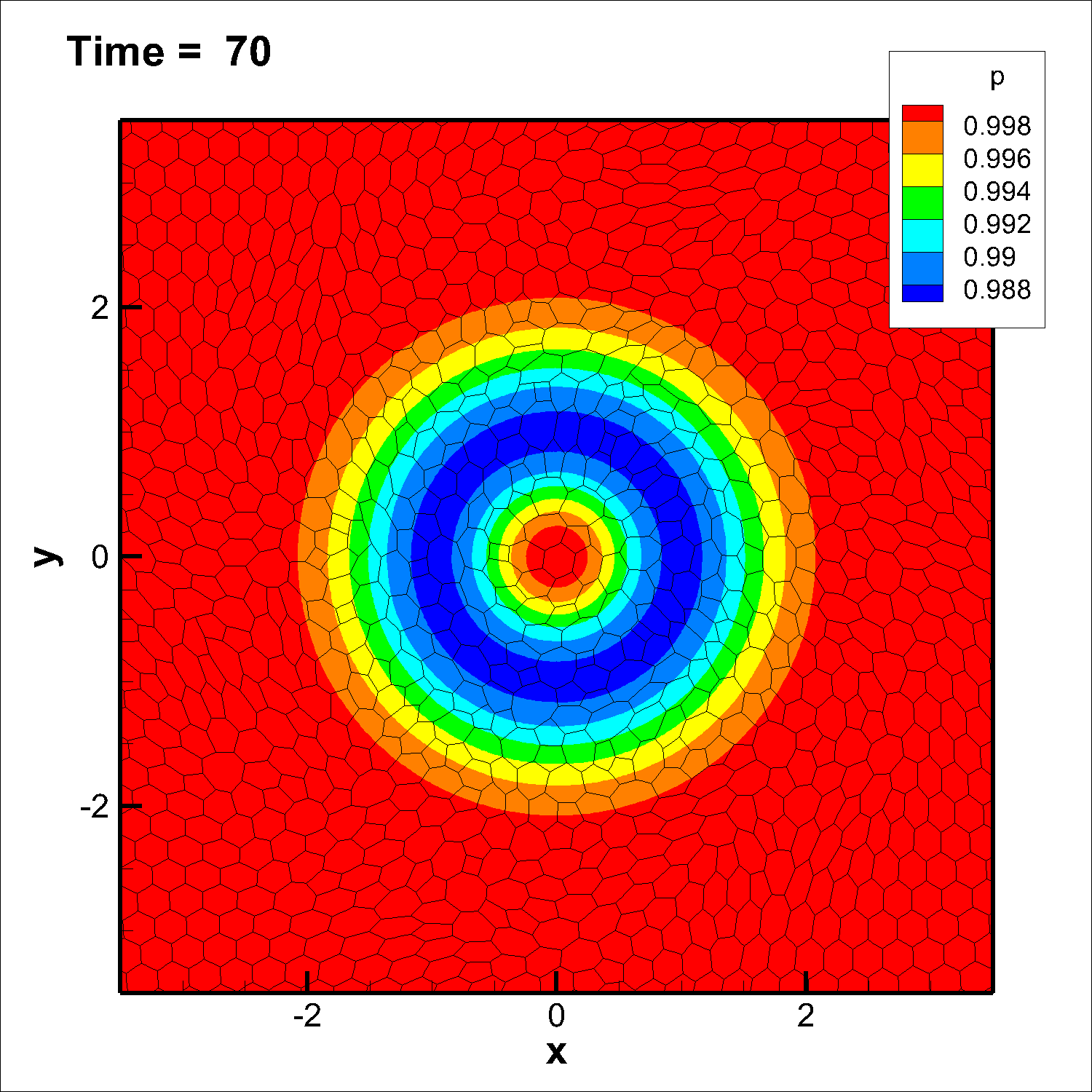}%
    \includegraphics[width=0.249\linewidth]{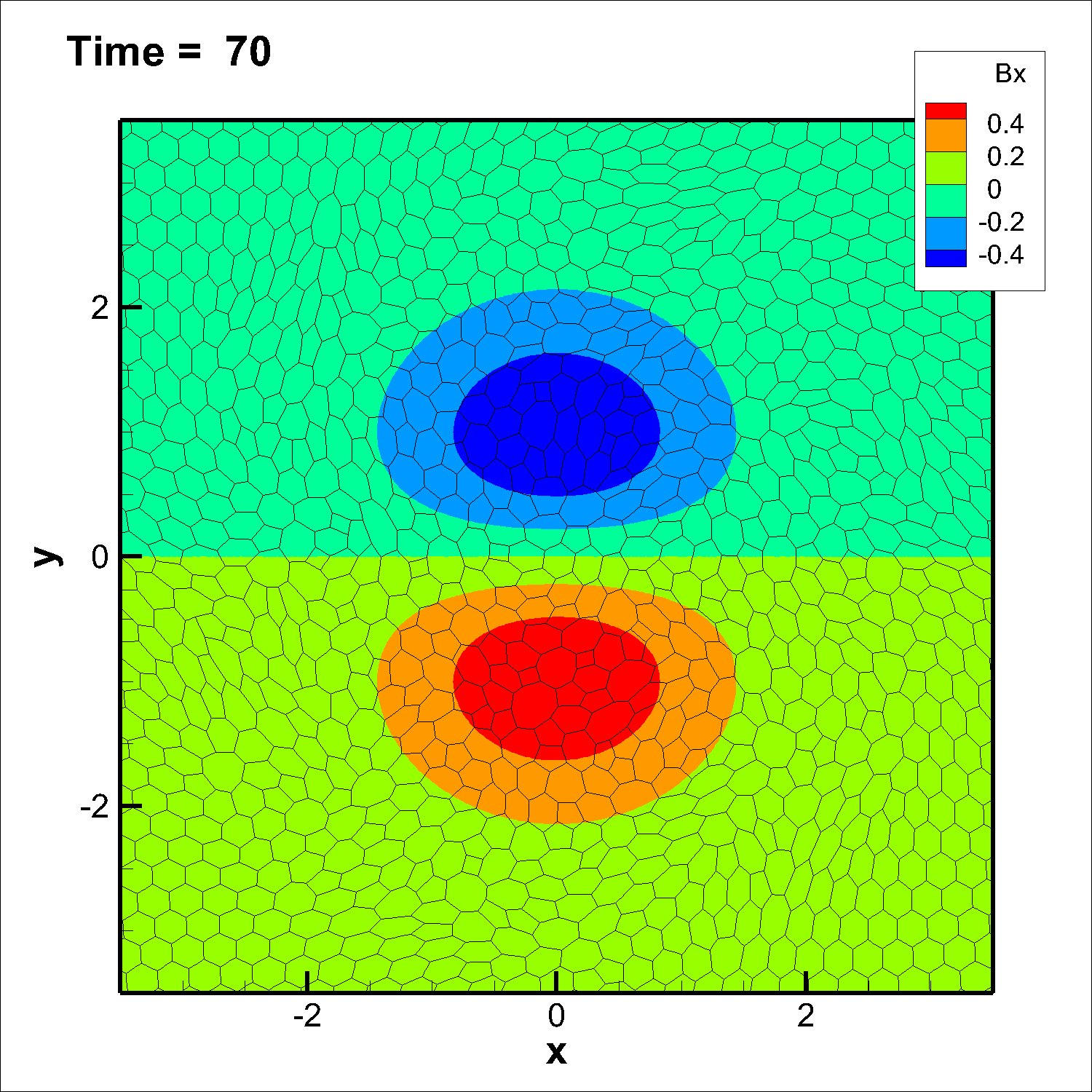}%
    \includegraphics[width=0.249\linewidth]{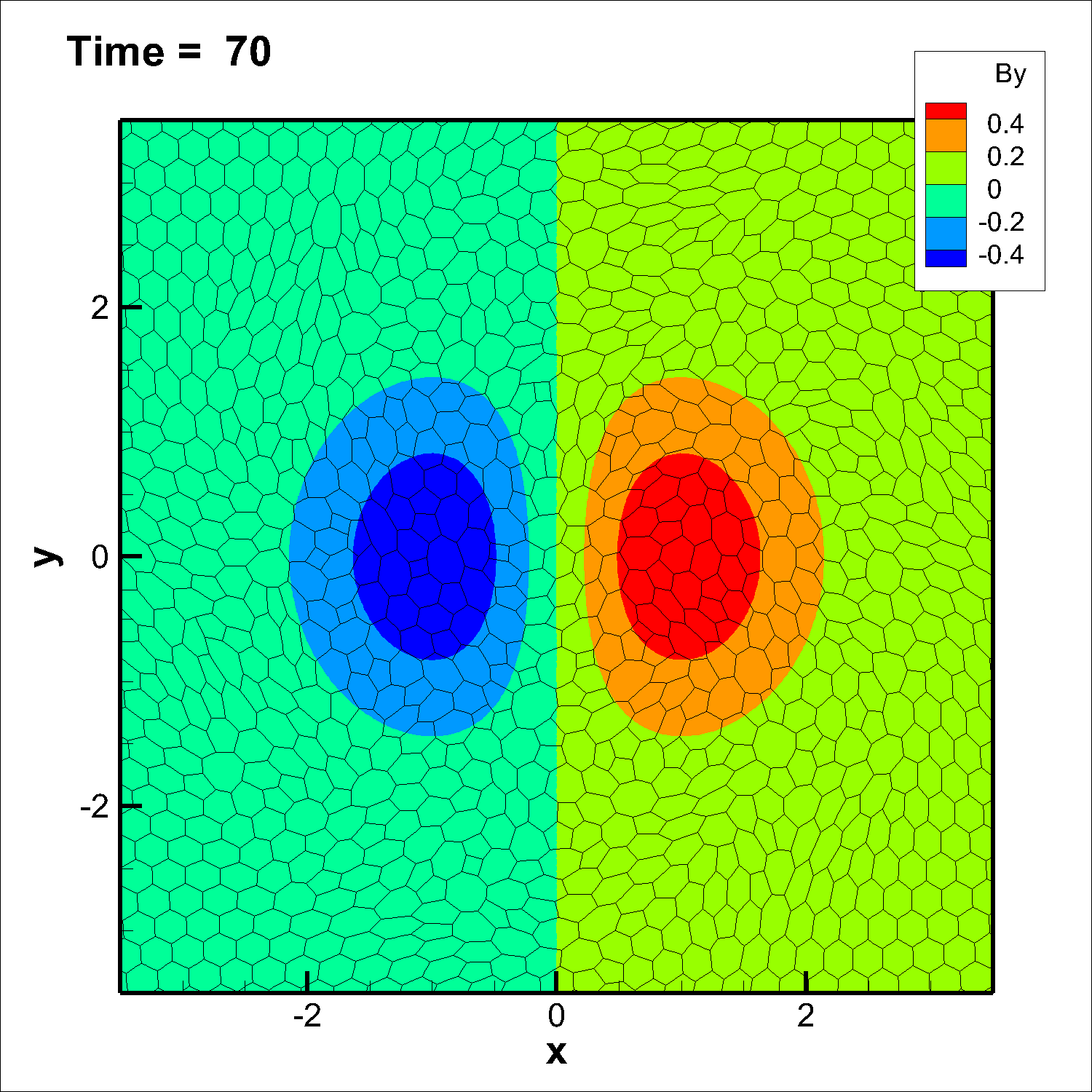}%
    \includegraphics[width=0.249\linewidth]{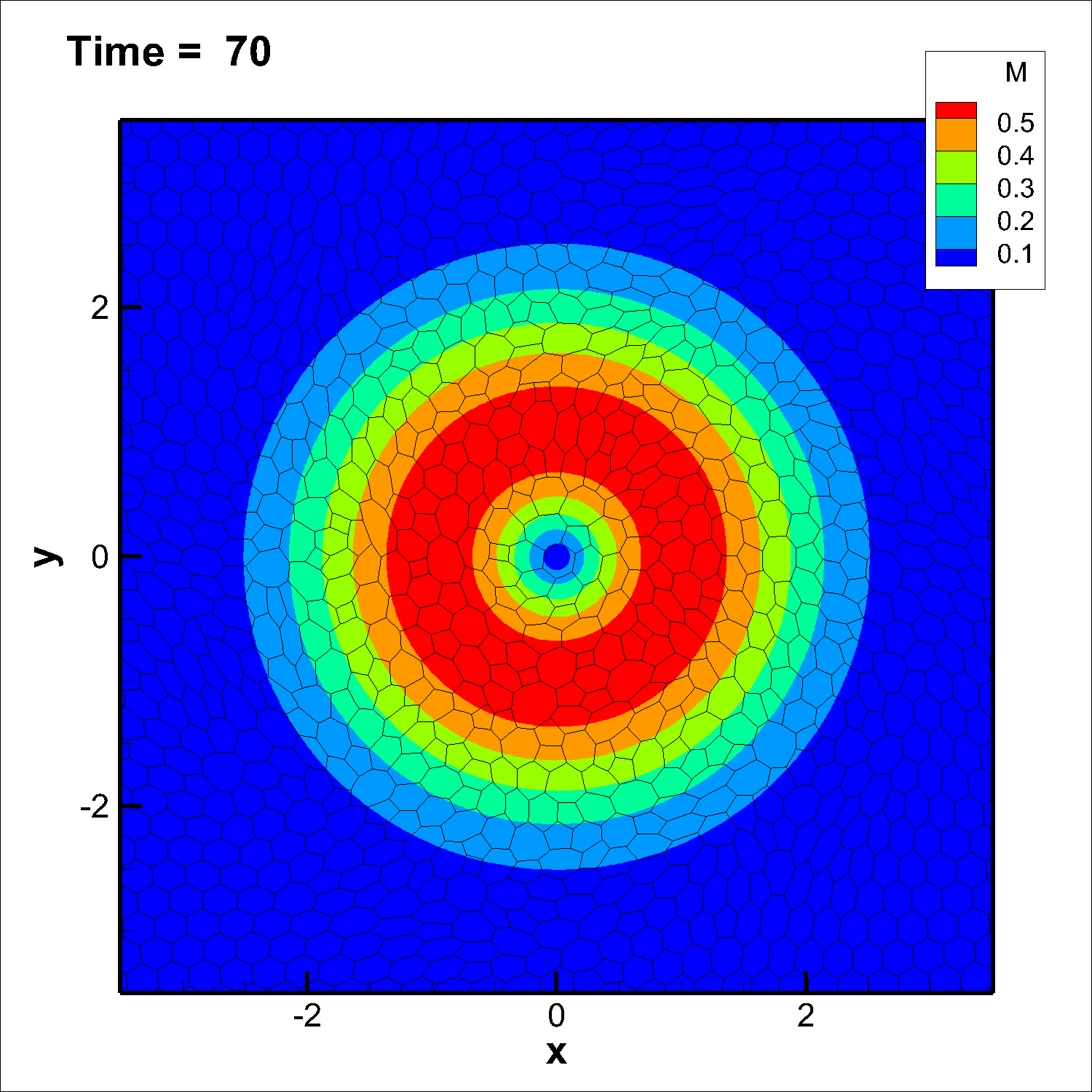}
    \caption{\RShashkovcolor{MHD vortex solved with our $P_3P_3$ DG scheme on a moving Voronoi mesh of $1900$ elements using the Lloyd-like smoothing algorithm with $\mathcal{F} = 10^{-5}$. We depict, from left to right, the pressure profile, the $x-$ and $y-$ components of the magnetic field, and the value of $M = \sqrt{B_x^2+B_y^2+B_z^2}$ at the initial time $t=0$ and at the final time of $t_f=70$, corresponding to 2 complete loops of the elements located at $r = 1$. 
            The connectivity changes (see also Figure~\ref{fig.MHDvortex_P3P3connectivity}) together with the high order methods allow to preserve the stationary MHD vortex for long times.}}
    \label{fig.MHDvortex_P3P3}
\end{figure*}

\paragraph{Convergence}
Tables~\ref{tab.orderOfconvergenceFV_MHD} and~\ref{tab.orderOfconvergenceDG_MHD} report the convergence rates from second up to fifth order of accuracy for the MHD vortex test problem run on a sequence of successively refined meshes up to the final time $t=1.0$.
The optimal order of accuracy is achieved both in space and time for the FV schemes as well as for the DG schemes.

\begin{figure}[!bp]
    \centering
    \includegraphics[width=0.495\linewidth]{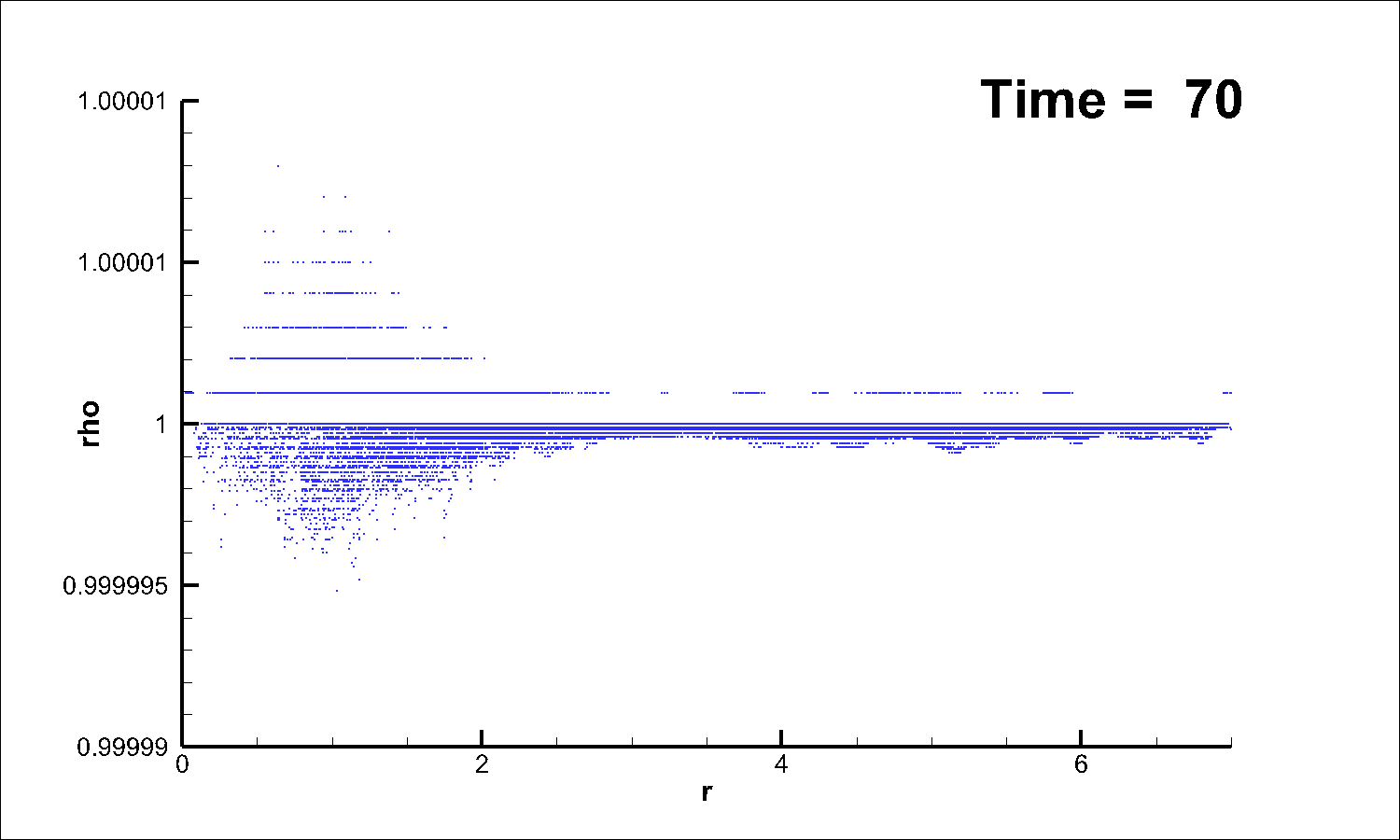}%
    \includegraphics[width=0.495\linewidth]{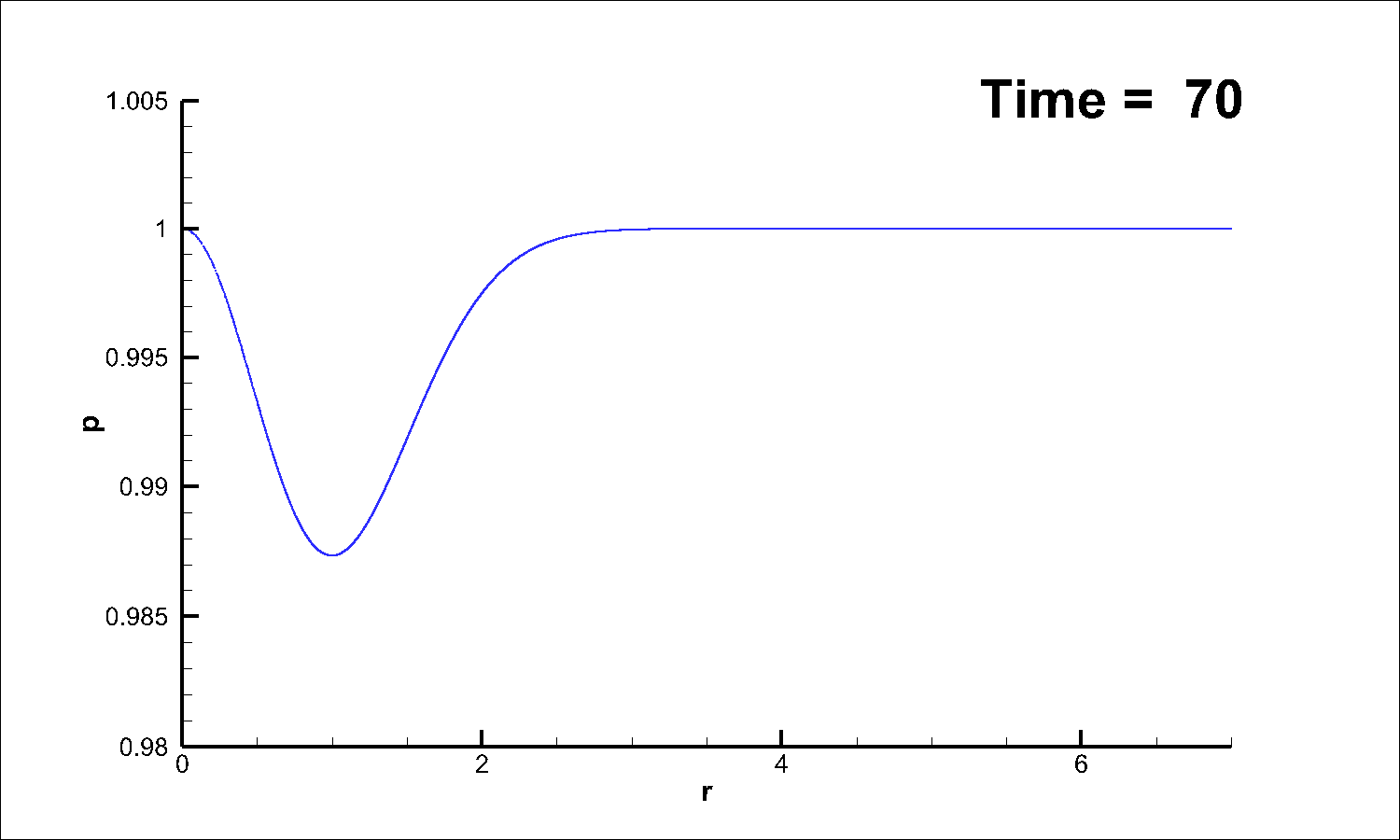}
    \caption{\RIIIcolor{MHD vortex solved with our $P_3P_3$ DG scheme on a moving Voronoi mesh of $1900$ elements using the Lloyd-like smoothing algorithm with $\mathcal{F} = 10^{-5}$. We report the scatter plot of the constant density profile (left) and of the pressure profile (right) at the final time $t_f=70$.}}
    \label{fig.MHDConstantRhoAndScatterP}
\end{figure}

\paragraph{Quality} In Figure~\ref{fig.MHDvortex_P3P3} we show the pressure profile and the magnetic 
field obtained with our fourth order $P_3P_3$ DG scheme at the initial time and after a long simulation with $t_f = 65$.
Once again, the profile of the vortex is simulated and conserved for a longer time with respect to 
standard conforming ALE scheme, for which mesh tangling would occur and stop the simulation earlier. 
\RIIIcolor{Moreover, scatter plots of the constant density profile and of the pressure profile are 
reported in Figure~\ref{fig.MHDConstantRhoAndScatterP}. One can observe that the errors in the density profile are very low. } 

\begin{figure*}[!bp]
    \centering
    \includegraphics[width=0.249\linewidth]{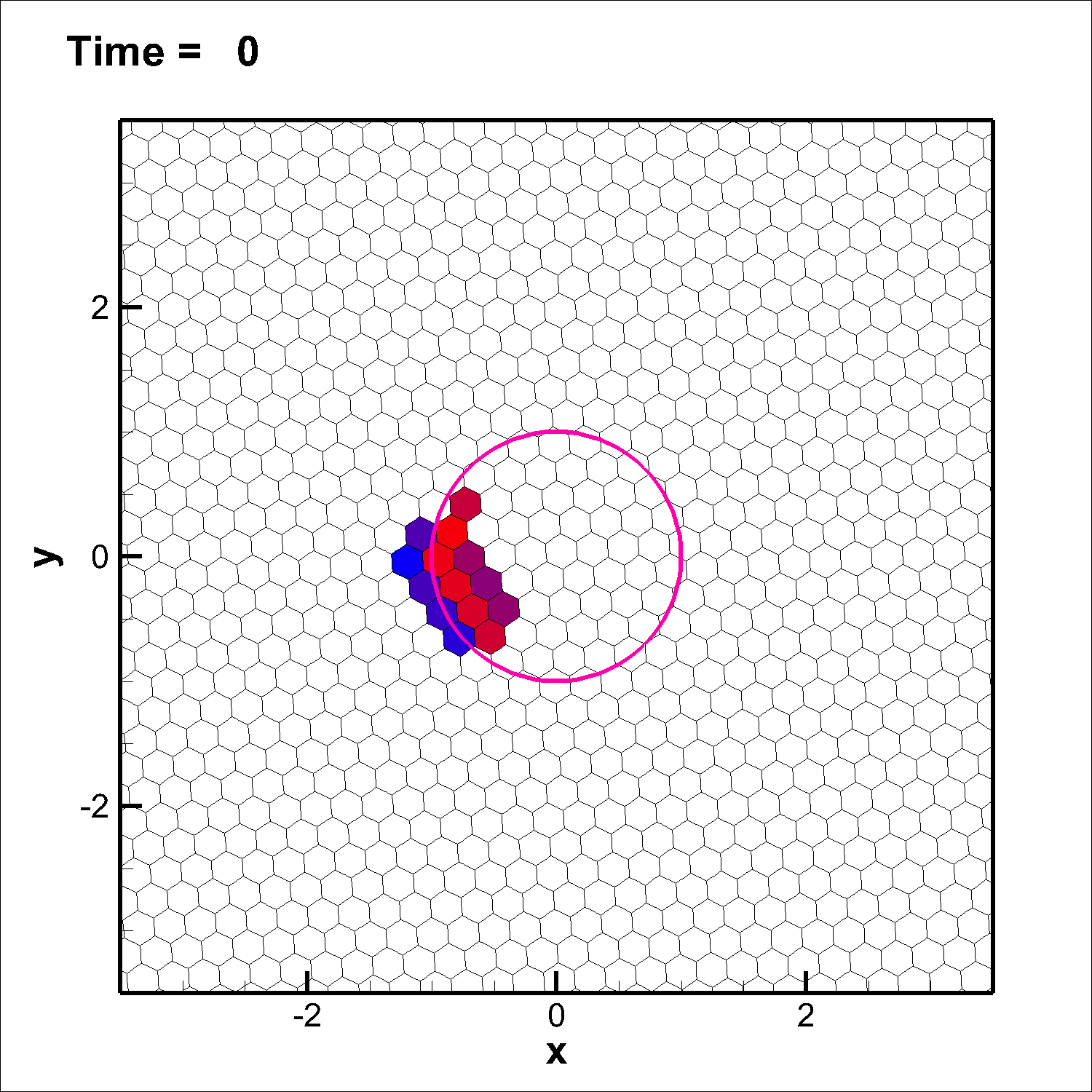}%
    \includegraphics[width=0.249\linewidth]{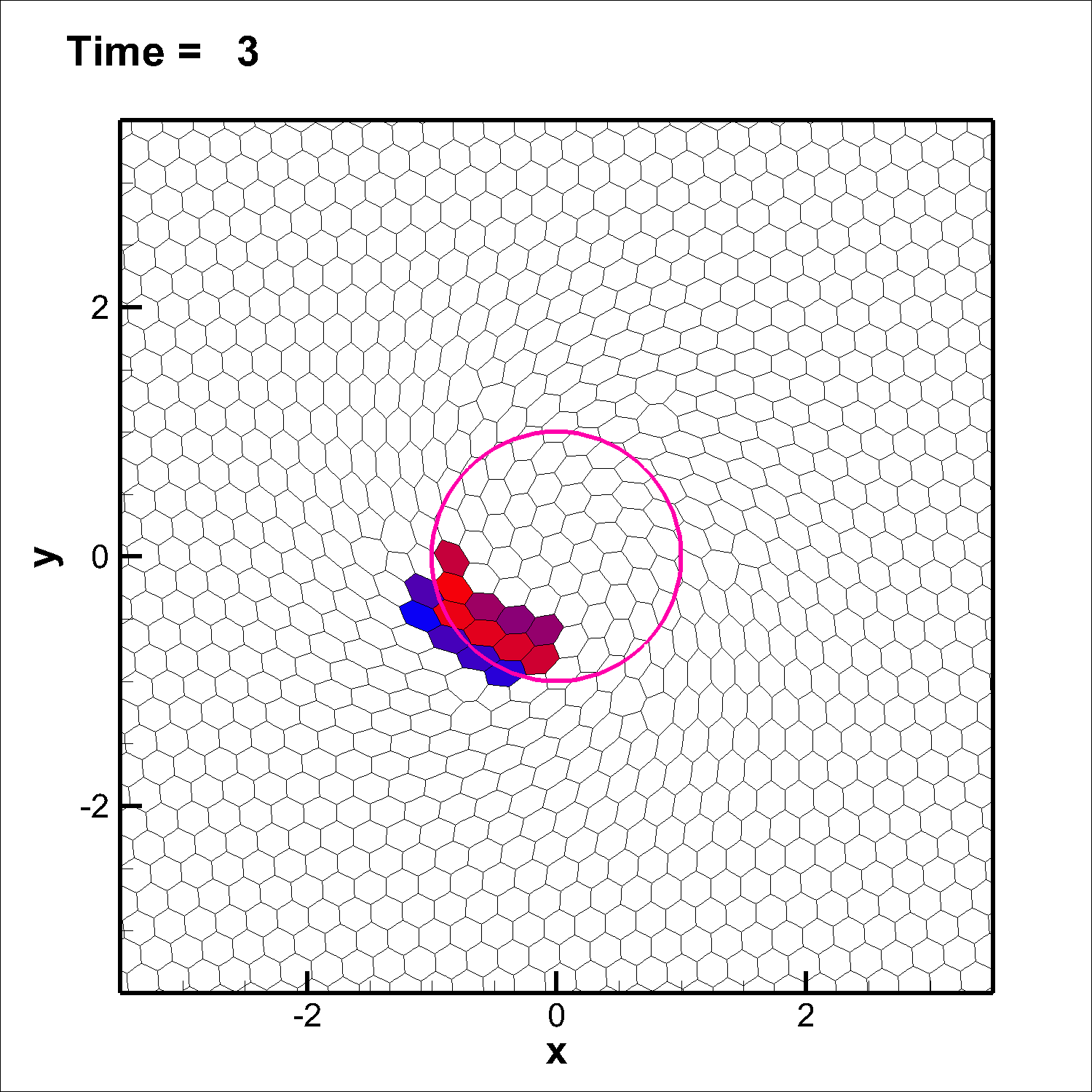}%
    \includegraphics[width=0.249\linewidth]{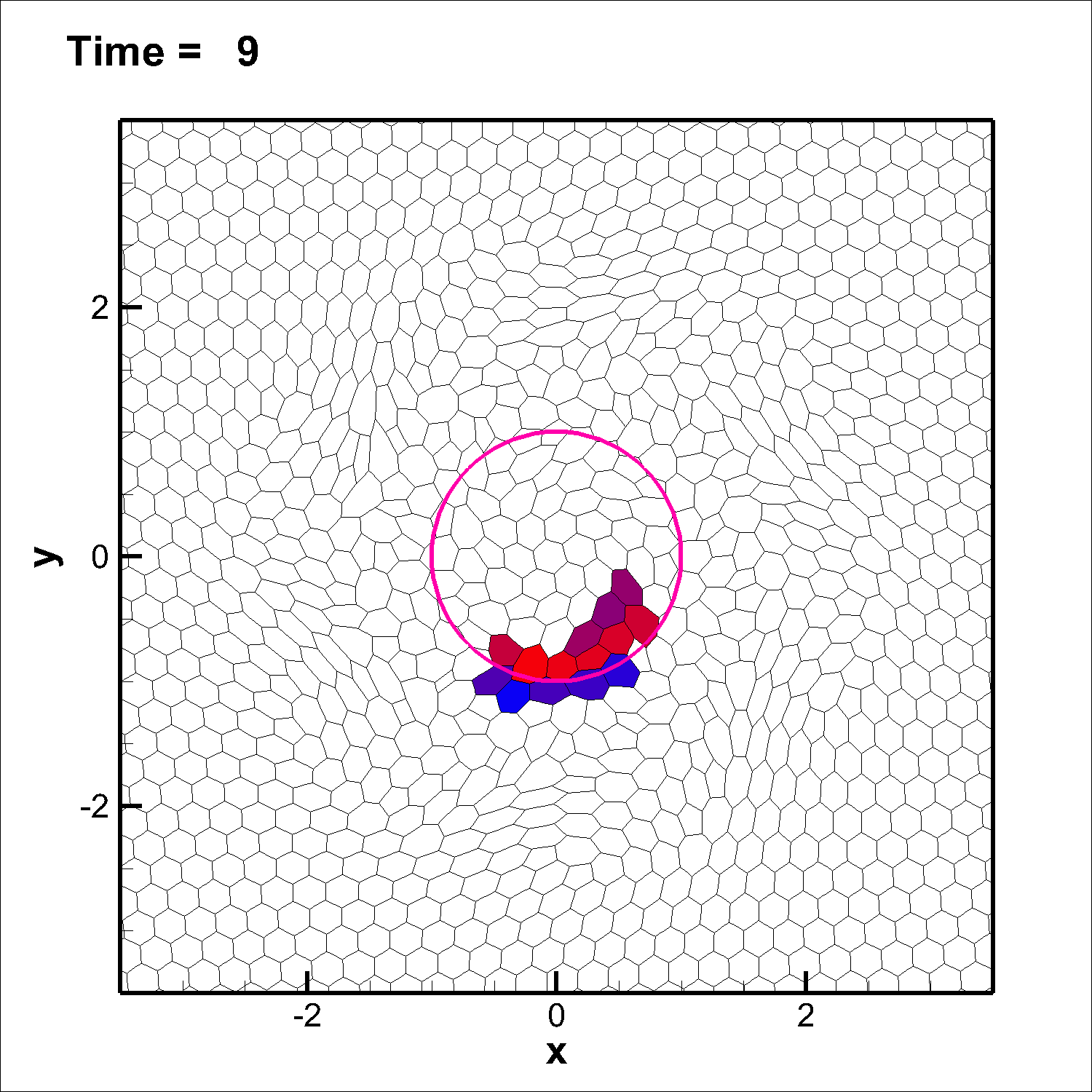}%
    \includegraphics[width=0.249\linewidth]{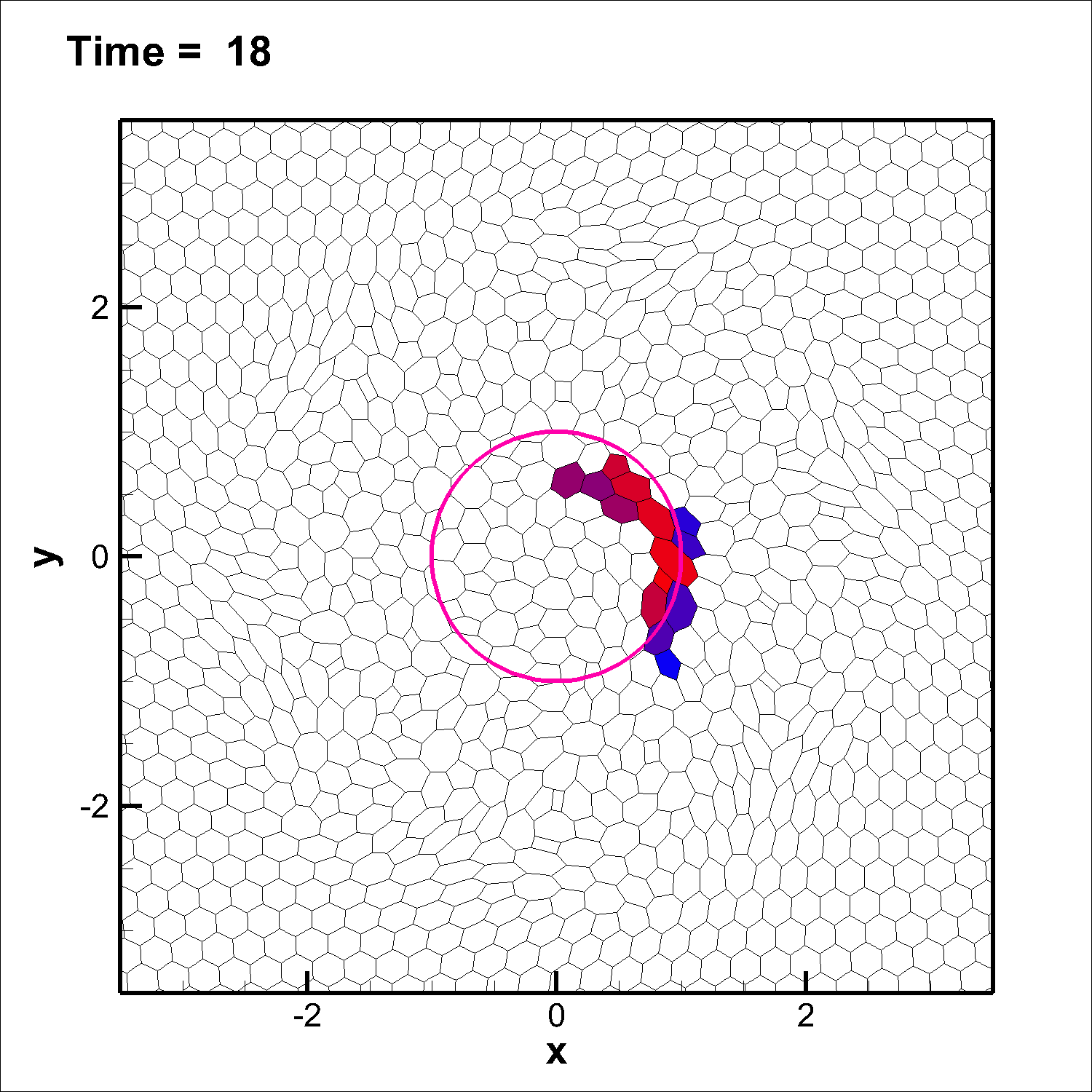}\\[-1pt]
    \includegraphics[width=0.249\linewidth]{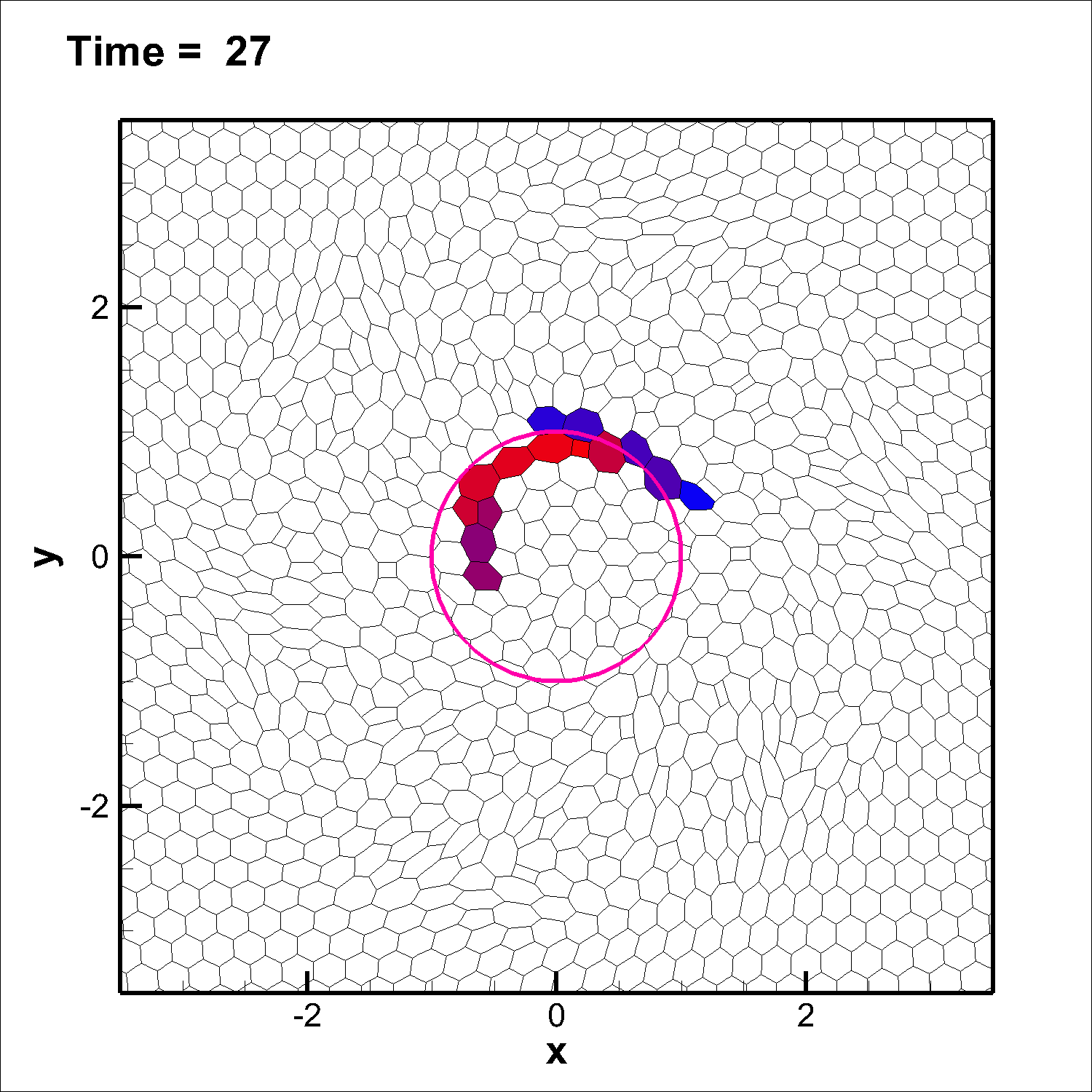}%
    \includegraphics[width=0.249\linewidth]{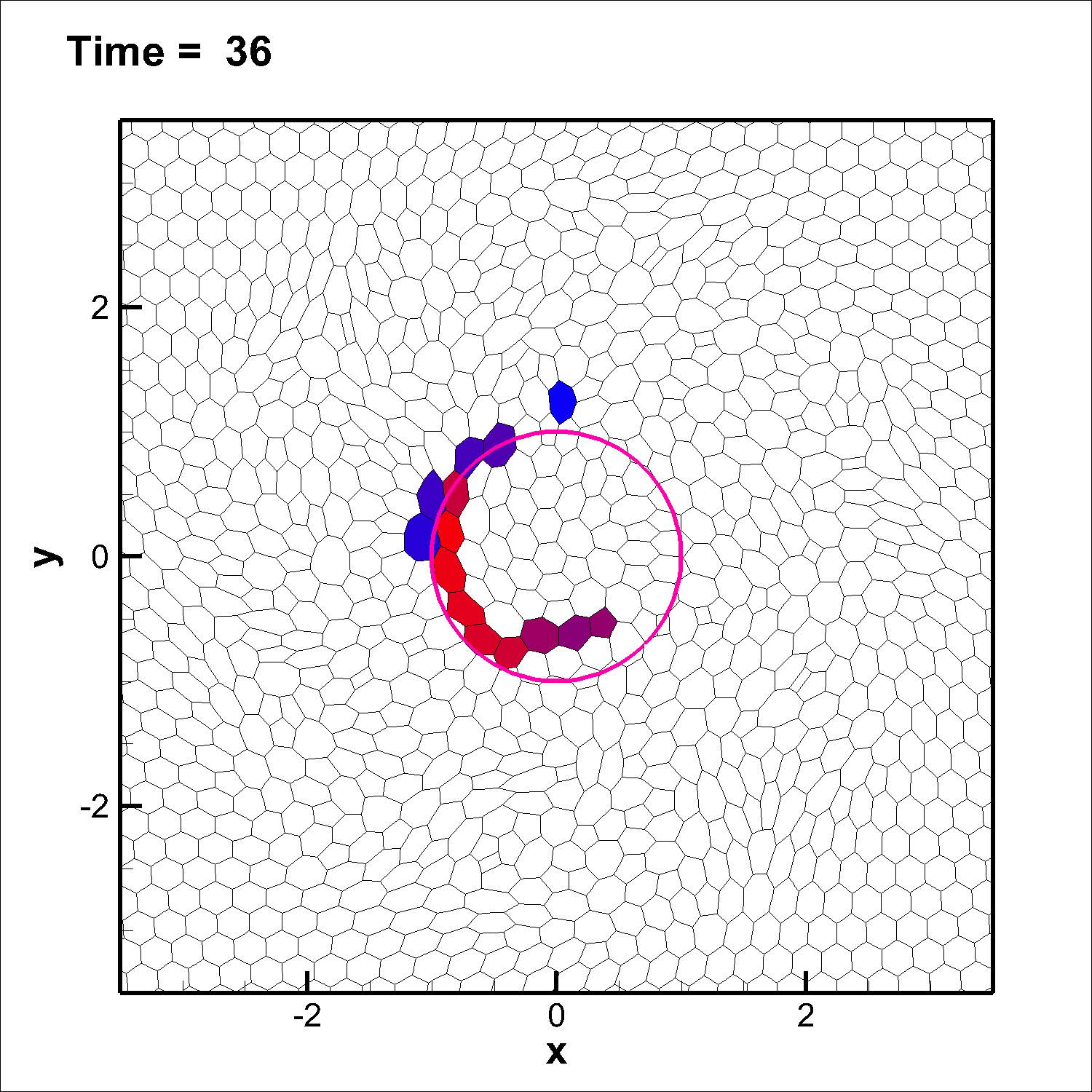}%
    \includegraphics[width=0.249\linewidth]{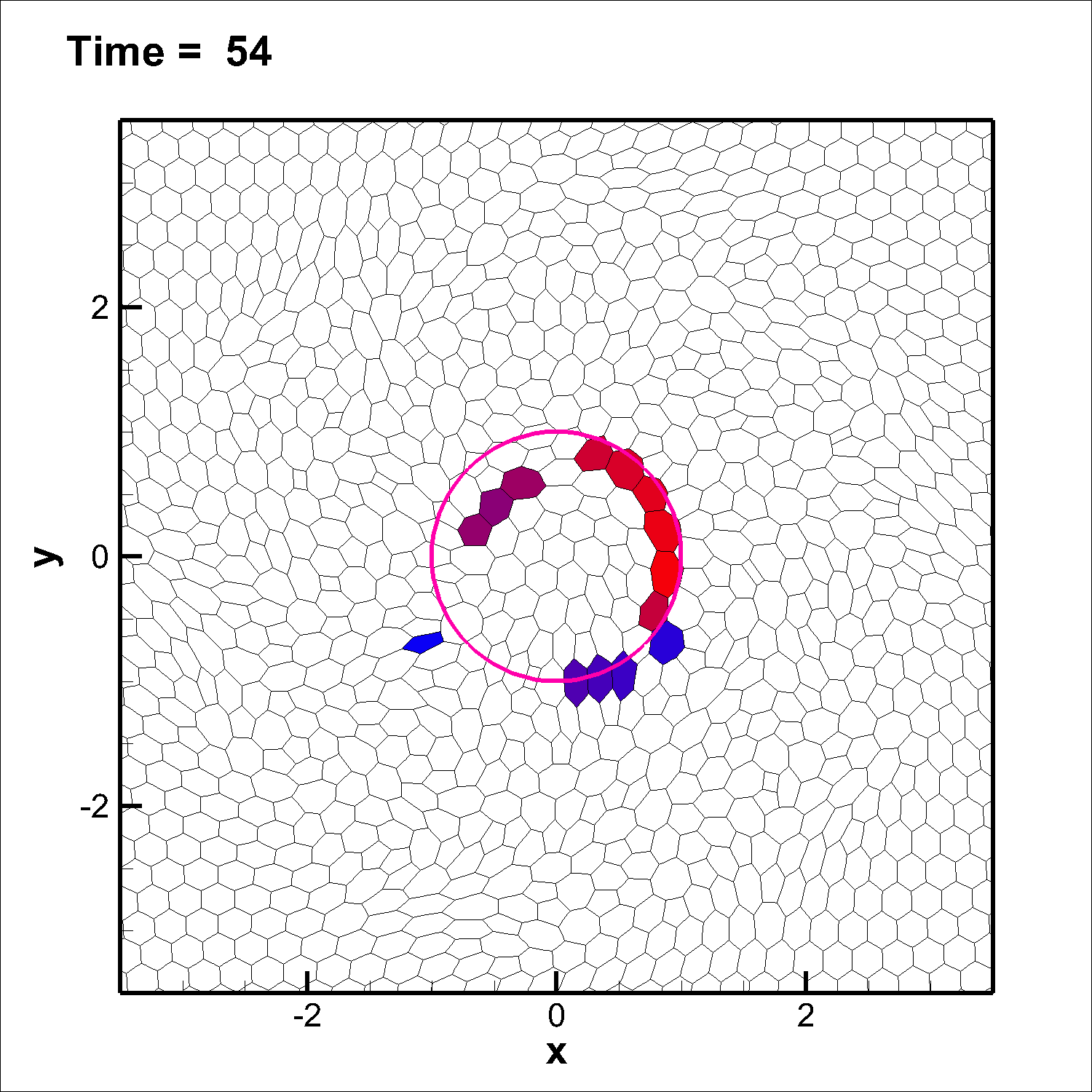}%
    \includegraphics[width=0.249\linewidth]{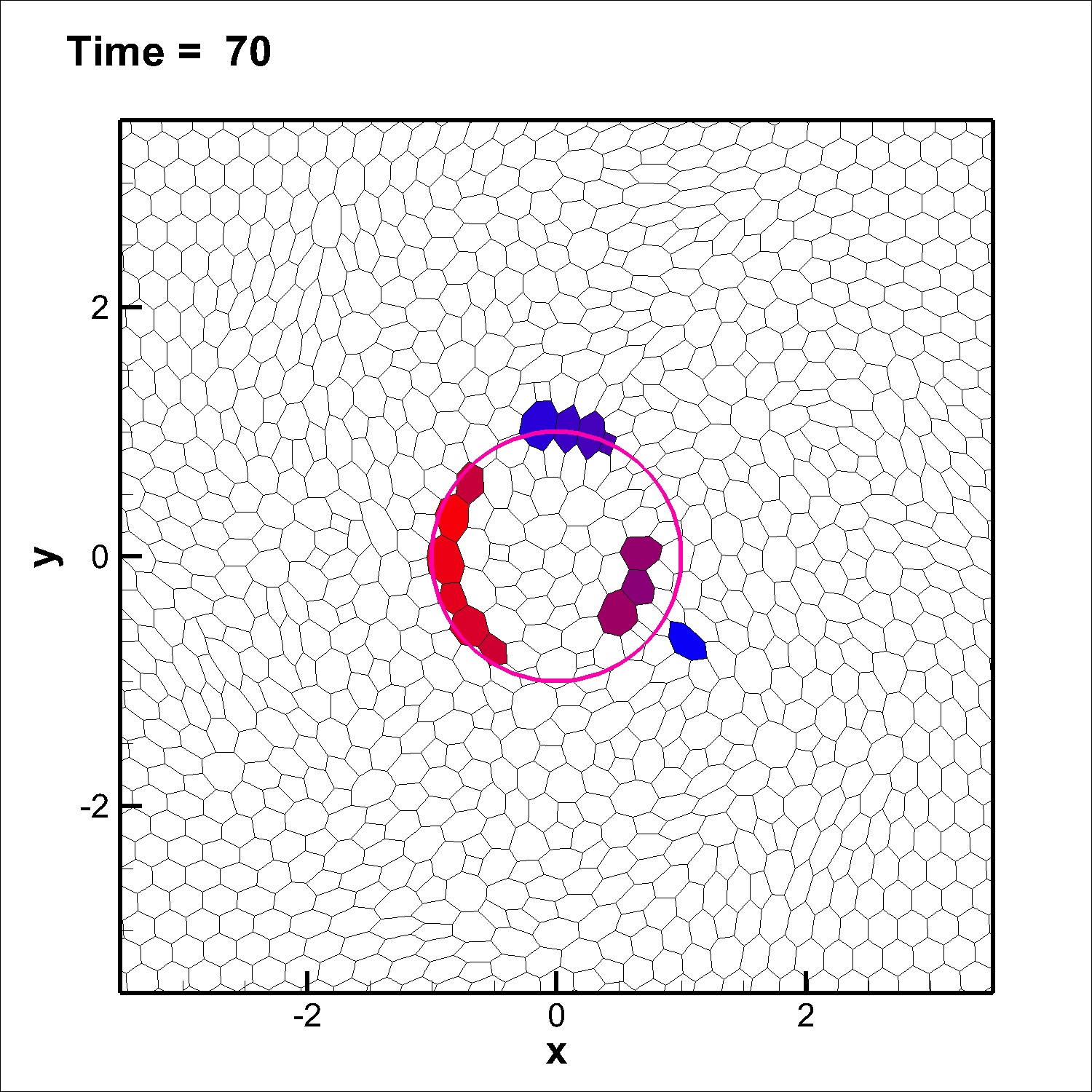}%
    \caption{\RShashkovcolor{MHD vortex solved with our $P_3P_3$ DG scheme on a moving Voronoi mesh of $1900$ elements using the Lloyd-like smoothing algorithm with $\mathcal{F} = 10^{-5}$.
            In this Figure we show the position of a bunch of initially neighbors elements at different times until $t=70$. 
            In this way one can notice the evolution of the grid topology during time and the necessity of allowing the mesh changing its 
            connectivity in order to correctly follow the fluid motion without distortion.}}
    \label{fig.MHDvortex_P3P3connectivity}
\end{figure*}

In Figure~\ref{fig.MHDvortex_P3P3connectivity} we report the position of a bunch of initially adjacent elements at different times: 
this clearly shows how strong the shear is, to which the mesh elements are subject, and why
the changes in the mesh topology are necessary. 

\RIVcolor{Finally, we report some statistics on the number of sliver elements created over the total number of time steps, 
and on the percentage of computational time required for the geometrical part of the code and for the $P_NP_M$ algorithm.
}
\RIVcolor{
With our $P_3P_3$ 
    DG scheme, the number of sliver elements that have been originated during the simulation
    until the final time $t_f=70$ is 21369, on a total of 62741 time steps. Three of these 
    time steps have been repeated through the MOOD loop described in Section~\ref{ssec.MOOD}. 
    The percentage of computational time employed by mesh regeneration and space time connectivity generation is  0.17\%, while the 
    fourth order
    predictor-corrector step on standard elements and on sliver elements accounts for 97.39\% and 7.5$\times10^{-4}$\% of the total wallclock time respectively.
    It turns out that the cost due to mesh rearrangement 
    and sliver computations is minimal, that MOOD restart activates very rarely 
    and that sliver elements are an essential ingredient to carry out the computation on moving
    meshes with topology change.}

\subsubsection{MHD rotor problem}
\label{test.MHDRotor}

\begin{figure}[!bp]
    \centering
    \includegraphics[width=0.33\linewidth]{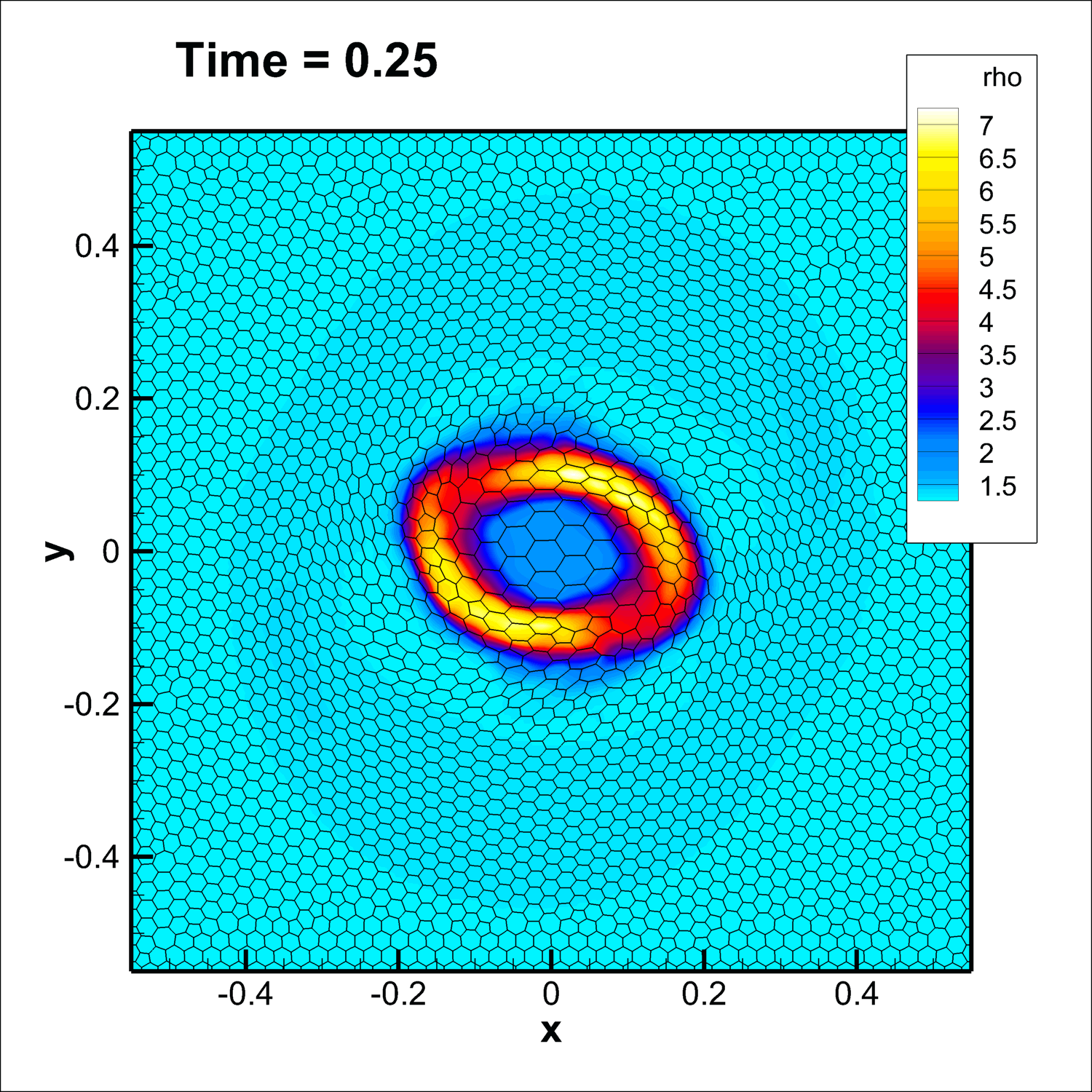}%
    \includegraphics[width=0.33\linewidth]{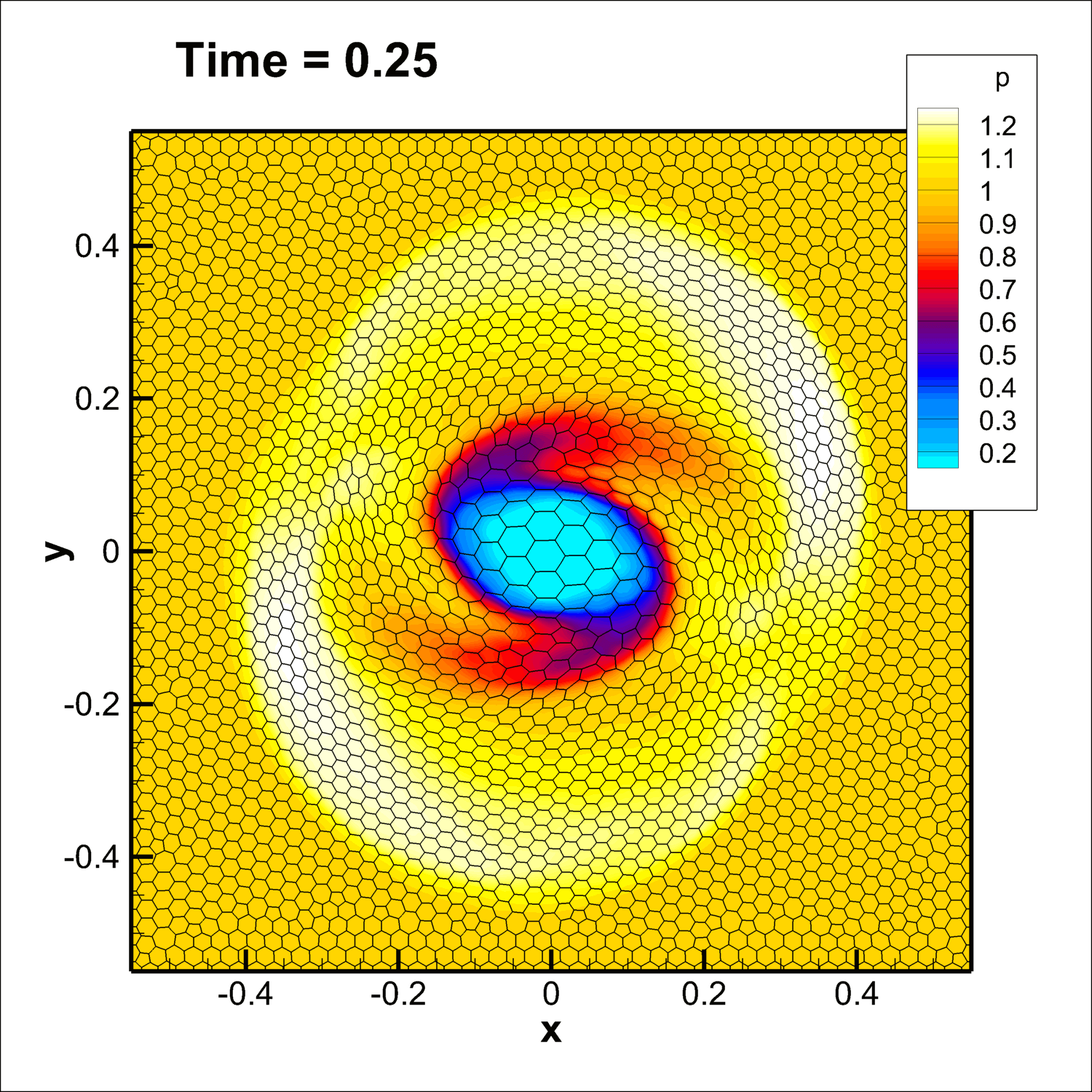}%
    \includegraphics[width=0.33\linewidth]{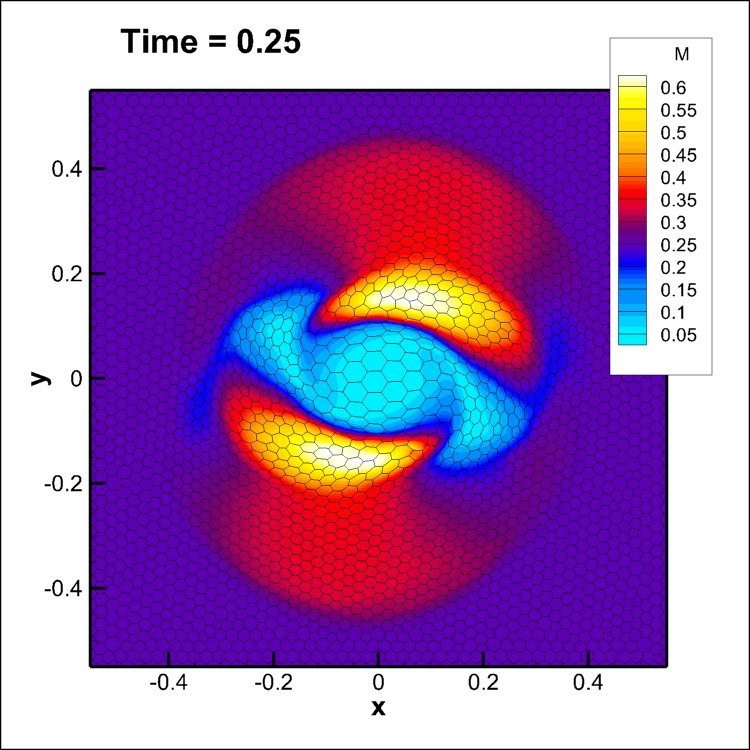}\\[-1pt]
    \includegraphics[width=0.33\linewidth]{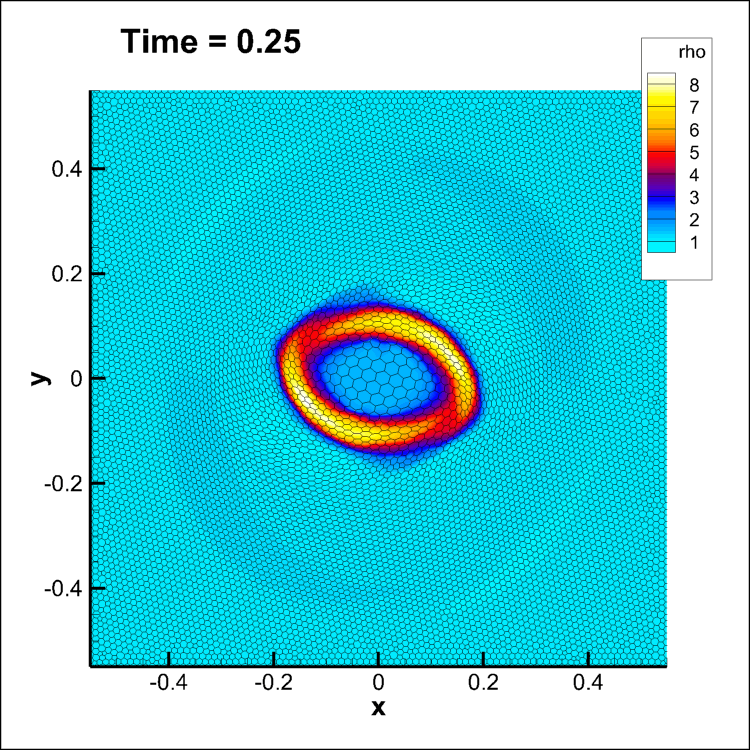}%
    \includegraphics[width=0.33\linewidth]{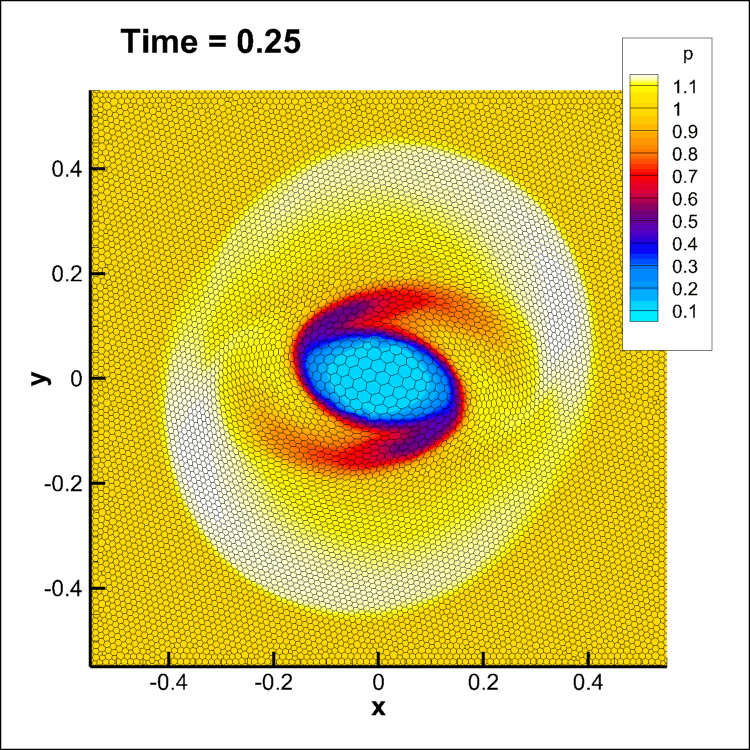}%
    \includegraphics[width=0.33\linewidth]{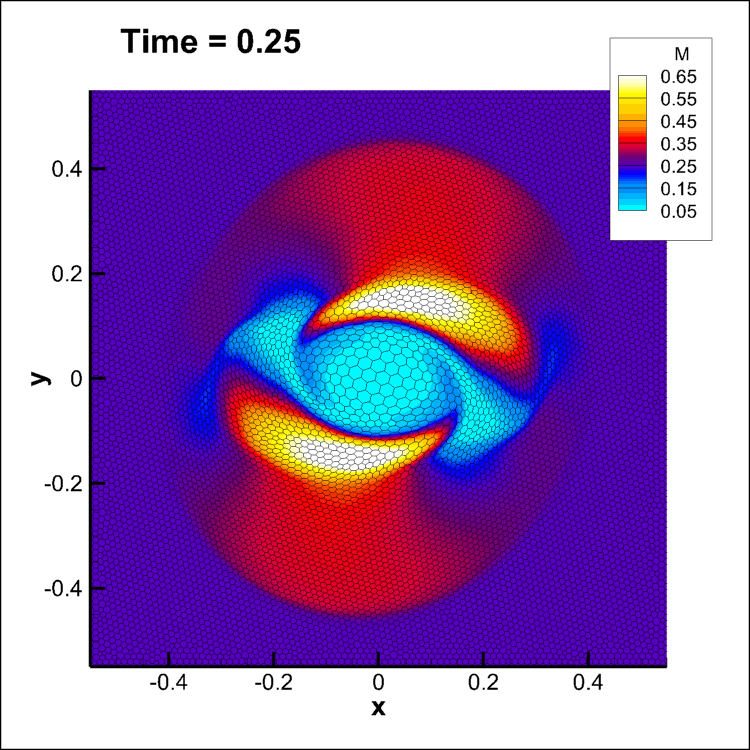}%
    \caption{MHD rotor problem solved with our $P_2P_2$ DG scheme on two moving Voronoi meshes respectively a coarse one made by $2727$
     and a finer one made by $10394$ elements. This test uses the Lloyd-like smoothing algorithm with $\mathcal{F} = 10^{-3}$.
    We depict the density profile (left) the pressure profile (middle) and the magnetic density profile $M = (B_x^2 + B_y^2+B_z^2)/(8\pi)$ (right) on the coarse mesh (top) and the fine mesh (bottom).}
    \label{fig:MHDrotor_P2P2_moving}
\end{figure}

This last MHD test case is the classical MHD rotor problem proposed by Balsara and Spicer in~\cite{BalsaraSpicer1999}. 
It consists of a rapidly rotating fluid of high density embedded in a fluid at rest with low density. Both fluids 
are subject to an initially constant magnetic field. 
\begin{figure}[!bp]
    \centering
    \includegraphics[width=0.33\linewidth]{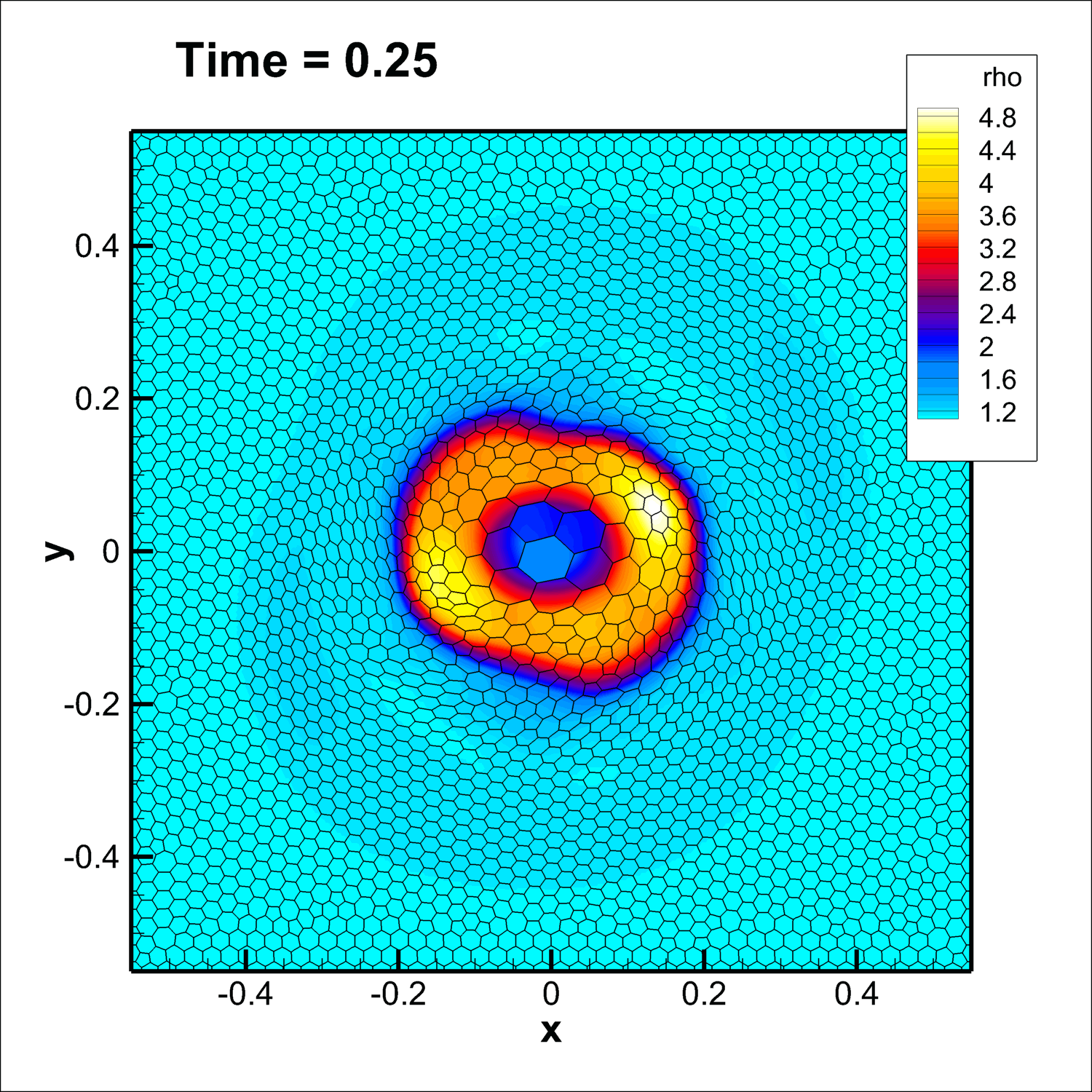}%
    \includegraphics[width=0.33\linewidth]{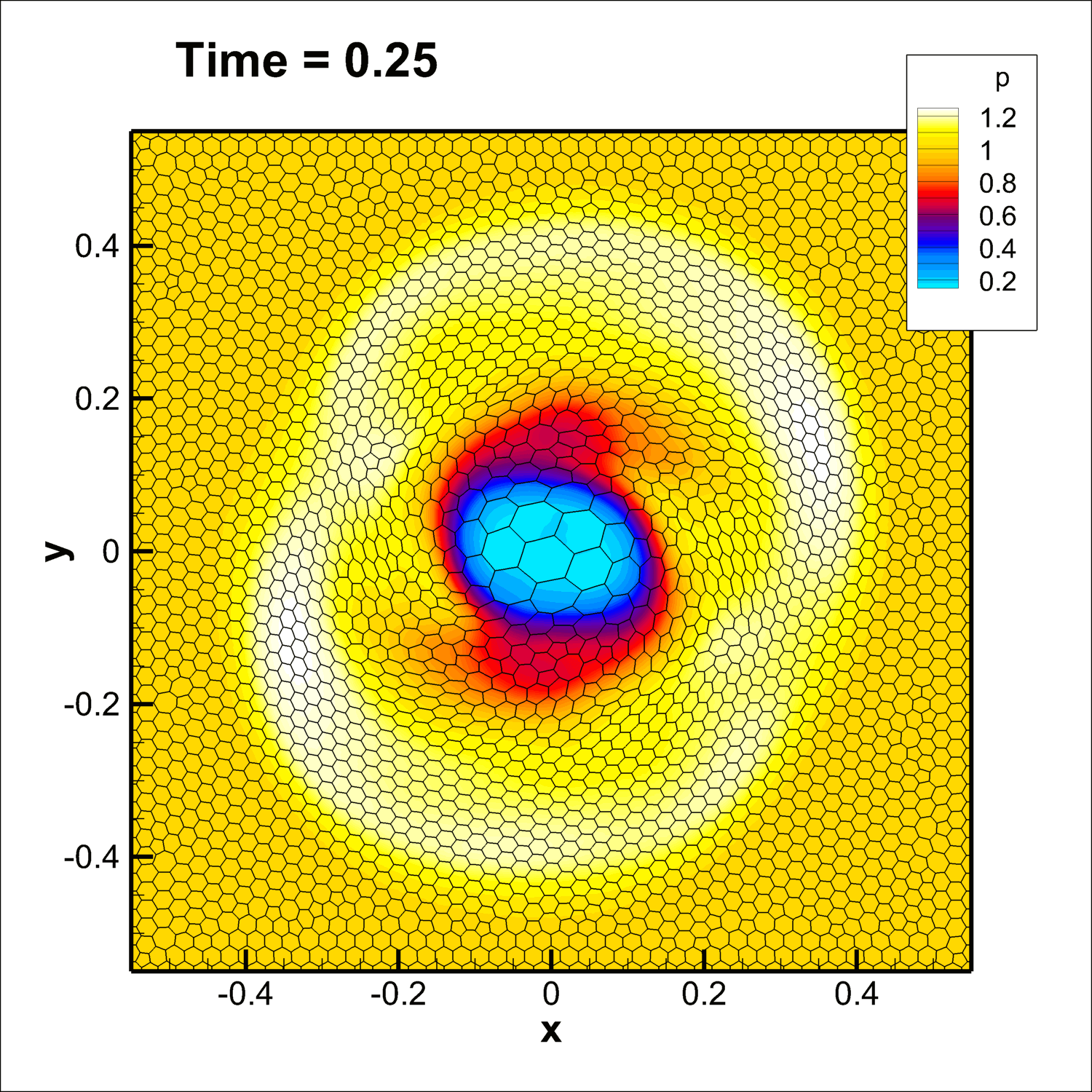}%
    \includegraphics[width=0.33\linewidth]{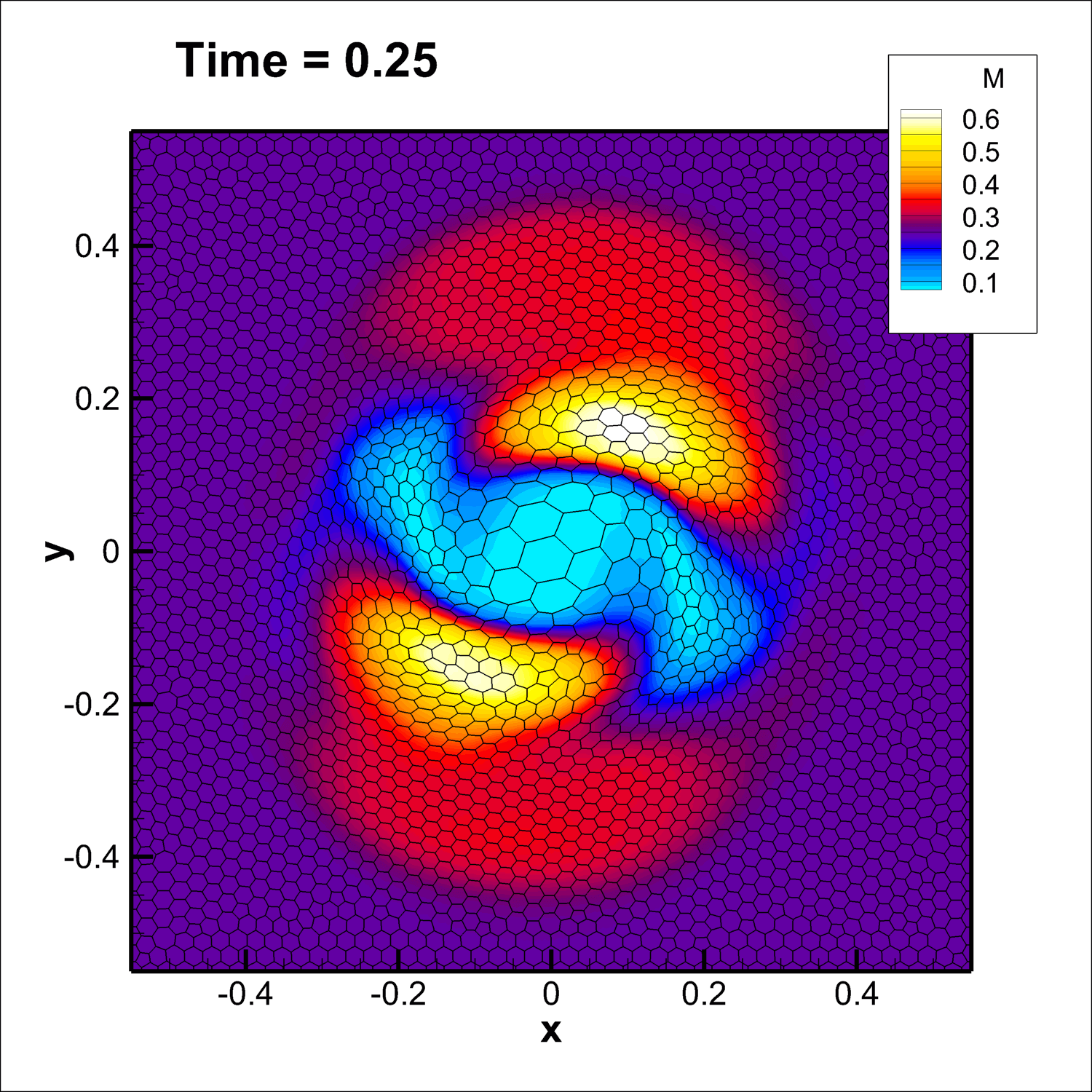}\\[-1pt] 
    \includegraphics[width=0.33\linewidth]{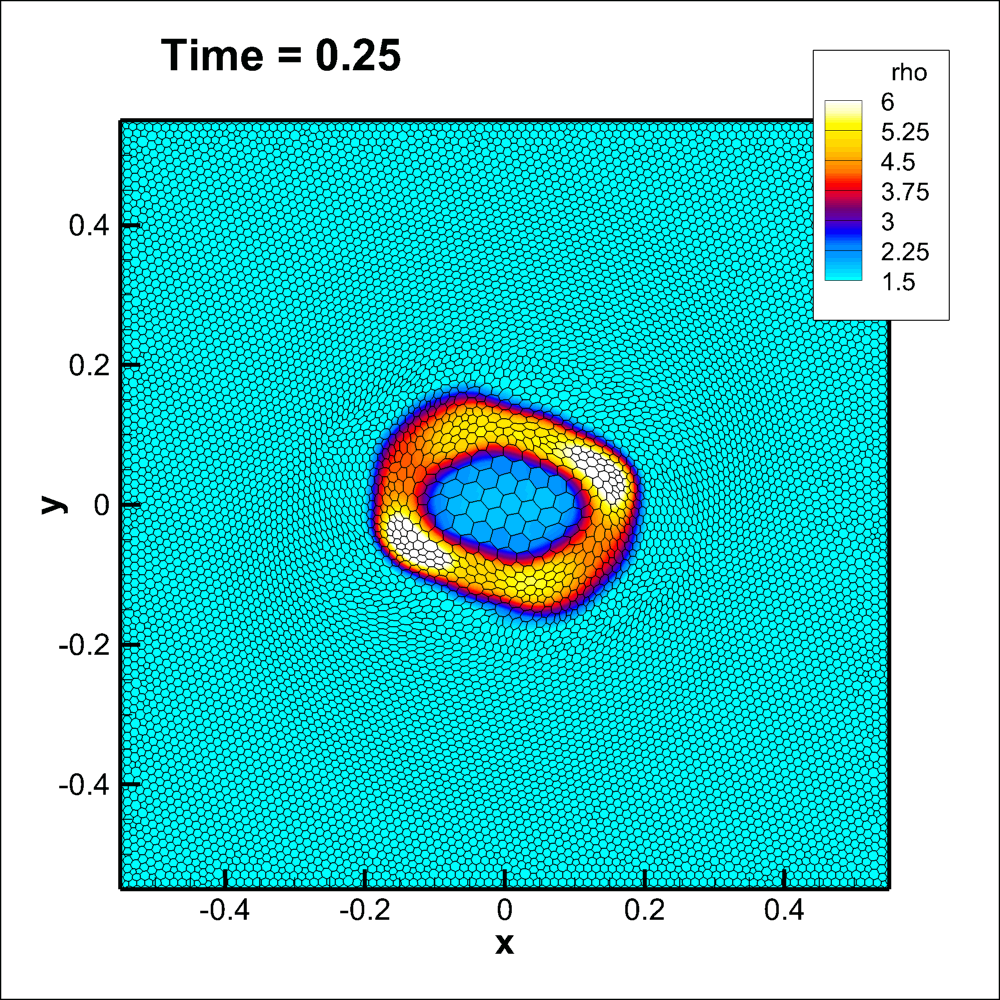}%
    \includegraphics[width=0.33\linewidth]{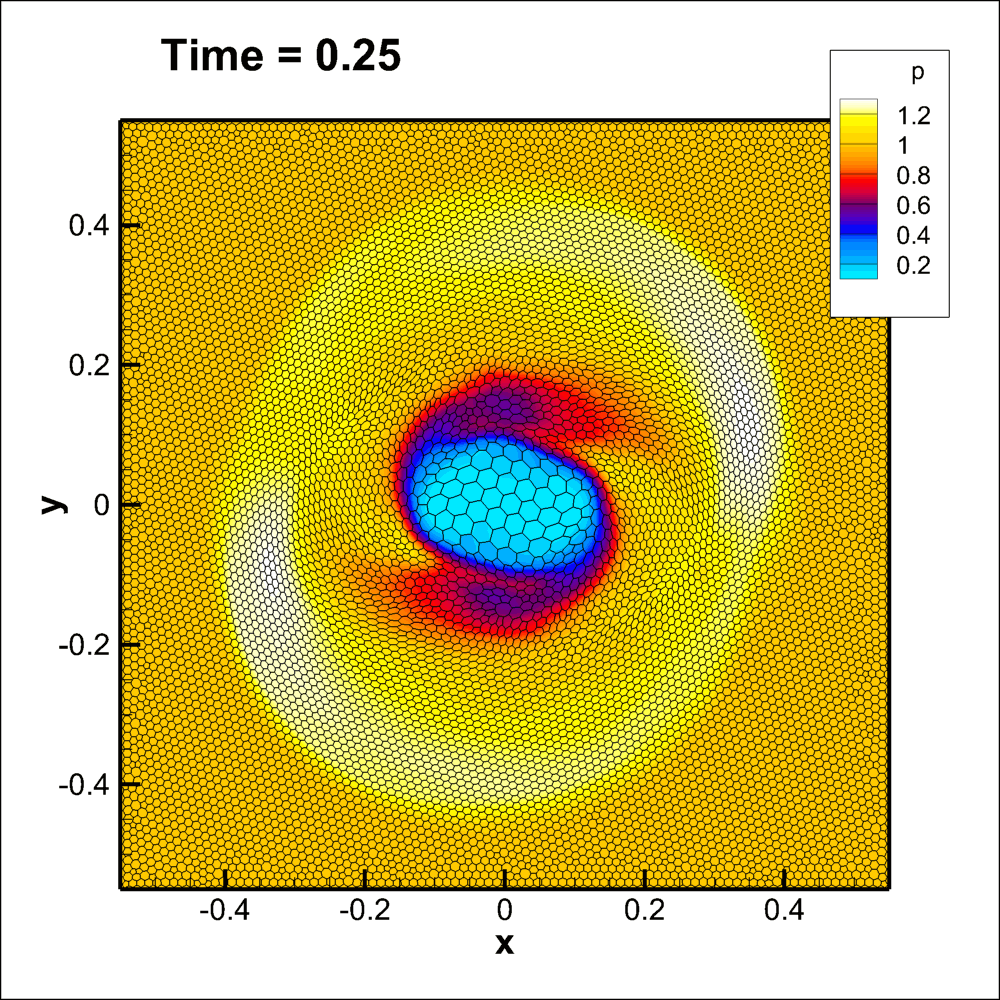}%
    \includegraphics[width=0.33\linewidth]{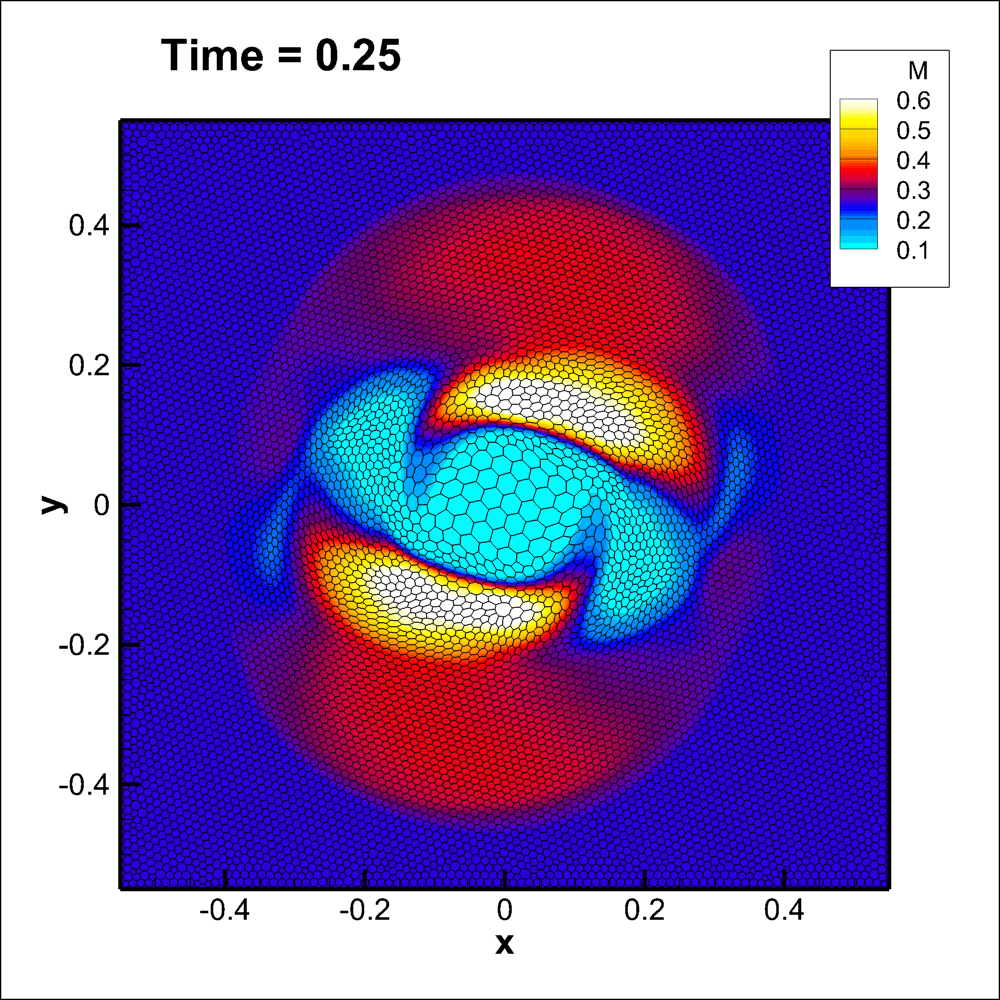}\\ 
    \caption{MHD rotor problem solved with our $P_0P_3$ fourth order FV scheme on on two moving Voronoi
     meshes respectively a coarse one made by $2727$ and a finer one made by $10394$ elements. This test uses the Lloyd-like smoothing algorithm with $\mathcal{F} = 10^{-3}$.  
     We depict the density profile (left) the pressure profile (middle) and the magnetic density profile $M = (B_x^2 + B_y^2+B_z^2)/(8\pi)$ (right).}
    \label{fig:MHDrotor_P0P3_moving}
\end{figure}
The rotor produces torsional Alfv\'en waves that are launched into the outer fluid at rest, resulting in a decrease
of angular momentum of the spinning rotor. The computational domain is taken to be $\Omega=[-0.55, 0.55]\times[-0.55, 0.55]$.
The density inside is $\rho=10$ for $0 \leq r \leq 0.1$ while the density of the ambient fluid at rest is
set to $\rho=1$. 
The rotor has an angular velocity of $\omega=10$. 
The pressure is $p=1$ and the magnetic field vector is set to $\B = (2.5, 0, 0)^T$ in the entire domain.
As proposed by Balsara and Spicer we apply a linear taper to the velocity and to the density in the range from 
$0.1 \leq r \leq 0.12$ so that density and velocity match those of the ambient fluid at rest at a radius of $r=0.12$.  
The speed for the hyperbolic divergence cleaning is set to $c_h=2$ and $\gamma=1.4$ is used. Wall boundary conditions are applied everywhere. 
We run this problem on a very coarse mesh $M_1$ made of $2727$ Voronoi elements and a finer one $M_2$ made of $10394$ Voronoi elements, and in two different configurations
\begin{description}
	\item[(a)]
	For the first test case we have applied our fourth order $P_0P_3$ Finite Volume scheme, see Figure~\ref{fig:MHDrotor_P0P3_moving}.
	\item[(b)]
	Then we have employed our third order accurate $P_2P_2$ DG scheme, see the results in Figure~\ref{fig:MHDrotor_P2P2_moving}.
\end{description} 
In all the cases, we can observe a good agreement between the obtained numerical results and those available in the literature. 
Moreover also a visual convergence can be deduced from the presented results on a coarse and a finer mesh.
Finally, by comparing with the reference results found in literature, it is clear that 
\RIIcolor{the results obtained with the DG scheme of order \textit{three} 
are sharper than those of the FV scheme of order \textit{four}}. 

\section{Conclusion}
\label{sec.conclusion}

In this work we have developed the first unified framework for explicit and arbitrary high order accurate direct  Arbitrary-Lagrangian-Eulerian FV and DG schemes on moving unstructured Voronoi meshes with \textit{topology change}, in order to benefit simultaneously from high order methods, high quality grids and substantially reduced numerical dissipation. 
Indeed, we would like to \RIIIcolor{emphasize} that in the current literature at least one of the previous ingredients is always missing: Lagrangian methods are usually affected by large mesh distortions, and available algorithms which are able to avoid it are usually only at most second order accurate; Eulerian methods are in general high order accurate, but exhibit significant dissipation errors due to the convective terms. 
In particular, the results on vortex flows give evidence of the advantages conveyed by the  proposed algorithm, and a large set of different numerical tests shows its robustness and
accuracy. 

We recall that the key ingredient of our novel algorithm is the generalization of the $P_NP_M$ scheme~\cite{Dumbser2008, Lagrange3D} to Voronoi and sliver space--time elements, which has required the investigation of several intricate steps.  
First, an automatic procedure to construct space-time meshes resulting from the connection of moving Voronoi meshes with different topologies at two consecutive time levels has not been used before.  
Next, computations on Voronoi elements have required their subdivision into triangular prisms, 
the adaptation of the basis functions, the neighbors search, the projection and reconstruction algorithms, 
and also a change in the notions of areas, volumes and characteristic mesh sizes. 
Finally, the presence of space--time sliver elements forced us to revisit the core of the $P_NP_M$ scheme, 
i.e. the space--time predictor and the update of the solution through flux computations, in order to maintain the property of mass, momentum and energy conservation, essential for solving nonlinear hyperbolic equations.


Future work may regard the improvement of the present algorithm in \textit{three} different directions.  
First, we plan to incorporate a path-conservative method to treat non conservative products, so that also a well balanced treatment of sources and a proper well-balanced preservation of stationary equilibria of the PDE system will be possible, following the ideas outlined in~\cite{Pares2006, Castro2008, gaburro2018diffuse, gaburro2018well, del2018asymptotic, grosheintz2019high, berberich2019high}. \RVcolor{Furthermore, non-conservative products would also allow the straightforward numerical discretization of diffuse interface models for compressible multi-phase flows, see e.g.  \cite{BaerNunziato1986,SaurelAbgrall,AbgrallSaurel,USFORCE2,DIM2D,DIM3D,RomenskiTwoPhase2007,RomenskiTwoPhase2010}}. Future applications of our new algorithm will then concern the unified first order hyperbolic formulation of continuum physics recently proposed in \cite{PeshRom2014,HPRmodel,HPRmodelMHD}.  
Above all, we plan to incorporate the presented high order techniques inside the massively parallel second order accurate ALE-FV code \textit{AREPO}~\cite{Springel}, 
which currently includes one of the most advanced moving Voronoi mesh generators in 2D and 3D. In  this way, \RIcolor{we will add the possibility of refining or coarsening our mesh by adding and deleting generators,} and we will gain a very efficient parallel environment which also redistributes the moving elements among the CPU cores in a dynamic load balancing approach. 
At this point, even challenging astrophysical simulations would be feasible.
Finally, the extension to three-dimensional domains is also envisaged. Although the  \texttt{AREPO} code is already available in three space dimensions, it is currently still low order accurate and does not yet provide any information about the space--time connectivity of the Voronoi meshes between two consecutive time levels, which is, however, needed by our high order DG and FV schemes. In our opinion, the realization of a coherent 4D space--time connection will be complex, but feasible (a first hint in this direction could be taken from~\cite{re2017interpolation}), and formally the direct ALE $P_NP_M$ schemes would require the same adaptations here  introduced in order to deal with degenerate four dimensional space--time control volumes.

\section*{Acknowledgments}

The research presented in this paper has been partially financed by the European Research Council (ERC) under the 
European Union's Seventh Framework Programme (FP7/2007-2013) with the research project \textit{STiMulUs},  
ERC Grant agreement no. 278267. 

E.~G. and W. B. have been also supported by two national mobility grants for young researchers in Italy, funded by GNCS-INdAM. 
E.~G acknowledges the support given by the University of Trento in Italy through the \textit{UniTN Starting Grant}.
W.~B. acknowledges support via INdAM (Italian National Institute of High Mathematics) under the program \textit{Young researchers funding 2018}. 

M.~D. has been funded by the European Union's Horizon 2020 Research and Innovation  Programme under the project \textit{ExaHyPE}, grant no. 671698 (call FETHPC-1-2014). M.D. also acknowledges the financial support received from 
the Italian Ministry of Education, University and Research (MIUR) in the frame of the Departments of Excellence  Initiative 2018--2022 attributed to DICAM of the University of Trento (grant L. 232/2016) and in the frame of the 
PRIN 2017 project \textit{Innovative numerical methods for evolutionary partial differential equations and  applications}. Furthermore, M.~D. has also received funding from the University of Trento via the Strategic 
Initiative \textit{Modeling and Simulation}. 

Moreover, C.~K. and V.~S. acknowledge the German Science Foundation (DFG) grant `Exascale simulations of the evolution of the universe including magnetic fields' within the priority program SPP1648 `Software for Exascale Computing'. \\

Last but not least, the authors would also like to thank all the seven anonymous referees of this paper for their  constructive comments and suggestions.


\clearpage 

\bibliographystyle{plain}
\bibliography{references}

\end{document}